%% file: 2002-23.tex
\documentclass{gtart}

\input gtoutput

\lognumber{293}

\volumenumber{6}
\papernumber{23} 
\volumeyear{2002}
\pagenumbers{657}{814} 
\received{21 November 2002}
\accepted{18 December 2002}
\published{18 December 2002}
\proposed{Robion Kirby}
\seconded{Yasha Eliashberg, John Morgan}

\usepackage{amssymb,amsmath}

\makeatletter
\def\@listii{\leftmargin12pt\parsep 0pt\topsep 5pt 
 \itemsep2pt plus3pt minus2pt}
\def\@listiii{\leftmargin12pt\parsep 0pt\topsep 5pt 
 \itemsep2pt plus3pt minus2pt}
\makeatother

\let\bigskip\medskip

\def\autonum{\addtocounter {equation}{1}\hfill{\rm(\theequation)}}
\def\autonumm{\addtocounter {equation}{1}\hfill\rlap{\hbox to
35pt{\hss\rm(\theequation)}}}
\def\narrower{\bgroup\par\medskip\leftskip28.5pt\rightskip35pt}
\def\endnarrower{\par\medskip\egroup}
\def\itemm{\noindent\llap{$\bullet$\qua}}

\newtheorem{lemma*}{Lemma \emph{}}

\numberwithin{equation}{section}

\def\lam{\lambda}
\def\opart{\overline{\partial}}
\def\oz{\overline{z}}
\def\olam{\overline{\lambda}}
\def\arg{{\rm \;arg}}
\newcommand{\Lam}{\Lambda}

\def\cN{\mathcal{N}}
\def\cH{\mathcal{H}}
\def\cM{\mathcal{M}}
\def\fM{\mathfrak{M}}
\def\dist{{\rm dist}}
\def\fX{\aleph}
\def\cL{\mathfrak{L}}
\def\cR{\mathfrak{R}}

\def\fM{\mathfrak{M}}
\def\cP{\mathcal{P}}
\def\uD{\underline{D}}
\def\uSigma{\underline{\Sigma}}
\def\cO{\mathcal{O}}
\def\uk{\underline{k}}
\def\uC{\underline{C}}

\def\dim{\text{dim}}

\def\cE{\mathcal{E}}

\def\mr{\mathbb{R}}

\def\mz{\mathbb{Z}}
\def\mc{\mathbb{C}}
\def\mP{\mathbb{P}}

\def\uZ{\underline{Z}}
\def\ux{\underline{x}}
\def\uv{\underline{v}}
\def\uc{\underline{c}}
\def\uw{\underline{w}}

\def\uvar{\underline{\varphi}}
\def\uh{\underline{h}}
\def\GF{\underline{G}}
\def\uvar{\underline{\varphi}}
\def\uth{\underline{\theta}}
\def\la{\langle}
\def\ra{\rangle}

\def\oww{\overline {w}}

\def\Pic{\text{\Pic}}
\def\uD{\underline{D}}

\begin{document}

\title{A compendium of pseudoholomorphic beasts\\in $\mr \times (S^1\times
S^2)$}
\asciititle{A compendium of pseudoholomorphic beasts in R x (S^1 x S^2)}

\covertitle{A compendium of pseudoholomorphic beasts\\in 
${\noexpand\bf R} \times (S^1\times S^2)$}

\author{Clifford Henry Taubes}
\address{Department of Mathematics,
Harvard University\\Cambridge, MA 02138, USA} 
\email{chtaubes@math.harvard.edu}

\begin{abstract}
This article describes various moduli spaces of pseudoholomorphic curves
on the symplectization of a particular overtwisted
contact structure on $S^1\times S^2$.  This contact structure appears when one
considers a closed self dual form on a 4--manifold as a symplectic form on
the complement of its zero locus.  The article is focussed mainly on
disks, cylinders and three-holed spheres, but it also supplies groundwork
for a description of moduli spaces of curves with more punctures and
non-zero genus.
\end{abstract}

\asciiabstract{This article describes various moduli spaces of
pseudoholomorphic curves on the symplectization of a particular
overtwisted contact structure on S^1 x S^2.  This contact structure
appears when one considers a closed self dual form on a 4-manifold as
a symplectic form on the complement of its zero locus.  The article is
focussed mainly on disks, cylinders and three-holed spheres, but it
also supplies groundwork for a description of moduli spaces of curves
with more punctures and non-zero genus.}



\maketitlepage

\section{Introduction}

The purpose of this article is to describe various moduli spaces of
pseudoholomorphic subvarieties in the symplectization of a certain
over twisted contact 1--form on $S^1 \times S^2$.  This said, the
motivation for such a study comes from 4--manifold differential
topology using three key observations.  Here is the first: As
explained in \cite{T1} and \cite{T2}, every compact, oriented 4--manifold with
positive sum of second Betti number and signature has a closed 2--form
that is symplectic where non-zero and whose zero set is a finite,
disjoint union of embedded circles.  Moreover, as explained in \cite{T2},
this closed 2--form restricts to a certain $[0, \infty) \times (S^1
\times S^2)$ neighborhood of each of its vanishing circles as the
symplectization of a particular contact 1--form, the one of interest
here.

Here is the second key observation stemming from \cite{T2}: Numerical
invariants of the moduli spaces of pseudoholomorphic subvarieties in
the complement of the zero circles most probably contain the
4--manifold's Seiberg--Witten invariants.  An optimist would hope to
find novel 4--manifold invariants here as well \cite{T3}.

Granted the second key observation, here is the third: Hofer \cite{H1,H2,H3};
Hofer, Wysocki and Zehnder \cite{HWZ1,HWZ2,HWZ3} (see also the references in \cite{H3});
Eliashberg \cite{E}, and Eliashberg with Hofer \cite{EH} have studied the
salient issues that confront the construction of numerical invariants
from moduli spaces of pseudoholomorphic subvarieties on non-compact
symplectic manifolds with symplectization type ends.  In particular,
they teach that such constructions require an understanding of the
analogous moduli spaces on the corresponding symplectizations.  In any
event, given the 4--manifold circumstances just described, the
symplectization is that of the contact form in question on $\mr\times
(S^1\times S^2)$.

With the preceding understood, it is time to be precise about the
relevant geometry.  For this purpose, introduce coordinates $(s, t,
\theta, \varphi)$ for $\mr\times (S^1\times S^2)$ where $s$ is the
coordinate for the $\mr$ factor in $\mr\times (S^1\times S^2)$, $t \in
\mr/2\pi\mz$ is the coordinate for the $S^1$ factor and
$(\theta,\varphi) \in [0, \pi]\times (\mr/2\pi\mz)$ are standard
spherical angle coordinates for the $S^2$ factor.  This done, the
contact form in question is
\begin{equation}
\alpha\equiv-(1-3\cos^2\theta) dt-\sqrt{6}\cos\theta\sin^2\theta
d\varphi.
\end{equation}
The resulting symplectic form on $\mr\times (S^1\times S^2)$ is
\begin{equation}
\omega=d(e^{-\sqrt{6}s}\alpha).
\end{equation}
In this regard, note that the convention here is such that the $s
\rightarrow\infty$ end of $\mr\times (S^1\times S^2)$ is the concave
side end in that $|\omega|$ drops to zero as $e^{-\sqrt{6}s}$ in this
direction.  Conversely, the end where $s \rightarrow -\infty$ is the
convex end.  Said differently, the contact form $\alpha$ is of concave
type with $S^1\times S^2$ viewed as the boundary of $[0, \infty)
\times (S^1 \times S^2)$.

By the way, the factor of $\sqrt{6}$ that enters above and
subsequently propagates throughout this article is a consequence of a
desire to have $\omega$ define a self-dual 2--form with respect to the
standard product metric, $ds^2 + dt^2 + d\theta^2 + \sin^2 \theta
d\varphi^2$, on $\mr \times (S^1 \times S^2)$.
  
It proves convenient in the ensuing discussion to have introduced
functions $f$ and $h$ on $\mr \times (S^1 \times S^2)$ defined as
follows:
\begin{equation}
f \equiv e^{-\sqrt{6}s} (1 - 3 \cos^2 \theta) \quad \text{and}\quad h
\equiv\sqrt{6}\; e^{-\sqrt{6} s} \cos \theta \sin^2\theta .
\end{equation}
This done, we have
\begin{equation}
\omega = dt \wedge df + d\varphi\wedge dh .
\end{equation}

The almost complex structure $J$ used here to define the term
`pseudoholomorphic' is specified by the relations
\begin{equation}
J\partial_t = g\partial_f\quad\text{and}\quad J\partial_\varphi = g
\sin^2 \theta\partial_h,
\end{equation}
where $g \equiv\sqrt{6} \;e^{-\sqrt{6} s} (1 + 3 \cos^4
\theta)^{1/2}$.  This almost complex structure is $\omega$--compatible.
This is to say that the bilinear form
\begin{equation}
g^{-1}\omega(\cdot, J(\cdot))
\end{equation}
defines a smooth metric on $\mr \times (S^1 \times S^2)$.  Infact, the
metric in (1.6) is the standard product metric, $ds^2 + dt^2 +
d\theta^2 +\sin^2\theta d\varphi^2$ but written in terms of $t$, $f$,
$\varphi$ and $h$ as
\begin{equation}
dt^2 + g^{-2} (df^2 + \sin^{-2}\theta dh^2) + \sin^2 \theta d\varphi^2
.
\end{equation}
Note that $J$ is not integrable.  By the way, $J$ sends the vector
field $\partial_s$ to a multiple of the Reeb vector field, $\hat{v} =
-g^{-1}[(1 - 3 \cos^2\theta) \partial_t + \sqrt{6} \cos\theta
\partial_\varphi]$, the unique vector field that contracts with
$\alpha$ to give 1 and is annihilated by $d\alpha$.  In addition, $J$
is invariant under translations of the coordinate $s$ on $\mr \times
(S^1 \times S^2)$.  Thus, $J$ is a standard almost complex structure
for the `symplectization' of the contact structure defined by
$\alpha$.

As remarked, the almost complex structure in (1.5) defines the notion
used here of a pseudoholomorphic subvariety.  A certain subset of the
latter, called here HWZ subvarieties, are of particular interest.
Here is the definition:
\bigskip

\noindent
{\bf Definition 1.1}\qua {An \emph{HWZ subvariety}, $C\subset\mr
\times (S^1 \times S^2)$, is a non-empty, closed subset with the
following properties:}

\begin{itemize}
\item The complement in $C$ of a countable, nowhere accumulating
subset is a two-dimensional submanifold whose tangent space is
$J$--invariant.

\item $\int_{C\cap K}\omega<\infty$ if $K\subset \mr \times (S^1
\times S^2)$ is an open set with compact closure.
 
\item $\int_C d\alpha<\infty$.
\end{itemize}

The HWZ subvarieties are remarkably well behaved.  As explained in the
next section, each intersects the large $\pm s$ portions of $\mr\times
(S^1\times S^2)$ as a finite, disjoint union of cylinders.  Moreover,
each such cylinder intersects the appropriate component of each large
and constant $|s|$ slice of $\mr\times (S^1\times S^2)$ transversely
and the resulting $s$--parameterized family of circles in $S^1\times
S^2$ converges pointwise to multiply cover an embedded circle whose
tangent lines are annihilated by $d\alpha$.

As indicated by the preceding remarks, the closed, integral curves in
$S^1\times S^2$ of the kernel of $d\alpha$ play a prominent role in
this story.  They are called `closed Reeb orbits'.  Here is the full
list of such circles:

\begin{itemize}

\item There are two distinguished ones, labeled $(+)$ and $(-)$, these
being the respective loci where $\theta = 0$ and where $\theta=\pi$.

\item The others are labeled by data $((p, p'), \upsilon)$ where
$\upsilon\in\mr/2\pi\mz$ and where $(p, p')$ are integers subject to
three constraints:

\begin{itemize}

\item[(a)]At least one is non-zero, if $p = 0$, then $p' = \pm 1$, if
$p'= 0$, then $p = 1$ and if both are non-zero, then they are
relatively prime.

\item[(b)]$|p'|/|p| > \sqrt{3}/\sqrt{2}$ when $p < 0$.

\item[(c)] $p > 0$ when $|p'|/|p| < \sqrt{3}/\sqrt{2}$.
\end{itemize}
\end{itemize}
The Reeb orbit labeled by these data are the loci where

\begin{itemize} 
\item[(1)] $p' t - p\varphi = \upsilon$

\item[(2)] $\theta$ is constant and such that $p'(1 - 3 \cos^2 \theta)
=p \sqrt{6 }\cos \theta$ with\nl $p'\cos \theta\ge 0$.\autonum
\end{itemize}

As it turns out, the second constraint determines $\theta$ in terms of
$(p, p')$, and vice-versa.  What follows is a schematic drawing
of the possible values of $(p, p')$ as a function of the azimuthal
angle $\theta$ on the sphere.

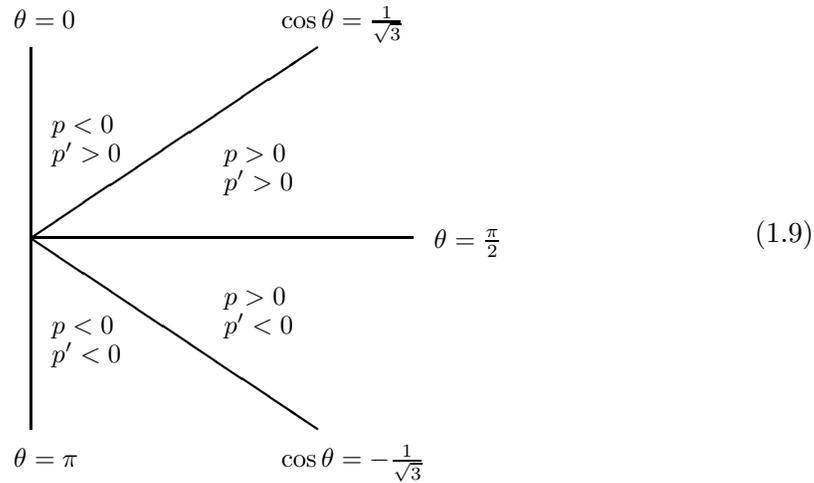
\vbox{\vglue 1.4in \cl{\small\unitlength 1 in
\begin{picture}(3,0) \thicklines
\put(0,0){\line(3,2){1.5}} \put(0,0){\line(3,-2){1.5}}
\put(0,0){\line(0,1){1}} \put(0,0){\line(0,-1){1}}
\put(0,0){\line(1,0){2}} \put(-.1,1.1){$\theta=0$}
\put(-.1,-1.2){$\theta=\pi$} \put(1.3,1.1){${\rm
cos}\,\theta=\frac{1}{\sqrt3}$} \put(1.3,-1.2){${\rm
cos}\,\theta=-\frac{1}{\sqrt3}$}
\put(2.1,-.05){$\theta=\frac{\pi}{2}$} \put(.1,.55){$p<0$}
\put(.1,.4){$p'>0$} \put(.1,-.5){$p<0$} \put(.1,-.65){$p'<0$}
\put(1,.4){$p>0$} \put(1,.25){$p'>0$} \put(1,-.35){$p>0$}
\put(1,-.5){$p'<0$}
\end{picture}}
\vglue-.45in \hbox{}\autonum \vglue 1.3in}

A momentary digression is in order here to comment on the fact that
the closed Reeb orbits are not isolated in $S^1\times S^2$.  In
particular, this is a direct consequence of the fact that the contact
form $\alpha$ is invariant under an $S^1\times S^1$ subgroup, $T$, of
the isometry group of $S^1\times S^2$.  The convention used here takes
the first factor of $S^1$ in $T$ to rotate the $S^1$ factor in
$S^1\times S^2$ via translation of the coordinate $t$; the second
factor of $S^1$ in $T$ rotates the $S^2$ factor of $S^1\times S^2$ via
translation of the spherical angle $\varphi$.  This understood, the
$\theta=0$ and $\pi$ Reeb orbits are the only closed Reeb orbits that
are invariant under the whole of $T$.  Each of the others is preserved
by no more than the product of a finite subgroup with a 1--parameter
subgroup of $T$.

As indicated above, each end of an HWZ subvariety $C$ determines a
closed Reeb orbit by the asymptotics of its constant $s$ slices as
either $s\rightarrow\infty$ or $s\rightarrow-\infty$.  This said,
introduce the number, $\fX_C$, of such convex side
$(s\rightarrow-\infty)$ ends with limit Reeb orbit where $\theta$ is
neither 0 nor $\pi$.  This integer $\fX_C$ plays a key role in the
subsequent discussion.

A second integer, $I_C$, also plays a key role here.  What follows is
an informal definition of $I_C$; the somewhat technical formalities
are relegated to Section 3.  The discussion starts with the
observation that the full set of HWZ subvarieties has a reasonable
structure of its own.  As explained later in Section 3, this set
enjoys a topology whereby a subvariety, $C$, has a neighborhood that
is homeomorphic to the zero set of a smooth map from a Euclidean ball
of some dimension $N_C$ to Euclidean space of a possibly different
dimension, $n_C$.  In this regard, the difference, $I_C=N_C-n_C$
defines a locally constant function and should be thought of as formal
dimension at $C$ of the space of HWZ subvarieties.  For future
reference, note that the \noindent assignment of the integer $\fX_C$
to an HWZ subvariety $C$ also defines a locally constant function in
this topology.

A second integer, $I_C$, also plays a key role here.  What follows is
an informal definition of $I_C$; the somewhat technical formalities
are relegated to Section 3.  The discussion starts with the
observation that the full set of HWZ subvarieties has a reasonable
structure of its own.  As explained later in Section 3, this set
enjoys a topology whereby a subvariety, $C$, has a neighborhood that
is homeomorphic to the zero set of a smooth map from a Euclidean ball
of some dimension $N_C$ to Euclidean space of a possibly different
dimension, $n_C$.  In this regard, the difference, $I_C=N_C-n_C$
defines a locally constant function and should be thought of as formal
dimension at $C$ of the space of HWZ subvarieties.  For future
reference, note that the \noindent assignment of the integer $\fX_C$
to an HWZ subvariety $C$ also defines a locally constant function in
this topology.

By the way, the interpretation of $I_C$ as a dimension is based on the
following observation: When $n_C$ can be taken to be zero, then $n_C$
can be taken zero on a neighborhood of $C$.  This said, the subset of
those HWZ subvarieties $C$ with $n_C = 0$ forms an open subset with a
natural smooth manifold structure and $I_C$ gives the dimension of
this manifold at $C$.

The set of HWZ subvarieties with the topology just described is called
the ``moduli space'' of HWZ subvarieties and is denoted by $\fM$.  Of
interest in this article are the components of $\fM$ where
$I_C\le\fX_C+1$.  In particular, the most prominent result of this
article is a complete description of these $I_C\le\fX_C+1$ components.
Theorems A.1--4 summarize many of the salient conclusions.  The author
currently plans to discuss some of the $I_C>\fX_C+1$ moduli spaces in
a sequel to this article.

The focus here on the $I_C\le\fX_C+1$ components of $\cM$ has its
ultimate justification in the proposed use of the constructs of Hofer
and his coworkers to study invariants of smooth 4--manifolds using
two-forms that are symplectic where non-zero.  However, as such
applications are dreams for now (and $I_C\le\fX_C+1$ cases surely
enter such dreams) take the results of the theorems here as simply
answers to the following question: What do the small dimensional, HWZ
moduli spaces look like?

The statements of Theorems A.1, A.2 and A.3 make explicit reference to
the fact that the group $T$ acts on $\fM$ as does the group $\mr$.
The former via its action on the $S^1\times S^2$ factor in
$\mr\times(S^1\times S^2)$ and the latter by its action on the $\mr$
factor via translation of the coordinate $s$.  As these group actions
commute with each other, so $\mr\times T$ acts on $\fM$.  Moreover,
this action is continuous and smooth on the smooth manifold parts of
$\fM$.
  
Theorems A.1, A.2 and A.3 refer to `irreducible' subvarieties.  The
term here denotes a subvariety that cannot be disconnected by the
removal of any finite set of points.  These theorems also use as
notation $\fM_C$ to denote the component in $\fM$ of a particular HWZ
subvariety $C$.
\bigskip

{\bf Theorem A.1}\qua\sl Let $C\subset X$ be an irreducible, HWZ
pseudoholomorphic subvariety.  Then the following are true:

\begin{itemize}
\item $I_C\ge\fX_C$.

\item If $I_C=\fX_C$, then $\fX_C= 0, 1$, or 2.
 
\begin{itemize}
\item[\rm(a)] $\fX_C=0$ if and only if $C=\mr\times\gamma$ where
$\gamma$ is a $\theta\in\{0,\pi\}$ closed Reeb orbit from the first
point in (1.8).  In this case, $\fM_C\subset\fM$ consists of the point
$C$.

\item[\rm(b)]$\fX_C=1$ if and only if $C = \mr\times\gamma$ where
$\gamma$ is a closed Reeb orbit from the second point in (1.8).  In
this case, $\fM_C$ is a smooth manifold, diffeomorphic to $S^1$,
invariant under the action of $\mr$ and a single orbit under the
action of $T$.

\item[\rm(c)] $\fX_C=2$ if and only if $C$ is a cylinder that is
invariant under a 1--parameter subgroup of $T$.  Here, $C$ is embedded
and the limit Reeb orbits from $C$ are characterized as in (1.8) by
data $(\pm(p,p'),\upsilon)$ where $|p'|/|p| > \sqrt{3}/\sqrt{2}$.  In
this case, $\fM_C$ is a smooth manifold, one orbit under $\mr\times T$
and diffeomorphic to $\mr\times S^1$.
\end{itemize}

\item If $I_C=\fX_C+1$, then $\fX_C=1, 2$ or 3.

\begin{itemize}

\item[\rm(a)] $\fX_C=1$ if and only if one of the following three
scenarios prevail:

\begin{itemize}
\item[\rm(1)]$C$ is an embedded disk invariant under the second factor
of $S^1$ in $T$.  Here, $C$ is embedded and its limit Reeb orbit is a
$((0, \pm 1), \upsilon)$ case from (1.8).  In this case, $\fM_C$ is a
smooth manifold, one orbit under $\mr\times T$, and diffeomorphic to
$\mr\times S^1$.

\item[\rm(2)]$C$ is an embedded cylinder, invariant under an $S^1$
subgroup of $T$ with only one convex side limiting Reeb orbit. Here,
$C$ is embedded and the convex side Reeb orbit is characterized by
$((p, p'),\upsilon)$ with $p < 0$ and $|p'|$ is the least integer that
is greater than $p \sqrt{3}/\sqrt{2}$.  Meanwhile, the concave side
Reeb orbit is characterized, as in (1.8), by either $(+)$ or $(-)$
with the sign in question that of $p'$.  In this case, $\fM_C$ is a
smooth manifold, one orbit under $\mr\times T$, and diffeomorphic to
$\mr\times S^1$.

\item[\rm(3)]$C$ is an embedded cylinder, invariant under an $S^1$
subgroup of $T$ with two convex side limiting Reeb orbits.  Here, one
is characterized by $((p, p'), \nu)$ with $p > 0$, and with $ |p'|$
equal to the greatest integer that is less than $p \sqrt{3}/\sqrt{2}$.
Meanwhile, the second convex side limit Reeb orbit is characterized,
as in (1.8), by either $(+)$ or $(-)$ with the sign in question
opposite that of $p'$.  In this case, $C$ is also embedded.  Again,
$\fM_C$ is a smooth manifold, one orbit under $\mr\times T$, and
diffeomorphic to $\mr\times S^1$.
\end{itemize}  
 
\item[\rm(b)] If $\fX_C = 2$ or 3, then $C$ is an immersed, thrice
punctured sphere with no limit Reeb orbits described, as in (1.8), by
either $(+)$ or $(-)$.  In each case, $\fM_C$ is a smooth manifold.
Theorems A.2 and A.3, below, describe the classification and structure
of $\fM_C$.

\end{itemize}      
\item All other cases have $I_C\ge\fX_C+2$.
\end{itemize}\rm
\medskip

Note that Theorem A.4, below, describes the number of singular points
in the immersed subvarieties that appear Part b of the third point.

Theorem A.1 makes reference to HWZ subvarieties that are preserved by
subgroups of $\mr\times T$.  In this regard, an irreducible subvariety
is preserved by an $\mr\times S^1$ subgroup of $\mr\times T$ if and
only if it is the cylinder $\mr\times\gamma$ where $\gamma$ is a
closed Reeb orbit.  Such a cylinder is preserved by the whole of
$\mr\times T$ if and only if $\gamma$ is either the $(+)$ or $(-)$
Reeb orbit in the notation of (1.8).  Meanwhile, an irreducible HWZ
subvariety $C$ that is moved by the $\mr$ action but preserved by an
$S^1$ subgroup of $T$ is either a cylinder or a disk.  All such are
described in great detail in Section 4 of this paper.  Here is a table
that summarizes the sorts of subvarieties that appear in Theorem A.1:

\narrower \itemm $\mr\times T$ invariant cylinders: The $\theta = 0$
and $\theta=\pi$ loci.

\itemm $\mr\times S^1$ invariant cylinders: $\mr\;\times$~(Closed Reeb
orbit).

\itemm $S^1$ invariant disks: These have $t=$~constant.  The one end
is on the convex side with constant $s $ slices that converge as
$s\rightarrow-\infty$ onto a closed Reeb orbit where $\cos^2\theta =
1/3$.

\itemm $S^1$ invariant cylinders: These can have either one or two
convex side ends.  In any event, each $\theta\not\in\{0,\pi\}$, closed
Reeb orbit determines exactly two $S^1$ invariant cylinders with
constant $s$ slices that converge to it as $s\rightarrow-\infty$.

\itemm Three-holed spheres with one concave side and two convex side
ends.

\itemm Three-holed spheres with no concave side and three convex side
ends.  \autonumm \endnarrower

No thrice-punctured sphere from Theorem A.1 is fixed by a non-trivial
element in $T$.  Nonetheless, much is known about these punctured
spheres and the moduli space components that contain them.  The next
theorem summarizes what is known about the $\fX_C=2$ cases.

\bigskip
\noindent
{\bf Theorem A.2}\qua\sl The $\fX_C=2$ components of $\fM$ that appear
in Part b of the third point of Theorem A.1 consist solely of
immersed, thriced punctured spheres.  Each such component is a smooth
manifold.  Moreover, these components have the following
classification and structure:

\begin{description}
\item[Classification] The components are classified by ordered sets of
four integers having the form $((p, p'), (q, q'))$ subject to the
following constraints:

\begin{itemize}
\item[\rm(a)] $\Delta\equiv pq'-p'q>0$.

\item[\rm(b)] $q'-p'>0$ {unless} $p'q'>0$

\item[\rm(c)] If $(m, m')$ denotes either $(p, p')$ or $(q, q')$, then
$|m'|/|m| > \sqrt{3}/\sqrt{2}$ when $m < 0$, and $m > 0$ when
$|m'|/|m| < \sqrt{3}/\sqrt{2}$.
\end{itemize}

\item[Structure] The component of $\fM$ that corresponds to $((p,
p')$, $(q, q'))$ is diffeomorphic to $\mr\times T$.  Moreover, this
diffeomorphism is $\mr\times T$ equivariant.
\end{description}
\rm\bigskip

The story on the $\fX_C = 3$ cases from Theorem A.1 is provided by the
next result.

\bigskip
\noindent
{\bf Theorem A.3}\qua\sl The $\fX_C = 3$ components of $\fM$ that
appear in Part b of the third point of Theorem A.1 consist solely of
immersed, thrice punctured spheres.  Each such component is a smooth
manifold.  Moreover, these components have the following
classification and structure:

\begin{description}
\item[Classification] The components are in 1--1 correspondence with
the unordered sets of three pair of integers that are constrained in
the following way: Such a set, $L$, can be ordered as $\{(p, p'), (q,
q'), (k, k')\}$ with

\begin{itemize}
\item[\rm(a)] $p + q + k = 0$ { and} $p' + q' + k' = 0$.

\item[\rm(b)]$|k'/k| > \sqrt{3}/\sqrt{2}$.

\item[\rm(c)] $\{(p, p'), (q, q')\}$ obey the constraints in the first
point of Theorem A.2.  In this regard, a set $L$ with an ordering that
satisfies these three conditions has precisely two distinct orderings
that satisfy the conditions.
\end{itemize} 

\item[Structure] The component of $\fM$ that is labeled by $L$ is a
smooth manifold that is diffeomorphic to $(0, 1) \times\mr\times T$.
Moreover, this diffeomorphism is $\mr\times T$ equivariant for the
$\mr\times T$ action on $(0, 1) \times\mr\times T$ that fixes the $(0,
1)$ factor and acts in the canonical fashion on the $\mr\times T$
factor.  Finally, the quotient, $(0, 1)$, of this moduli space
component by the $\mr\times T$ action has a natural compactification
as $[0, 1]$, where the two added points label the $\mr\times T$
quotient of two components of the $\fX = 2$ moduli space components
described in Theorem A.2.  In this regard, the relevant components are
labeled by the first two pairs from the two possible orderings of $L$
that obey the three constraints given in the preceding point.
\end{description}\rm
\bigskip

The final theorem in this section describes the number of singular
points of the thrice punctured spheres described in Theorems A.2 and
A.3.  For a subvariety $C$, this number, $m_C$, `counts' the number of
singular points in the immersion.  A precise definition of this count
is given in Section 3a.  In any event, $m_C = 0$ if and only if $C$ is
embedded, and $m_C$ is the number of double points when all
singularities are locally transversal intersections of pairs of disks.
  
The statement of Theorem A.4 implicitly views $S^1$ as the unit radius
circle about the origin in $\mc$.
\bigskip

\noindent
{\bf Theorem A.4}\qua\sl Suppose that $C$ is either described by
Theorem A.2 and its moduli space component is classified by the data
$\{(p, p'), (q, q')\}$, or else $C$ is described by Theorem A.3 and
its moduli space component is described by the data set $\{(p, p'),
(q, q'), (k, k')\}$.  In either case, the integer $m_C$ is one half of
the number of pairs $(\eta,\eta') \in S^1\times S^1$ such that
$\eta\not=\eta'$, neither $\eta$ nor $\eta'$ equals 1, and
$\eta^p\eta^q= \eta^{p'} \eta^{q'}=1$.  Thus,
$$
2 m_C = \Delta - gcd(p, p') - gcd(q, q') - gcd(p + q, p'+ q') + 2
$$
where $gcd(m, m')$ denotes the greatest common divisor of $m$ and
$m'$. For example, $m_C = 0$ and so $C$ is embedded if and only if one
of the following conditions holds:

\begin{itemize}
\item $|p q' - p'q| $ is either 1 or 2.

\item $|p q' - p'q| $ divides with integer remainder both members of
at least one of the pairs of integers $(p, p'), (q, q')$ and $(k = - p
- q, k' = - p'- q')$.
\end{itemize}\rm\bigskip

The remainder of this article has five more sections that are
organized along the following lines:
  
Section 2 states and proves a theorem that describes the behavior of a
natural class of pseudoholomorphic subvarieties on a non-compact
symplectic 4--manifold whose ends are symplectomorphic to either the $s
> 0$ or $s < 0$ portions of $\mr\times(S^1\times S^2)$.  In this
regard, the symplectomorphism is required to identify the almost
complex structure with that depicted in (1.5).  Propositions 2.2 and
2.3 summarize the principle results of Section 2.

Section 3 considers the structure of the moduli spaces of the
subvarieties from Section 2.  In particular, this section defines the
topology for $\fM$ and provides, in Proposition 3.2, a local model for
neighborhoods of points in $\fM$.  In addition, Proposition 3.6
provides an `index theorem' that computes the analog of $I_C$ in
explicitly geometric terms.

Section 4 focuses on the explicit example of $\mr\times(S^1\times
S^2)$ and provides a proof of Theorem A.1.  By the way, the proof of
Theorem A.1 derives additional constraints on the possibilities for
$I_C$.  The latter are summarized in Proposition 4.3.

Section 5 focuses on $\mr\times(S^1\times S^2)$ and proves Theorem
A.2.  This section also proves the assertions of Theorem A.4 about
Theorem A.2's subvarieties.

Section 6 focuses on $\mr\times(S^1\times S^2)$ and contains the proof
of Theorem A.3.  This section also contains the proofs of the
assertions in Theorem A.4 about the subvarieties from Theorem A.3.

\section{Regularity}

The discussion here and in the third section concerns an oriented,
symplectic 4--manifold $X$ that can be described in the following way:
Start with a smooth, oriented 4--manifold with boundary, $X_0$, where
the boundary of $X_0$ is a disjoint union of copies of $S^1 \times
S^2$.  Suppose that this boundary can be written as $\partial^-X_0
\cup \partial^+X_0$, where each component of $\partial^-X_0$ has a
neighborhood with an orientation preserving diffeomorphism to $[0, 1)
\times (S^1 \times S^2)$, and each component of $\partial^+X_0$ has
one to $(-1, 0] \times (S^1 \times S^2)$.  In this regard, view the
latter as subsets of $\mr \times (S^1 \times S^2)$ with its
orientation defined by the symplectic form in (1.2).

Given $X_0$, then $X$ is obtained by attaching $(-\infty, 0] \times
(S^1 \times S^2)$ to each component of $\partial^-X_0$ and attaching
$[0, \infty) \times (S^1 \times S^2)$ to each component of $\partial
^+X_0$.  Meanwhile, the symplectic form, $\omega$, on $X$ is required
to restrict to an open neighborhood of each component of $X-X_0$ as
either the form in (1.2) or else as this form after passing to a
suitable 2--fold cover.  In this regard, the deck transformation for
this cover is the fixed point-free involution, $\sigma$: $S^1 \times\
S^2 \rightarrow S^1 \times S^2$ that sends
\begin{equation}
(t, \theta , \varphi)\overset{\sigma}{\rightarrow} (t + \pi,\pi -
\theta, -\varphi)
\end{equation}
Note that $\alpha$ in (1.1) is invariant under the action of $\sigma$,
and thus $\omega$ is invariant under the induced involution on $\mr
\times (S^1 \times S^2)$.  Therefore, both descend to the associated
quotient.

Here is some terminology used below: An end of $X$ is a component of
$X-X_0$.  A component of $X-X_0$ that comes from a component of
$\partial^-X_0$ is called a \emph{convex end} of $X$, while one that
comes from a component of $\partial^+X_0$ is called a \emph{concave
end}.  Also, a component of $X-X_0$ where $\omega$ restricts directly
as the form in (1.2) is said to have \emph{orientable} $z$\emph{-axis
line bundle}.  A component where passage to the double cover is
required is said to have \emph{non-orientable} $z$\emph{-axis line
bundle}.

Of course, the prime example in this paper of such an $X$ is
$\mr\times (S^1 \times S^2)$ with the symplectic form in (1.2).  In
this case, there is one convex end and one concave end.  Moreover both
ends of $X$ have an orientable $z$--axis line bundle. As explained in
\cite{T2}, other examples come from compact 4--manifolds with 2--forms that
are symplectic where non-zero and vanish on an embedded union of
circles.  In the latter examples, all of the ends are concave, but
there can be some with a non-orientable $z$--axis line bundle.
 
By the way, when $X$ comes from a compact 4--manifold as just
described, a result of Gompf \cite{Go} asserts that the parity of the
number of ends of $X$ with orientable $z$--axis line bundle is opposite
that of the sum of the first Betti number, the second Betti number and
the signature of the original compact 4--manifold.  On the other hand,
start with such a compact 4--manifold and, according to Luttinger \cite{Lu},
the closed form can be manipulated so that the resulting manifold $X$
has only orientable $z$--axis ends.
  
Given $X$ with its symplectic form as just described, there are almost
complex structures on $X$ that are $\omega$--compatible and restrict to
some open neighborhood of each end of $X$ as follows: If the end has
orientable $z$--axis line bundle, then $J$ restricts as the almost
complex structure in (1.5).  On the other hand, if the end has
non-orientable $z$--axis line bundle, then $J$ should restrict as the
push-forward of (1.5) via the covering map induced by $\sigma$.  In
this regard, note that (1.5) is $\sigma$--invariant.

Unless explicitly noted otherwise, assume that every almost complex
structure that appears below has the properties just described. This
said, fix such a $J$ and Definition 1.1 has the following analog:
\bigskip

{\bf Definition 2.1}\qua With $X$, $\omega$ and $J$ as just described,
a subset $C \subset X$ is an \emph{HWZ subvariety} when the following
conditions are met:

\begin{itemize}
\item $C$ is closed, and the complement of a countable,
non-accumulating set is a smooth submanifold with a $J$--invariant
tangent space.

\item Let $K \subset X$ be any open set with compact closure.  Then
$\int_{C\cap E}\;\omega<\infty$.

\item Let $E \subset X-X_0$ be any component.  Then $\int_{C\cap
E}\;d\alpha<\infty$ where $\alpha$ is the contact form in (1.1) when
$E$ has orientable $z$--axis line bundle; otherwise, $\alpha$ is the
push-forward of the form in (1.1) via the covering map defined by
$\sigma$ in (2.1).\end{itemize}
\bigskip

With the preceding understood, it can now be said that the purpose of
this section is to state and then prove the two propositions that
follow that describe the ends of an HWZ subvariety in $X$.  With
regards to the proofs, note that they introduce various constructions
that are used in later portions of this article.
\bigskip

\noindent
{\bf Proposition 2.2}\qua\sl Let $C \subset X$ be an {\rm HWZ}
subvariety.  Then:
\begin{itemize}  
\item $C$ has a finite number of singular points.
\item $C$ intersects each sufficiently large and constant $|s|$ slice
of $X-X_0$ transversely.
\item There is a finite union, $\Gamma$, of disjoint closed Reeb
orbits in $\partial X_0$, thus integral curves of the distribution
kernel $(d\alpha)$, with the following significance:
\begin{itemize}
\item[\rm(a)]Let $U\subset\partial X_0$ be any tubular neighborhood of
$\Gamma$.  Then $C$'s intersection with each sufficiently large and
constant $|s|$ slice of $X-X_0$ lies in $U$.
\item[\rm(b)]Fix a tubular neighborhood projection from $U$ to
$\Gamma$ and then $C$'s intersection with each sufficiently large and
constant $|s|$ slice of $X-X_0$ projects to each component of $\Gamma$
as a finite to one covering map.
\end{itemize}
\item There exists a complex curve $C_0$ with a finite set of
cylindrical ends together with a proper, pseudoholomorphic map into
$X$ whose image is $C$ and which embeds the complement of a finite
set.
\end{itemize}\rm
\medskip

The curve $C_0$ will be called `the model curve' for $C$.  The set
$\Gamma$ will be called the `limit set' for $C$.  Note that the closed
Reeb orbits are listed in (1.8) for the components of $\partial X_0$
with orientable $z$--axis line bundle.  The closed Reeb orbits in the
unorientable $z$--axis line bundle components of $\partial X_0$ are the
images of the orbits listed in (1.8) under the 2--1 covering map
induced by the map $\sigma$ in (2.1).
 
It follows from Proposition 2.2 that the $|s| \rightarrow \infty$
limit of the constant $|s|$ slices of $C$ on $X-X_0$ converge to some
union of closed Reeb orbits.  The following proposition gives a more
detailed picture of this convergence:
\bigskip

\noindent
{\bf Proposition 2.3}\qua\sl Let $C \subset X$ be an {\rm HWZ}
subvariety.  There exists a finite union, $\Gamma\subset\partial X_0$,
of closed Reeb orbits, and, after a tubular neighborhood of $\Gamma$
is identified via an exponential map with $\Gamma\times D$ with
$D\subset \mr^2$ an open disk, there exists $s_0 > 0$ such that
\begin{itemize}
\item The intersection of $C$ with each constant $|s| \ge s_0 $ slice
of $X-X_0$ lies in $\Gamma\times D$.

\item Each component of the $|s| \ge s_0$ portion of $C$'s
intersection with $X-X_0$ can be parameterized by $[s_0, \infty)
\times S^1$ via a map which sends $(s, \tau)$ to $(\pm
s,\gamma(m_\gamma \tau), \eta(s, \tau))$ where $m_\gamma$ is a
positive integer, $\gamma\in\Gamma$ and $\eta$ is a smooth map from
$[s_0, \infty) \times S^1$ to $D$.  Here, the $+$ sign is used with a
component in a concave end of $X$ and the $-$ sign with a component in
a convex end.

\item There exists $\epsilon>0$ and, for each integer $k \ge 0$, a
constant $\zeta_k \ge 0$ such that the $C^k$--norm of $\eta$ is bounded
by $\zeta_k \;e^{-\epsilon |s|}$.
\end{itemize}\rm
\medskip

Note that these two propositions would follow directly from Theorems
1.2 and 1.4 in \cite{HWZ1} but for the fact that the latter assume a
non-degeneracy condition on the closed Reeb orbits that is not obeyed
here.

The remainder of this section is occupied with the proofs of
Propositions 2.2 and 2.3.

\sh{(a)\qua Proof of Proposition 2.2}
  
Before getting to specifics, note that the assertions of the
proposition are local to the ends of $X$ and so no generality is lost
by assuming that the ends have orientable $z$--axis line bundle. This
is because the almost complex structure in (1.5) is $\sigma$--invariant
and thus $\sigma$ preserves the conditions for the appellation HWZ
subvariety.  Thus, the setting for a non-orientable $z$--axis line
bundle end can be pulled up via the double cover map induced by
$\sigma$ in (2.1) and viewed as a $\sigma$--equivariant example of the
orientable $z$--axis line bundle case.  This understood, all ends in
the subsequent discussion are implicitly assumed to have orientable
$z$--axis line bundle.

To begin the proof, remark first that the first and fourth points
follow directly from the second and third.  Thus, the argument below
focuses on the latter two points.  This argument is broken into ten
steps.

{\bf Step 1}\qua For each integer $n\ge 2$, let $C_n$ denote the
intersection of $C$ with the portion of $X-X_0$ where $|s| \in [n - 2,
n + 2]$.  Via the evident identification of this cylinder with $W
\equiv [-2, 2] \times\partial X_0$, each $C_n$ can be viewed as a
proper, pseudoholomorphic subvariety of $W$.  In this regard, note
that $\omega=d(e^{-\sqrt{6}\;s}\alpha)$ is a symplectic form on $W$.
The first claim here is that there is an $n$--independent upper bound
to the symplectic area of $C_n$.  This follows from the finiteness of
$d\alpha$'s integral over $C \cap (X-X_0)$. Indeed the fact that
$d\alpha$ has finite integral has two key implications:

\narrower\sl
\itemm $\lim_{n\rightarrow\infty}\int_{C_n}\;d\alpha=0$.

\itemm Fix a component of $\partial X_0$ and then the sequence of
numbers (indexed by $n$) obtained by integrating $\alpha$ over the
intersection of $C$ with the $|s| = n$ slice of the corresponding
component of $X-X_0$ is convergent.\autonumm
\endnarrower

Note that the second point implies the assertion about the upper bound
for the symplectic area of $C_n$.
\bigskip

{\bf Step 2}\qua Use the compactness theorem of Proposition 3.3 in
\cite{T4} to conclude that any subsequence of $\{C_n\}$ has inside it, a
subsequence (hence renumbered consecutively from 1) which converges to
a proper, pseudoholomorphic subvariety $C'\subset W$.  This
convergence is in the following sense: 

First, the sequence
\begin{equation}
\left\{\text{sup}_{x\in C'} \text{dist}(x, C_n) + \text{sup}_{x\in
C_n} \text{dist}(x, C')\right\}
\end{equation}
converges with limit zero.  

Second, $C'$ can be written as a finite
union $\{C'_\alpha\}_{1\le\alpha\le N}$ for some integer $N\ge 1$,
where each $C'_\alpha$ is a proper, pseudo-holomorphic subvariety and
the intersection of $C'_\alpha$ with $C'_{\alpha'}$ is finite when
$\alpha\not=\alpha'$.  Moreover, there is a corresponding sequence of
positve integers $\{m_\alpha\}$ such that the set pairs $c
\equiv\{(C'_\alpha, m_\alpha)\}$ has the property that for any 2--form
$\mu$ on $W$, the sequence
\begin{equation}
\left\{\int_{C_n}\;\mu\right\}
\end{equation}
converges with limit
\begin{equation}
\sum_{(C',m) \in c} m \int_{C'}\;\mu\;.
\end{equation}

{\bf Step 3}\qua It follows from the preceding step that $d\alpha$
vanishes on $C'$, and this implies that $C'$ has the form $[-2, 2]
\times\Gamma$, where $\Gamma$ is a finite union of closed Reeb orbits.
These closed Reeb orbits are listed in (1.8).  Were they are isolated,
it would follow from the discussion in Step 1 that the data $c =
\left\{([-2, 2]\times\gamma, m)\right. : \gamma\in\Gamma$ and
$\left.\left. m \in\{1, 2,\cdots\right\}\right\}$ describing the limit
of $\{C_n\}$ is uniquely defined from $C$.  However, in the present
case, the closed Reeb orbits are not isolated.  Even so, one can still
draw this same conclusion:
\bigskip
 
\noindent
{\bf Lemma 2.4}\qua {\sl Let $C \subset X$ be an {\rm HWZ} subvariety.
Then all limits of the corresponding sequence $\{C_n\}$ produce the
same data set $c$ as limit.  In particular, there is a finite union,
$\Gamma$, of closed Reeb orbits with the following significance: Given
$\epsilon > 0$, there exists $s_\epsilon$ such that each point of
$C'$s intersection with any constant $|s| \ge s_\epsilon$ portion of
$X-X_0$ has distance $\epsilon$ or less from $\Gamma$.  Conversely,
each point of $\Gamma$ has distance $\epsilon$ or less from a point in
$C'$s intersection with any constant $|s| \ge s_\epsilon$ portion of
$X-X_0$.}
\bigskip

This lemma is proved below in the final step of the proof of
Proposition 2.3, so accept it for the time being.

Before continuing to Step 4, note that the remainder of the
proof of Proposition 2.2 (including that of Lemma 2.4) will assume
that the component of $\partial X_0$ in question is concave.  In this
regard, the arguments for the convex case are identical save for some
notational changes.
\bigskip

{\bf Step 4}\qua To begin, re-introduce the coordinate functions $(t, f,
h, \varphi)$ for $\mr \times (S^1 \times S^2)$ as defined in (1.3).
In this regard, note that the submanifolds where $(t, f)$ are constant
are pseudoholomorphic cylinders or pairs of disks (when $f = 0$),
while those where $(\varphi, h)$ are constant are pseudoholomorphic
cylinders.

Fix attention on a component, $\gamma$ of $\Gamma$ and remember from
(1.8) that $\gamma$ can be labeled either $(+)$, $(-)$ or $((p, p'),
\nu)$.  In all cases, the spherical angle $\theta$ is a constant,
$\theta_0$, on $\gamma$.
 
The value of $\theta_0$ forces two cases to be distinguished: Case 1
has $\gamma\not= ((0, \pm1),\nu)$, which is equivalent to the
condition that $\cos^2 \theta_0 \not= 1/3$.  Meanwhile, Case 2 has
$\gamma=((0, \pm 1), \nu)$ and therefore $\cos^2 \theta_0 = 1/3$.  The
discussion below focus on the case where $\cos^2 \theta_0 \not= 1/3$
and this condition will be assumed implicitly.  The discussion for the
$cos^2 \theta_0 = 1/3$ case is essentially identical to that below
after switching the roles which are played by the coordinates $t$ and
$\varphi$ and also by $f$ and $h$.  For this reasons, the discussion
in the $\cos^2\theta_0 = 1/3$ case is omitted.
	
With the preceding understood, return to the loop $\gamma$, and
observe that in the case where $\theta_0\notin\{0,\pi\}$, this loop
has a parameterization by a periodic variable
$\tau\in\mr/(2\pi|p|\mz)$ via

\begin{itemize}
\item $t=\tau$,
\item $\theta=\theta_0$,
\item $\varphi=\varphi_0+\tau p'/p$.\autonum
\end{itemize}

Here, $\varphi_0\in\mr/(2\pi\mz)$ is a constant.  In the case where
$\theta_0\in\{0,\pi\}$, the parameterization of $\gamma$ is also given
by (2.6) with the last two lines absent.  Even so, take $p = -1$ when
$\theta_0\in\{0,\pi\}$.  In any case, the projection from $S^1 \times
S^2$ to the $S^1$ factor restricts to $\gamma$ as an $|p|$ to 1
covering map.
 
Note that the functions $f$ and $h$ in (1.3) restricts to $\gamma$ to
obey
\begin{equation}
h = (p'/p) \sin^2(\theta_0) f.
%(2.7)
\end{equation}

{\bf Step 5}\qua Lemma 2.4 implies that at all large values of $s$,
the intersection of $C$ with some fixed radius tubular neighborhood of
$\{s\} \times\gamma$ lies very close to $\{s\} \times\gamma$.  In
particular, take the tubular neighborhood in question to have disjoint
closure from the other components of $\Gamma$.  With this noted, there
must exist $f_1$ such that $C$ has empty intersection with the
boundary of the closure of this tubular neighborhood where $|f| \le
f_1$.  In particular, there is an unambiguous component of the
intersection between $C$ and $\mr\times(S^1 \times S^2)$ where $|f|\le
f_1$ which lies in the given tubular neighborhood.  Given the
preceding, agree to restrict attention to the just mentioned $|f| \le
f_1$ portion of $C$.  By the way, note that on this portion of $C$,
the limits $|f|\rightarrow 0$ and $s \rightarrow\infty$ can be assumed
equal by constraining the tubular neighborhood of $\gamma$ so $|\ln((1
- 3 \cos^2 \theta)/(1 - 3 \cos^2 \theta_0))| < 1$ on it.

In the subsequent steps, this particular $|f| \le f_1$ portion of $C$
will still be denoted by $C$.
	
With these last remarks understood, it then follows from Lemma 2.4 and
(2.6) that $C$ intersects each constant $(t, f)$ pseudoholomorphic
subvariety in $m|p|$ points counting multiplicities.  Here, $m$ is the
weight which appears with $\gamma$ in the set c from Step 3.  Note
that all of the multiplicities here are positive.  Moreover, because
of (2.7), these intersection points are grouped in subsets of $m$
points (counting multiplicities), where each of the $|p|$ points in
the set
\begin{equation}
\left\{(\varphi=\varphi_0 + t + 2\pi k p'/p, h = (p'/p)
\sin^2(\theta_0) f) : 1 \le k \le |p|\right\}
\end{equation}
%(2.8)
is very close to a unique such subset.  Indeed, the following
assertion is a direct consequence of Lemma 2.4:
\medskip

\narrower
\sl Fix $\epsilon > 0$, and $\delta > 0$ exists such that if $|f| <
\delta$, $t \in\mr/(2\pi|p|\mz)$ and $\eta$ is a point in (2.8),
then the $m (\varphi,h)$ values near $\eta$ coming from points in
$C$ with coordinates $(t, f)$ have distance $\epsilon$ or less from
$\eta$.\autonumm
\endnarrower

Agree to use $\pi$ to denote the map which assigns the coordinates
$(t, f) \in S^1 \times \mr$ to the points in $\mr\times(S^1 \times
S^2)$.  The restriction of $\pi$ to $C$ has the following key
property: There exists a countable, non-accumulating set
$\Lambda\subset C$ such that $\pi$ restricts to $C-\Lambda$ as an $m
|p|$ to 1 covering map.  That is, if $x\in\pi(C-\Lam)$, then the
multiplicity of each point in $\pi^{-1}(x)$ is precisely 1.  On the
other hand, if $x \in\Lam$, then
\narrower
\sl $\pi^{-1}(\pi(x))$ contains less than
$ m|p|$ points; thus some of them have
multiplicty greater than 1.\autonumm
\endnarrower
%(2.10)

By the way, it is important to note that each point in $\Lam$ is
either a singular point of $C$ or a smooth point of $C$ but a critical
point of $\pi$.  (These last assertions and (2.10) all follow more or
less directly from the fact that the constant $(t, f)$ surfaces are
pseudoholomorphic.  A detailed argument can be had by mimicking,
almost verbatim, the discussion in Part a of the Appendix to \cite{T4}.)

With regards to $\Lam$, the key observation now is that if $\Lam$ is
finite, then any small, positive and constant $|f|$ slice of $C$
projects as finite to one covering map over the analogous slice of the
$(t, f)$ cylinder.  This said, Steps 6--8 of the proof demonstrate
that $\Lam$ is finite; the ninth step shows how the points of the
proposition then follow.

{\bf Step 6}\qua This step introduces the space $\pi^*C$ which
consists of the triples $(w, \eta_-,\eta _+)$ with $(w,\eta_\pm )\in
C$ and $\pi(w, \eta _\pm) = w$.  Here, $\eta_\pm$ consists of the
2--tuple $(\varphi_\pm, h_\pm)$.

The local structure of $\pi^*C$ can be analyzed by mimicking the
discussion in Part a of the Appendix to \cite{T4}.  What follows is a
summary of some of the important features.  First, $C$ embeds in
$\pi^*C$ as the diagonal where $\eta_- = \eta_+$.  In what follows,
$C$ and its image in $\pi^*C$ will not be notationally distinguished.
Second, $\pi^*C$ is a smooth manifold except near points $(w, \eta_-,
\eta_+)$ where one or both of $(w, \eta_\pm) \in\Lam$.  Meanwhile, the
structure of $\pi^*C$ near one of these singular points can be
described in detail as Part a of the Appendix in \cite{T4}.  In any event,
the projection $\pi$ mapping $(w,\eta_-,\eta_+)$ to $w$ restricts over
the complement of $\pi(\Lam)$ as an $(m |p|)^2$ to 1 covering map.
Third, let $\pi^*C«$ denote the closure in $\pi^*C $ of $\pi^*C-C$.
Then $\pi^*C' \cap C = \Lam\subset C$.  Indeed, this follows from
(2.10).

There is one more crucial point to make about $\pi^*C$, namely: If
$f_1 > 0$ is sufficiently small, then:

\narrower\sl
$\varphi_+-\varphi_-$ is a bonafide, real-valued function on $\pi^*C$;
in fact, given $\epsilon>0$, there exists $\delta > 0$ such that when
$|f| <\delta$, then $|\varphi_+-\varphi_-|\le\epsilon.$\autonumm
\endnarrower
%(2.11)

\noindent 
This assertion is an immediate consequence of (2.8) and (2.9).  Thus,
it is a corollary to Lemma 2.4, and this is essentially the only place
in the argument for Proposition 2.2 that requires Lemma 2.4.  However,
the conclusion expressed in (2.11) is absolutely crucial for the
subsequent arguments.

\medskip
{\bf Step 7}\qua Introduce the set $G\subset\pi^*C«$ which consists of
the points $(w, \eta_-, \eta_+)$ where the function $\underline{h
}\equiv h_+-h_-$ is zero.  This set is an example of an `embedded
graph'.  This is to say that $G$ is a locally compact subset with the
following additional properties: First, $G$ has a distinguished
subset, $G_v$, which is a locally finite collection of points.
Elements in $G_v$ are called vertices.  Meanwhile, $G - G_v$ is a
locally finite set of properly embedded, open intervals in $\pi^*C'- G
_v$.  The closure of each component of $G-G_v$ is called an edge.
Finally, each vertex has a neighborhood in $\pi^*C«$ whose
intersection with $G$ consists of a finite union of properly embedded
images of the half open interval $[0, 1)$ by an embedding which sends
0 to the vertex in question.  Moreover, these embedded intervals
intersect pairwise only at the given vertex.
 
In the present case, $G$ has some additional properties which are
summarized below:
\begin{itemize}\sl
\item {Each vertex of $G$ is either an $\underline{h} = 0$
critical point of $h$ in the smooth part of $\pi^*C«$ or else a point
$(w, \eta_+, \eta_-)\in \pi^*C$ where one or both of $(w, \eta_\pm)
\in\Lam$.}
\item {Let $\underline{\varphi}\equiv\varphi_+-\varphi_-$.  The
1--form $d\uvar$ is non-zero on the tangent space of $G-G_v$.  Infact,
at all points on the interior of each edge,the 2--form $d\uvar\wedge
d\uh$ orients $\pi^*C'$ so that $\pi$ is an orientation preserving
map.}
\item {The intersection of $G$ with some open neighborhood of
each vertex is a finite union of embedded, half open arcs with
endpoints lying on the vertex, but disjoint otherwise.  Moreover, the
tangent lines to the arcs at the vertex are well defined and
disjoint. The interior of each arc is part of an edge of the graph.
The number of such arcs is non-zero and even. Exactly half of the arcs
are oriented by $d\uvar$ so $\uvar$ increases towards the vertex while
half are oriented by $d\uvar$ so that $\uvar$ decreases towards the
vertex.}
\item {If $f_2 \in (0, f_1)$ is chosen to be sufficiently
generic, then the $|f| = f_2$ locus in $\pi^*C$ intersects only the
smooth part of $\pi^*C$, and where $df \not= 0$.  In addition,
$\Gamma$ intersects this $|f| = f_2$ submanifold of $\pi^*C$ as a
finite set of points, all in the interior of edges and this
intersection is transverse.}\nl
\hbox{}\autonum%(2.12)
\end{itemize}

These points are proved by copying the arguments in Steps 3--5 of Part
b in the Appendix to \cite{T4}.  In this regard, note that the observation
that $d\uvar$ is non-zero on $G-G_v$ (where $\uh = 0$) plays a
starring role in the next step.
 
By the way, the conclusion of the second point about $d\uvar$
orienting $G-G_v$ is a specific consequence of the fact that $J$ from
(1.5) maps $dh$ to a $\varphi$--independent multiple of $d\varphi$.
Indeed, because $C$ is pseudoholomorphic, the pair $(\varphi, h)$ obey
a Cauchy Riemann like equation on $C-\Lam$ as functions of the
variables $(t, f)$.  This equation can be written schematically as
$d\varphi = j\cdot dh$.  Here, $j$ is a function of $f$ and $h$ only;
in particular, $j$ is independent of $\varphi$.  One consequence of
this equation $d\varphi=j\;dh$, is that the function
$\uvar=\varphi_+-\varphi_-$ on $\pi^*C-\pi^{-1}(\pi(\Lam))$, when
viewed as a function of $(t, f)$, obeys
$d(\varphi_+-\varphi_-)=j(f,h_+)dh_+-j(f,h_-)dh_-$.  In particular,
where $h_+ = h_-$, this reads $d\uvar=j(f,h)d\uh$; and this last
equation directly implies the assertion in the second point of (2.12).
\bigskip

{\bf Step 8}\qua To complete the proof that $\Lam$ is finite (modulo
the proof of Lemma 2.4), note first that each point of $C \cap\pi^*C'$
is a vertex of $G$.  Moreover, these are precisely the vertices of $G$
where $\uvar=0$.  Thus, it is sufficient to show that there are at
most a finite number of $\uvar = 0$ vertices of $G$.  In fact, the
claim is that the number of $\uvar = 0$ vertices of $G$ where $|f| <
f_2$ is no greater than the number of points where $|f| = f_2$ on $G$.
This last number is finite by virtue of the final point in (2.12)).

To prove the preceding claim, choose a $\uvar = 0$ vertex of $G$ and
then follow some edge out from this vertex where $\uvar$ is
increasing.  The existence of such an edge is guaranteed by the third
point in (2.12).  Continue to travel along this edge.  According to
the second point of (2.12), the function $\uvar$ continues to
increase.  Either this edge eventually hits another vertex of $G$
where $|f|\le f_2$, or else it hits the $|f| = f_2$ locus of $G$.  In
this regard, note that $|f|$ is bounded away from zero on such an
edge, since $|\uvar|$ converges to zero along any path in $\pi^*C$
where $|f|$ limits to zero because of (2.11).  Indeed, remember that
$|\uvar|$ started at zero and then increases along the edge.

If the edge ends in a second vertex of $G$, then $\uvar> 0$ at this
vertex, and there is another edge coming into this vertex on which
$\uvar$ is increasing in the outward pointing direction (by the fourth
point in (2.12).)  Continue out on this new edge.  Note that $\uvar$
still is increasing.  Iterate this procedure.  As $\uvar$ always
increases, the path so traced out remains in a compact subset of the
$|f| \le f_2$ portion of $\pi^*C'$.  Meanwhile, no vertices are hit by
this piecewise smooth path in $G$ more than once.  By compactness, the
path must end, and the only possible way to do so is to hit the $|f| =
f_2$ locus.

With the preceding understood, let $G^{(1)}\subset G$ denote the
compliment of the interiors of the edges which are traversed by the
path just described.  Note that $G^{(1)}$ is also described by (2.12)
except that $G^{(1)}$ may have some isolated vertices.  Agree to
ignore these as they play no role in what follows.  By the way,
observe that the intersection of $G^{(1)}$ with the $|f| = f_2$ locus
contains one fewer point then that of $G$.

Given $G^{(1)}$, repeat the procedure just described in the previous
three paragraphs, but with $G^{(1)}$ replacing $G$.  The result is a
$G^{(2)}\subset G^{(1)}$ which is described by (2.12) (except maybe
for some isolated vertices).  Note that $G^{(2)}$ has two fewer
intersection points with the $|f| = f_2$ locus as did $G$.  Of course,
one can continue in this vein, creating $G^{(3)}, \ldots,$ etc., each
time reducing by one the size of the set of intersections with the
$|f| = f_2$ locus.  Eventually, this finite set of $|f| = f_2$
intersections is exhausted, say for $G^{(k)}$, in which case $G^{(k)}$
has no non-isolated $\uvar = 0$ vertices.  In particular, this means
that the original number of $\uvar = 0$ vertices in $G$ is no greater
than the size of the $|f| = f_2$ locus in $G$.

\medskip
{\bf Step 9}\qua Given that the set $\Lam$ is finite, it follows that
there exists some $f_1 > 0$ such that when $|f| < f_1$, then the
constant $|f|$ slice of $C$ defines an $m|p|$ to 1 covering map over
the corresponding slice of the $(t, f)$ cylinder.  This understood,
then $C$ can be parameterized as a multiple cover of the appropriate
component of the $|f| < f_1$ portion of the $(t, f)$ cylinder.  One
such parameterization uses coordinates $(\tau, \rho)$ where
$\tau\in\mr/(2\pi m|p|\mz)$ and $\rho\in (\rho_0, \infty)$.  In the
case where $\theta\notin \{0, \pi\}$, the latter parameterize the $|f|
< f_1$ portion of $C$ via the map that sends $(\tau, \rho)$ to
\begin{itemize}
\item $t=\tau$,
\item $f=\text{sign}(p)e^{-\kappa\rho}$,
\item $\varphi=\varphi_0+\tau p'/p+x(\tau,\rho)$,
\item
$h=e^{-\kappa\rho}\sin^2\theta_0(p'/p+\kappa^{-1}\text{sign}(p)w(\tau,\rho))$,
\autonum%(2.13)
\end{itemize}

where $\kappa = 6^{-1/2} (1 + 3 \cos^4 \theta_0)^{-1/2}|
1-3\cos^2\theta_0|$.  Meanwhile, as $\theta$ and $\varphi$ are not
good coordinates near the poles of $S^2$, the parameterization in case
when $\theta_0 \in \{0,\pi\}$ replaces the latter by the functions
\begin{equation}
a_1 \equiv 6^{-1/4 }|f|^{-1/2} |h|^{1/2} \cos(\varphi)\quad \text{
and}\quad a_2 \equiv 6^{-1/4} |f|^{-1/2} |h|^{1/2} \sin(\varphi) .
%(2.14)
\end{equation}
This done, the parameterization sends $(\tau,\rho)$ to
$(t=\tau$,$f=-e^{-\sqrt{6}\;\rho}$, $a_1(\tau,\rho)$,
$a_2(\tau,\rho))$.
  
Here is the point of such a parameterization: Introduce the
2--component, column vector $\lam$, with either top entry $x$ and
bottom entry $w$ or top entry $a_1$ and bottom entry $a_2$ as the case
may be.  By virtue of (1.5) and the fact that $C$ is
pseudoholomorphic, this vector obeys a differential equation with the
schematic form
\begin{equation}
\partial_\rho\lam+L_0\lam+\cR(\lam,\lam_\tau)=0.
\end{equation}
%(2.15)
Here, $L_0$ in (2.15) denotes the operator
\begin{equation}
L_0= \left(
\begin{matrix}
-\zeta'&-\partial_\tau\\ \partial_\tau&-\zeta
\end{matrix}
\right)
\end{equation}
%(2.16)
where $\zeta'$ and $\zeta$ are constants.  Meanwhile,
$\cR(\cdot,\cdot)$ is affine linear in the second factor and obeys
\begin{equation}
|\cR(a, b)|\le\xi (|a|^2 + |a| |b|)\; ,
\end{equation}
%(2.17)
where $\xi$ is independent of $\rho$, $a$ and $b$ where $\rho$ is
large and $|a|$ is small.  Equation (2.15) also exhibits the notation,
used subsequently, where the partial derivative of a function by a
parameter is denoted by the function's symbol adorned with the
parameter as subscript.

As a solution to (2.15), the vector $\lam$ also obeys the asymptotic
condition
\begin{equation}
\lim_{\rho\rightarrow\infty}|\lam|= 0
\end{equation}
%(2.18)
by virtue of the fact that the constant and large $\rho$ slices of $C$
converge to the Reeb orbit in question.

Observe now that the equation in (2.15) is, for small $|\lam|$,
uniformly elliptic.  This noted, and given (2.18), standard elliptic
regularity arguments as in Chapter 6 of \cite{Mo} apply to (2.15) and find
that derivatives of $\lam$ to all orders converge to zero as $\rho$
tends to infinity.  For future reference, this conclusion is stated
formally as the following lemma.
\bigskip

\noindent
{\bf Lemma 2.5}\qua {\sl Given a non-negative integer $k$ and
$\epsilon> 0$, there exists $\rho_\epsilon$ such that
$|\nabla^{\otimes k}\lam|<\epsilon$ at all points with
$\rho>\rho_\epsilon$.}
\bigskip

This observation of Lemma 2.5 implies the assertions of Proposition
2.2.  For example, the transversal intersection of $C$ with the large
and constant $|s|$ slices follows after first pulling the differential
$ds$ back to $C$ while noting that $ds = -\sqrt{6} (f \;df + h \;dh)$.
This done, write $h$ in terms of the components of $\lam$ and the
$\rho\rightarrow\infty$ limit of zero for $|\nabla\lam|$ implies that
$\partial_\rho s\not= 0$ at large $\rho$ on $C$.  The assertions of
the third point of the proposition follow in a similar vein.

\bigskip

{\bf Step 10}\qua This step is devoted to the following proof.
\bigskip

\noindent
{\bf Proof of Lemma 2.4}\qua A digression comes first to summarize
certain observations made prior to the statement of the lemma.  To
begin the digression, recall that the possible limits are described by
sets of the form $\{(m_\gamma, \mr\times\gamma): \gamma\in\Gamma\}$,
where $m_\gamma$ is a positive integer weight, and where $\Gamma$ is a
finite set of distinct closed Reeb orbits in $S^1 \times S^2$.
Moreover, the function $\theta$ on $S^2$ has some constant value,
$\theta_0$, on each circle.  In this regard, the possible $\theta_0$
values which can appear for some $\gamma\in\Gamma$, and the
corresponding weight $m_\gamma$ are uniquely determined as the set of
possibilities for these data is discrete.  In fact, these parameters
can be distinguished by intersection numbers of $C$ with various
pseudoholomorphic submanifolds of $\mr\times(S^1 \times S^2)$.  This
is to say that there might be different sets $\Gamma$ which appear as
limits, but each such set has the same collection of $\theta_0$
values, and the multiplicities $m_\gamma$ for the elements
$\gamma\in\Gamma$ depend only on these $\theta_0$ values and so are
the same for all of the possible sets $\Gamma$ which could appear.

With regard to these $\theta_0$ values, the argument for Lemma 2.4
given here are valid for the elements in $\Gamma$ with $\cos^2\theta_0
\not= 1/3$.  The argument for the elements with $\cos^2 \theta_0 =
1/3$ are left to the reader in as much as they are essentially the
same as those given below after switching the roles of $(\varphi, h)$
with those of $(t, f)$.

With the digression now complete, suppose that $\Gamma$ is as above
and $\gamma\in\Gamma$ has $\cos^2 \theta_0 \not= 1/3$.  Then
$\gamma\in\Gamma$ can be parameterized as in (2.6).  Of course, the
key point is that the $\theta_0$ value in (2.6) determines $\gamma$
only up to the constant $\varphi_0 \in\mr/(2\pi\mz)$.  Thus, the
different possible limits are distinguished by having different values
for $\varphi_0$, and it is this possibility that will now be ruled
out.
  
To begin this task, first focus on a point in $C-\Lam$ which is very
close to $\mr\times\gamma$ for a particular $\gamma\in\Gamma$.  Here,
suppose that $\gamma$ is parameterized by (2.6) for some choice
$\theta_0$ and $\varphi_0$. With respect to this question of
$\varphi_0$, remember that for any fixed $\delta > 0$, but small, the
points of $C$ where $|f| < \delta$ that have identical $(t, f)$
coordinates form some $|p|$ subsets, each with $m\equiv m_\gamma$
members.  All members of the same subset have $\varphi$ coordinates
which differ by a very small amount, while members of different
subsets have $\varphi$ coordinates which differ by $2\pi k/p$ where $k
\not= 0\; \text{mod}(p)$.

In any event, some neighborhood of the chosen point in $C-\Lam$ has a
parameterization by coordinates $(\tau, u)$ with $\tau$ periodic and
$u$ small in absolute value and with sign that of $p$ via
\begin{itemize}
\item $t=\tau$,
\item $f=u$,
\item $\varphi=\varphi_0+\tau\; p'/p+x(\tau,u)$,
\item $h=(p'/p)\sin^2(\theta_0)u+y(\tau,u)$.
\autonum%(2.19)
\end{itemize}
Here, $|x|$ and $|y|/|u|$ are small and, in particular, much less than
$2\pi/|p|$.  By virtue of (1.5), the fact that $C$ is
pseudoholomorphic manifests itself in the fact that $x$ and $y$ obey
the differential equation
\begin{itemize}
\item $ x_u = -g^{-2 }\sin^{-2}(\theta) y_\tau$ ,
\item $ x_\tau = \sin^{-2}(\theta) y_u +
(p'/p)(\sin^2(\theta_0)/\sin^2(\theta) - 1)$.
\autonum%(2.20)
\end{itemize}
Here, $g = \sqrt{6}\; e^{-\sqrt{6}s} (1 + 3 \cos^4 \theta)^{1/2.}$
Note that the subscripts `$\tau$' and `$u$' on the variables $x$ and
$y$ indicate the partial derivative by the corresponding coordinate.
This notation is used frequently in the subsequent discussions.
 
In this last equation, both $\theta$ and $g$ are functions of the
variables $\tau$ and $u$, but the $\tau$ dependence is only implicit,
through the dependence of $y$ in (2.19) on $\tau$.  Thus, view
$\theta$ and $g$ as functions of $f$ and $h$, the former via $h/f =
\sqrt{6} \cos \theta \sin^2\theta (1 - 3 \cos^2 \theta)^{-1}$, and the
latter via the relation $g = \sqrt{6} (f^2 + h^2
\sin^{-2}\theta)^{1/2}$.  This understood, the right side of the first
line in (2.19) can by written as the $\tau$ derivative of the
restriction to $C$ of a function $Q$ on $\mr\times(S^1 \times S^2)$.
Indeed, $Q$ can be any function of the variables $f $ and $h$ whose
partial derivative in $h$ is $-g^{-2} \sin^{-2}(\theta)$.  For
example,
\begin{equation}
Q = - \sqrt{6}\; f^{-1} \ln(csc \theta+ \cot \theta) .
\end{equation}
%(2.21)
In any event, with $Q$ chosen, the first line in (2.21) reads
\begin{equation}
x_u = Q_\tau .
\end{equation}
%(2.22)
With (2.22) understood, remember that there are some $m$ members of
$C$ which are very close to the chosen point and have the same $(t,
f)$ value.  Thus, $C$ determines not just one pair of functions $(x,
y)$ as in (2.19), but a set, $\{(x_j, y_j)_{1\le j\le m}\}$ of $m$
such pairs.  Note that this set can only be ordered locally near each
point in the compliment of $\pi(\Lam)$.  The ordering may be permuted
around circles which enclose points of $\pi(\Lam)$ or around the $u =
$~constant circle.  In any event, $\Sigma_j x_j$ is a bonafide
function of $(\tau, u)$, and
\begin{equation} 
\underline{x}(u) \equiv \int_{ 0\le\tau\le 2\pi|p|} \Sigma_jx_j(\tau,
u)\;d\tau
\end{equation}
%(2.23)
is a function just of the coordinate $u$.
 
The claim here is that $\underline{x}$ is the constant function.
Indeed, this claim is a consequence of the observation that $\Sigma_j
Q(y_j)$ is also a bonafide function of $(\tau, u)$; and thus it
follows from (2.22) that the derivative of $\underline{x}$ is zero.

The fact that $\underline{x}$ is constant implies that the group of
$m$ points under observation must keep the same $\varphi_0$ value no
matter the size of $u$.  Indeed, a change in $\varphi_0$ must change
$\underline{x}$, as can be seen from (2.19).  As argued at the outset,
the constancy of $\varphi_0$ implies the assertion of Lemma 2.4.

\sh{(b)\qua Proof of Proposition 2.3}
  
The first two points of Proposition 2.3 are simply restatements from
Proposition 2.2.  The only point at issue here is the last one.  As
remarked at the outset, the last point in the proposition would follow
directly from Theorem 1.4 of \cite{HWZ1} were the closed orbits of the Reeb
vector field non-degenerate in a certain technical sense which is
satisfied here only for the closed orbits with $\theta_0\in\{0,\pi\}$
and, in the non-orientable $z$--axis line bundle case, that for which
the fundamental class generates the first homology over $\mz$.
However, the degeneracies here are due entirely to the fact that the
contact form $\alpha$ is invariant under the $T = S^1 \times S^1$
subgroup of isometries of $S^1 \times S^2$, those that fix the poles
of $S^2$.  This fact can be used to modify the arguments in \cite{HWZ1} to
apply here.  Although such a modification is straightforward, the
presentation of the details would be lengthy, and thus an alternate
proof is offered below.

Before embarking on the details of the proof, there are some
preliminary comments to be made.  First, the proof of the last point
in Proposition 2.3 is given below only for the case where the end in
question is concave and has orientable $z$--axis line bundle.  As in
the proof of Proposition 2.2, the argument for the convex end case is
identical in all essential aspects and thus left to the reader.
Meanwhile, as the assertions are local to the ends of $X$, the
non-orientable $z$--axis line bundle case can be treated as a
$\sigma$--equivariant example of the orientable $z$--axis line bundle
case.  Moreover, as the assertions in the third point of Proposition
2.3 are local to each end of $C$, attention can be restricted to a
single end.  This said, there is but one limit closed Reeb orbit
involved.
   
The final comment here is that the proof of the last point of
Proposition 2.3 is broken into six steps, and all but the final step
assume that the limit Reeb orbit for the end in question has neither
$\cos^2 \theta_0 = 1/3$ nor $\theta_0 \in \{0,\pi\}$.  The case where
the Reeb orbit has $\cos^2 \theta_0 = 1/3$ can be dealt with using
simple modifications of the arguments given below; basically, the
modifications involve the switching of the roles played by the pair
$(t, f)$ with those of $(\varphi, h)$.  As the details add nothing
novel, they will be left to the reader.  Likewise, the essentials of
the argument in the case where $\theta_0\in\{0,\pi\}$ are the same as
those given below, however, there are some specific differences which
deserve comment.  In particular, these comments constitute the final
step of the proof.  (As remarked previously, the case of Proposition
2.3 where the Reeb orbit has $\theta_0 \in \{0, \pi\}$ is also
directly a consequence of some general results in \cite{HWZ1}.)
  
With the preceding understood, what follows are the details of the
proof of the third point of Proposition 2.3.
\bigskip

{\bf Step 1}\qua To set the stage, return to the notation used in Step
9 of the proof of Proposition 2.2.  Thus, the end in question of $C$
is referred to henceforth as $C$; and its small $|f|$ portion is
parameterized as in (2.13).  This done, re-introduce the two component
column vector $\lam$ from (2.15) and (2.18).  Concerning the equation
in (2.15), note that the constant $\zeta'$ that appears in $L_0$ from
(2.16) is zero while $\zeta > 0$.  In fact,
\begin{equation}
\zeta = 6^{1/2} \sin^2 \theta_0 (1 + 3 \cos^2 \theta_0) (1 + 3 \cos^4
\theta_0)^{-1/2} |1 - 3 \cos^2 \theta_0|^{-1}.
\end{equation}
%(2.24)
The operator $L_0$ is a formally self-adjoint operator on the $\mr^2$
valued functions on $S^1$ and so has a complete set of eigenvectors.
Having constant coefficients, the eigenvalues and eigenvectors can be
readily found.  In particular, the corresponding eigenvalue set is
\begin{equation}
\left\{2^{-1 }(-\zeta \pm (\zeta^2 + 4 n^2/(m|p|)^2)^{1/2}): n =
 0, 1, 2, \ldots\right\}.
\end{equation}
%(2.25)
Note that there is a single zero eigenvalue, one of the $n = 0$ cases
in (2.25).  The corresponding eigenvector is the constant column
vector $e_-$ with top entry 1 and lower entry 0.  The other $n = 0$
eigenvalue is $-\zeta$ and it also has multiplicity one with a
constant eigenvector.  The $n \not= 0$ eigenvalues have multiplicity 2
and the components of the eigenvectors are linear combinations of the
functions $\sin(n\tau/(m|p|))$ and $\cos(n\tau/(m|p|))$.
	
\medskip
{\bf Step 2}\qua Introduce the $L^2$--orthogonal projections,
$\lam^{_+,-,0}$, of the vector $\lam$ in (2.15) onto the respective
spans of the eigenvectors of $L_0$ with eigenvalues where are positive
$(\lam^+)$, negative $(\lam^-)$ and zero $(\lam^0)$.
  
Here is a key fact:
\narrower\sl
The $L^2$--orthogonal projection, $\lam^0$, of the vector
$\lam$ in (2.15) onto the kernel of
$L_0$ is zero.\autonumm
\endnarrower
%(2.26)

Indeed, $\lam^0$ is constant by virtue of (2.22) as its bottom entry
is zero and its top entry is the average value of $x(\tau, \rho)$
around the $\rho$ = constant circles.  This understood, the constant
in question is zero by virture of (2.18).

{\bf Step 3}\qua The purpose of this step is to prove that the
function $f^+ (\rho) = \int d\tau |\lam^+|^2(\rho, \tau)$ and the
analogously defined function $f^-(\rho)$ both decay exponentially fast
to zero as $\rho$ tends to infinity.  To see that such is the case,
first let $E$ denote the smallest of the absolute values of the
non-zero eigenvalues of $L_0$.  Second, use Lemma 2.5 to find $\rho_0$
so that the $\rho > \rho_0$ versions of the $\cR$ term in (2.15) obey
$|\cR(\lam,\nabla\lam)| \le 10^{-2} E |\lam|$.  Third, for fixed $\rho
> \rho_0$, consider (2.15) to be an equality between $\mr^2$--valued
functions on the circle.  This done, then the $L^2$ inner product on
the $\tau$ parameterized circle of both sides of (2.15) with
$\lam^+(\rho, \cdot)$ leads to the inequality
\begin{equation}
2^{-1} \partial_\rho f^+ + Ef^+ - 10^{-2 }E (f^+ + f^-) \le 0 .
%(2.27)
\end{equation}
Meanwhile, the analogous inner product of both sides of (2.15) with
$\lam^-$ leads to
\begin{equation}
2^{-1} \partial_\rho f^- - Ef^- + 10^{-2 }E (f^+ + f^-) \ge 0 .
%(2.28)
\end{equation}
It now follows from these last two equations that $f\equiv f^-
- 0.2 f^+$ obeys the differential inequality
\begin{equation}
2^{-1 }\partial_pf- 0.97 E f \ge 0 .
\end{equation}
This last equation can be integrated to find that when $\rho'$ is
large and $\rho>\rho'$, then
\begin{equation}
f(\rho) \ge e^{2\delta(\rho-\rho')}f(\rho') ,
\end{equation}
where $\delta = 0.97 E$.
  
There is one immediate conclusion to draw from (2.30) which is this:
As $f$ is supposed to have zero for its $\rho\rightarrow\infty$ limit,
it follows from (2.30) that $f(\rho)$ is nowhere positive.  This is to
say that for all sufficiently large $\rho$,
\begin{equation}
f^-\le 0.2 f^+.
\end{equation}
The preceding inequality can now be inserted into (2.27) to yield
\begin{equation}
2^{-1}\partial_\rho f^+ +0.97 Ef^+\le 0.
\end{equation}
This equation can be readily integrated to see that
\begin{equation}
f^+(\rho)\le e^{-2\delta(\rho- \rho')} f^+(\rho')
\end{equation}
whenever $\rho'$ is large and $\rho>\rho'$.  Here, $\delta= 0.97 E$.
Together, (2.26), (2.31) and (2.33) assert that the function $g(\rho)
= \int d\tau |\lam(\rho, \tau)|^2$ has exponential decay to zero as
$\rho$ tends to infinity.

\medskip

{\bf Step 4}\qua This step proves the following assertion: For any
integer $k \ge 0$, the function $g_k(\rho)
\equiv\int\;d\tau|(\nabla^{\otimes k}\lam)(\rho,\tau)|^2$ has
exponential decay to zero as $\rho$ tends to infinity.
  
Here is the argument: First, remark that $\nabla$ commutes with both
$\partial_\rho$ and $L$.  Second, remark that $\nabla^{\otimes
k}\lam^+$ is in the span of the eigenvectors of $L_0$ with positive
eigenvalue while $\nabla^{\otimes k}\lam^-$ is in the span of those
with negative eigenvalue.  Third, differentiate both sides of (2.27)
and (2.28) $k$--times, and then take the respective $L^2$ inner
products with $\nabla^{\otimes k}\lam^+$ and $\nabla^{\otimes
k}\lam^-$.  Fourth, let $f^\pm(\rho)$ now denote $\int \;
d\tau|\nabla^{\otimes k}\lam^\pm|^2$.  Fifth, invoke Lemma 2.5 to
conclude that when $\rho$ is large (with lower bound depending on
$k$), this new $f^\pm$ obeys (2.27) and the new $f^-$ obeys (2.28).
Fifth, repeat the argument in the preceding step to obtain the desired
conclusion.
\bigskip

{\bf Step 5}\qua The third point in Proposition 2.3 follows from the
conclusions of the previous step using standard Sobolev inequalities.

\medskip
{\bf Step 6}\qua This step assumes that the closed Reeb orbit under
consideration has $\theta_0\in\{0,\pi\}$, and infact, $\theta_0 = 0$
as the discussion for the $\theta_0 =\pi$ case is identical save for
some innocuous sign changes.  The purpose of this step is to point out
the two places in the argument for the $\theta_0 = 0$ case where the
modifications to the just concluded argument are more than cosmetic.
In particular, the argument here requires the parameterization of $C$
using the functions $(a_1, a_2)$ as in (2.14).  This done, the only
other significant modification to the argument involves the constants
$\zeta'$ and $\zeta$ that appear in (2.16).  In this $\theta_0=0$
case, these are $\zeta'=\zeta=\frac{\sqrt{3}}{\sqrt{2}}$.  This said,
then the spectrum of the operator $L_0 $ in (2.16) is the set
\begin{equation}
\left\{-\sqrt{3}/\sqrt{2} + n/m : n \in\mz\right\} .
\end{equation}
Here, each eigenspace is two-dimensional; and the components of the
eigenvector for any given $n$ are linear combinations of
$\sin(n\tau/m)$ and $\cos(n\tau/m)$.
  
With the preceding understood, the rest of the argument for the
$\theta_0 = 0$ case is even simpler than that for the cases considered
previously because the $\theta_0 = 0$ version of $L_0$ has no zero
eigenvalue.

\section{Deformations}

Let $C$ be an irreducible, HWZ pseudoholomorphic subvariety in $X.$
Here, $X$ is as described in the introduction to the previous section,
with some ends concave and others convex, some with orientable
$z$--axis line bundle and others with the latter non-orientable.  After
a preliminary discussion on $C$Õs topology, the focus here is on
deformations of $C$ which preserve both the topology and its status as
an HWZ pseudoholomophic subvariety.  In this regard, note that the
conclusions of the subsequent discussions are summarized by
Propositions 3.1, 3.2 and 3.6.  In particular, Proposition 3.1 asserts
an adjunction formula that generalizes a formula for compact
pseudoholomorphic subvarieties that relates the self-intersection
number and intersection number with the canonical divisor to the Euler
characteristic.  Next, Proposition 3.2 asserts that the set of HWZ
subvarieties in $X $ has a natural topology that gives $C$ a
neighborhood that is homeomorphic to the zero locus of a smooth map
from a neighborhood of the origin in $\mr^N$ to $\mr^n$ for a suitable
choice of integers $N$ and $n$.  Here, $\mr^N$ and $\mr^n$ naturally
appear as kernel and cokernel of a Fredholm operator on $C$.  Finally,
Proposition 3.6 provides a geometric formula for the index, $I_C
\equiv N - n$, of this Fredholm operator.

\sh{(a)\qua The Euler characteristic of $C$}

The Euler characteristic of any embedded, compact, connected,
pseudoholomorphic submanifold $C \subset X$ is determined via the
adjunction formula from the class, $[C]$ of $C$ in $H_2(X; \mz)$.  In
this regard, remember that the symplectic form pulls back without
zeros to $C$ and so endows $C$ with a canonical orientation.  In any
event, the adjunction formula reads:
\begin{equation}
-\chi(C) = \la e, [C]\ra + \la c_1, [C]\ra,
\end{equation}
%(3.1)
where, $\la , \ra$ denotes the pairing between cohomology and
homology, $e\in H^2(X)$ is the image of the Poincar\'e dual to $[C]$,
and $c_1\in H^2(X)$ is the first Chern class of the canonical line
bundle $K \equiv\Lam^2T^{1,0}X$.  Here, $T^{1,0}X$ is the
$J$--holomorphic part of the complexified cotangent bundle of $X$.

Now, suppose that $C$ is required only to be an irreducible subvariety
whose singularities are purely transversal double points with local
intersection number 1. In particular, $C$ is the image of a connected
surface, $C_0$, via an immersion.  Let $m_C$ denote the number of
double points in $X$.  The corresponding adjunction formula in this
case reads:
\begin{equation}
-\chi(C_0) = \la e, [C]\ra + \la c_1, [C]\ra - 2 m_C .
\end{equation}
%(3.2)
In the general case where $C$ is irreducible and has singularities
other than transverse double points, there is still an adjunction
formula which gives the Euler characteristic of the smooth model for
$C$.  This is to say that $C$ is the image in $X$ of a compact,
complex curve, $C_0$, via a pseudoholomorphic map which is an
embedding off of a finite set in $C_0$.  And, the Euler characteristic
of $C_0 $ is given by the left-hand side of (3.2) where $m_C$ now
denotes the number of double points in a symplectic deformation of the
map from $C_0$ into $X$ which immerses $C_0$ with only transversal
double points self intersections that have local intersection number
1.  Indeed, the existence of such a deformation follows from the fact
that the singularities of a pseudoholomorphic variety are essentially
those of a complex curve in $\mc^2$ as shown by \cite{Mc}.  In any event,
make such a deformation and then (3.2) applies.

The next order of business is to introduce a version of the adjunction
formula that applies to a non-compact pseudoholomorphic subvariety
$C\subset X$.  In this case, the right-hand sides of (3.1) and (3.2)
are presently meaningless as $[C]$ is in $H_2(X, X-X_0)$ while both
$e$ and $c_1$ are in $H^2(X)$.  However, a pairing of a class
$\kappa\in H^2(X; \mz)$ with $[C]$ can be unambiguously defined with
the choice of a suitably constrained section of the restriction to $C$
of the complex line bundle $E_\kappa\rightarrow H$ with first Chern
class $\kappa$.  In particular, the section must have a compact zero
set on $C$.  With such a section chosen, a generic, compactly
supported perturbation produces a section which vanishes transversely
at a finite set of points in the smooth part of $C$.  Then, the count
of these points with the standard $\pm 1$ weights defines the pairing
$\la\kappa, [C]\ra$.  Note that this number is unchanged when the
chosen section of $E_\kappa$ is deformed through sections whose zero
sets all lie in a fixed, compact subset of $C$.

To put the just described count definition for $\la\kappa, [C]\ra$ in
a slightly larger context, note first that the chosen section of
$E_\kappa|_C$ can be extended as a section of $E_\kappa$ over the
whole of $X$.  The zero set of this extended section then carries a
fiducial two-dimensional, relative homology class that represents the
Poincar\'e dual of $\kappa$.  For example, if the original section
over $C $ and its extension to $X$ are chosen to have transversal zero
set, then this relative two-dimensional homology class is the
fundamental class of the zero set.  In this case, the pairing
$\la\kappa, [C]\ra$ is simply the intersection number of $C$ with the
zero set of the extended section.

This definition of $\la\kappa, [C]\ra$ needs some elaboration when
$\kappa\equiv e$ and $[C]$ is $e$'s Poincar\'e dual.  To proceed with
this elaboration, remark first that when $C$ is embedded in $X$, then
$E_e$ is defined by its transition function over the intersection of
two coordinate patches.  The first coordinate patch is $X-C$.  The
second is the image in $X$ via the metric's exponential map of a
certain open disk bundle, $N_0$, in the normal bundle of $C$.  Here,
the fiber radius of $N_0$ varies smoothly over $C$ to ensure that the
exponential map's restriction is an embedding.  This understood,
identify $N_0$ with its exponential map image.  Now, declare
$E|_{X-C}$ to be the trivial bundle $(X - C) \times\mc$ and declare
$E|_{N_0}$ to be the pull-back via projection to $C$ of $C$'s normal
bundle.  Thus, $E|_{N_0-C}$ has a canonical section and thus a
canonical trivialization with which to identify it with the
restriction to $(N_0 - C)$ of $(X - C) \times\mc$.

In the case where $C$ is not embedded, there are perturbations of $C$
in any given neighborhood its singular points that result in an
embedded, oriented submanifold.  Choose such a perturbation and use
the resulting submanifold to define $E_e $ as a proxy for $C$.

With the preceding understood, a three part digression is in order to
define the appropriate sections over $C$ of $E_e$ and the canonical
line bundle $K$.
\bigskip

{\bf Part 1}\qua This first part of the digression simply recalls a
definition from Proposition 2.2.  According to Propositions 2.2 and
2.3, each constant and large $|s|$ slices of $C \cap (X-X_0)$ consists
of a finite union of embedded circles that converges in the $C^\infty$
topology as $s\rightarrow\infty$ to a finite union of closed Reeb
orbits on $\partial X_0$.  Remember that this set of orbits, $\Gamma$,
is called \emph{the limit set} for $C$.  In particular, each end of
$C$ defines an element in $\Gamma$ and each element in $\Gamma$
corresponds via the aforementioned limit to one or more ends of $C$.
\bigskip

{\bf Part 2}\qua This part of the digression specifies the section of
$E_e|_C$ to be used when defining $\la e,[C]\ra$.  In this regard, it
is enough to specify the section over $C$'s intersection with the $|s|
\ge s_0$ portion of $X-X_0$ for any $s_0 \ge 0$ and so this is the
purpose of the subsequent discussion.  As previously noted, once such
a definition is made, then all extensions of this section to the
remainder of $C$ give the same count for the $\la e,[C]\ra$.  Thus,
the focus below is the definition of a nowhere zero section of the
restriction of $C$'s normal bundle to the large $|s|$ portion of $C$.

The task at hand begins with three remarks.  The first remark is that
Proposition 2.2 asserts that when $s_0$ is large, then $C$'s
intersection with $|s| \ge s_0$ portion of $X-X_0$ is an embedded
submanifold with boundary.  The second remark is that with $s_0$ so
chosen, then the restriction of $E_e|_C$ to this intersection is to be
identified with $C$'s normal bundle.  Here is the final remark: It
follows from Propositions 2.2 and 2.3 that when $s_0$ is large, then a
specification of a nowhere zero section of the restriction to $C$'s
intersection with the $|s| \ge s_0$ portion of $X-X_0$ of $C$Õs normal
bundle is determined, up to homotopy through non-vanishing sections,
by a nowhere zero section of the normal bundle in $\partial X_0$ to
$C$'s limit set.
    
To elaborate on these last remarks, note first that a component, $M =
S^1 \times S^2$, of $\partial X_0$ has two -dimensional homology whose
generator is the image of the fundamental class, $[S^2]$, of any
2--sphere of the form $\{\text{point}\} \times S^2 \subset S^1 \times
S^2$.  In this regard, these 2--spheres are oriented by the form $\sin
\theta d\theta\wedge d\varphi$ .  This said, then the pairing of $e$
with $[S^2]$ determines, up to bundle equivalence, the restriction of
$E_e$ to $M$ and thus to $ [s_0, \infty)\times M$.  In the context at
hand, this pairing is equal to the intersection number between $C$ and
any copy of $\{\text{point}\} \times S^2$ in $[s_0, \infty) \times M$.
 
A more explicit formula for this intersection number is available.
However, a digression is required before stating this formula.  This
digression explains how the large $|s|$ asymptotics of any given end,
$\cE$, of $C$ can be characterized in part by a suitable pair, $(p,
p')$, where $p$ is an integer and $p'$ is either another integer, or
one of the symbols $+$ or $ -$.  Here, the pair $(p, p')$ is deemed
suitable when the following requirement is met: If $p'= +$ or $-$,
require $p < 0$.  If $p'$ is an integer and $p < 0$, require that
$|p'/p| > \sqrt{3}/\sqrt{2}$; and if $p'$ is an integer with $|p'/p| <
\sqrt{3}/\sqrt{2}$, require that $p > 0$.
 
The meaning of $(p, p')$ is as follows: If $p' = +$ or $-$, then the
large and constant $|s|$ slices of the end in question converge to the
Reeb orbit with $\theta = 0$ or $\theta= \pi$, respectively.  This
understood, then $m \equiv |p|$ gives the multiplicity by which these
large $|s|$ circles cover the Reeb orbit.  If $p'$ is an integer, let
$m\ge 1$ denote the greatest common divisor of $p$ and $p'$.  Then,
the pair $(p/m, p'/m)$ determine a circleÕs worth of closed Reeb
orbits as dictated in (1.8), and the large $|s|$ slices of the end in
question are asymptotic to one of the latter.  Meanwhile, the integer
$m$ gives the multiplicity by which these large $|s|$ slices of the
end cover the Reeb orbit.

With the digression complete, remark that the intersection number in
question is the sum of the first components from the integer pairs
$(p, p')$ that come from the ends of $C$ that lie in the component
given by $M$ of $X-X_0$.

Now, suppose that the large $|s|$ slices of the end, $\cE$, is
asymptotic, in the manner described by Propositions 3.2 and 3.3, to a
Reeb orbit $\gamma\subset M$.  What follows describes how nowhere zero
sections of $\gamma$'s normal bundle in $M$ produce nowhere zero
sections of the large $|s|$ portion of $\cE$'s normal bundle.  To
begin the story, let $N_\gamma$ denote a small radius disk bundle in
the normal bundle to $\gamma$ in $M$.  Here, the radius should be such
that the metric's exponential map embeds $N_\gamma$ in $M$.  This
understood, identify $N_\gamma$ with its exponential map image.  With
these identifications made, then, Propositions 3.2 and 3.3 provide
$s_1 > s_0$ and a description of the $|s| > s_1$ portion of $\cE$ as
the image over $(s_1,\infty) \times\gamma$ of a multi-valued section
of $N_\gamma$.  This is to say that there is a degree $m$ covering map
$\pi\co S^1\rightarrow\gamma$ and a section, $\eta$, over $(s_1,
\infty) \times S^1$ of $\pi^*N_\gamma$ such that the composition of
$\eta$ with the tautological map $\hat{\pi}\co\pi^*N_\gamma\rightarrow
N_ \gamma$ maps $(s_1, \infty) \times S^1$ diffeomorphically onto the
$|s| > s_1$ portion of $\cE$.  This composition, of $\eta$ with
$\hat{\pi}$ identifies the normal bundle to the large $|s|$ portion of
$\cE$ with that of the image of $\eta$ in $\pi^*N \rightarrow (s_1,
\infty)\times S^1$.  Meanwhile, the latter is canonically isomorphic
to $\pi^*NÕ$s pullback over itself via its defining projection to
$(s_1, \infty) \times\gamma$.  Following this chain of bundle
isomorphisms produces a nowhere zero section of the large $|s|$
portion of $\cE$'s normal bundle from a nowhere zero section of the
Reeb orbit's normal bundle in $M$.
   
Now, to provide such a section of $\gamma$'s normal bundle, consider
first the case where the component $M$ has orientable $z$--axis line
bundle.  Also, suppose that $\theta = \theta_0 \notin \{0, \pi\}$ on
$\gamma$.  Then, a nowhere zero section of $\gamma$'s normal bundle is
defined by the vector field $\partial_\theta$ along $\gamma$.
 
If $\theta_0 = 0$ or $\pi$ on $\gamma$, the functions $(x_1
\equiv\sin\theta \cos\varphi$, $ x_2\equiv \sin\theta \sin\varphi)$
are zero on $\gamma$ and the triple $(t, x_1, x_2)$ are good
coordinates for $M$ near $\gamma$.  With this understood, the vector
field tangent to any line through the origin in the $(x_1, x_2)$ plane
defines along $\gamma$ a nowhere vanishing section of $\gamma$'s
normal bundle.

Now consider the case where $M \subset\partial X_0$ defines an end of
$X$ whose $z$--axis line bundle is non-orientable.  The first order of
business is to write out the 2--1 covering map, $\hat{\sigma}\co S^1
\times S^2 \rightarrow M$ $(= S^1 \times S^2)$ whose deck
transformations are generated by $\sigma$ in (2.1).  To do so, view
$S^2$ as the unit sphere in $\mr^3$ and introduce Cartesian
coordinates $(t, x_1 = \sin \theta \cos \varphi$, $x_2 = \sin \theta
\sin \varphi$, $x_3 = \cos \theta)$ on the domain $S^1 \times S^2$.
Let $(t', x'_1, x'_2, x'_3)$ denote the analogous coordinates for the
range.  Then, the map is defined so that Cartesian coordinates of
$\hat{\sigma}(t, x_1, x_2, x_3)$ are

\begin{itemize}
\item $t' = 2 t$,
\item $x'_1 = x_1$,
\item $x'_2 = \cos(t) x_2 + \sin(t) x_3$,
\item $x'_3 = - \sin(t) x_2 + \cos(t) x_3$. \autonum%(3.3)
\end{itemize}

With $\hat{\sigma}$ defined, turn back to the task at hand.  In this
regard, each component of the inverse image of a closed Reeb orbit
under the map $\hat{\sigma}$ is either described as in Step 4 of the
proof of Proposition 2.2, or else lies in a constant $t$ slice where
$\cos^2 \theta = 1/3$.  In any event, fix attention on a closed Reeb
orbit $\gamma\subset M$.  Here, there are two cases to consider,
depending on whether $[\gamma]$ is or is not a generator of $H_1(M;
\mz)$.

Consider first the case where $\gamma$ is not an integral generator of
$H_1$.  Then $\gamma$ will necessarily be an even multiple of a
generator. (This follows from the form of $\sigma$ in (3.3) and the
fact that $\theta$ is constant on any inverse image of $\gamma$.)  As
an even multiple of a generator of $H_1(M; \mz)$, the inverse image
under the map in (3.3) of $\gamma$ must have two components where each
is mapped diffeomorphically by (3.3) onto $\gamma$.  With this said,
it follows that a nowhere vanishing section of $\gamma$'s normal
bundle in $M$ is defined via the following two-step procedure: First,
take a section of the normal bundle in $S^1 \times S^2$ of one of
$\gamma$'s inverse images under (3.3) by the rules which were just
described for the orientable $z$--axis line bundle case.  Second, use
(3.3) to push this section down to $M$ as a section of the normal
bundle to $\gamma$.

Next, consider the case where $\gamma$ generates $H_1(M; \mz)$.  In
this case, the inverse image of $\gamma$ under (3.3) must be connected
and (3.3) maps this inverse image as a 2--1 covering map onto
$\gamma$.  Moreover, as $\theta$ = constant on such an inverse image,
it follows that the functions $x'_2$ and $x'_3$ from (3.3) vanish on
$\gamma$ and together with $t'$ define good coordinates for $M$ near
$\gamma$.  With this understood, the vector field tangent to any line
through the origin in the $(x'_2, x'_3)$ plane defines along $\gamma$
a nowhere zero section of $\gamma$'s normal bundle.

\bigskip

{\bf Part 3}\qua This part of the digression defines the desired
section over $C$ of the canonical bundle $K$ and thus gives meaning to
$\la c_1, [C]\ra$.  In this regard, note that it is again sufficient
to specify the section only over $C$'s intersection with each
component of the $|s| \ge s_0$ part of $X-X_0$ for very large $s_0$.
In the discussion below, only the case of a concave component will be
considered, as the convex case can be directly obtained from the
latter.

With the preceding understood, consider first the promised section of
$K$ over $C$'s intersection with $[s_0, \infty) \times M$ where $M
\subset\partial X_0$ is a concave component with orientable $z$--axis
line bundle.  In this regard, take $s$ so large that $C \cap ([s_0,
\infty) \times M)$ is a union of cylindrical components with each
defining in the limit a particular loop in $C$'s limit set.  The
specification of a section of $K$ over a given component now depends
on whether the corresponding loop in the limit set has
$\theta\in\{0,\pi\}$ or not.  Consider first the case where
$\theta\notin\{0,\pi\}$ on the loop.  Near such a component, the
complex valued 1--forms $dt + i g^{-1} df$ and $\sin^2 \theta d\varphi
+ i g^{-1} dh$ span $T^{1,0}X$ and so their wedge product, $(dt + i
g^{-1} df) \wedge (\sin^2 \theta d\varphi + i g^{-1} dh)$ gives a
section of $K$.  (Note that this section is defined over the whole of
$[0, \infty) \times M$, and its zero set is the locus where
$\theta\in\{0, \pi\}$.)
 
To consider the case where a limit set loop of a component of $C \cap
([s_0, \infty) \times M)$ has $\theta = 0$ or $\pi$, first introduce
the `Cartesian coordinates' $(x_1 \equiv \sin\theta \cos\varphi,$ $x_2
\equiv \sin\theta \sin\varphi)$.  These are smooth coordinates near
the $\theta = 0$ and $\theta = \pi$ loci.  Furthermore, $(dt + i ds)
\wedge (dx_1 - i dx_2)$ specifies a non-vanishing section of $K$ over
the $\theta = 0$ locus while $(dt + i ds) \wedge (dx_1 + i dx_2)$
plays the same role for the $\theta = \pi$ locus.  Now, the required
section of $K$ over a component of $C \cap ([s_0,\infty) \times M)$
with limit set loop in the $\theta = 0$ locus is obtained by
restriction of any smooth extension of the section $(dt + i ds) \wedge
(dx_1 - i dx_2)$ to a tubular neighborhood of this locus.  The
required section of $K$ in the case where the limit set loop lies in
the $\theta = \pi$ locus is defined analogously.

Now consider the promised section of $K$ over $C$'s intersection with
$[s_0, \infty) \times M$ where $M \subset\partial X_0$ is a concave
component with non-orientable $z$--axis line bundle.  For this purpose,
it is important to keep in mind that the canonical line bundle of $[0,
\infty) \times M$ pulls up via the map in (3.3) to the canonical
bundle just studied for the orientable $z$--axis case.
  
With the preceding point taken, fix $s$ large enough to insure that $C
\cap ([s_0, \infty) \times M)$ is a union of cylindrical components
with each defining in the limit a particular loop in $C$'s limit set.
Now, the inverse image of this loop has $\theta$ locally constant and
there are two cases, as before, depending on whether $\theta\in\{0,
\pi\}$ or not on this inverse image.  Consider the latter case first.
If the constant $s$ slices of the end in $C$ converge as
$s\rightarrow\infty$ as a multiple cover of a closed Reeb orbit with
two inverse images under the map in (3.3), use the push-forward from
one of them of the section $(dt + i g^{-1} df) \wedge (\sin^2 \theta
d\varphi + i g^{-1} dh)$. If the constant $s$ slices of the end in
question converge as $s\rightarrow \infty$ as a multiple cover of a
closed Reeb orbit with a single inverse image under (3.3), proceed as
follows to obtain the appropriate section: Note first that the section
$(dt + i g^{-1} df) \wedge (\sin^2 \theta d\varphi + i g^{-1} dh)$
changes sign under the involution in (2.1) but its product with
$e^{it}$ does not. Therefore, $e^{it}(dt+ig^{-1}df)\wedge(\sin^2\theta
d\varphi+ig^{-1}dh)$ is the pull-back via (3.3) of a section over
$[0,\infty)\times M$ of the canonical bundle of $X$.  Use the latter
section for the definition of $\la c_1,[C]\ra$.

In the case where $\theta\in\{0, \pi\}$ on the inverse image of the
limit loop, then the inverse image under the map in (3.3) of the
component of $C \cap([s_0, \infty) \times M)$ in question has two
components, one with $\theta$ very nearly zero and the other with
$\theta$ very nearly equal to $\pi$.  With this point understood, a
section on the component in question of $C \cap ([s_0, \infty) \times
M)$ is obtained as follows: Use $\sigma$ to push forward any nowhere
zero section on a tubular neighborhood of the $\theta = 0$ locus whose
restriction to the $\theta = 0$ locus equals $(dt + i ds) \wedge (dx_1
- i dx_2)$.

With the digression now over, return to the original subject, which is
an adjunction formula for the smooth model for $C$.  In this regard,
note that Proposition 2.2 asserts that $C$ is the image of a complex
curve, $C_0$, with cylindrical ends via a pseudoholomorphic map which
is an embedding off of a finite set of points.  Moreover, Propositions
2.2 and 2.3 imply that the Euler characteristic of $C_0$ can be
computed as the usual algebraic count of the zeros of a section of
$C_0$'s tangent bundle if, on the large $|s|$ portion of $C_0$, this
same tangent vector pushes forward to $\mr$ via the composition of the
pseudoholomorphic map with the projection $\mr\times M \rightarrow\mr$
as a nowhere zero multiple of the vector $\partial_s$.  With this last
point understood, it then follows from the given choice of sections of
$E_e|_C$ and $K$ for the respective definitions of $\la e, [C]\ra$ and
$\la c_1, [C]\ra$ that $\chi(C_0)$ is given by (3.2) where $m_C$ is
defined as in the case where $C$ was compact.  Indeed, the argument
for (3.2) in this case amounts to an essentially verbatim
recapitulation of a standard proof for the compact case as given in
\cite{GH}.

The following proposition summarizes the conclusions of the preceding
adjunction formula discussion:
\bigskip

\noindent
{\bf Proposition 3.1}\qua\sl Let $X$ be as described in the
introduction to Section 1.  Let $C \subset X$ be an irreducible,
pseudoholomorphic subvariety, and let $C_0$ denote the smooth model
curve for $C$.  Then $\chi(C_0)$ is given by
$$-\chi(C_0) = \la e, [C]\ra + \la c_1, [C]\ra- 2 m_C,$$
where:
\begin{itemize}  
\item The cohomology-homology pairings are defined as described
above.
\item $m_C$ denotes the number of double points of any
perturbation of the defining map from $C_0$ onto $C$ which is
symplectic, an immersion with only transversal double point
singularities with locally positive self-intersection number, and
which agrees with the original on the complement of some compact set.
\end{itemize}
If $C$ is immersed, then this formula can be rewritten as
$$-\chi(C_0) = \text{degree}(N) + \la c_1, [C]\ra,$$
where degree$(N)$ is the degree of the normal bundle to the
immersion as defined using a section that is non-vanishing at large
$|s|$ and is homotopic at large $|s|$ through non-vanishing sections
to that described in Part 2, above.\rm
\bigskip

Some readers may prefer the formula that writes $-\chi(C_0)$ in terms
of degree$(N)$ because the pairing $\la e, [C]\ra$ is not invariant
under deformations of $C$ unless $C$'s large $|s|$ slices converge
with multiplicity 1 to a set of distinct, limit Reeb orbits.

\sh{(b)\qua A topology on the set of pseudoholomorphic subvarieties}

The first task is the introduction of a certain topology on the set
me, $\fM_{e,\chi}$ of irreducible, pseudoholomorphic subvarieties in
$X$ with fundamental class Poincar\'e dual to the given class $e \in
H^2(X; \mz)$ and with the given number $\chi$ equal to the Euler
characteristic of the model curve $C_0$.  The topology in question
comes from the metric for which the distance between a pair $C$, $C'
\subset \fM_{ e,\chi}$ is
\begin{equation}
\text{sup}_{x\in C} \text{distance}(x, C') + \text{sup}_{x'\in
C'}\text{distance}(C, x')\; .
\end{equation}
%(3.4)
Given this definition, then the next order of business is a structure
theorem for a neighborhood in $\fM_{e,\chi}$ of any given subvariety.
The following proposition summarizes the story:
\bigskip

\noindent
{\bf Proposition 3.2}\qua {\sl Let $C \subset \fM_{e,\chi}$.  There
exists a Fredholm operator $D$ and a homeomorphism from a neighborhood
of $C$ in $\fM_{e,\chi}$ to the zero set of a smooth map from a ball
in the kernel of $D$ to the cokernel of $D$.}
\bigskip

The proof of this proposition is given in the next subsection.

The description of $D$ is simplest when $C$ is compact and the
associated pseudoholomorphic map $\varphi$: $C_0 \rightarrow X$ is an
immersion.  In this case, there is a well defined `normal bundle', $N
\rightarrow C_0$ that is a real 2--plane bundle whose restriction to
any open $K \subset C_0$ embedded by $\varphi$ is the
$\varphi$--pullback of the normal bundle to $\varphi(K)$.  The almost
complex structure on $X$ endows $N$ with the structure of a complex
line bundle over $C_0$, and the associated Riemannian metric from $X$
can then be used to give $N$ the structure of a Hermitian line bundle
with a holomorphic structure.  The induced $d$--bar operator on the
space of sections of $N$ will be denoted by $\opart$.  (The $\opart$
operator used here is twice the usual $\partial/\partial\oz$.)
  
In the case at hand, the operator in Proposition 3.2 is the first
order, $\mr$--linear operator from $C^\infty(N)$ to $C^\infty(N\otimes
T^{0,1}C_0)$ that sends a section $\lam$ of $N$ to
\begin{equation}
D\lam\equiv\opart\lam+\nu\lam+\mu\olam .
\end{equation}
%(3.5)
Here $\nu$ and $\mu$ are respective sections of $T^{0,1}C_0$ and
$N^2\otimes T^{0,1}C_0$ that are determined by the 1--jet of the almost
complex structure $J$ along $C$.  Although the kernel dimension may
depend on $\nu$ and $\mu$, the index~$\equiv$~dim(kernel) --
dim(cokernel) of $D$ does not.  In fact, the index is the same as that
of (as an $\mr$--linear operator) namely: index$(D) = 2 $ degree$(N) +
\chi(C_0$).  As $\chi C_0)$ obeys (3.2) and degree$(N) = \la e, [C]\ra
- 2 m_C$, this index can also be written in various equivalent ways,
for example:

\begin{itemize}
\item index$(D)=\la e,[C]\ra-\la c_1,[C]\ra-2m_C$.
\item index$(D) = -\chi(C_0) - 2 \la c_1, [C]\ra$.\autonum%(3.6)
\end{itemize}

When $C \subset \fM_{e,\chi}$, still compact, is not immersed,
Proposition 3.2's operator $D$ is more complicated.  What follows is a
brief description of this new operator.  The definition of $D$
requires, as a preliminary step, the introduction of a first-order
differential operator, $\uD$, which sends a section of
$\varphi^*T_{1,0}X$ to one of $\varphi^*T_{1,0}X \otimes T^{0,1}C_0$.
Here and below, $T_{1,0}$ denotes the holomorphic part of the
corresponding complexified tangent bundle.  The operator $\uD$ differs
from the corresponding $\opart$ by a zeroÕth order, $\mr$--linear
multiplication operator and thus has the same schematic form as
depicted on the right-hand side of (3.5).  This operator $\uD$ is
defined so that its kernel provides the vector space of deformations
of the map $\varphi$ which remain, to first-order, pseudoholomorphic
as maps from $C_0$ into $X$.  In this regard, $\uD$ is not quite the
sought after operator as its use in Proposition 3.2 would allow only
those deformations which preserve the induced complex structure on the
image curve.  The point is that the cokernel of $\uD$ is too large
when $\chi \le 0$.  Meanwhile, when $\chi \ge 0$, the kernel of $\uD$
is too big as it contains deformations which come from holomorphic
automorphisms of $C_0$.

To address these problems, introduce, first of all, the usual $\opart$
operator which sends sections of $T_{1,0}C_0$ to $T_{1,0}C_0 \otimes
T^{0,1}C_0$.  The kernel of this operator, $V$, is trivial when
$\chi(C_0) < 0$, but not trivial otherwise.  (Its dimension over $\mc$
is $\chi + 1$.)  Meanwhile, let $V_\bullet$ denote the cokernel of
this same version of $\opart$.  The complex dimension of $V_\bullet$
is $-3\chi/2$ when $\chi < 0$, one for a torus and zero for a sphere.
Fix some favorite subspace of smooth and compactly supported sections
of $T_{1,0}C_0 \otimes T^{0,1}C_0$ that projects isomorphically to the
cokernel of $\opart$ and identify the latter with $V_\bullet$.
  
Next, remark that as $\varphi$ is pseudoholomorphic, its differential
provides a $\mc$--linear map, $\partial\varphi$: $T_{1,0}C_0
\rightarrow\varphi^*T_{1,0}X$ and thus one, also denoted by
$\partial\varphi$, from $T_{1,0}C_0 \otimes T^{0,1}C_0$ to
$\varphi^*T_{1,0}X \otimes T^{0,1}C_0.$ In particular, note that the
appropriate version of $\partial\varphi$ sends $V$ and $V_\bullet$
injectively into the respective kernel and target space of $\uD$.
With this understood, then $\uD$ induces an operator, the desired $D$,
that maps $C^\infty(C_0; \varphi^*T_{1,0}X)/\partial\varphi(V)$ to the
$L^2$ complement in $C^\infty(C_0; \varphi^*T_{1,0}X \otimes
T^{0,1}C_0)$ of $\partial\varphi(V_\bullet)$.  The conclusions of
Proposition 3.2 hold for this $D$ when $C$ is not immersed.
  
By the way, note that the index of this new $D$ is still given by the
formula in (3.6) where $m_C$ is now interpreted as in Proposition 3.1.

The remainder of this subsection describes Proposition 3.2's $D$ in
the case where $C \subset X$ is a non-compact, pseudoholomorphic
subvariety.  The discussion has been divided into three parts.
\bigskip

{\bf Part 1}\qua As the story is simplest when $C$ is immersed, this
condition will be assumed until the final part.  In this regard, note
that the removal of the immersion assumption requires no new
technology since a pseudoholomorphic subvariety is, in any event,
embedded where $|s|$ is large on $X-X_0$.
  
The first remark is that the operator $D$ is formally the same as that
which is described in (3.5).  In particular, $N$ is defined as before,
while the hermitian structure on $N$ and the Riemannian structure on
$TC_0$ are both induced by the Riemannian metric on $X$.  However,
there is some subtlety with the range and domain of $D$.  In
particular, these are defined as follows: First, a very small $\delta
> 0$ must be chosen.  There is an upper bound to the choice which is
determined by the properties of $C'$s limit set.  In any event, with
$\delta$ chosen, the domain of $D$ is the (Hilbert space) completion
of the set of smooth sections of $N$ for which
\begin{equation}
\int e^{-\delta |s|} (|\nabla\lam|^2 + |\lam|^2)
\end{equation}
%(3.7)
is finite; moreover, (3.7) defines the square of the relevant norm.
Here, and below, the function $s$ which was originally specified only
on $X-X_0 $ has been extended to the remainder of $X$ as a smooth
function.  By the way, the integration measure in (3.7) and in
subsequent integrals is the area form from the Riemannian metric on
$C$ that is induced by the metric on $X$.  Of course, given that $C$
is pseudoholomorphic, this measure is the same as that defined by the
restriction to $C$ of $|\omega|^{-1}\omega$ .
  
Meanwhile, the range space for $D$ is the completion of the set of
smooth sections of $N$ for which
\begin{equation}
\int e^{-\delta|s|} |\lam|^2
\end{equation} 
%(3.8)
is finite; in this case (3.8) gives the square of the relevant norm.
Let $\cL_1$ denote the just defined domain Hilbert space for $D$ and
let $\cL_0$ denote the range.
  
With the preceding understood, here is the key lemma:

\bigskip
\noindent
{\bf Lemma 3.3}\qua\sl {If $\delta$ is positive, but sufficiently
small, then the operator $D$ as just described extends as a bounded,
Fredholm operator from $\cL_1$ to $\cL_0$.  Moreover,}
\begin{itemize} 
\item {The index of $D$ as well as the dimensions of the kernel
and cokernel of $D$ are independent of $\delta$.}
\item {In fact, the kernel of $D$ is the vector space of sections
of $\lam$ of $N$ with $D\lam = 0$ and ${\rm sup} |\lam| < \infty$.}
\item {Polarize the norm square in (3.8) to obtain an inner
product on $\cL_0$ and represent the cokernel of $D$ as the orthogonal
complement to $D$'s image.  Then, multiplication by $e^{-\delta|s|}$
identifies {\rm cokernel}$(D)$ with the space of sections $\lam$ of $N
\otimes T^{0,1}C_0$ with $\int|\lam|^2 < \infty$ and which are
annihilated by the formal $L^2$ adjoint, $D^*$, of $D$.}
\item {There exists $\delta_1$ which is independent of $\delta$
and which has the following significance: Let $\cP_1$ be an
$\mr$--linear bundle homomorphism from $T^*C_0$ to ${\rm Hom}(N;
T^{0,1}C_0)$ with norm $|\cP_1| < \delta_1$.  Also, let $\cP_0$ be an
$\mr$--linear bundle homomorphism from $N$ to $N \otimes T^{0,1}C_0$.
Moreover, suppose that $e^{2\delta|s|} (|\cP_0| + |\cP_1|)$ is
bounded.  Then, $D + \cP_1(\nabla) + \cP_0$ extends as a bounded,
Fredholm operator from $\cL_1$ to $\cL_0$ whose index is the same as
$D$'s. }
\end{itemize}\rm
\medskip

This lemma is also proved below.
\bigskip

{\bf Part 2}\qua This part of the discussion describes the operator
$D$ in the case when $C$ is not immersed. But keep in mind that there
exists $s_0$ such that $C$'s intersection with the $|s| > s_0$ portion
of $X-X_0$ is nonetheless a disjoint union of embedded cylinders which
intersect the constant $|s|$ slices of $X-X_0$ transversely.

As in the case where $C$ is compact, the operator $D$ is constructed
by first introducing the operator $\uD$ that is defined just as in the
compact case.  Thus, $\uD$ maps sections of $\varphi^*T_{1,0}X$ to
those of $\varphi^*T_{1,0}X \otimes T^{0,1}C_0$.  By analogy with the
case where $C_0$ is immersed, the domain for $\uD$ is the completion
of the space of sections of $\varphi^*T_{1,0}X$ for which the
expression in (3.7) is finite; and (3.7) defines the square of the
norm for this completion.  Meanwhile, the range of $D$ is the
completion of the space of sections of $\varphi^*T_{1,0}X \otimes
T^{0,1}C_0$ for which (3.8) is finite, and (3.8) defines the square of
the norm for the completion in this case.  The domain Hilbert space
will be denoted by $\underline{\cL}_1$ and the range by
$\underline{\cL}_0$.
  
Also needed in this discussion are vector spaces which play the role
here that the kernel and cokernel of $\opart$: $T_{1,0}C_0\rightarrow
T_{1,0}C_0 \otimes T^{0,1}C_0$ play in the compact case.  In this
regard, the two versions of the linear map $\partial\varphi$ are still
available as $\varphi$ is, in any event, pseudoholomorphic.  With this
understood, let $V$ denote the vector space of sections $T_{1,0}C_0$
which are annihilated by $\opart$ and for which (3.7) is finite.
Meanwhile, let $V_\bullet$ denote the space of sections of $T_{1,0}C_0
\otimes T^{0,1}C_0$ for which (3.8) is finite and whose product with
$e^{-\delta |s|}$ is annihilated by the formal, $L^2$ adjoint of
$\opart$.
 
The linear map $\partial\varphi$ maps $V$ injectively into $\uD$'s
kernel in $\underline{\cL}_1$, and it maps $V_\bullet$ injectively
into $\underline{\cL}_0$.  Thus, $\uD$ induces an operator, $D$, from
$\cL_1 \equiv \underline{\cL}_1/ \partial\varphi(V)$ to $\cL_0
\equiv\underline{\cL}_0^T$ where $\underline{\cL}_0^T$ denotes the
orthogonal complement to $\partial\varphi(V_\bullet)$.  One then has:
\bigskip

\noindent
{\bf Lemma 3.4}\qua {\sl Whether or not $C$ is immersed, the assertions
of Lemma 3.3 hold with the operator $D$: $\cL_1\rightarrow\cL_0$ as
just described in the preceding paragraph.}
\bigskip

\noindent
The proof of this lemma will be left to the reader in as much as its
proof is a straightforward marriage of the arguments given below for
Lemma 3.3 with those given for the compact case.  Note that were all
closed Reeb orbits in $C$'s limit set isolated, then this lemma more
or less restates results from \cite{HWZ3}.
\bigskip

{\bf Part 3}\qua The remainder of this subsection contains the
following proof.\bigskip

\noindent
{\bf Proof of Lemma 3.3}\qua Although the lemma can be proved by
modifying the analysis in \cite{HWZ3}, an abbreviated argument will be
provided.  The argument below has two steps.

\medskip
{\bf Step 1}\qua The argument presented here for the proof is based on
some general facts about elliptic differential operators on manifolds
with cylindrical ends.  These facts were originally established in
\cite{LM}, but see also \cite{T5} which analyses a first order elliptic operator
in a context that has many formal analogies with the context here.
What follows here is simply a summary of those facts which are
relevant to the case at hand.

To set the stage for the subsequent discussion, consider a non-compact
manifold, $Y$, which comes with an open set having compact closure,
and a diffeomorphism from the complement of this open set to
$[-1,\infty) \times Z$.  This diffeomorphism will be used to identify
$[-1, \infty)\times Z$ as a subset of $Y$.  With this understood, let
$\rho$: $Y \rightarrow [-2, \infty)$ be a smooth function which
restricts to $[0, \infty) \times Z$ as the projection to the first
factor.

Now, suppose that $D$ is a first-order, elliptic operator on $Y$
taking sections of one vector bundle, $E$, to those of a second,
$E_\bullet$.  Suppose further that both of these bundles are provided
with fiber metrics and metric compatible connections.  Parallel
transport via these connections along the paths $[0, \infty)\times
\{\text{point}\} \subset [0, \infty) \times Z$ then identifies these
bundles with their respective restrictions to $\{0\} \times Z$, and
this identification will be explicit in what follows.  Note that the
latter identifications make the fiber metrics $\rho$--independent on
$[0, \infty) \times Z$.  Use $\nabla$ to denote the covariant
derivative of either connection.

Next, suppose that $\rho$ is the Euclidean coordinate on the
$[0,\infty)$ factor of $[0, \infty) \times Z$ and that the restriction
of $D$ here has the form:
\begin{equation}
D = A_0 \partial_\rho + L_0 + \alpha_1(\nabla) + \alpha_0
\end{equation}
%(3.9)
where:

\begin{itemize}
\item $A_0$ {is a $\rho$--independent, isometric isomorphism
between $E$ and $E_\bullet$.}
\item $L_0$ {is a $\rho$--independent, first-order operator which
differentiates only along vectors which are tangent to $Z$.  Require
that $A_0^\dagger L_0 = L_0^\dagger A_0$.}
\item $\alpha_1$ {is vector bundle homorphism from $T^*([0,
\infty) \times Z)$ into ${\rm Hom}(E, E_\bullet)$.  Here, it is
required that $|\alpha_1|$ should be small, where a precise upper
bound, $\delta_1$, depends on the norm of the symbol of $L_0$.
However, $\delta_1 < 1/10$ in any event.}
\item $\alpha_0$ {is a section of ${\rm Hom}(E, E_\bullet)$.}
\item {There exists $\epsilon > 0$ such that $e^{\epsilon\rho}
(|\alpha_1| + |\alpha_0|)$ is bounded on\nl
 $[0, \infty) \times Z$.}\autonum%(3.10)
\end{itemize}
 
To define appropriate domain and range spaces for $D$, it is necessary
to first choose a Riemannian metric on $X$ whose restriction to $[0,
\infty)\times Z$ is the product metric constructed using $d\rho^2$ on
the first factor and an $\rho$--independent metric on the second.

The definition of the domain and range for $D$ also requires the
choice of $\delta\in\mr$.  With $\delta$ chosen, the domain Banach
space, $\cL_1$, is obtained by completing the set of sections of $E$
for which (3.7) is finite using the expression in (3.7) for the square
of the norm after substituting $p$ for $s$.  Meanwhile, the range
Banach space is obtained by completing the set of section of
$E_\bullet$ for which the corresponding version of (3.8) is finite
using the latter for the square of the norm.

Under, the preceding assumptions, here is the fundamental conclusion
from \cite{LM}:
\bigskip

\noindent
{\bf Lemma 3.5}\qua\sl {If $\delta \not= 0$, but $|\delta|$ is
sufficiently small, then $D$ defines a Fredholm map from $\cL_1 $to
$\cL_0$.  Moreover, the following are true: }
\begin{itemize} 
\item {Each element in {\rm kernel}($D$) has a well defined
$\rho\rightarrow\infty$ limit on $[0, \infty) \times Z$.  In
addition:}
\begin{itemize}
\item[\rm(a)] {This limit is zero when $\delta < 0$, and it lies in
the kernel of $L_0$ when $\delta > 0$.}  
\item[\rm(b)] {When $\delta
> 0$, these limits define a homomorphism $\iota_\delta$: ${\rm
kernel}(D) \rightarrow {\rm kernel}(L_0)$.}
\item[\rm(c)] {Let $\lam\in {\rm kernel}(D)$ and let $k \in \{1, 2,
\ldots\}$.  Then $e^{|\delta|\rho/2} |\nabla^{\otimes k}\lam|$ is
square integrable on $Y$; and if $|\lam| \rightarrow 0$ as
$\rho\rightarrow\infty$, then $e^{|\delta|\rho/2 }|\lam|$ is also
square integrable on $Y$.}
\end{itemize}  
\item {Identify the cokernel of $D$ with the orthogonal
complement to the image of $D$ as defined by the inner product induced
by the norm in (3.8) on $\cL_0$.  With this identification understood,
let $\eta\in\lam\; {\rm cokernel}(D)$ and set $\lam\equiv
e^{-\delta\rho}\eta$.  Then $\lam$ is annihilated by the formal $L^2$
adjoint of $D$ it has a well defined $\rho\rightarrow\infty$ limit on
$[0, \infty) \times Z$.  In addition:}
\begin{itemize} 
\item[\rm(a)] {This limit is zero when $\delta > 0$, and it lies in
the kernel of $L_0$ when $\delta < 0$.  }
\item[\rm(b)] {When $\delta < 0$, these limits define a homomorphism
$\underline{\iota}_\delta$: ${\rm cokernel}(D)\rightarrow {\rm
kernel}(L_0)$. }
\item[\rm(c)] {Let $\lam\equiv e^{-\delta\rho}\eta$ with $\eta\in
{\rm cokernel}(D)$ and let $k \in \{1, 2, \ldots\}$.  Then, $
e^{|\delta|\rho/2 }|\nabla^{\otimes k}\lam|$ is square integrable on
$Y$; and if $|\lam|\rightarrow 0 $ as $\rho\rightarrow\infty$, then
$e^{|\delta|\rho/2 }|\lam|$ is also square integrable on $Y$.}
\end{itemize}
\item {If $\delta > 0$, then the images of
$\underline{\iota}_\delta$ and $\underline{\iota}_{-\delta}$ are
orthogonal, complementary subspaces in ${\rm kernel}(L_0)$.}
\end{itemize}\rm

\bigskip
With this last lemma understood, the proof of Lemma 3.3 is reduced to
verifying that the conditions in (3.9) and (3.10) are satisfied in the
present case.
\bigskip

{\bf Step 2}\qua This step in the proof verifies for certain cases
that the operator $D$ which appears in Lemma 3.3 obeys the conditions
in (3.10).  This is accomplished using the fact that when $s_0$ is
large, then $C$'s intersection with the $|s| \ge s_0$ part of $X-X_0$
is a disjoint union of cylinders where the constant $s$ slice of each
cylinder is very close to some closed Reeb orbit.  Moreover, according
to Proposition 2.3, as $|s| \rightarrow \infty$, this constant $|s|$
slice converges in the $C^k$ topology for any $k$ exponentially fast
in $|s|$ to some multiple covering of the closed Reeb orbit.  The
implications of this observation are slightly different depending on
whether the cylinder is in a component of $X-X_0$ which is concave or
convex and has or does not have an orientable $z$--axis line bundle.
In this regard, only the concave case is presented below as the
argument for the convex case is identical save for changing various
signs.  With $M$ now taken to define a concave end of $X$, the case
with orientable $z$--axis line bundle is considered in this step, and
the non-orientable $z$--axis line bundle case is considered in the
next.

Until directed otherwise, assume that $M = S^1 \times S^2
\subset\partial X_0$ defines a concave end of $X$ with orientable
$z$--axis line bundle.  Now, suppose that $\gamma\subset M$ is an
element in the limit set for $C$.  Here, there are three cases which
are treated separately.  The first case has $\theta = \theta_0$ on
$\gamma$ with $\theta_0 \notin \{0,\pi\}$ and with $\cos^2\theta_0
\not= 1/3$.  In this case, a component of $C$ which lies in a small
radius tubular neighborhood of $\gamma$ can be parameterized as
(2.13).  It then follows from (2.15) that $D$ has the form of (3.9)
where $A_0$ is the identity $2 \times 2$ matrix and where $L_0$ is the
operator $L$ in (2.16).  Moreover, the fact that $\lam$ in (2.15)
decays exponentially fast to zero implies that the requirements in
(3.10) are met in this case.

Next, consider the case where $\gamma$ is such that $\cos^2 \theta_0 =
1/3$.  In this case, the roles of $(t,f)$ and $(\varphi, h)$ can be
switched in the discussion of Step 9 of Proposition 2.3 to find a
parameterization of the part of $C$ near $[s_0, \infty)\times\gamma$
by coordinates $(\rho, \tau)$ for which $D$ has the form in (3.9)
where $A_0$ is the identity and $L_0$ is again given by (2.16), but
with $\zeta' = 0$ and $\zeta = 2$.  In this case, the coodinate $\rho$
takes values in an interval of the form $[\rho_0, \infty)$ while
$\tau$ takes values in $[0, 2\pi m]$ where $m$ is the number of times
each large, but constant $s$ slice of $C$ wraps around a given tubular
neighborhood of $\gamma$.  In this parameterization, $h = \epsilon
e^{-2\rho}$ with $\epsilon$ either 1 or $-1$, and $\varphi = \tau$.
Meanwhile, the coordinates $t$ and $f$ are parameterized as $t =
x(\rho,\tau)$, $f = 3 \epsilon e^{-2\rho}w(\rho,\tau)$ and where $x$
and $w$ are functions of $\rho$ and $\tau$ which are periodic in
$\tau$ and which decay to zero exponentially fast as $\rho$ tends to
infinity.  Note that a column vector with $x$ the top entry and $w$
the bottom obeys a differential equation with the schematic form of
(2.15).

The third case to consider is where $\theta_0$ is either 0 or $\pi$ on
$\gamma$.  The discussion here concerns solely the $\theta_0 = 0$ case
as the $\theta_0 = \pi$ case has an identical story modulo some
inconsequential sign changes.  In this case, use of the
parameterization in (2.14) finds that the operator $D$ has the form in
(3.9) using the version of $L_0$ in (2.16) that has $\zeta' = \zeta =
\sqrt{3}/\sqrt{2}$.  By the way, note that in this case, kernel$(L_0)
= 0$ because $\sqrt{6}$ is irrational.
\bigskip
   
{\bf Step 3}\qua This step considers the form of $D$ when the end of
$C$ in question lies in a concave end of $X$ with non-orientable
$z$--axis line bundle.  Here, there are two cases to consider; they
depend on whether the fundamental class of the corresponding element,
$\gamma$, of $C$'s limit set does not or does generate $H_1(S^1 \times
S^2; \mz)$.

In the case where the fundamental class of $\gamma$ does not generate
$H_1(S^1 \times S^2; \mz)$, then the inverse image of $\gamma$ via
the map in (3.3) has two distinct components which differ in the sign
of $\cos(\theta_0)$.  Likewise, the inverse image of the relevant
portion of $C$ via the map in (3.3) has two components, one near each
component of the inverse image of $[s_0, \infty) \times\gamma$.
Choose one such component and parameterize it as in the previous step.
The result gives an operator $D$ of the form in (3.9) which obeys the
constraints in (3.10).

Now consider the case where the fundamental class of $\gamma$ does
generate $H_1(S^1 \times S^2; \mz)$.  In this case, the inverse image
of $\gamma$ under the map in (3.3) is the circle, $\gamma_0$, with
$\theta_0 =\pi/2$ and $\varphi = 0$.  This is to say that the
coordinates $(x'_1, x'_2, x'_3)$ for the $S^2$ factor which appear in
(3.3) are either $(1, 0, 0)$ or $(-1, 0, 0)$ on $\gamma_0$.
Meanwhile, the corresponding inverse image of $C$ will be very close
to $[s_0, \infty) \times\gamma_0$.  In the case where this inverse
image has two components, choose one and parameterize it as in Step 2.
The resulting expression for $D$ will then have the form of (3.9) and
obey the constraints in (3.10).  In this regard, note that the inverse
image here will have two components precisely when the constant $s >
s_0$ circles in this end of $C$ are even multiples of a generator of
$H^1(M; \mz)$.

Finally, suppose that the inverse image of the relevant part of $C$
has just one component.  In this case, $C$ can be parameterized by
coordinates $\rho\in [\rho_0,\infty)$ and $\tau'\in [0, 2\pi m]$ where
$m$ is an odd, positive integer.  This parameterization writes the
coordinates $(s, t', x'_2, x'_3)$ for this end of $C$ as $(\rho_,
\tau', a', b')$ with $a'$ and $b'$ functions of $\rho$ and $\tau'$
which are periodic in $\tau'$ with period $2\pi m$.

With these last points understood, it follows that the inverse image
of $C$ under the map in (3.3) can be parameterized by the same
function $\rho$ and a function $\tau\in[0, 2\pi m]$ which writes the
coordinates $(s, t, x_2, x_3)$ as $(\rho, \tau, a, b)$ where $a$ and
$b$ are functions of $\rho$ and $\tau$ which obey
\begin{equation}
a(\rho, \tau + \pi m) = -a(\rho,\tau) \text{and}\quad b(\rho,
\tau + \pi m) = -b(\rho, \tau)\;.
\end{equation}
%(3.11)
Moreover, the column vector $\lam$ with top entry $a$ and bottom entry
$b$ obeys an equation with the schematic form given by (2.15) with
$\zeta' = 0$ and $\zeta = \sqrt{6}$.  Finally, $\lam$ decays
exponentially fast to zero as $\rho$ tends to infinity.  By the way,
note that the involution in (2.1) is realized on the inverse image
curve by sending $\tau$ to $\tau + \pi m$.
 
In any event, these last remarks imply that the operator $D$ on $C$
has the form given in (3.9) and satisfies the constraints in (3.10).
Indeed, this is because the pullback of $D$ to $C$'s inverse image
curve is as described in Step 3 for the $\theta_0 = \pi/2$
case. However, be forwarned that the domain and range Hilbert spaces
on $C$, the spaces $\cL_0$ and $\cL_1$ in Lemma 3.5, pullback to the
inverse image curve as subspaces of $\mr^2$ valued functions over the
inverse image curve which change sign when $\tau$ is changed to $\tau
+ \pi\; m$.  Here, $D$ on $C$ should be viewed as an operator on
$\mr^2$ valued functions by trivializing $C$'s normal bundle using the
restriction to $C$ of the vector fields which are tangent to the
$x'_2$ and $x'_3$ axis.

\sh{(c)\qua The proof of Proposition 3.2}

The strategy and most of the technical details for the proof follow
those for the proof in the case where $C$ is compact (see, for
example, \cite{MS}).  Certain special cases of the proposition also follow
from the analysis in \cite{HWZ3}.
  
The argument given below for Proposition 3.2 has two parts.  The first
finds a ball $B \subset \text{kernel}(D)$, a smooth map, $f$, from $B$
into the cokernel of $D$ and an embedding from $f^{-1}(0)$ into a
neighborhood of $C$ in $\fM_{e,\chi}$.  The second part of the
argument proves that this embedding of $f^{-1}(0)$ into $\fM_{e,\chi}$
is onto an open set.  Both parts consist mostly of fairly
straightforward generalization of arguments that are used in the
compact case and in \cite{HWZ3}.  Thus, the discussion below will be brief,
with many of the details left for the reader.  In this regard, the
discussion that follows will consider only the case where $C$ is
immersed with purely transversal double point self intersections; the
general case is left as an exercise.  In any event, these two parts to
the proof of Proposition 3.2 constitute Steps 1--6 and Step 7 of the
seven steps into which the proof below is divided.
\bigskip

{\bf Step 1}\qua What follows is a brief summary of the formal set up
for the first part of the proof of Proposition 3.2 (as described
above).  In this regard, note that the basic conclusions here in Step
1 follow more or less automatically from an application of the
implicit function theorem.  However, there are two subtle points in
this application.  The first such point involves the choice of the
appropriate spaces and the map between them.  The second subtle point
involves a reference to a particular regularity theorem in \cite{Mo}.

To start the summary, the basic observation is that a constant
$\epsilon_0 > 0$ and an `exponential' map, $q$, from the radius
$\epsilon_0$ disk bundle $N_0\subset N$ to $X$ can be found with the
following properties:
\begin{itemize}\sl
\item $q$ {restricts to the zero section as $\varphi$ and embeds
each fiber of $N_0$ as a pseudoholomorphic disk.}
\item {At each point along the zero section, $q$'s differential
maps $T_N$ isomorphically to $T_X$.}
\item {If $\lam$ is a section of $N_0$, then the image of
$q(\lam(\cdot))$ is a pseudoholomorphic submanifold of $X$ if and only
if $\lam$ obeys an equation with the schematic form}
$$
D\lam + \cR(\lam, \nabla\lam) = 0 \;,
$$ 
{where $\cR$ is a smooth, fiber preserving map from $N_0\oplus
(N\oplus T^*C)$ to $N\oplus T^{0,1}C$ that is affine in the second
factor and obeys}
$$
|\cR(a, b)| \le \zeta (|a|^2 + |a| \;|b|)
$$
{for some constant $\zeta$.}
\item {In fact, $q$ can be constructed so that $\cR(a, b) =
\cR_0(a) + \cR_1(a) b^{1,0}$ where $\cR_0$ and $\cR_1$ respectively
map $N_0$ to $N \otimes T^{0,1}C$ and to ${\rm Hom}_\mc(T^{1,0}C; N
\otimes T^{0,1}C)$, and where $b^{1,0}$ denotes the projection of $b$
onto the $(1,0)$ summand in $T^*C_{\mc}$.}\autonum
\end{itemize}

The map $q$ can be constructed as described in Section 2 of \cite{T6}.
Indeed, the basic point is that $q$ maps the fibers of $N_0$ into $X$
as pseudoholomorphic disks.  This understood, then a local,
$\mc$--valued fiber coordinate $\eta$ on $N_0$ can be found with
respect to which the $q$ pullback of $T^{1,0}X$ is spanned by a form
that annihilates vertical vectors and $e = d\eta + \alpha$ where
$\alpha$ vanishes on the zero section and also annihilates vertical
vectors.  The various points in (3.12) then follow by exploiting these
last remarks.

With (3.12) understood, the goal is to use the implicit function
theorem coupled with various standard elliptic regularity theorems to
describe the small solutions to (3.12).  To broaden the subsequent
discussion so as to cover certain generalizations of (3.12), it proves
useful to consider the form of $\cR$ in the third point of (3.12)
without assuming the validity of the final point.  Thus, the
subsequent discussion makes no reference to this final point of
(3.12).  It also proves useful to set up the implicit function theorem
in a weighted Sobolev space of sections of $N$ where finite norm
demands local square integrability of the section and its covariant
derivative, but makes no demand for an $L^\infty$ bound.  Were the
purpose solely that of proving Proposition 3.2, a norm with $L^\infty$
implications can be used.  In any event, with the lack of an
$L^\infty$ bound from the norm, a trick is employed to handle the
nonlinearity in $\cR$.
  
Here is the trick: First, introduce a smooth function $\kappa$: $[0,
\infty) \rightarrow [0, 2)$ whose value at $t \in [0, \infty)$ is 1
when $t < 1$, 1$/t$ when $t > 2$.  Given $\epsilon > 0$, set
$\kappa_\epsilon$: $[0, \infty) \rightarrow [0, 2\epsilon)$ to denote
the function whose value at $t$ is $\kappa(t/\epsilon)$.  Now, fix
$\epsilon \ll 1$ and the trick is to replace the given $\cR$ in the
third point of (2.14) with $\cR_\epsilon\co N_0 \oplus(N\otimes T^*C)
\rightarrow N\otimes T^{0,1}C$ whose value on $(a, b)$ is
$\cR(\kappa_\epsilon(|a|) a, b)$.  Solutions to the $\cR_\epsilon$
version of (3.12) are then found via a straightforward implicit
function theorem argument.  This understood, and given that $\epsilon$
is small, Theorem 5.4.1 in \cite{Mo} guarantees that these solutions are
pointwise smaller than $\epsilon$ over the whole of $C$, and so they
solve the desired $\cR$ version of (3.12).

Solutions to the $\cR_\epsilon$ version of (3.12) are found in the
following way: A certain space of sections of $N$ is split as
kernel$(D) \oplus\cL$ while a corresponding space of sections of
$N\otimes T^{0,1}C_0$ is split as cokernel($D)\oplus\cL'$.  Then, for
each $\lam$ in kernel$(D) \oplus\cL$, the section $D\lam +
\cR_\epsilon(\lam, \nabla\lam)$ of $N\otimes T^{0,1}C_0$ is projected
into $\cL'$ to define a map from kernel$(D) \oplus \cL$ to $\cL'$.
The differential of this last map at $\lam = 0$ (which is formally
$D$) can be seen to identify $\cL$ with $\cL'$.  Thus, the implicit
function theorem gives a ball $B$ about the origin in kernel($D$) and
a smooth map, $\Phi\co B\rightarrow \cL$ with the following properties:
\begin{itemize}\sl
\item {When $\eta\in B$, then $\Phi(\eta) = \cO(|| \eta ||^2)$.}
\item {When $\eta\in B$, set $\lam\equiv\eta + \Phi(\eta)$ and
then $D\lam + \cR_\epsilon(\lam, \nabla\lam)$ projects to zero in
$\cL'$.}
\autonum%(3.13)
\end{itemize}

This ball $B$ is the one that used in Proposition 3.2.

With (3.13) understood, define a map $f\co B\rightarrow 
\text{cokernel}(D)$ by taking the projection onto cokernel($D$) of
$\cR_\epsilon(\eta + \Phi(\eta)$, $\nabla(\eta + \Phi(\eta))$.  This
is the map to be used in Proposition 3.2.  Indeed, if $|\eta +
\Phi(\eta)|$ is everywhere less than $\epsilon$, then it follows from
the third point in (3.12) and the second in (3.13) that the image in
$X$ of the map $q(\eta + \Phi(\eta))$: $C_0 \rightarrow X$ is a
pseudoholomorphic submanifold in $\cM_{e,\chi}$ near to $C$ if and
only if $f(\eta) = 0$.  Meanwhile, standard elliptic estimates for $D$
find some $\epsilon_D > 0$ such that the following is true: If
$\epsilon < \epsilon_D$ and if $B$ is defined so that its elements
have supremum norm less than $\epsilon/2$, then $|\eta + \Phi(\eta)|$
is guaranteed to have supremum norm less than $\epsilon$ and so $\lam=
\eta + \Phi(\eta)$ solves $D\lam +\cR(\lam, \nabla\lam) = 0$.  Thus,
$f^{-1}(0)$ is truly mapped to a neighborhood of $C$ in
$\cM_{e,\chi}$.  Conversely, as $\Phi$ is smooth, some relatively
standard $s \rightarrow\infty$ estimates on the sections of $N$ in
kernel$(D) \oplus\cL$ guarantee that the map just described embeds
$f^{-1}(0)$ in $\cM_{e,\chi}$.
\bigskip

{\bf Step 2}\qua This step introduces the important function spaces
involved.  For this purpose, fix $\delta > 0$, but very small and fix
$k \in \{0, 1,\ldots\}$.  Let $E$ denote either $N$, $N\otimes
T^{0,1}C_0$ or the tensor product of $N$ with some multiple tensor
product of $T^*C_0$.  Let $L_{k,\delta}^2(E)$ denote the completion of
the set of smooth, compactly supported sections of $E$ using the norm
whose square is
\begin{equation}
|| \lam ||_{k,\delta}^2 \equiv \int e^{\delta |s|} \sum_{0\le p\le k}
|\nabla^k\lam|^2 \; .
\end{equation}
%(3.14)
Note that this norm uses the growing exponential $e^{\delta s}$ while
those in (2.8) and (2.9) use the shrinking exponential $e^{-\delta
s}$.

The following is a basic fact about these spaces which is left to the
reader to verify:

\narrower\sl The covariant derivative extends to a bounded map from\nl 
$L_{1,\delta}^2$ {\it to} $L_{0,\delta}^2$.  \autonumm
\endnarrower

{\bf Step 3}\qua Reintroduce the operator $D$ from (3.5).  On each end
of $C$, this operator has the schematic form in (3.9) and thus
associated to each end of $C$ is the kernel of the relevant version of
the operator $L_0$.  Note that each of these kernels is either zero or
one-dimensional.  In any event, let $W$ denote the direct sum (indexed
by the ends of $C$) of these kernel($L_0$) vector spaces.

Fix $\delta > 0$, but very small and, for the moment, consider $D$
with its range and domain as described in Lemma 3.3.  Lemma 3.5
provides the homomorphism $\iota_\delta$: kernel$(D)\rightarrow W$.
Note that the kernel of $\iota_\delta$ is the intersection of
kernel($D$) with $L_{0,\delta}^2(N)$ and thus equal to the kernel of
this same operator $D$ but viewed as a map from $L_{1,\delta}^2(N)$ to
$L_{0,\delta}^2 (N\otimes T^{0,1}C_0)$.  Use $K$ to denote the kernel
of the latter version of $D$.  Thus, kernel($D$) can be split as
\begin{equation}
\text{kernel}(D) = \text{image}(\iota_\delta) \oplus K \; .
\end{equation}
%(3.16)  
Viewed as mapping $L_{1,\delta}^2(N)$ to $L_{0,\delta}^2(N\otimes
T^{0,1}C_0)$, the operator $D$ is Fredholm, a fact implied by Lemma
3.5.  Use $K_*\subset L_{0,\delta}^2(N\otimes T^{0,1}C_0)$ to denote
the corresponding cokernel.  In this regard, note that the domain for
Lemma 3.5's map $\underline{\iota}_\delta$ is this vector space $K_*$.
And, according to the second point of Lemma 3.5, the space $K_*$ can
be split as
\begin{equation}
K_* = \text{image}(\underline{\iota}_\delta) \oplus e^{-2\delta
|s|}\text{ kernel}(D^*)\; .
\end{equation}
%(3.17)
The latter splitting is made in a canonical way by requiring the two
summands to be orthogonal with respect to the inner product on
$L_{0,\delta}^2(N\otimes T^{0,1}C_0)$.
\bigskip

{\bf Step 4}\qua The purpose of this step is to identify $W$ as a
subspace of sections of $N$ over $C_0$.  For this purpose, use Lemma
3.5 to decompose $W = \text{image}(\iota_\delta) \oplus
\text{image}(\underline{\iota}_\delta)$.  Then, identify the
image$(\iota_\delta)$ summand of $W$ with the corresponding subspace
of kernel($D$) in (3.16).  (Note that the splitting in (2.18) is not
canonical.)

To identify the image$(\underline{\iota}_\delta)$ summand of $W$ as
space of sections of $N$, first use the coordinates on each end of $C$
which are described in the proof of Lemma 3.3 to view
image$(\underline{\iota}_\delta)$ as a subspace of sections of $N$
over the $s > s_0$ portion of $C$.  In this regard, note that each
such section of $N$ is pointwise bounded.  Next, choose a smooth
function, $\beta$, on $\mr$ which has values between 0 and 1, and
which is zero on $(-\infty, 2s_0)$ and which is one on $(3s_0,
\infty)$.  Compose $\beta$ with the function $s$ on $C_0$ to view
$\beta$ as a function on $C_0$.  Finally, embed
image$(\underline{\iota}_\delta)$ as a subspace of sections of $N$
over the whole of $C_0$ by sending $w \in
\text{image}(\underline{\iota}_\delta)$ to $\beta w$.  Note that with
the latter embedding understood, the vector space
image$(\underline{\iota}_\delta)\subset W$ has now been identified
both as a subspace of sections of $N$ and also, via (3.17), as a
subspace of sections of $N\otimes T^{0,1}C_0$.  These two versions of
image$(\underline{\iota}_\delta)$ will not be notationally
distinguished.
\bigskip

{\bf Step 5}\qua Let $K^1\subset L_{1,\delta}^2 (N)$ denote the
orthogonal complement of $K$ as defined using the inner product on
$L_{0,\delta}^2(N)$.  (Note that this is meant to be the analog of the
$L^2$ inner product as opposed to that of the $L_1^2$ inner product.)
Likewise, define $K_*^\perp\subset L_{0,\delta}^2(N\otimes
T^{0,1}(C_0)$.  Also, introduce the corresponding
$L_{0,\delta}^2$--orthogonal projection $\prod$:
$L_{0,\delta}^2(N\otimes T^{0,1}(C_0) \rightarrow K_*^\perp$.

The discussion in Step 1 referred to spaces $\cL$ and $\cL'$.  In this
regard, define $\cL$ to be
$\cL\equiv\text{image}(\underline{\iota}_\delta)\otimes K^\perp$ and
define $\cL'$ to be
$\cL'\equiv\text{image}(\underline{\iota}_\delta)\oplus K_*^\perp$.
Note that $\cL'$ sits in $L_{0,\delta}^2(N\otimes T^{0,1}(C_0)$, and
let $\prod$: $L_{0,\delta}^2(N\otimes T^{0,1}(C_0)\rightarrow \cL'$
denote the $L_{0,\delta}^2$--orthogonal projection.
\bigskip
  
{\bf Step 6}\qua With $\epsilon > 0$, but small, define a map from
kernel$(D) \oplus\cL$ to $\cL'$ by sending $(\eta,w)$ in the former
space to
\begin{equation}
\prod (Dw + \cR_\epsilon(\eta+ w,\nabla(\eta+ w))
\end{equation}
%(3.18)
It is a straightforward task (which is left to the reader) to check
that this map has surjective differential at $(0, 0)$ along the $\cL$
summand of kernel$(D)\oplus \cL$.  This understood, then the implicit
function theorem finds a ball $B \subset \text{kernel}(D)$ and a
smooth map $\Phi$: $B \rightarrow \cL$ such that when $\eta\in B$,
then $(\eta, \Phi(\eta))$ solves (3.18).  Moreover, the implicit
function theorem also finds $\epsilon_1>0$ such that any pair $(\eta,
w)$ solving (3.18) with $\eta\in B$ and $||w||_{1,\delta} <
\epsilon_1$ has $w =\Phi(\eta)$.  This understood, define $f$:
$B\rightarrow \text{cokernel}(D)$ by sending $\eta$ to $f(\eta) \equiv
(1 - \prod)(D\Phi(\eta) + \cR_\epsilon(\eta + \Phi(\eta), \eta +
\Phi(\eta))$.  Thus, $\lam\equiv\eta + \Phi(\eta)$ solves the
$\cR_\epsilon$ version of (3.12) when $\eta\in f^{-1}(0)$.
  
Here is one other automatic consequence of the implicit function
theorem: The graph of $\Phi$ in kernel$(D) \oplus\cL$ is homeomorphic
to a neighborhood of $(0, 0)$ in the space of solutions to (3.18).
Thus, $f^{-1}(0)$ is homeomorphic to a neighborhood of 0 in the space
of $(\eta, w)\in \text{kernel}(D) \oplus\cL$ for which
$\lam\equiv\eta+\Phi(\eta)$ obeys the $\cR_\epsilon$ version of
(3.18).
\bigskip

{\bf Step 7}\qua Theorem 5.4.1 in \cite{Mo} now finds some $\epsilon_3 \in
(0, \epsilon)$ such that when $\eta$ lies in the centered, radius
$\epsilon_3$ ball in $B$, then $\lam\equiv\eta +\Phi(\eta)$ is
everywhere bounded in norm by $\epsilon$.  Thus, $\lam$ solves the
$\cR$ version of (3.12) and, as a consequence, the map that sends
$\eta\in f^{-1}(0)$ to $\eta+\Phi(\eta)$ properly embeds $f^{-1}(0)$
in a neighborhood of $C$ in $\fM_{e,\chi}$.
 	
Meanwhile, an application of Proposition 2.3 finds $\epsilon_2 \in(0,
\epsilon)$ such that if $\lam\in\cL_1$ obeys the original, $\cR$
version of (3.12) and has $|\lam| < \epsilon_2$ everywhere, then
$\lam= \eta + \Phi(\eta)$ with $\eta\in f^{-1}(0)$.  Thus, the
aforementioned map from $f^{-1}(0)$ to $\fM_{e,\chi}$ is a
homeomorphism onto an open neighborhood of $C$.

\sh{(d)\qua The index of $D$}

It is useful to have a formula for the index of $D$ that generalizes
those in (3.6).  Such a formula is given in Proposition 3.6 at the end
of this subsection.  The derivation of the analog of (3.6) is a six
step affair which follows this preamble.  However, before starting,
note that the concave and convex end discussions are not identical,
but are completely analogous.  Thus, the discussion below focuses on
the concave case while only summarizing the modifications that are
required when convex ends are present.  In particular, the first four
steps below involve only the case where $X$ has solely concave ends.
Note here that this assumption in Steps 1--4 about $X$ is mostly
implicit.  Also, except for Step 6, the operating assumption is that
$C$ is immersed in $X$ with only transversal, double point
self-intersections.

Before starting, it is pertinent to remark on the absence in the
subsequent discussion of reference to a Maslov or Conley-Zehnder index
as in \cite{HWZ3}.  The first point is that the index formula for a
Cauchy-Riemann operator on some given completion of the space of
sections of a bundle over a surface must have the following schematic
form: Index = 2 degree + $\chi$ + `boundary correction'.  Here, $\chi$
is the Euler characteristic of the surface, degree is a first Chern
number of the bundle in question and `boundary correction' is just
what it says.  In this regard, the definition of the degree requires
some choice of trivialization of the bundle along the ends of the
surface and a different choice might change the degree.  Of course, if
it does, it will also change the `boundary correction' term to
compensate.  In this regard, an index formula with a formal
Conley-Zehnder index term correction to the basic `2 degree + $\chi$'
simply reflects a particularly natural choice for the trivialization
of the normal bundle of the pseudoholomorphic surface on its ends.  In
particular, these Conley-Zehnder trivializations are trivializations
that are induced from certain trivializations of the normal bundles to
the limiting closed Reeb orbits.
  
The preceding understood, remark that the formula given below in
Proposition 3.6 also uses certain natural trivializations of the
normal bundle to the surface on its ends that are induced from those
of the limiting closed Reeb orbits.  In fact, the latter were
described previously in Part 2 of Section 3a.  However, the formal
introduction of a `Conley-Zehnder' index in this case is not done here
for three reasons.  First, an integer valued Conley-Zehnder index is
awkward to define in the present circumstances because the first Chern
class of $T^{1,0}X$ evaluates with absolute value 2 on the $S^2$
factors of the ends of $X$.  Second, the index formula for the
non-isolated Reeb orbits necessarily has terms with no analog in the
index formula \cite{HWZ3}.  Third, the formula given in Proposition 3.6
are, in any event, reasonably easy to use without the additional
burden of a Conley-Zehnder index computation.
\bigskip

{\bf Step 1}\qua Consider first an end of $C$ which lies in an end of
$X$ with orientable $z$--axis line bundle.  By virtue of Proposition
2.3, the large $s$ part of the corresponding end of $C_0$ has a
product structure with special coordinates $(\tau,\rho)$ such that $J$
maps $\partial_\tau$ to $\pm\partial_\rho$.  Here, $\tau\in\mr/(2\pi
m|p|\mz)$ while $\rho\in [\rho_0,\infty)$ with $m > 0$, $p \not= 0$
and both integers.  In this regard, the integer $m$ should be viewed
as the multiplicity by which the large and constant $|s|$ slices of
this end of $C$ wrap around the limiting Reeb orbit.  This
interpretation of $m$ requires taking $p = -1$ when the limiting Reeb
orbit for this end is characterized as in (1.8) by one of the symbols
$+$ or $-$.  Note that these coordinates differ from those in (2.13)
or its $\cos^2 \theta_0 = 1/3$ analog; however, they differ only by
functions whose derivatives to any order decay to zero at an
exponential rate as $s \rightarrow\infty$.
 
There is also a special trivialization of the normal bundle over the
end in question as the product $\mr^2$ bundle for which $J$ sends the
column vector $\left(\begin{smallmatrix} 1\\0\end{smallmatrix}\right)$
to $\pm\left(\begin{smallmatrix}0\\1\end{smallmatrix}\right)$ With
respect to this trivialization, the operator $D$ defines an operator
on the $\mr^2$ valued functions that has the form
\begin{equation}
D\lam = \partial_\rho\lam + L_0\lam + \cP\lam\; .
\end{equation}
%(3.19)
Here, $\lam$ is a 2--component column vector, $L_0$ is described by
(2.16) and $\cP$ is a $2 \times 2$ matrix valued function that is
small and drops to zero exponentially fast as $\rho\rightarrow\infty$.
Meanwhile, $\zeta'$ and $\zeta$ are the constants in (2.16), and thus
one of the following hold:
\begin{itemize}
\item $\zeta' = 0$ {and} $\zeta > 0$.
\item $\zeta' = \zeta > 0$.
\autonum%(3.20)
\end{itemize}
In this regard, the second case occurs only when the given end has its
corresponding $\theta_0$ in $\{0, \pi\}$.
 
Because the coefficients of $\cP$ drop to zero exponentially fast at
large $\rho$, the index of $D$ is the same as that of $D_0 \equiv D -
\beta(\rho)\cP$, where $\beta(\rho)$: $[\rho_0, \infty) \rightarrow
[0, 1]$ is 0 on a certain interval of the form $[\rho_0, \rho_1]$ and
$\beta(\rho)$ is 1 on $[2\rho_1, \infty)$.  Here, $\rho_1 > \rho_0$
can be taken to be very large.

This operator $D_0$ is introduced because its coefficients are
constant on $[2\rho_1,\infty)$ $\times S^1$.  In particular, this
implies (see \cite{APS}) that $D_0$ has the same index as the identical
operator on the domain in $C_0$ where $\rho \le 4\rho_1$ but with
appropriate spectral boundary conditions imposed on the $\rho =
4\rho_1$ boundary of this domain.  By way of motivation, this
translation to a spectral boundary condition problem is made so as to
facilitate the comparison of the index of $D_0$ with that of a certain
complex linear operator with known index.

In any event, there is one key observation that is used to determine
the appropriate spectral boundary conditions on the $\rho = 4\rho_1$
circle: Any bounded element in the kernel of $D_0$ must restrict to
the half cylinder $S^1 \times [2\rho_1, \infty)$ in the span of the
set
\begin{equation}
\left\{e^{-E\rho}\eta_E(\tau):E\ge 0\quad\text{and}\quad
L_0\eta_E=E\eta_E\right\}\;.
\end{equation}
%(3.21)
The eigenvalues $E$ which appear above have the form
\begin{equation}
E = -2^{-1} (\zeta' + \zeta - [(\zeta'« - \zeta)^2 + 4
n^2/(m|p|)^2]^{1/2}) \;.
\end{equation}
%(3.22)
Here, $n \in \{0, 1, \ldots\}$ is such as to make $E \ge 0$.  In
particular, with regard to this last concern, note that when the first
point of (3.20) holds, then $E > 0$ for all positive $n$ and $E = 0$
for $n = 0$.  On the other hand, when the second point of (3.20)
holds, then $E > 0$ if and only if $n > m \sqrt{3}/\sqrt{2}$ because
$\zeta'= \zeta = \sqrt{3}/\sqrt{2}$ and $|p| = 1$ in this case.

Concerning the eigenspace for the eigenvalue $E$, it is important to
note that this is a two-dimensional vector space over $\mr$ if the
integer $n$ in (3.22) is strictly positive.  Such is also the case
when $n = 0$ and the second point in (3.20) holds.  However, when
$\zeta' = 0$ in (3.20), then each of the $n = 0$ values of $E$ has
1--dimensional eigenspace.  In any case, the components of the
eigenfunctions $\eta_E$ in (3.21) are certain linear combinations of
$\cos(n\tau/(m|p|))$ and $\sin(n\tau/(m|p|))$.
 
It now follows directly from (3.21) that the index of $D_0$ is the
same as that of the identical operator on the domain where $\rho\le
4\rho_1$ in $C_0$ with spectral boundary conditions on the $\rho =
4\rho_1$ boundary of this domain which restrict the sections of the
normal bundle under consideration to lie in the set spanned by
[$\eta_E: E \ge 0$ is given by (3.22)].
\bigskip
	
{\bf Step 2}\qua Here and in Step 3, assume that $C$ only intersects
components of $X-X_0$ that have orientable $z$--axis line bundle.  With
this assumption understood, the remainder of this step considers not
$D_0$, but a particular $\mc$--linear operator version of the $\opart$
operator on the sections of the normal bundle.  The reason for this
digression should be apparent by the end of the next step.

The first order of business is to specify the version of $\opart$ to
be used.  In particular, the operator is constrained only on the ends
of $C_0$; and on an end, in terms of the coordinates $(\rho, \tau)$
and the aforementioned trivialization of the normal bundle, this
operator sends an $\mr^2$ valued function $\lam$ to
\begin{equation}
\opart\lam\equiv\lam_\rho + \left(
\begin{matrix}
0&-\partial_\tau \\ \partial_\tau&0
\end{matrix}
\right) \lam\;.
\end{equation}  
%(3.23)

By the way, the difference of two $\opart$ operators for the same
complex line bundle is, as a $\mc$--linear operator, the tensor product
with a section of $T^{0,1}C_0$.  Were, $C_0$ compact, then the two
operators would have the same index, but such an event is not
guaranteed when $C_0$ is non-compact.
	
The index of the version of $\opart$ in (3.23) is readily computed
with the help of the Riemann-Roch theorem, and the computation yields
the identity: index$(\opart) = 2\cdot \text{degree}(N) + \chi(C_0) +
\aleph_T$.  Here, degree$(N)$, $\chi(C_0)$ and $\aleph_T$ denote the
following: First, degree$(N)$ is the degree of the line bundle $N$ as
determined via the sections over the end of $C_0$ which are described
in Part 2 of Section 3a; meanwhile $\chi$ and $\aleph_T$ respectively
denote the Euler characteristic and the number of ends of $C_0$.
  
The formula just described for the index of $\opart$ can be rewritten
using Proposition 2.1 and the fact that degree$(N) = \la e, [C]\ra - 2
m_C$.  For example, here is an equivalent formula: index($\opart ) =
\la e, [C]\ra - \la C_1, [C]\ra - 2 m_C + \aleph_T$.

Note that the kernel of $\opart$ on an end of $C_0$ can also be
written as in (3.21) except that $E = n/(m|p|)$ replaces the formula
for $E$ in (3.22).  Here, $n$ can be any non-negative integer.  For
the case of $\opart$, all of the eigenspaces are two-dimensional.

As a constant coefficient operator on the ends of $C_0$, the index of
the $\opart$ operator can also be identified with the index of the
same operator on the $\rho \le 4\rho_1$ domain in $C_0$ with spectral
boundary conditions on each $\rho = 4\rho_1$ circle.  In particular,
these boundary conditions restrict the sections of $N$ on the $\rho =
4\rho_1$ circle to be linear combinations from those column vectors
with top component $\cos(n\tau/(m|p|))$ and bottom component
$\sin(n\tau/m|p|)$ or else top component $-\sin(n\tau/(m|p|))$ and
bottom component $\cos(n\tau/(m|p|))$.  Here, $n\in \{0, 1, \ldots\}$.
\bigskip

{\bf Step 3}\qua This step compares the spectral boundary conditions
for $D_0$ and for $\opart$.  In this regard, note that for both
operators, the relevant boundary condition for each boundary component
restricts the sections of $N$ to lie in a certain direct sum of finite
dimensional spaces.  In the case of $\opart$ , these summands are
indexed by integers $n \in\{0, 1, 2, \ldots\}$ and each is
two-dimensional.  Meanwhile, the summands in the case of $D_0$ are
subspaces of those for $\opart$.

In particular, on an end of $C_0$ for which the corresponding element
in $C$'s limit set has $\theta_0\notin\{0,\pi\}$, {the $n > 0$
summands for $D_0$ are the same as those for $\opart$, while the $n =
0$ summand for $D_0$ is a 1--dimensional subspace of that for $\opart$.
On the other hand, if the corresponding element in $C$'s limit set has
$\theta_0 \in \{0, \pi\}$, then there is some integer $m_0 \ge 1$ such
that the $n \ge m_0$ summands for $D_0$ and for $\opart$ agree, while
the $n < m_0$ summands for $D_0$ are all trivial (zero-dimensional).
Here, $m_0$ is the least integer that is greater than $m
\sqrt{3}/\sqrt{2}$ with $m$ defined as follows: The large and constant
$s$ slices of the end in question define a class in $H_1(S^1 \times
S^2; \mz)$.  Then $m$ is the absolute value of the pairing between
this class and a generator of $H^1(S^1 \times S^2; \mz)$.
 
With this last paragraph understood, it now follows from the analysis
in \cite{APS} that the index of $D_0$ is less than that of $\opart$, with
the deficit, $\Theta$, accounted for by a contribution from each end
of $C_0$.  Here, an end $\cE\subset C_0$ with $\theta_0 \notin\{0,
\pi\}$ contributes 1 to $\Theta$, while one with $\theta_0 \in \{0,
\pi\}$ has the non-zero contribution $2 m_0(\cE)$.  Thus, with
$\aleph_0$ denoting the sum, indexed by the ends $\cE\subset C_0$ with
$\theta_0 \in \{0, \pi\}$, of the corresponding quantities $(1 - 2
m_0(\cE))$, the following formulae hold: index$(D) =
2\cdot\text{degree}(N) + \chi(C_0) +\aleph_0$.  Equivalently,
index$(D) = \la e_e, [C]\ra - \la c_1, [C]\ra - 2 m_C +\aleph_0$, and
also index$(D) = -\chi(C_0) - 2 \la c_1, [C]\ra + \aleph_0$.
\bigskip

{\bf Step 4}\qua Now consider the possibility that $C_0$ has some ends
which lie in $[0, \infty) \times M$ where $M \subset\partial X_0$ is a
component with non-orientable $z$--axis line bundle.  In this case, the
index of $D$ is computed via the route just used, first through a
deformation on the ends of $C_0$ to a constant coefficient operator,
next through a reinterpretation of the latter as an operator with
spectral boundary conditions, and finally via a comparison with the
index of the analogous spectral boundary condition interpretation of a
certain $\opart$ operator.
 
With regard to this $\opart$ operator, note that the latter will be
Fredholm with index$(\opart ) = 2\cdot \text{degree}(N) + \chi(C_0) +
\aleph_T$.  (As before, index$(\opart) = \la e, [C]\ra - \la C_1,
[C]\ra- 2 m_C + \aleph_T$ is an equivalent formula.)  The $\opart$
operator used here is, once again, restricted on the ends of $C_0$.
The restrictions described above for $\opart$ are used for the ends of
$C_0$ which map to components of $X-X_0$ where the corresponding
$z$--axis line is orientable.  The restrictions on the remaining ends
of $C_0$ are described below.

In any event, the resulting analysis finds index($D$) = index$(\opart)
- \Theta$, where $\Theta$ is again a sum of contributions from each
end of $C_0$.  An end of $C_0$ which maps to a component of $X-X_0$
where the corresponding $z$--axis line bundle is orientable makes the
same contribution to $\Theta$ as before.  Meanwhile, the description
of this contribution for an end of $C_0$ which lies where the $z$--axis
line bundle is non-orientable requires the distinction between two
cases.  The first case discusses those ends where the constant $s >
s_0$ circles define even multiples of a generator of $H^1(M; \mz)$.
Note that this case occurs automatically when the corresponding
element of $C$'s limit set is not a generator of $H^1(M; \mz)$.  The
second case discusses the situation when these circles are odd
multiples of a generator.

In the first case, the inverse image in $[s_0,\infty) \times (S^1
\times S^2)$ of the end in question under the map in (2.3) has two
components, and so determines a particular value of $\theta_0$ up to
the replacement $\theta_0 \rightarrow \pi - \theta_0$.  Moreover,
after the identification of the end with either of these components,
the analysis for the orientable $z$--axis line bundle case repeats here
hitch free.  In this regard, please note that the $\opart$ operator
for use on $C_0$ is restricted so as to yield (3.23) after the
identification of the given end with one of its inverse images in
$[s_0, \infty) \times (S^1\times S^2)$
  
With the preceding understood, here is the result of the analysis: The
end in question contributes 1 to the deficit $\Theta$ if the
corresponding $\theta_0 \notin \{0, \pi\}$, and otherwise the
contribution is twice the value of the associated integer $m_0$.

Now turn to the case of an end of $C_0$ where the corresponding large,
constant $|s|$ circles is an odd multiple of a generator of $H_1(M;
\mz)$.  Let $C_{01} \subset [s_0, \infty) \times (S^1 \times S^2)$
denote the inverse image of $C_0$ via the map in (3.3).  As noted
earlier, $C_{01}$ is a pseudoholomorphic cylinder with
$\theta_0=\pi/2$ and with $\varphi$ asymptotic as $s
\rightarrow\infty$ to either 0 or $\pi$.  Thus, $C_{01}$ has the
previously described parameterization by coordinates $(\rho, \tau) \in
[\rho_0,\infty) \times [0, 2\pi m]$ where $m$ is an odd, positive
integer.  In this regard, note that the assumptions in play here force
the relevant element in the limit set of $C_{01}$ to have $\theta_0 =
\pi/2$, and this implies that the integer $p$ is equal to one.  Use
coordinates $(\rho', \tau') \in [\rho_0, \infty) \times[0, 2\pi m]$
for $C_0$ so that the map in (3.3) restricts to $C_{01}$ as that which
writes $\rho'= \rho$ and $\tau' = 2 \tau$.
  
Under the previously described identification between the normal
bundle of $C_{01}$ and the pull back via (3.3) of the normal bundle of
$C_0$, sections of $C_0$'s normal bundle pull up as $\mr^2$ valued
functions on $C'_0$ which obey $\lam(\rho, \tau + m \pi) = -\lam(\rho,
\tau)$.  This follows from (3.3) and the comments in the final
paragraph of Part 2 of Section 3a.  The operator $D$ on $C_0$ then
pulls up to be the corresponding $D$ on $C_{01}$, but with the domain
and range restricted to the functions which are odd under the
involution $\tau\rightarrow\tau + m \pi$.
 
Given these last points, the analysis done for the other ends of $C_0$
can be repeated to find that the index of $D$ on $C_0$ is the same as
that of this differential operator on the $\rho\le 4\rho_1$ portion of
$C_0$ with spectral boundary conditions given as follows: The sections
on the $\rho = 4\rho_1$ boundary should pull up to $C_{01}$ as a
linear combination of two column vectors where one has top component
$\cos(n\tau/m)$ and bottom component $\sin(n\tau/m)$ while the other
has top component $-\sin(n \tau/m)$, and bottom component $\cos(n
\tau/m)$.  However, here $n$ must be a positive, odd integer.

Now consider the constraint for the $\opart$ operator on this end of
$C_0$.  To guarantee that this operator has index equal to 2
degree$(N) + \chi(C_0) - \aleph_T$, the constraint is as follows: Use
the coordinates $(\rho', \tau')$ on the end of $C_0$ and trivialize
the normal bundle over this end using the vectors tangent to the
$x'_2$ and $x'_3$ coordinates which appear in (3.3).  With these
coordinates and the bundle trivialization understood, the $\opart$
operator acts on $\mr^2$ valued functions.  In this guise, the
operator should be that which sends a column vector function $\lam$ to
$2^{-1}\lam_{\rho'} + L'$ where $L'$ is the matrix operator given by
the second term in (3.23) with $\tau'$ derivates replacing those of
$\tau$.

To compare this $\opart$ operator with $D$, it proves convenient to
pull the former up to $C_{01}$ and express it as an operator on the
normal bundle to the latter.  In particular, the resulting operator
differs from the operator in (3.23) by the zero'th order operator
which sends $\lam$ to the operator which is given by (3.23) but with
$\partial_\tau$ replaced by $\partial_\tau + 1$.  As with the pullback
to $C_{01}$ of $D$, the range and domain of this operator must be
restricted to those $\mr^2$ valued functions $\lam$ which obey
$\lam(\rho, \tau + m \pi) = - \lam(\rho, \tau)$.
  
With the preceding understood, the analysis previously done can be
repeated to find that the index of $\opart$ on $C_0$ is the same as
that of the same operator on the $\rho\le 4\rho_1$ portion of $C_0$
with spectral boundary conditions given as follows: The sections on
the $\rho = 4\rho_1$ boundary should pull up to $C_{01}$ as linear
combinations of two column vectors, one whose top component is
$\cos(n\tau/m)$ and whose bottom component is $\sin(n\tau/m)$ and the
other whose top component is $-\sin(n\tau/m)$ and whose bottom
component is $\cos(n \tau/m)$.  In this case, $n$ here must be an odd
integer and no less than $-m$.

Given now all of the preceding, a comparison of the spectral boundary
conditions just described for $D$ and for $\opart$ finds that the
given end of $C_0$ contributes $m+1$ to the deficit $\Theta$.
  }
Now, each end $\cE\subset C_0$ where the $z$--axis line bundle is
non-orientable and where the constant $s > s_0$ circles define odd
multiples of a generator of $H_1(M; \mz)$ supplies a positive, odd
integer $m \equiv m(\cE)$ and the sum, indexed by such ends of $C_0$,
of $ - m(\cE)$ defines an integer $\aleph_1$.  Meanwhile, each end
$\cE\subset C_0$ where the corresponding $\theta_0$ is in $\{0, \pi\}$
(whether or not the ambient $z$--axis line bundle is orientable)
defines an non-zero integer $m_0(\cE)$ (as described above).  Use
$\aleph _0$ to denote the sum, indexed by this last set of ends, of
$(1 - 2 m_0(\cE))$.  Then, the index of $D$ is given by
index($D$)~=~2$\cdot$degree$(N) + \chi(C_0) +\aleph_0 + \aleph_1$,
which is no different than saying that index($D$)~=~$\la e, [C]\ra -
\la c_1, [C]\ra - 2 m_C +\aleph_0 + \aleph_1$ or index($D$)~=
~$-\chi(C_0) - 2 \la c_1, [C]\ra + \aleph_0 + \aleph_1$.
\bigskip

{\bf Step 5}\qua This step considers the possibility that $X$ has some
convex ends that are intersected by $C$.  The first point to make is
that the analysis for the convex end case is completely analagous to
that for the concave ends.  Indeed, the difference amounts to
essentially considering the eigenvalues $E$ of the operator $L_0$ in
(3.19) which are non-positive rather than non-negative.  The result of
this analysis (whose details are left to the reader) follows.

To begin, let $\cE\subset $C be an end which lies in a convex end of
$X$ with orientable $z$--axis line bundle.  Then $\cE$ adds a factor 1
to the previous index formula, 2 degree$(N) + \chi(C_0) + \aleph_0 +
\aleph_1$, when the corresponding element in $C$'s limit set has
$\theta_0 \notin \{0, \pi\}$.  On the other hand, when $\theta_0 \in
\{0, \pi\}$, then $\cE$ contributes $2m_0(\cE) - 1$ to the integer
$\aleph_0$, where $m_0(\cE)$ is defined as in Step~3.

Next, let $\cE\subset C$ denote an end which lies in a convex end of
$X$ with non-orientable $z$--axis line bundle.  In this case, the
contribution of $\cE$ to the index formula again depends on whether
the large, but constant $s$ circles in $\cE$ give an even or an odd
multiple of a generator of $H_1(S^1 \times S^2; \mz)$.  In the case
where these circles give an even multiple of a generator, then $\cE$
contributes 1 to the previous index formual when the components of the
inverse image of the end under the map in (3.3) define limiting Reeb
orbits with $\theta_0 \in \{0, \pi\}$.  On the other hand, where such
lifts produce limiting Reeb orbits with $\theta_0 \in \{0, \pi\}$,
then $\cE$ contributes $2m_0(\cE) - 1$ to the integer $\aleph_0$,
where $m_0(\cE)$ is defined as in Step 3 by either of the lifts.
  
In the case where the constant $s$ circles in $\cE$ give an odd
multiple of a generator of the first homology of $S^1 \times S^2$,
then $\cE$ contributes $m(\cE) $ to the integer $\aleph_1$.  Here,
$m(\cE)$ is the absolute value of this multiple.
\bigskip      

{\bf Step 6}\qua This step considers the general case and so drops the
assumption that $C$ is immersed in $X$.  For such $C$, the operator
$D$ is the one that Lemma 3.4 considers.  But for a change of
interpretation, its index is again index$(D) = \la e, [C]\ra - \la
c_1, [C]\ra - 2 m_C +\aleph_0 + \aleph_1$.  Here, the change concerns
the integer $m_C$ which now must be interpreted as in the second point
of Proposition 3.1.
  	
To justify the preceding formula, first remember that the given map
from $C_0$ to $X$ has arbitrarily small, $C^\infty$ perturbations for
which the result satisfies the requirements in the second point of
Proposition 3.1.  The perturbation can be done so that the perturbed
map is pseudoholomorphic for an appropriate almost complex structure
on $X$ which differs from the original on a compact subset of $X$.
The operator $D$ of Lemma 3.4 is defined for this new map, and if the
perturbation is small in the $C^2$ topology, it will differ from the
original by an operator with small norm; thus, the new and old
versions of Lemma 3.4's $D$ will have the same index.  Meanwhile, the
new map, being an immersion, has an associated operator, call it $D'$,
which is given by (3.5) and is the one discussed in Steps 1--4, above.
It is left as an exercise to check that the operator $D$ of Lemma 3.4
and the operator in (3.5) have the same index when $C$ is immersed.

The following proposition summarizes the results of the preceding
discussion:
\bigskip

\noindent
{\bf Proposition 3.6}\sl\qua {Let $C\subset \fM_{e,\chi}$ and define
the operator $D$ either by (3.5) when $C$ is immersed, or as in Lemma
3.4.  Let $\cL_0$ and $\cL_1$ be as defined in either Lemma 3.3 or
Lemma 3.4 as the case may be and view $D$ as a Fredholm operator from
$\cL_0$ to $\cL_1$.  Then, the following formulae hold:}
\begin{itemize}
\item ${\rm index}(D) = \la e, [C]\ra - \la c_1, [C]\ra - 2 m_C
+\aleph_0 + \aleph_1 + \aleph$.
\item ${\rm index}(D) = -\chi(C_0) - 2 \la c_1, [C]\ra +\aleph_0 +
\aleph_1 + \aleph$.
\end{itemize}

\noindent
{Here, the integers $m_C$, $\aleph_0$, $ \aleph_1$ and $\aleph$
are defined as follows:}
\begin{itemize}
\item {$m_C$ denotes the number of double points of any
perturbation of the defining map from $C_0$ onto $C$ which is
symplectic, an immersion with only transversal double-point
singularities with locally positive self-inter\-sec\-tion number, and
which agrees with the original on the complement of some compact set.}
\item { $\aleph_0$ is the sum of the contributions of the form
$\epsilon (1 - 2m_0(\cE))$ from each end $\cE$ of $C_0$ for which the
corresponding limiting Reeb orbit has $\theta_0 \in \{0, \pi\}$.  In
this regard, $\epsilon = +1$ when $\cE$ lies in a concave end of $X$
and $\epsilon = -1$ otherwise.  Meanwhile, $m_0(\cE)$ is a positive
integer, and here is its definition when the corresponding end of $X$
has orientable $z$--axis line bundle: Let $m(\cE)$ denote the absolute
value of the pairing between a generator of $H_1(S^1 \times S^2; \mz)$
and any constant, $s > s_0$ circle in $\cE$.  Then, $m_0(\cE)$ is the
least integer which is greater than $m(\cE) \sqrt{3}/\sqrt{2}$.  In
the case when the corresponding component of $X$ has non-orientable
$z$--axis line bundle, use this last formula but with $\cE$ replaced by
either component of its inverse image under the map in (3.3).}
\item {$\aleph_1$ is the sum of the contributions of the form
$\epsilon m(\cE)$ from each end $\cE$ of $C_0$ that satisfies the
following criteria:}
\begin{itemize}
\item[\rm(a)] {The end $\cE$ lies in an end of $X$ with
non-orientable $z$--axis line bundle.}
\item[\rm(b)] {The absolute value of the pairing between a generator
of the first cohomology of the corresponding end of $X$ and any
constant, $s > s_0$ circle in $\cE$ is odd.  The preceding understood,
take $\epsilon = -1$ when $\cE$ lies in a concave end of $X$ and take
$\epsilon = ++1$ otherwise. Meanwhile, take $m(\cE)$ to be the
absolute value of the cohomology pairing just described.}
\end{itemize}
\item {$\aleph$ is the number of ends of $C_0$ that satisfy one
of the following two criteria:}
\begin{itemize}
\item[\rm(a)] {The end $\cE\subset C_0$ lies in a convex end of $X$
with orientable $z$--axis line bundle.  In addition, $\cE$'s
corresponding limit Reeb orbit has $\theta_0\notin\{0, \pi\}$.}
\item[\rm(b)] {The end $\cE\subset C_0$ lies in a convex end of $X$
that has non-orientable $z$--axis line bundle.  In addition, the
inverse image via (3.3) of $\cE$ has 2--components and neither defines
a limiting Reeb orbit with $\theta_0 \in \{0, \pi\}$.}
\end{itemize}
\end{itemize}

\noindent 
{Finally, when $C$ is immersed, the index formula can also be
written as}
$$
{\rm index}(D) = {\rm degree}(N) - \la c_1, [C]\ra + \aleph_0 +
\aleph_1 + \aleph,
$$
{where $N$ denotes the normal bundle of $C$, and where its degree
is defined using a section that is non-vanishing on the ends of $C$
and there homotopic through non-vanishing sections to the sections
described in Part 2 of Section 3a.}\rm

\section{Subvarieties in $\mr \times (S^1 \times S^2)$}
	
This section turns away from the general discussion of the previous
two sections to concentrate on the special case where $X = \mr \times
(S^1 \times S^2)$ with the symplectic structure and associated almost
complex structure as described in (1.1)--(1.4).  In particular, the
discussion here highlights general features of the moduli space of HWZ
subvarieties in $\mr \times (S^1 \times S^2)$ and results in a proof
of Theorem A.1.  In fact, given Proposition 3.2, the latter follows
directly from Proposition 4.2, below.  Other features deemed of
particular interest are summarized in Propositions 4.3 and 4.8.

\sh{(a)\qua Examples}

Before diving into generalities, it proves useful to introduce various
explicit examples of HWZ pseudoholomorphic subvarieties as these cases
will be referred to in the later subsections.  The examples below
exhaust the set of HWZ subvarieties that are invariant under some
1--parameter subgroup of the subgroup $T = S^1 \times S^1$.  In this
regard, remember that $T$ is generated by the vector fields
$\partial_t$ and $\partial_\varphi$ and it is the component of the
identity of the subgroup of the isometry group of $S^1 \times S^2$
that preserves the contact form $\alpha$.
\bigskip

\textbf{Example 1}\qua The simplest examples are the cylinders $C = \mr
\times \gamma$, where $\gamma \subset S^1 \times S^2$ is a closed Reeb
orbit.  Remember that these closed Reeb orbits are described in (1.8).
In particular, each such $\gamma$ is characterized, in part, by some
constant value, $\theta_0$, for the polar angle $\theta$ on $S^2$.
For example, if $\theta_0 \notin \{0, \pi\}$, then the closed Reeb
orbit determines a pair of relatively prime integers $(p, p')$ and if
$p \not= 0$, then the corresponding cylinder is parameterized by a
periodic coordinate $\tau\in \mr/(2\pi|p|\mz)$ and $u$ with $|u| \in
(0, \infty)$ and $\text{sign}(u) = \text{sign}(p)$ by the rule
\begin{equation}
(t = \tau, \;\;f = u,\;\; \varphi = \varphi_0 + p'/p \tau,\;\; h =
p'/p \sin^2(\theta_0) u)\;.
\label{(4.1)}
\end{equation}
Note that the condition $p \not= 0$ is equivalent to the assertion
that $\cos^2 \theta_0 \not= 1/3$.
 
Meanwhile, if $p = 0$, then $p' = \pm 1$ and $\cos \theta_0 =
p'/\sqrt{3}$.  In this case, the parameterization of the corresponding
cylinder is by $\tau \in \mr/(2\pi\mz)$ and $u$ with $|u| \in
(0,\infty)$ and $\text{sign}(u) = p'$ according to the formula
\begin{equation}
(t = t_0,\;\; f = 0,\;\; \varphi = \tau,\;\; h = u) ,
\label{(4.2)}
\end{equation}
where $t_0$ is a constant.

Finally, in the case where $\theta_0$ is 0 or $\pi$, the
parameterization of the corresponding cylinder is by $\tau \in
\mr/(2\pi \mz)$ and $u \in (0, \infty)$ according to the formula
\begin{equation}
(t = - \tau, \;\;f = - u)\; .
\label{(4.3)}
\end{equation}
Note that in this case, the function $s$ along $C$ is given by $s = -
6^{-1/2} \ln(u/2)$.
\bigskip

\textbf{Example 2}\qua This example has $t = \text{constant}$ and $f =
-\kappa < 0$ with $\kappa$ being constant.  Each such example is a
once punctured sphere---a plane---whose limit set lies in the $(s
\rightarrow-\infty)$ and thus convex end of $\mr \times (S^1 \times
S^2)$.  Here, the limit set is a $ (p = 0$, $p' = \pm 1)$ and so
$\cos^2 \theta_0 = 1/3$ closed Reeb orbit.  Given that the limit
closed Reeb orbit in question has $\cos \theta_0 = p'/\sqrt{3}$, such
a subvariety can parameterized by polar coordinates on the plane,
$(\tau, u) \in \mr/(2\pi \mz) \times [0, \infty)$ according to the
rule
\begin{equation}
(t = t_0, \;\;f = -\kappa,\;\; \varphi = \text{sign}(p') \tau,\;\; h =
\text{sign}(p') u).
\label{(4.4)}
\end{equation}
The maximum value of $s$ achieved by this plane occurs at the plane's
origin, where $u = 0$.  At this point, $s_{\max} = - 6^{-1/2}
\ln(\kappa/2)$ and $\theta$ is either 0 or $\pi$ depending on whether
$p' = +1$ or $Ð1$, respectively.
\bigskip

\textbf{Example 3}\qua In this case, $t = \text{constant}$ and $f =
\kappa > 0$ with $\kappa$ again being constant.  Each such variety is
an embedded cylinder whose limit set lies in the convex end of
$\mr\times (S^1 \times S^2)$.  In this regard, the limit set consists
of both $\cos^2 \theta_0 = 1/3$ closed Reeb orbits.  Such a subvariety
can be parameterized by $\tau \in \mr/(2\pi \mz)$ and $u \in \mr$
according to the formula
\begin{equation}
(t = t_0,\;\; f = \kappa, \;\;\varphi = \tau,\;\; h = u)\; .
\label{(4.5)}
\end{equation}
In this case, the maximum value of $s$ is achieved when $u = 0$ with
value $s_{\max} = -6^{-1/2} \ln(\kappa/2)$.  This maximum value occurs
where $\theta = \pi/2$.
\bigskip

\textbf{Example 4}\qua These next examples are all cylinders where
$\varphi$ and $h$ are both constant, but with $h \not= 0$.  Consider
first the case where $\varphi = \varphi_0$ and $h = \kappa \not= 0$.
The limit set for such a cylinder consists of one $(p = 1,$ $p' = 0)$
closed Reeb orbit (so $\theta_0 = \pi/2$, $\varphi = \varphi_0$) and
either the closed Reeb orbit with $\theta = 0$ or that with $\theta =
\pi$.  The $\theta = 0$ closed Reeb orbit appears when $\kappa > 0$
and the $\theta = \pi$ closed Reeb orbit appears when $\kappa < 0$.
In any case, both of ends of the cylinder lie in the convex end of
$\mr \times (S^1 \times S^2)$.  The cylinder is parameterized by
coordinates $\tau\subset \mr/(2\pi \mz)$ and $u \subset \mr$ by the
rule
\begin{equation}
(t = \tau,\;\; f = u,\;\; \varphi = \varphi_0,\;\; h = \kappa)\; .
\label{(4.6)}
\end{equation}
The maximum value of $s$ on this orbit is $s_{\max} = - 6^{-1/2} \ln(3
\kappa/(2\sqrt{2}))$; it occurs where $u = 0$ and so $\cos^2 \theta =
1/3$.
\bigskip

The next three examples describe 1--parameter families of embedded,
pseudoholomorphic cylinders. However, the subsequent discussion
requires the preliminary digression that follows.

To start the digression, remark that the families considered below are
labeled in part by a pair $(p, p')$ of relatively prime integers with
$p > 0$ and an angle $\varphi_0 \subset \mr/(2\pi \mz)$.  In this
regard, the combination $p\varphi - p't$ is constant on a cylinder
with label $((p, p')$, $\varphi_0)$ and so such a cylinder is
invariant under the subgroup of $T$ generated by the vector field
$p\partial_t + p'\partial_\varphi$.

A cylinder labeled by $((p, p'), \varphi_0)$ can be parameterized by
$\tau\int \mr/(2\pi|p|\mz)$ and $u \in \mr$.  The parameterization is
given as
\begin{equation}
(t = \tau,\; f = u, \;\varphi = \varphi_0 + (p'/p) \tau,\; h = h(u))
\label{(4.7)}
\end{equation}
where the function $h(\cdot)$ is constrained to obey the differential
equation
\begin{equation}
h_u = (p'/p) \sin^2 \theta\; .
\label{(4.8)}
\end{equation}
Here, the subscript Ô$u$Õ signifies differentiation with respect to
$u$.  Note that $\sin^2 \theta$ in (4.8) should be viewed as an
implicit function of $u$ and $h$ through its dependence on $f$ and
$h$.  Indeed, the identity
\begin{equation}
h/f = \sqrt{6} \cos \theta \sin^2 \theta/(1- 3 \cos^2 \theta)
\label{(4.9)}
\end{equation}
can be inverted locally to write $\theta = \theta(h/f)$.

Note that (4.8), being a first-order, ordinary differential equation
for one unknown function, has a single constant of integration which
parameterizes its solution set.  This constant of integration
determines, and is completely determined by fixing the value of the
coordinate $s$ in some fiducial manner.

It is illuminating to use (4.8) and the condition $f = u$ to derive
equivalent differential equations for both $\theta$ and $s$.  In
particular, the latter read:
\bigskip

$\bullet$\qua $\theta_u = 6^{-1/2} e^{\sqrt{6} s} (1 + 3 \cos^4
\theta)^{-1} \{\sqrt{6} \cos \theta - p'/p (1 - 3 \cos^2 \theta)\}
\sin \theta$ .

$\bullet$\qua  $s_u = - 6^{-1/2} e^{\sqrt{6}s} (1 + 3 \cos^4 \theta)^{-1}
\{1 - 3 \cos^2 \theta+ \sqrt{6} p'/p \cos \theta \sin^2 \theta\} $.\autonum

For any given functional dependence of $s$ on $u$, the top equation in
(4.10) has constant $\theta$ solutions only for $\theta = 0$, $\pi$,
$\theta_0$ and, in the case when $|p'|/p > \sqrt{3}/\sqrt{2}$, also at
a fourth angle, $\underline{\theta}_0$.  Here, $\theta_0$ and
$\underline{\theta}_0$ are determined as follows: In the case where
$p' = 0$, set $\theta_0 = \pi/2$.  For the other cases, digress
momentarily to note that
\begin{equation}
(\sqrt{6} \alpha)^{-1} \{- 1 + (1 + 2 \alpha^2)^{1/2}\}
\label{(4.11)}
\end{equation}
has absolute value less than 1/3 for any choice of $\alpha ­ 0$.  In
particular, for any such $\alpha$, the expression in (4.11) can be
written as cos $\theta_0$ for precisely one $\theta_0 \in (0, \pi)$.
Meanwhile, if $\alpha$ is such that $|\alpha| > \sqrt{3}/\sqrt{2}$,
then
\begin{equation}
(\sqrt{6} \alpha)-1 \{- 1 - (1 + 2 \alpha^2)^{1/2}\}
\label{(4.12)}
\end{equation}
has absolute value less than 1.  Thus, this last expression equals
$\cos \uth_0$ for precisely one angle $\uth_0 \in (0, \pi)$.  Note
that the expressions in (4.11) and (4.12) have opposite signs, and
thus $\cos \theta_0$ and $\cos \uth_0$, have opposite signs.
   
With the preceding understood, end the digression and define
$\theta_0$ and $\uth_0$ for $(p, p')$ when $p'$ is non-zero by
replacing $\alpha$ in the preceding digression by $p'/p$.  As
remarked, the angles $\theta_0$, $\uth_0$ (when $|p'|/p >
\sqrt{3}/\sqrt{2})$ and the angles $\theta = 0$ and $\theta = \pi$
provide the only $h = \text{constant}$ solutions to (4.9).  In this
regard, the constant $\theta$ solutions to (4.9) provide the
pseudoholomorphic cylinders that are already described in Example 1
above.
  
With the constant solutions to the top equation in (4.10) understood,
consider next the non-constant solutions.  These solutions have the
property that $\theta_u$ is either always positive or always negative.
Thus, in the case where $|p'|/p < \sqrt{3}/\sqrt{2}$, either $\theta$
ranges without critical points between 0 and $\theta_0$ or between
$\theta_0$ and $\pi$, with these particular values giving the inf and
sup of $\theta$ as the case may be.  On the other hand, when $p'/p >
\sqrt{3}/\sqrt{2}$, then there are three possible ranges for $\theta$.
The first has $\theta$ ranging without critical values between 0 and
$\theta_0$, the second between $\theta_0$ and $\uth_0$ and the third
between $\uth_0$ and $\pi$.  The analgous situation exists when $p'/p
< -\sqrt{3}/\sqrt{2}$, for here $\theta$ ranges without critical
points either between 0 and $\uth_0$, or between $\uth_0$ and
$\theta_0$, or between $\theta_0$ and $\pi$.

\eject
What ever the range of $\theta$, as long as the latter is not
constant, then it is a simple matter to change variables from $u$ to
$\theta$ in the second line of (4.10) and thus see $s = s(\theta)$ as
an anti-derivative of the function

$- \{1 - 3 \cos^2 \theta + \sqrt{6} \alpha \cos \theta
\sin^2 \theta\}/ \{(\sqrt{6} \cos\theta - \alpha (1 - 3 \cos^2
\theta)) \sin \theta\}$.\autonum

\noindent
Note that any pair of anti-derivatives for (4.13) differ by an
additive constant, and the latter can be taken as another parameter
which distinguishes the cylinders in any the families that are
considered below.
 
With $s = s(\theta)$ now viewed as an anti-derivative of (4.13), the
$\theta$ dependence of $u$ can be obtained by either solving the
algebraic equation $u = e^{-\sqrt{6}s} (1 - \cos^2 \theta)$, or by
using $s(\theta)$ in the top line of (4.10) to view the latter
autonomous equation for $\theta$ as a function of $u$.  In any event,
the second approach does freely identify the sign of $\theta_u$ along
the cylinder.

With the digression now over, consider the examples.
\bigskip

\textbf{Example 5}\qua In this example, all of the cylinders have both
components of the limit set in the convex end of $\mr \times(S^1
\times S^2)$.  To describe these cylinders, first choose $(p, p')$ and
$\varphi_0$ as constrained in the preceding digression.  Then, each of
the cylinders under consideration is parameterized by (4.7) as
described in this same digression.

If $|p'|/p > \sqrt{3}/\sqrt{2}$, then $p'/p$ and $\varphi_0$
determines a 1--parameter family of examples where the parameter can be
taken to be the maximum value that the function $s$ achieves on the
given example.  If $p'/p > \sqrt{3}/\sqrt{2}$, then the limit set
consists of the $\theta = 0$ closed Reeb orbit and the closed Reeb
orbit with $(\theta = \theta_0, \varphi = \varphi_0 + p'/p t)$.
Meanwhile, the coordinate $\theta$ varies along the pseudoholomorphic
cylinder without critical points between 0 and $\theta_0$.  If $p'/p <
-\sqrt{3}/\sqrt{2}$, then the limit set consists of the $\theta = \pi$
closed Reeb orbit and the $(\theta = \theta_0, \;\varphi = \varphi_0 +
p'/p t)$ closed Reeb orbit while $\theta$ varies on the cylinder
without critical points between these two extremes.

In the case where $|p'|/p < \sqrt{3}/\sqrt{2}$, then each such pair
$((p, p'), \varphi_0)$ determines two single parameter families of
pseudoholomorphic cylinders.  On each such family, the parameter can
still be taken to be the maximum value achieved by the function $s$.
For both of these families, the closed Reeb orbit with $(\theta =
\theta_0, \varphi = \varphi_0 + p'/p t)$ comprises one of the
components of the limit set.  But, the families are distinguished by
the other component of the limit set, which is either the closed Reeb
orbit where $\theta = 0$ or $\theta = \pi$.  As before, the function
$\theta$ varies without critical points with its supremum and infimum
given by the limit set values.
\bigskip	

\textbf{Example 6}\qua These examples consider the cases of (4.7) where
$|p'|'/p > \sqrt{3}/\sqrt{2}$ and where $\theta$ ranges between the
extreme values $\theta_0$ and $\uth_0$.  In this regard, note that
$\cos \theta_0$ has the same sign as $p'$ while $\cos \uth_0$ has the
opposite sign.  For these cylinders, both components of the limit set
again lie in the convex end of $\mr\times (S^1 \times S^2)$.  Here,
the limit set consists of the two closed Reeb orbits where $(\theta =
\theta_0, \; \varphi = \varphi_0 + p'/p t)$ and $(\theta = \uth_0,
\;\varphi = \varphi_0 + p'/p t)$.  Note that in this example,
$\theta_u$ is negative when $p'/p > \sqrt{3}/\sqrt{2}$ and $\theta_u$
is positive when $p'/p < -|sqrt{3}/\sqrt{2}$.  The parameter which
distinguishes the elements in any one family can again be taken to be
the value of $s$ at its maximum.
\bigskip

\textbf{Example 7}\qua These examples also consider the cases of (4.7)
where $|p'|/p > \sqrt{3}/\sqrt{2}$, but here $\theta$ ranges between
$\uth_0$ and $\pi$ in the positive $p'$ case, and between 0 and
$\uth_0$ in the negative $p'$ case.  These examples are embedded
cylinders with one component of the limit set in the convex end of $R
\times (S^1\times S^2)$ and the other in the concave end.  In this
regard, the orbit in the convex end is the closed Reeb orbit with
$(\theta = \uth_0, \;\varphi = \varphi + p'/p t)$.  Meanwhile, the
component of the limit set in the concave end is, depending on the
sign of $p'$, either the $\theta = \pi$ or $\theta = 0$ closed Reeb
orbit.

In any event, for fixed $(p, p')$ and $\varphi_0$, there is, once
again, a 1--parameter family of such examples.  However, in this case,
the function s restricts to the cylinder with neither maxima nor
minima, and so the value of $s$ at some specified $\theta$ value can
be taken as the parameter.
	      
\sh{(b)\qua The Index of the Operator $D$}

The purpose of this subsection is to describe certain aspects of the
kernel, cokernel and index for the operator $D$ of Propositions 3.2
and 3.6 for an HWZ subvariety in $X = \mr \times (S^1 \times S^2)$.
  
The next proposition summarizes the index story by restating
Proposition 3.6 in this special case.  With regard to the statement of
the subsequent proposition, remember that $\chi(C_0)$ denotes the
Euler characteristic of the smooth model curve, $C_0$, for a given
pseudoholomorphic subvariety $C$.  Also, note that the integer $\la
c_1, [C]\ra$ is to be defined as described in Section 3a.
\bigskip

\noindent
{\bf Proposition 4.1}\sl\qua {Let $C \subset X$ be an irreducible,
HWZ pseudoholomorphic subvariety, and use $C$ to define the operator
$D$ as described in Section 3b.  Then,}
\begin{equation}
{\rm Index}(D) = - \chi(C_0) - 2 \la c_1, [C]\ra + \aleph + \aleph_+ +
\aleph_-\;,
\label{(4.14)}
\end{equation}
{where $\aleph$, $\aleph_+$ and $\aleph_-$ are defined as
follows:}
\begin{itemize}
\item $\aleph$ {denotes the number of ends of $C$ which lie in
the convex end of $X$ and which approach an element of $C$Õs limit set
with $\theta_0 \not= \{0, \pi\}$.}
\item $\aleph_+$ {is the sum of contributions of the form $(1-
2m_0(\cE))$ from each end of $C_0$ which lies in the concave end of
$X$ and for which the corresponding element in $C$Õs limit set has
$\theta_0 \in \{0, \pi\}$.  Here, $m_0(\cE)$ is the positive integer
which is defined as follows: Let $m(\cE)$ denote the absolute value of
the pairing between a generator of $H^1(S^1 \times S^2; \mz)$ and any
sufficiently large, but constant $s$ circle in $\cE$.  Then,
$m_0(\cE)$ is the least integer which is greater than $m(\cE)
\sqrt{3}/\sqrt{2}$.}
\item $\aleph_-$ { is the sum of contributions of the form
$2m_0(\cE) - 1$ from each end of $C_0$ for which the corresponding
element in $C$Õs limit set has $\theta_0 \in \{0, \pi\}$ and lies in
the convex end of $X$.  Here, $m_0(e)$ is defined as above.}
\end{itemize}\rm
\bigskip

It is important to note that (4.14) places serious constraints on the
subvarieties with a given index for the operator $D$.  The following
proposition lists the constraints on the subvarieties with
$\text{index}(D) \le\aleph + 1$.
\bigskip

\noindent
{\bf Proposition 4.2}\sl\qua {Let $C \subset X$ be an irreducible,
HWZ subvariety.  Then the following is true:}
\begin{itemize}  
\item $\text{Index}(D) \ge \aleph$.
\item {If ${\rm index}(D) = \aleph$, then $\aleph = 0$, 1, or 2.}
\begin{itemize}  
\item[\rm(a)] {If $\aleph = 0$, then $C$ is a $\theta_0 \in \{0,
\pi\}$ case from Example 1.}
\item[\rm(b)] {If $\aleph = 1$, then $C$ is a $\theta_0 \notin \{0,
\pi\}$ case from Example 1.}
\item[\rm(c)] {If $\aleph = 2$, then $C$ comes either from Example 3
or 6.}
\end{itemize}
\item {If ${\rm index}(D) = \aleph + 1$, then $\aleph = 1$, 2 or
3.}
\begin{itemize}
\item[\rm(a)] $C$ {is the plane from Example 2, so} $\aleph = 1$.
\item[\rm(b)] $C$ {is a cylinder from Example 5 where one closed Reeb
orbit in the limit set is characterized by a pair $(p, p')$ where
$|p'|$ is the greatest integer that is less than
$(\sqrt{3}/\sqrt{2})|p|$. The other closed Reeb orbit has $\theta_0 =
0$ if $p' < 0$ and $\theta_0 = \pi$ if $p' > 0$.  Here, $\aleph = 1$.}
\item[\rm(c)] $C$ {is a cylinder from Example 7 where one closed Reeb
orbit in the limit set is characterized by a pair $(p, p')$ where
$|p'|$ is the least integer that is greater than $(\sqrt{3}/\sqrt{2})
|p|$.  Here also, $\aleph = 1$.}
\item[\rm(d)] $C$ {is an immersed, thrice punctured sphere with
$\aleph = 2$, or 3.  Moreover, $C$ has no intersections with the
$\theta\in \{0, \pi\}$ locus and none of its limit set closed Reeb
orbits have $\theta_0 \in \{0, \pi\}$.}
\end{itemize}
\item {All other cases have ${\rm index}(D) \ge \aleph + 2$.}
\end{itemize}\rm
\bigskip

The question of existence and the classification of the thrice
punctured spheres from Part d of the propositionÕs third point is
deferred to Sections 5 and 6.  Sections 5 and 6 also provide a formula
for the number of double points for these immersed spheres.  Note that
Theorem A.1 in the Introduction follows from Propositions 3.2 and 4.2.
 
The following proposition elaborates on the inequality in the final
point of Proposition 4.2.
\bigskip

\noindent
{\bf Proposition 4.3}\sl\qua {Let $C \subset X$ be an irreducible,
HWZ subvariety.  Then}
$$
{\rm Index}(D) \ge 2 (-1 + g + Q + \aleph + \aleph_0^{cc} +
\aleph_0^{cv}) + \aleph^c
$$
{where}
\begin{itemize}
\item $g$ {is the genus of $C_0$.}
\item $Q$ {is the number, counted with multiplicity, of
intersections of $C$ with the $\theta \in \{0, \pi\}$ locus.}
\item $\aleph_0^{cc}$ {is the number of concave side ends of $C$
where the $s \rightarrow\infty$ limit of $\theta$ is either 0 or
$\pi$.
\item $\aleph_0^{cv}$ {is the number of convex side ends of $C$
where the $s\rightarrow\infty$ limit of $A\theta$ is either 0 or
$\pi$.}
\item $\aleph^c$ is the number of concave side ends of $C$ where the
$s\rightarrow\infty$ limit of $\theta$ is not 0 nor $\pi$.}
\end{itemize}\rm

With regard to the organization of the remainder of this section, the
next subsection, 4c, contains the first of three parts to the proof of
Proposition 4.2.  Subsection 4d constitutes a digression that proves
Proposition 4.3.  With aspects of the latter proof then available,
Subsection 4e resumes the proof of Proposition 4.2 and contains the
latterÕs second part.  The final part of the proof of Proposition 4.2
is in Subsection 4f.  There is an extra subsection that considers the
cokernel dimension of the operator $D$ when the relevant subvariety is
any from an example in Section 4a or any thrice punctured sphere from
Part d of the third point of Proposition 4.2.  In particular,
Proposition 4.8 in this subsection asserts that this cokernel is
trivial in all these cases.

\sh{(c)\qua Proof of Proposition 4.2, Part 1}

Consider in this part of the proof solely the case where no closed
Reeb orbit in $C$Õs limit has $\theta_0 \notin \{0, \pi\}$.  The first
point to make is that the preceding assumption implies that
\begin{equation}
\la c_1, [C]\ra \le 0 \;.
\label{(4.15)}
\end{equation}
Indeed, $\la c_1, [C]\ra$ computes a sum of non-zero integer terms,
where each corresponds to a zero along $C$ of the section over $X$ of
the line bundle $K$ given by $(dt + i g^{-1} df) \wedge (\sin^2 \theta
d\varphi + i g^{-1} dh)$.  Since neither $dt$ nor $df$ vanish on $X$,
this section vanishes only where $\sin^2 \theta d\varphi$ and $dh$
vanish, which is along the pair of cylinders where $\theta\in \{0,
\pi\}$.  Thus, $\la c_1, [C]\ra$ counts, with the appropriate integer
weight, the intersections of $C$ with these half cylinders.  It is
left to the reader to verify that the integer weight here is negative
in all cases.
  
Now, under the given assumptions,$\aleph+ = 0$, $\aleph_- = 0$ and
$\aleph \ge 1$.  Moreover, as $C_0$ is connected, $\chi(C_0) \le 1$ so
(4.15) implies that Index$(D) \ge \aleph -1$ with equality if and only
if $C_0$ is a plane, (and thus $\aleph = 1$) and $\la c_1, [C]\ra =
0$.  To rule out this possibility, remember that a plane has one end,
and so there is just one closed Reeb orbit in the limit set.  The
latter either has $cos^2 \theta_0 = 1/3$ or not.  If not, then, as
will be argued momentarily, the zeros of the 1--form $dt$ count with
appropriate weights to give $\chi(C_0)$, and all of these weights are
negative.  In particular, this means that $\chi(C_0) \le 0$ which
rules out the $\cos^2 \theta_0 \not= 1/3$ possibility.

The claim that the zeros of $dt$ count with negative weights to
compute $\chi(C_0)$ is valid for any model curve $C_0$ for an
irreducible, HWZ subvariety with no $\cos^2 \theta_0 = 1/3$ closed
Reeb orbits in its limit set.  Here is a digression to explain why:
First, because no limit set closed Reeb orbit has $\cos^2 \theta_0 =
1/3$, the gradient on $C_0$ of the pull-back of the function $f$ is
not tangent to any constant, but sufficiently large $|s|$ circle in
$C_0$.  This implies that the zeros of $df$ count, with the usual
weights, $\chi(C_0)$.  Here, one should be careful to count degenerate
zeros appropriately.  However, the zeros of $df$ are the same as those
of $dt$ and all are isolated and all count with negative weight.  To
see that such is the case, introduce the real and imaginary parts,
$(x_1, x_2)$, of a complex parameter on a plane in $C_0$.  Then, by
virtue of (1.5), the functions $(t, f)$ obey a version of the
Cauchy-Riemann equations,

\begin{itemize}
\item $g t_1 = f_2$,
\item $g t_2 = - f_1$.\autonum
\end{itemize}

\noindent
This last equation implies via fairly standard elliptic equation
techniques that $dt$ and $df$ have the same zeros, that $t$ and $f$
are real analytic functions on $C_0$, that their zeros are isolated
and that all count with negative weights to give $\chi(C_0)$.

By the way, note that the pair $(\varphi, h)$ obey an analogous
equation:
\begin{itemize}
\item $g \sin^2 \theta\varphi_1 = h_2$,
\item $g \sin^2 \theta\varphi_2 = - h_1 $.\autonum
\end{itemize}

In particular, (4.17) implies that the zeros of $d\varphi$ or $dh$
count with negative weights to give $\chi(C_0)$ when $C_0$ is the
model curve for an irreducible, HWZ subvariety which lacks both
intersections with the $\theta_0 \in \{0, \pi\}$ cylinders and
$\theta_0 \in \{0,\pi/2,\pi\}$ closed Reeb orbits in its limit set.
  
With the digression complete, consider now the only remaining
possibility in the case under consideration, which is for the closed
Reeb orbit in $C$'s limit set to have $\cos^2 \theta_0 = 1/3$.  In
this case, the 1--form $dj$ pulls back without zeros on the constant
but large $|s|$ circles in $C_0$.  Moreover, as $C$ does not intersect
either the $\theta = 0$ or $\theta = \pi$ locus, this 1--form pulls
back as a smooth 1--form on $C_0$.  As remarked above, this pull-back
also count with negative weights to give $\chi(C_0)$.  Thus,
$\chi(C_0) \le 0$ and in no case can $C_0$ be a plane.

Now consider the possibilities when ${\rm Index}(D) = \aleph$.  Since
$\aleph_\pm= 0$, this can happen only if $\chi(C_0) = 0$ and $\la c_1,
[C]\ra = 0$.  As $\chi(C_0) = 0$, the surface $C_0$ is a cylinder and
so has two ends.  Thus, either $\aleph = 1$ or $\aleph = 2$.
 
Meanwhile, $C$ has no intersections with the $\theta = 0$ and $\theta
= \pi$ loci because $\la c_1, [C]\ra = 0$.  As before, this implies
that the 1--form $dj$ pulls back to $C_0$ as a smooth 1--form.  Of
course, dt always pulls back to $C_0$ as a smooth 1--form.  Since each
is closed, the integral of each along the constant $s$ circles in
$C_0$ must be independent of $s$.  The latter assertion demands that
the closed Reeb orbits at the ends of $C$ have the same value for the
integers $p$ and $p'$ in (1.8).  In particular, $p' t + p \varphi$
pulls back to $C_0$ as a bonafide function, moreover, one which
approaches a constant value asymptotically on the ends of $C_0$.

Now, $dt$ and also $d\varphi$ are either identically zero on $C_0$ or
nowhere zero as $\chi(C_0)$ would be negative otherwise.  The argument
is the same as given previously.  Moreover, both cannot vanish
identically, so at least one is nowhere zero.  Either being nowhere
zero implies that $C$ is immersed.

If $dt \equiv 0$, then $C_0$ is a $cos^2 \theta_0 = 1/3$ case from
Example 1 and index$(D) = \aleph = 1$, or else$C$ is described in
Example 3 and index$(D) = \aleph = 2$.  If $dt\not= 0$, then $C$ can
be parametrized by a periodic coordinate $\tau\in\mr/(2\pi|p|\mz)$ and
a linear coordinate $u$ as in (2.19).  Here, the pair of functions
$(x, y)$ obey (2.20).  Note also that $x$ is just the restriction of
$p' t + p \varphi$ to $C$.  In any event, as argued subsquently to
(2.20), the function $x$ in (2.23) must be constant, and this implies
that the asymptotic values of $p' t + p \varphi$ on the ends of $C$
are identical.  Meanwhile, as is demonstrated momentarily, (2.20) can
be employed with the maximum principle to prove that the function $x$
has neither local maxima nor local minima on $C$.  Thus, $p' t + p
\varphi$ is constant on $C$.  This implies that $C$ is either a case
from Example 1 with $\theta_0 \notin \{0, \pi\}$, or else a case from
Example 6.  The former has index$(D) = \aleph = 1$ and the latter has
index$(D) = \aleph = 2$.
  
To complete the index$(D) = \aleph$ story, return now to the postponed
part of discussion of local maxima and minima for the function $x$ in
(2.19).  For this purpose, note that (2.20) implies that $x$ obeys the
second-order differential equation
\begin{equation}
(g^2 \sin^2 \theta x_u)_u + (\sin^2 \theta x_\tau)_\tau + 2
\frac{p'}{p} (\sin \theta \cos \theta) \theta_\tau = 0 .
\label{(4.18)}
\end{equation}
Meanwhile, the identity $h/f = \sqrt{6} \cos \theta \sin^2 \theta/ (1
- 3 \cos^2 \theta)$ can be inverted where $\theta \not= \{0, \pi\}$ to
write q as a function of $h/f$, and thus the restriction of $\theta$
to $C$ can be viewed as a function of $h/u$.  Hence, the $\tau$
dependence of this function comes via the dependence of $\theta$ on
$h$.  In particular, one can write $\theta_\tau = \theta_h h_\tau$ and
then employ the top line in (2.20) to rewrite (4.18) as
\begin{equation}
(g^2 \sin^2 \theta x_u)_u + (\sin^2 \theta x_\tau)_\tau - 2
\frac{p'}{p} (g^2 \sin^3 \theta \cos \theta) \theta_hx_u = 0 .
\label{(4.19)}
\end{equation}
The strong form of the maximum principle applies directly to the
latter equation and precludes $x$ from having local maxima and minima.

Finally, this part of Proposition 4.2's proof ends by considering the
possibilities when index$(D) = \aleph + 1$.  Here, (4.14) allows only
two possibilities for the pair $\chi$ and $\la c_1, [C]\ra$; the first
has $\chi(C_0) = 1$ and $\la c_1, [C]\ra = -1$, while the second has
$\chi(C_0) = -1$ and $\la c_1, [C]\ra = 0$.  In the first case, $C_0$
is a plane which intersects the union of the $\theta = 0$ and $\theta
= \pi$ loci exactly once.  Moreover, the fact that $C_0$ is a plane
implies that the coordinate $t$ must restrict to $C_0$ as an
$\mr$--valued function, and thus $t$ must restrict as an $\mr$--valued
function to the closed Reeb orbit which comprises $C$'s limit set.
The only closed Reeb orbits with this property have $\cos^2 \theta_0 =
1/3$, and thus $t$ is constant.  Therefore, $t$ approaches a constant
value on the end of $C_0$ which implies, via the maximum principle,
that $t$ is constant on the whole of $C_0$.  Hence, $C$ is described
by Example 2 in Section 4a.

The other possibility has $\la c_1, [C]\ra = 0$ and $\chi(C_0) = -1$.
Now, by virtue of the definition of the pairing $\la c_1, \cdot \ra$
in Section 3a, this first condition makes $C$ disjoint from the
$\theta\in\{0, \pi\}$ locus.  Meanwhile, with Euler characteristic
$Ð1$, $C_0$ is either a once punctured torus with $\aleph = 1$ or else
a thrice punctured sphere in which case the possibilities for $\aleph$
are 1, 2 or 3.  The torus case is ruled out by the following argument:
If the one end has $p \not= 0$, then $|f|$ increases uniformly with
increasing $|s|$ when the latter is sufficiently large.  Thus, if $p >
0$, the function $f$ would have a global minimum on $C$ and if $p <
0$, then $f$ would have a global maximum.  As neither can happen, no
such torus exists.  In the case $p = 0$, then $p' \not=­ 0$ and the
same argument applies with $h$ substituted for $f$.
 
Consider now the case for punctured spheres.  As Proposition 4.7,
below, asserts that the cases that arise in Proposition 4.2 are
immersed, all that is left to say about these $\chi(C_0) = -1$ and
$\la c_1, [C]\ra = 0$ cases is summarized by the following assertion:
\begin{equation}
\text{{\sl There are no cases with}}\; \aleph = 1.
\label{(4.20)}
\end{equation}

The argument for this claim is a somewhat more sophisticated version
of the preceding argument that ruled out punctured tori. In fact, the
argument that follows proves (4.20) with no preconditions on the
number of concave side ends or the genus.  The only precondition is
that $C$ has no intersections with the $\theta = 0$ and $\theta = \pi$
loci.  The claim in (4.20) is an immediate corollary to the following
assertion:

\begin{itemize}\sl
\item {If $C$ has no intersections with the $\theta=0$ and
$\theta=\pi$ loci, and $C$ is not an $\mr$--invariant cylinder then the
pullback of the function $\theta$ to $C_0$ has neither local maxima
nor local minima.}
\item {In addition, $\theta$'s restriction to any concave side
end takes values at arbitrarily large $s$ that are both larger and
smaller than its $s\rightarrow\infty$ limit on the end.}\autonum
\end{itemize}

To prove the first point, use (4.16) and (4.17) to derive a
second-order differential equation for $\theta$ to which the maximum
principle can be applied.  In this regard, it is important to note
that the latter equation has the form
\begin{equation}
\sigma_1 \Delta\theta + \sigma_2\cdot d\theta = 0 ,
\label{(4.22)}
\end{equation}
where $\Delta$ is the Laplacian on $C_0$, and where $\sigma_1$ and
$\sigma_2$ are well defined provided that $\theta\notin\{0, \pi\}$.

To argue for the second point, it proves useful to introduce the
function $\lam\equiv h/f$ which is defined where $\cos^2\theta\not=
1/3$.  In particular, the latter is a monotonic function of $\theta$
on each of the three $\cos^2\theta\not= 1/3$ components of $(0,\pi)$
and it obeys a version of (4.22). To discuss ends where the
$s\rightarrow\infty$ limit of $\theta$ satisfies $\cos^2\theta=1/3$,
the function $\lam^{-1}=f/h$ will be used instead. The latter is a
monotonic function of $\theta$ on each $\theta\not=\pi/2$ component of
$(0,\pi)$ and also obeys a version of (4.22).

To proceed, suppose that $E$ is a concave side end of $C_0$ whose
corresponding closed Reeb orbit has angle $\theta=\theta_0$ with
$\cos^2\theta_0\not=1/3$. The argument for the case where
$\cos^2\theta_0=1/3$ is left to the reader in as much as it is
essentially identical to the one given below after changing the roles
of the pair $(f,t)$ and $(h,\varphi)$. With the
$\cos^2\theta_0\not=1/3$ assumption understood, note that the second
point in (4.21) follows if the assertion in question holds for the
function $\lam$ instead of $\theta$. In this regard, note that the
$s\rightarrow\infty$ limit of $\lam$ is $\sqrt{6}
\cos\theta_0\sin^2\theta_0 (1-3\cos^2\theta_0)^{-1}$.  Use $\lam_0$ to
denote the latter. To continue, remark that it is enough to prove the
following: Given $s_0>0$, then $\lam_0$ is neither the infinum nor
supremum of the values of $\lam$ on the $s\ge s_0$ portion of the end
$E$.  Since the arguments in either case are essentially the same,
only the argument .for the infimum will be given. To start the latter,
suppose that $\lam_0$ were $\lam$'s infimum. The argument that follows
derives a contradiction from this assumption.

To obtain the contradiction, note first that for large values of the
coordinate $s$, the end of $C$ in question can be parameterized as in
(2.13) where the functions $(x,w)$ obey (2.15).  Furthermore, $|x|$
and $|w|$ tend to zero as $\rho\rightarrow\infty$.  Also, as asserted
in Lemma 2.5, the derivatives of $x$ and $w$ tend to zero as
$\rho\rightarrow\infty$. Keep all of this in mind. The focus here is
on the function $w$ since up to a positive multiple, $w$ is
$\lam-\lam_0$. Thus, under the given assumption, $w\ge 0$ for large
$s$ on $E$ and thus $w>0$ at large $s$ as $\lam$ cannot have a local
minimum. However, this possibility is precluded by (1.23). Indeed,
were $w$ positive everywhere, then at all large values of the prameter
$\rho$ in (2.13), the bottom component of (2.15) would force the
differential inequality
\begin{equation}
w_\rho + x_\tau - 2^{-1} \zeta w \ge 0 \;.
\label{(4.23)}
\end{equation}
Here, $\zeta> 0$ is the constant that appears in (2.16).  To see why
(4.23) holds, first note that the term $\cR(a, b)$ as it appears in
(2.17) depends on its first entry, $a$, only through $a$'s bottom
component.  Indeed, this follows from the lack of $\varphi$ dependence
in the complex structure in (1.5).  This point understood, then the
$\cR$ term in (2.15) is bounded by a constant multiple of $|w| (|w| +
|w_\tau| + |x_\tau|)$ and then (4.23) follows from Lemma 2.5's
guarantee that the derivatives of $x$ and $w$ vanish in the limit as
$\rho\rightarrow\infty$.

Now, (4.23) implies that the function, $\uw(\cdot)$, of $\rho$ which
is obtained by averaging $w$ over the $\rho = $constant circles would,
per force, be greater than zero and obey
\begin{equation}
\uw_\rho - 2^{-1} \zeta \uw \ge 0
\label{(4.24)}
\end{equation}
at large values of $\rho$.  Of course, the latter inequality forces
the growth of $w$ as $\rho$ gets large; and this last conclusion
provides the promised contradiction.

\eject
\sh{(d)\qua Proof of Proposition 4.3}

Given that the Euler characteristic $\chi(C_0)$ is equal to $2 - 2g -$
(the number of ends), the asserted inequality follows from the
inequality for $-2 \la c_1, [C]\ra + \aleph_+ + \aleph_-$ asserted by
the following result.
\bigskip

\noindent
{\bf Lemma 4.4}\sl\qua {Let $C$ be an irreducible, HWZ-subvariety.
Then}
\begin{equation}
-2 \la c_1, [C]\ra +\aleph_+ + \aleph_- \ge 2 Q + \aleph_0^{cc} +
\aleph_0^{cv}\;.
\label{(4.25)}
\end{equation}\rm

\noindent
The remainder of this subsection contains the following proof.
\bigskip

\noindent
{\bf Proof of Lemma 4.4}\qua As the proof of this lemma is long, it is
broken into seven steps.
\bigskip

\textbf{Step 1}\qua It is the computation of $\la c_1, [C]\ra$ that
complicates the proof; the complexity of this computation stems from
the fact that the defining section of the canonical bundle $K$ near
the $\theta_0 = \{0, \pi\}$ ends is different from that used near the
other ends.  This complication is addressed via a decomposition of
$\la c_1, [C]\ra$ as a sum of terms, one of which algebraically counts
the zeros on $C$ of the section $(dt + i g^{-1} df) \wedge (\sin^2
\theta d\varphi + i g^{-1} dh)$ of $K$, while the others, one for each
$\theta_0 \in \{0, \pi\}$ end of $C$, are `correction factors'.  This
decomposition of $\la c_1, [C]\ra$ is provided momentarily.  Coming
first is a lemma with proof which simplifies the definition of this
decomposition.
\bigskip

\noindent
{\bf Lemma 4.5}\qua {\sl An irreducible, HWZ-pseudoholomorphic
subvariety in $\mr\times (S^1\times S^2$) that is not a $\theta\in\{0,
\pi\}$ cylinder from Example 1 intersects such a cylinder a finite
number of times.}
\bigskip

\noindent
{\bf Proof of Lemma 4.5}\qua Since the $\theta = 0$ and $\theta = \pi$
loci are pseudoholomorphic cylinders, and isolated in the sense of
\cite{HWZ2}, this lemma is a version of Proposition 4.1 in \cite{HWZ2}.  Even
so, a proof is given below since various portions of it are used
subsequently.  In any event, the proof that follows is different from
that offered in \cite{HWZ2}.  The proof of this lemma constitues Steps 2
and 3 of the proof of Lemma 4.4.  In this regard, Step 2 establishes
that there are at most finitely many intersections on the concave side
of $\mr\times(S^1\times S^2)$, while Step 3 does the same for the
convex side.
\bigskip

\textbf{Step 2}\qua This step rules out the possibility of infinitely
many intersections between $C$ and the $\theta \in \{0, \pi\}$
subvarieties where $s > 0$.  To be more precise, only intersections
with the $\theta\equiv 0$ cylinder will be discussed here since the
analogous discussion for the intersections with the $\theta\equiv\pi$
cylinder are identical but for insignificant notational changes.
  
To start the story, pick some large and positive value, $s_0$, for $s$
with the property that $C$ is disjoint from both the $\theta\equiv 0$
cylinder and the $\theta = \pi/2$ and $\varphi = 0$ cylinder on the
slice $s = s_0$.  In this regard, note that either $C$ coincides with
one of these subvarieties or else the intersections with either have
no accumulation points.  By assumption, $C$ does not coincide with the
former, and if $C$ coincides with the latter, then there is nothing
further to discuss.  Thus, it is safe to assume that such an $s_0$
exists.  With $s_0$ chosen, let $\delta$ denote the minimum distance
between $C$'s intersection with the $s = s_0$ slice and the
intersection with this slice of the $\theta\equiv 0$ and $(\theta
\equiv \pi/2, \varphi(\varphi= 0)$ cylinders.

The discussion proceeds from here by assuming that $C$ does, in fact,
have an infinite number of intersections with the $\theta\equiv 0$
cylinder where $s > s_0$; some unacceptable foolishness is then
derived from this assumption.  In particular, a contradiction arises
by considering $C$'s $s > s_0$ intersection number with certain
pseudoholomorphic cylinders from Example 4 with $h \equiv\kappa > 0$
and $\varphi \equiv 0$.

In this regard, there are some preliminary facts to recall about
intersection numbers.  First, the local intersection numbers between
pairs of pseudoholomorphic subvarieties are strictly positive.
Second, the local intersection numbers between such subvarieties are
invariant under sufficiently small perturbations of the maps of the
model, smooth curves.  Third, if a subvariety has compact intersection
with the $s > s_0$ portion of $\mr\times (S^1 \times S^2)$ and no
intersections with $C$ in the $s = s_0$ slice, then it has a well
defined $s > s_0$ intersection number with $C$, and this intersection
number is invariant under compact deformations of the subvariety which
avoid $C$ on the $s = s_0$ slice.

Next, some comments are in order concerning the $(h \equiv\kappa,
\varphi\equiv 0)$ cylinders.  First, each such cylinder intersects the
$s > s_0$ portion of $\mr\times (S^1 \times S^2)$ in a compact set.
Second, there exists $\kappa_0$ such that if $0 < \kappa < \kappa_0$,
then the intersection of the $(h \equiv\kappa, \varphi\equiv 0)$
cylinder and the $s = s_0$ slice occurs in a radius $\delta/10$
tubular neighborhood of the union of the cylinders where $\theta\equiv
0$ and where $(\theta\equiv\pi /2, \varphi\equiv 0)$.  Thus, as long
as $0 < \kappa < \kappa_0$, no $(h \equiv\kappa, \varphi\equiv 0$
cylinder intersects $C$ where $s = s_0$.  Note that these last two
facts imply that each such $(h \equiv\kappa, \varphi\equivº 0)$
cylinder has a finite number of intersections with $C$ where $s >
s_0$, and thus finite $s > s_0$ intersection number with $C$.
Moreover, this intersection number is \emph{independent} of $\kappa$
as long as $\kappa\in(0, \kappa_0)$.  Third, fix $s_1 > s_0$ and a
tubular neighborhood of the $s_0 \le s \le s_1$ portion of the
$\theta\equiv 0$ cylinder.  Then, there exists some $\kappa_1$ such
that for $0 < \kappa < \kappa_1$, each $(h \equiv \kappa,
\varphi\equiv 0)$ cylinder intersects this tubular neighborhood as a
graph over the $s_0 \le s \le s_1$ portion of the $\theta\equiv 0$
cylinder, and as $\kappa\rightarrow 0$, these graphs converge in a
smooth manner to the trivial graph, the $\theta\equiv 0$ cylinder
itself.
  
Note that this third comment, plus the remarks about positivity of
local intersection numbers and their invariance under small
perturbations implies the following: Fix $n > 0$ and there exists
$\kappa[n] > 0$ such that when $\kappa\equiv (0, \kappa[n])$, then the
$(h \equiv \kappa, \varphi \equiv 0)$ cylinder has $s > s_0$
intersection number at least $n$ with $C$.  Of course, this conclusion
is ridiculus, because the $\kappa$--invariance of this intersection
number (for $0 < \kappa < \kappa_0)$ implies that the $(h =
\kappa_0/2, \varphi\equiv 0)$ cylinder has infinite intersection
number with $C$.
\bigskip

\textbf{Step 3}\qua This step rules out the possibility of $C$
intersecting the $\theta\equiv 0$ and $\theta\equiv\pi$ cylinders
infinitely many times where $s < 0$.  Here again, only intersections
with the $\theta\equiv 0$ cylinder will be discussed.  The strategy
here is similar to that used in Step 2: Assume that $C$ has infinitely
many negative $s$ intersections with the $\theta\equiv 0$ cylinder and
find a ridiculous conclusion.  In this case, the untenable conclusion
is that $C$ has infinitely many $s < 0$ intersections with certain
cylinders from Example 7.
  
To start the story, choose $s_0$ so that $C$ is disjoint from the
$\theta\equiv 0$ cylinder where $s = s_0$.  Now, let $\delta$ denote
the minimum distance between $C$'s intersection with the $s = s_0$
copy of $S^1 \times S^2$ and that of the $\theta\equiv 0$ cylinder.
Next, choose $p'$ to be an integer more negative than $-1$ and such
that with $\uth_0$ determined by $p'$ as in Example 7, no
$\theta\equiv \uth_0$ closed Reeb orbits lie $C$'s limit set.

Example 7 describes a 1--parameter family of pseudoholomorphic
cylinders all labeled by $\uth_0$ and some fixed choice for
an angle $\varphi_0 \in [0, 2\pi]$.  In this regard, fix $\varphi_0 =
0$ and then the resulting 1--parameter family of cylinders can be
labeled as $\{\Sigma_r\}_{r\in \mz}$ where the distinguishing feature
of $\Sigma_r$ is that its $s = r$ slice has distance $\delta/2$ from
the $s = r$ slice of the $\theta\equiv 0$ cylinder.  What follows are
some relevant facts to note about $\Sigma_r$.  First, if $r < s_0$,
then there are no $s = s_0$ intersections between $\Sigma_r$ and $C$,
and there are at most a finite number of such intersections where $s <
s_0$.  Thus, there is a well defined $s < s_0$ intersection number
between each $r \le s_0$ version of $\Sigma_r$ and $C$.  Moreover,
this intersection number is independent of $r \le s_0$.  In this
regard, note that there are no very negative $s$ intersections between
any $\Sigma_r$ and $C$ since the limit sets for $C$ and $\Sigma_r$ are
disjoint.

Here is a second crucial fact: Fix $s_1 < s_0$ and a tubular
neighborhood of the portion of the $\theta\equiv 0$ cylinder where $s
\in [s_1, s_0]$.  Then, there exists $r_1$ such that for $r < r_1$,
each $\Sigma_r$ intersects this tubular neighborhood as a graph over
the $s_1 \le s \le s_0$ portion of the $\theta\equiv 0$ cylinder, and
as $r \rightarrow-\infty$, these graphs converge in a smooth manner to
the trivial graph, the $\theta\equiv 0$ cylinder itself.
  
Now, as before, this last observation implies that given any positive
$n$, the $s < s_0$ intersection number between $\Sigma_r$ and $C$ is
at least $n$ if $r$ is sufficiently small.  However, the invariance of
this intersection number with variations of $r$ in $(-\infty, s_0]$
implies the silly conclusion that $C$ has infinite $s < s_0$
intersection number with any $\Sigma_r$.
\bigskip

\textbf{Step 4}\qua With the proof of Lemma 4.5 complete, here is the
advertised decomposition of $\la c_1, [C]\ra$:
\begin{equation}
\la c_1, [C]\ra = - \nu_0 + \Sigma_{\cE\in V}\nu(\cE)
\label{(4.26)}
\end{equation}
where $\nu_0$ and the sum are defined as follows: First, $\nu_0$ is
the intersection number between $C$ and the $\theta\equiv 0$ and
$\theta\equiv \pi$ cases from Example 1.  This is to say that $-\nu_0$
counts the zeros of the section $(dt + i g^{-1} df) \wedge (\sin^2
\theta d\varphi + i g^{-1} dh)$ of $K$ with each zero contributing the
usual weight to the count.  As explained above, these weights are all
negative.  By the way, Lemma 4.5 insures that there are at most a
finite number of terms which enter into the definition of $\nu_0$.
Meanwhile, the sum in (4.26) is indexed by the elements of the
collection, $V$, of ends of $C$ which correspond to the $\theta_0 \in
\{0, \pi\}$ closed Reeb orbits in $C$'s limit set.  The weight
$\nu(\cE)$ of an end $\cE\in V$ accounts for the fact that the section
in $(dt + i g^{-1} df) \wedge (\sin^2 \theta d\varphi + i g^{-1} dh)$
of $K$ is not the correct section to use on $\cE$.  The preferred
section over $\cE$ has the form
\begin{equation}
(dt + i ds)\wedge (dx_1 - i dx_2) + \cO(|\sin \theta|) ,
\label{(4.27)}
\end{equation}
where $x_1 = \sin \theta \cos \varphi$ and $x_2 = \sin \theta \sin
\varphi$.  Actually, any section of $K|_\cE$ can be used as long as it
is homotopic to that in (4.27) through sections which do not vanish at
large values of $|s|$ on $\cE$.  This freedom to use homotopic
sections simplifies the computation for $\nu(\cE)$.

Given (4.26) and with $\nu_0$ understood to be non-positive and zero
if and only if $C$ avoids the $\theta\in \{0, \pi\}$ loci, then Lemma
4.4 becomes an immediate corollary to the following lemma.
\bigskip

\noindent
{\bf Lemma 4.6}\sl\qua {The number $\nu(\cE)$ in (4.26) is
constrained to obey}
\begin{itemize}
\item $v(\cE) \le -m_0(\cE)$ {when $\cE$ is in the concave side
of} $\mr\times (S^1 \times S^2)$.
\item $v(\cE) \le m0(e) - 1$ {when $\cE$ is in the convex side of
$\mr\times (S^1 \times S^2)$.}
\end{itemize}
\noindent
{Here, $m_0(\cE)$ is defined as in Proposition 4.1.}\rm
\bigskip

The remaining Steps 5--7 of the proof of Lemma 4.4 are devoted to the
following proof.\bigskip

\noindent
{\bf Proof of Lemma 4.6}\qua Note that the discussions below consider
only the case where $\theta_0 = 0$, for the case $\theta_0 = \pi$ is
identical save for some notation and sign changes.
\bigskip
	
\textbf{Step 5}\qua The end $\cE$ can be parameterized by coordinates
$(\tau, \rho)$ as in (2.14).  Remember that this parametrization has
$\tau\in\mr/(2\pi m \mz)$ and $r$ is either in $(-\infty, 1]$ or $[-1,
\infty)$ depending on whether $\cE$ sits in the convex or concave side
of $\mr\times(S^1 \times S^2)$.  Here, $m$ is the same as the number
$m(\cE)$ that is defined in Proposition 4.1.  The parametrization then
writes
\begin{equation}
(t = \tau, f = -e-^{-\sqrt{6}\rho}, a_1 = a_1(\tau, \rho), a_2 =
a_2(\tau, \rho))
\label{(4.28)}
\end{equation}
where $(a_1, a_2) = 6^{-1/4} e^{-\sqrt{6}\rho/2} h^{1/2} (\cos
\varphi, \sin \varphi)$.  In this regard, note that the pair of
functions $(x_1, x_2)$ which appear in (4.27) is related to $(a_1,
a_2)$ via
\begin{equation}
(a_1, a_2) = (x_1, x_2) + \cO(x_1^2 + x_2^2) .
\label{(4.29)}
\end{equation}

\textbf{Step 6}\qua To continue with the definition of $\nu(\cE)$, it
proves useful to select a function $\beta$ on $\mr$ with the following
properties:
 
\begin{itemize}
\item $\beta(\rho)=1$ for $\rho\ge 1$
\item $\beta(\rho)=0$ for $\rho\le -1$
\item $\beta'\le 0$.\autonum
\end{itemize}

Then, for $R \in\mr$, set $\beta_R(\rho) \equiv \beta(\rho - R)$.

A section of $K$ on $\cE$ can be written using $\beta$ which
interpolates between the `wrong' section, $(dt + i g^{-1} df) \wedge
(\sin^2 \theta d\varphi + i g^{-1} dh)$, near $\rho = 0$ and the
preferred section in (4.27) where $|\rho| \gg 1$.  Doing so finds that
the number $\nu(\cE)$ in (4.26) is equal to the usual algebraic count
of the number of the zeros of a certain complex valued function.  In
particular, by virtue of (4.29) and the previously mentioned homotopy
flexibility, the following functions suffice:
\begin{itemize}
\item $\beta_R (a_1 - i a_2) - (1 - \beta_R)$,\; {with any $R \gg
1$ when $\cE$ is in the concave side of $\mr\times(S^1 \times S^2)$.}
\item $(1 - \beta_R) (a_1 - i a_2) - \beta_R$,\; { with any
$\mr\ll \ll -1$ when $\cE$ is in the convex side of $\mr\times(S^1
\times S^2)$.}\autonum
\end{itemize}

A straightforward homotopy argument will verify that the algebraic
counting of the zeros of the functions in (4.31) depends only on the
winding number of the large $|\rho|$ version of the map from the
circle $\mr/(2\pi m \mz)$ to $\mc-\{0\}$ which sends $\tau$ to $(a_1 -
i a_2)|_{(\tau, \rho)}$.  In particular, if this function winds like
$e^{-i k\tau/m}$, with $k\in\mz$, then
\begin{itemize}\sl
\item $v(\cE) = -k$ {when $\cE$ is in the concave side of
$\mr\times(S^1 \times S^2)$.}
\item $v(\cE) = k$ {when $\cE$ is in the convex side of
$\mr\times(S^1 \times S^2)$.}\autonum
\end{itemize}

\noindent
The verification of (4.32) is left as an exercise save for the warning
to remember that the orientation of $\cE$ is defined by the
restriction of $-d\tau \wedge d\rho$.  Note also that (4.32)
constitutes a special case of some more general conclusions in Section
4 of \cite{HWZ2}.
\bigskip

\textbf{Step 7}\qua This step reports on the possibilities for the
winding numbers of $a_1 - i a_2$ on the large, constant $|\rho|$
circles in $\cE$.  In particular, the possibilities for the winding
number of $a_1 - i a_2$ are constrained by the fact that $(a_1, a_2)$
satisfy (2.15) and in addition
\begin{equation}
\lim_{ |\rho|\rightarrow\infty} (a_1^2 + a_2^2) = 0 .
\label{(4.33)}
\end{equation}
Note also that Lemma 4.5 and equation (4.29) insure that $(a_1^2 +
a_2^2)$ is never zero when $|\rho|$ is sufficiently large.
  
To see how these constraints arise, use the positivity at large
$|\rho|$ of $a_1^2 + a_2^2$ to write the complex number $a_1 - i a_2$
as $a_1 - i a_2 = e^{-v-i(w+k\tau/ m)}$; here $v$ and $w$ are smooth,
real valued functions of $(\tau, \rho)$ where $|\rho|$ is large; while
$k\in\mz$ is the winding number in question.  Note that $v$ has no
limit as $|\rho|$ tends to infinity as (4.33) implies that
\begin{equation}
\lim_{|\rho|\rightarrow\infty}\inf_\tau v(\rho, \tau ) = \infty .
\label{(4.34)}
\end{equation}
In addition, the top component of (2.15) implies that
\begin{equation}
v_\rho = w_\tau + k/m - \sqrt{3}/\sqrt{2} + c ,
\label{(4.35)}
\end{equation}
where $|c| \le e^{-v}$ at large values of $|\rho|$ by virtue of (4.33)
and Lemma 2.5.
  
To apply (4.35), first introduce $\uv(\rho)$ to denote the average of
$v$ over a large, but constant $\rho$ circle.  Likewise, introduce
$\uc$.  Both are smooth function of $\rho$ where $|\rho|R$ is large,
and (4.35) asserts that
\begin{equation}
\uv_\rho = k/m - \sqrt{3/2} + \uc.
\label{(4.36)}
\end{equation}	
In order to use (4.36) to obtain the concave side constraints in Lemma
4.6, first pick $\epsilon > 0$ and then some large $\rho_0$ so that
both $\uv$ and $\uc$ are defined where $\rho\ge\rho_0$ and so that
$|\uc| \le \epsilon$ for such $\rho$.  Then, take $R \gg 1$ and
integrate both sides of (4.36) from $\rho_0$ to $\rho_0 + R$.  The
result provides the inequality
\begin{equation}
\uv(\rho_0+R) - \uv(\rho_0) \le (k/m - \sqrt{3}/\sqrt{2} + \epsilon)
R.
\label{(4.37)}
\end{equation}
Now, as k/m is rational, $\sqrt{3}/\sqrt{2}$ is irrational, $m$ is
fixed by the end $\cE$ and $\epsilon$ can be made as small as desired
by taking $\rho_0$ large, the inequality in (4.37) is compatible with
condition in (4.34) only when $k > \sqrt{3} m/\sqrt{2}$.  As
$m_0(\cE)$ is defined to be the smallest integer that obeys this
inequality, it follows that $k > m_0(\cE)$ which, together with the
first line in (4.32), gives the first assertion of Lemma 4.6.

To consider the convex side constraints in Lemma 4.6, pick $\epsilon >
0$ but small as before, and then some very negative $\rho_0$ so that
now both $\uv$ and $\uc$ are defined where $\rho\le\rho_0$ and so that
$|\uc| \le\epsilon$ for such $\rho$.  Then, take $R \gg 1$ and
integrate both sides of (4.36) from $\rho_0 - R$ to $\rho_0$.  The
result provides the inequality
\begin{equation}
\uv(\rho_0) - \uv(\rho_0-R)\ge (k/m - \sqrt{3}/\sqrt{2} - \epsilon) R.
\label{(4.38)}
\end{equation}
Arguing as before finds (4.38) compatible with (4.34) only when $k <
\sqrt{3} m/\sqrt{2}$.  Thus, $k$ is at most one less than $m_0(\cE)$
which, together with the second line in (4.32), implies the second
assertion of Lemma 4.6.

\sh{(e)\qua Proof of Proposition 4.2, Part 2}

This part of the proof of Proposition 4.2 considers the cases where
there is a closed Reeb orbit in $C$'s limit set which has $\theta_0
\in \{0, \pi\}$.  In this regard, note first that such a subvariety
has $\chi(C_0) \le 0$ since the pullback of $dt$ to $C_0$ is not exact
on an end which approaches a $\theta_0 \in \{0, \pi\}$ orbit.  The
same argument implies that $C$ must have one or more limit set closed
Reeb orbits with $\cos^2 \theta_0 \not= 1/3$.

The next restriction is simply that if $\aleph$ in (4.14) is zero,
then $C$ is a cylinder which is described by (4.2).  To show that such
is the case, note first that if all elements of $C$'s limit set have
$\theta_0 \in \{0, \pi\}$, then $f < 0$ on the ends of $C$ and the
maximum principle requires the $f < 0$ condition to hold on the whole
of $C$.  In particular, $f$ is not zero on $C$, and so $\lambda\equiv
h/f$ is well defined on $C$.  Moreover, $h/f\sim -\sqrt{3}/sqrt{2}
\sin^2 \theta$ as $\theta$ nears either 0 or $\pi$, so $\lambda$ tends
to zero on all ends of $C$.  Meanwhile, the maximum principle applies
to (4.21) and forbids non-zero local maxima or minima.  Thus $\lambda
\equiv 0$ and so $h \equiv 0$ and $C$ is described by (4.2). Now
consider that Lemma 4.4 and the observation that $\chi(C_0)\le 0$
imply the following:

\narrower\sl
An irreducible, HWZ subvariety with a $\theta_0\equiv\{0,\pi\}$ closed
Reeb orbit in its limit set has index $(D) > \aleph + 1$ unless it is
an immersed cylinder which avoids both the $\theta_0 \in \{0, \pi\}$
loci.\autonumm%4.39
\endnarrower

\noindent
With regard to (4.39), here is the explanation for the assertion that
$C$ is immersed: In the case at hand, $\chi(C_0)$ must be zero for
index$(D)$ to equal $\aleph + 1$.  Thus, $C_0$ is a cylinder and so
$C$ has exactly two ends.  Moreover, as $dt$ is homologically
nontrivial on the end where $\theta$ tends to either 0 or $\pi$, it
must be so on the other.  Therefore, the other end can not have $cos^2
\theta_0 = 1/3$, and so $dt$ pulls back to sufficiently large, but
constant $|s|$ circles without zeros.  It then follows from (4.16)
that the number of zeros of $dt$'s pullback to $C_0$ is equal to
$-\chi(C_0)$, and is thus zero.  This last point implies that $C$ is
immersed.

To summarize the preceding, the cylinder $C$ must have one of its
limiting closed Reeb orbit with $\theta_0 \notin \{0, \pi\}$ and
$cos^2 \theta_0 \not= 1/3$.  This `other' closed Reeb orbit is
therefore characterized in part by the condition that $\varphi -
\frac{p'}{p}t $ is constant, where $p$' and $p$ are relatively prime
integers and $p \not= 0$.  This last point is important because if
$\varphi - \frac{p'}{p}t $ is constant on $C$ itself, then $C$ comes
from Example 7 when the $\theta_0 \in \{0, \pi\}$ closed Reeb orbit
corresponds to end of $C$ from the concave side of $\mr\times (S^1
\times S^2)$; otherwise $C$ comes from Example 5.  In this regard,
only the cases which are described in Proposition 4.2 have index$(D) =
\aleph + 1 = 2$; all of the others have index$(D) > 2$.  This last
assertion follows from Propositions 4.1 and (4.32).

Thus, a demonstration that $\varphi - \frac{p'}{p} t$ is constant on
$C$ completes the proof of Proposition 4.2 but for the immersion
remark in Part d of the third point.  The proof of the latter is
deferred to the next subsection while the remainder of this subsection
demonstrates that $\varphi - \frac{p'}{p} t$ is indeed constant on
$C$.

The constancy of $\varphi - \frac{p'}{p} t$ is considered below only
for the case where $C$ has a $\theta_0 = 0$ closed Reeb orbit.  As
before, the considerations for the $\theta_0 = \pi$ case are identical
in all essential aspects.

The demonstration starts with the annunciation of a maximum principle:

\narrower\sl
The restriction to $C$ of the multivalued function
$\varphi-\frac{p'}{p} t$ has neither local maxima
nor minima.\autonumm\endnarrower

\noindent
To prove this claim, note first that $dt$ pulls back to $C$ without
zeros.  Indeed, this follows from a combination of three facts: First,
$C$ is a cylinder so has zero Euler characteristic.  Second, $dt$
pulls back without zeros to all sufficiently large and constant
$|\rho|$ circles so its zeros with the appropriate integer weight
count $\chi(C)$.  Finally, all such weights are negative by virtue of
(4.16).  To complete the proof of (4.40), note that as $dt$ pulls back
without zeros to $C$, so $df$ does too and this allows $C$ to be
parameterized as in (2.19) in terms of functions $(x, y)$.  In this
regard, the constant $\varphi_0$ together with $p'$ and $p$ label the
$\theta_0 \notin \{0, \pi\}$ closed Reeb orbit in $C$'s limit set.  In
particular, up to a constant, $\varphi-\frac{p'}{p} t$ is the function
$x$.  The latter, with $y$, satisfies (2.20) and (2.20) implies the
second-order differential equation for $x$ in (4.18) to which the
maximum principle applies.

Now, the function $x$ tends to zero on the end of $C$ that corresponds
to the $\theta_0 \notin \{0, \pi\}$ closed Reeb orbit, so if $x$ tends
to zero on the other end of $C$ as well, then (4.40) establishes that
the combination $\varphi- \frac{p'}{p} t$ is constant on $C$.  Thus,
the demonstration now focuses on the behavior of $x$ on the end of $C$
which is near the $\theta = 0$ cylinder.  In this regard, remember
that this end of $C$ is characterized in part by the pair of integers
$(k, m)$, where $m$ is the absolute value of the degree of the
fundamental class of a sufficiently large, but constant $|\rho|$
circle in $H_1(S^1\times S^2; \mz)$ and $k$ is the winding number
defined so at all sufficiently large and constant $|\rho$|, the
$C$--valued function $a_1 - i a_2$ is homotopic to $e^{-ik\tau/m}$ as a
map from $\mr/(2\pi m \mr)$ to $C-\{0\}$.  Here, $C$ is parametrized
as in (4.28).

The integers $(k, m)$ can be identified as follow: The fact that $C$
is a cylinder and $dt$ is closed implies that $m = |p|$; and the fact
that $C$ has no intersections with the $\theta\in\{0, \pi\}$ cylinders
and $d\varphi$ is closed implies that $k = \text{sign}(p) p'$.

Next, write $a_1 - i a_2 = e^{-v-i(w+k\tau/m)}$ so that $v$ and $w$
are smooth, real valued functions of $(\tau, \rho)$ where $|\rho|$ is
large.  Note that $w = x + \varphi_0$.  Also, $v$ satisfies (4.34).
Now introduce a new function, $z(\tau, \rho)$, by writing $v = (k/m -
\sqrt{3}/\sqrt{2}) \rho + z(\tau, \rho)$ and note that (2.15) implies
the Cauchy-Riemann like equations

\begin{itemize}
\item $z_\rho = w_\tau + c$,
\item $w_\rho = - z_\tau + c'$.\autonum
\end{itemize}

Here, $c$ appears already in (4.35) while $c'$ is a smooth function
which satisfies similar bounds as $c$.  In particular, by virtue of
(2.17), Lemma 2.5 and Proposition 2.3's insurance for the exponential
decay at large $|\rho|$ of $(a_1^2 + a_2^2)$, these functions at
sufficiently large $|\rho|$ obey
\begin{equation}
|c| + |c'| \le e^{-\delta|\rho|} ,
\label{(4.42)}
\end{equation}
with $\delta > 0$ some $\rho$ independent constant.

	As is demonstrated below, it is a moment's investment to
establish from (4.34), (4.41) and (4.42) that $w$ (and also $v$) have
constant limits as $|\rho|$ goes to infinity.  Given that such is the
case, it follows from the fact that $\ux$ in (2.23) is constant that
the limit of $w$ is $\varphi_0$.  As $w = x + \varphi_0$, this means
that $x$ vanishes asymptotically on the end of $C$ near the $\theta =
0$ cylinder, as required.

The argument that $(z, w)$ has a constant limit as $|\rho|$ goes to
infinity is given here for the concave side case only, as the other
case is settled with the identical argument up to a sign change or
two.  To proceed in this case, first note that the average values on
the constant r circles of $z$ and $w$ have limits as
$\rho\rightarrow\infty$ as can be seen by integrating both sides of
(4.41) over the constant $\rho$ circles.  Thus, the case is settled
with a proof that the remaining parts of $x$ and $w$ decay to zero at
large $\rho$.  For this purpose, introduce the $\zeta' = \zeta = 0$
version of the operator $L_0$ in (2.16) and let $f^+$ and $f^-$ denote
the functions of $\rho$ (for $\rho$ large) whose values are the
respective $L^2$ norm of the $L^2$ orthogonal projection of the column
vector with top component $z$ and bottow component $w$ onto the span
of the eigenvectors of $L_0$ with positive $(+)$ and negative $(-)$
eigenvalues.  Then (4.41) and (4.42) imply that
\begin{equation}
f^-_\rho \ge f^-/N - e^{-\delta\rho}\quad \text{and}\quad f^+_\rho \le
-f^+/N + e^{-\delta\rho}.
\label{(4.43)}
\end{equation}
Integration of the left most equation in (4.43) finds the dichotomy:
Either $f^-(\rho) \ge \nu^- e^{\rho/N}$ with $\nu^- > 0$ or else
$f^-(\rho) \le \nu^- e^{-\delta\rho}$.  Meanwhile, integration of the
right most equation in (4.43) finds that $f^+$ has no choice but to
decay exponentially fast as $\rho$ tends to infinity.  However,
exponential growth of $f^-$ is forbidden by (4.34) as the part of $(v,
w)$ that contributes to $f^-$ integrates to zero around any constant
$\rho$ circle and thus exponential growth of $f^-$ forces
exponentially large values on $-v$ at places on each large r circle.
	      
\sh{(f)\qua Proof of Proposition 4.2, Part 3}

The proof of Proposition 4.2 is completed here with a proof of the
following formal restatement of a portion of Part d of Proposition
4.2's third point:
\bigskip

\noindent
{\bf Proposition 4.7}\qua {\sl An HWZ pseudoholomorphic subvariety
$C$ whose model curve, $C_0$ is a thrice-punctured sphere with
index$(D) = \aleph + 1$ is the image of $C_0$ via an immersion.}
\bigskip

\noindent
{\bf Proof of Proposition 4.7}\qua To begin, consider the zeros of the
pullback to $C_0$ of $\kappa_t dt - \kappa_\varphi d\varphi$ in the
case where $\kappa_t$ and $\kappa_\varphi$ are constant real numbers.
Of course, this pullback has a zero at any local singular point of the
tautological map to $X$ since all pullbacks vanish at such points.  In
any event, each zero of the pullback of $\kappa_t dt - \kappa_\varphi
d\varphi$ counts with a negative weight when used in an `Euler class'
count.  A proof of this last assertion uses (4.16) and (4.17) to
derive a second order differential equation for the pullback to $C_0$
of $\kappa_t t - \kappa_\varphi\varphi$ to which the maximum principle
applies and rules out local extrema.  By the way, an argument near
points where $dt \not= 0$ can also be made directly from (4.18) since
the corresponding $x$ in the case of (4.18) where $p'/p =
\kappa_t/\kappa_\varphi$ differs locally by an additive constant from
a non-zero multiple of $\kappa_t t - \kappa_\varphi\varphi$.
   
Not all values of $\kappa_t$ and $\kappa_\varphi$ provide a form whose
zero count produces $\chi(C_0)$, hence the quotes in the preceding
paragraph around the words `Euler class.' Indeed, the equality with
$\chi(C_0)$ of the algebraic counting of the zeros of the pullback of
$\kappa_t dt - \kappa_\varphi d\varphi$ can be guaranteed only when
$\kappa_t/\kappa_\varphi \not= p'/p$ for all pairs $(p, p')$ which
come via (1.8) from the closed Reeb orbits in $C$'s limit set.  To
explain, this constraint on $\kappa_t/\kappa_\varphi$ arises precisely
because the form $p'/p dt - d\varphi$ pulls back as zero on the closed
Reeb orbit which supplied the pair $(p, p')$.  In particular, the
condition $\kappa_t/\kappa_\varphi \not= p'/p$ for all such pairs
guarantees that the pullback of $\kappa_t dt - \kappa_\varphi
d\varphi$ to all sufficiently large, constant $|s|$ circles in $C_0$
has no zeros.  When such is the case, a standard argument proves that
a multiplicity weighted count of the zeros of $\kappa_t dt -
\kappa_\varphi d\varphi$ yields $\chi(C_0)$.
 
Having digested the preceding, fix $\epsilon\not= 0$ but small and
take $(p'/p + \epsilon) dt - d\varphi$ where the pair $(p, p')$ comes
from one of the closed Reeb orbits in $C$'s limit set where $cos^2
\theta_0 \not= 1/3$.  The pullback of this form has norm
$\cO(\epsilon)$ on an end of $C$ which approaches a closed Reeb orbit
that supplies the pair $(p, p')$, and it is relatively large,
$\cO(1)$, on other ends of $C$.  In any event, the argument from the
preceding paragraph applies here and explains why the pullback of
$(p'/p + \epsilon) dt - d\varphi$ to $C_0$ has exactly one zero
provided that $\epsilon$ is small in absolute value but not zero.
Thus, if the tautological map from $C_0$ to $X$ has a singular point,
then the pullback of $(p'/p + \epsilon) dt - d\varphi$ vanishes only
at this point.  In particular, this pullback cannot have a zero on any
end of $C$.  However, just such a zero is exhibited below, and so the
tautological map from $C_0$ to $X$ lacks local singular points.

To exhibit the asserted zero, focus on an end $\cE\subset C$ which
approaches a closed Reeb orbit that supplies $(p, p')$, and
parameterize said end by coordinates $(\tau, u)$ as in (2.19) and
(2.20).  Thus, $\tau\in\mr/(2\pi m|p|)$ and either $|u| \in
[u_0,\infty)$ or $|u| \in (0, u_0]$ depending on whether $\cE$ is on
the convex or concave side of $\mr\times (S^1 \times S^2)$.  In terms
of this parameterization, $(p'/p + \epsilon) dt - d\varphi$ pulls back
as
\begin{equation}
(\epsilon - x_\tau) d\tau - x_u du ,
\label{(4.44)}
\end{equation}
and thus it vanishes on $\cE$ only at points where $x_u = 0$ and where
$x_\tau = \epsilon$.
  
The constraint $x_u = 0$ is satisfied at two or more points on every
constant $u$ circle since the first line in (2.20) identifies these
points with the critical points of $y$'s pullback to such a circle.
Meanwhile, $|x_\tau|$ limits to zero as $|s| \rightarrow\infty$ on
$\cE$.  Therefore, as $\epsilon$ can be as small as desired and chosen
either positive or negative as desired, the vanishing of $\epsilon +
x_\tau$ occurs on the $x_u = 0$ locus for arbitrarily small but
non-zero choices of $\epsilon$.  (Keep in mind here that the
simultaneous zeros of $x_\tau$ and $x_u$ are isolated, otherwise $x$
would be constant and $C_0$ would be a cylinder from one of the
examples in Section 4a.)  The preceding argument proves that there is
exactly one solution to (4.44).

\sh{(g)\qua The cokernel of $D$}

Future constructions with Proposition 4.2's subvarieties may simplify
with the knowledge that the corresponding operator $D$ has, in all
cases, trivial cokernel.  A formal statement of this assertion appears
below, and its proof occupies the remainder of this subsection.
\bigskip

\noindent
{\bf Proposition 4.8}\qua {\sl The cokernel of $D$ is trivial if $C$
is an irreducible, pseudoholomorphic HWZ subvariety with {\rm
index}$(D) \le \aleph + 1$.}
\bigskip

The remainder of this subsection is occupied with the following proof.
\bigskip

\noindent
{\bf Proof of Proposition 4.8}\qua All cases save the thrice-punctured
spheres mentioned in Part d of the third point follow directly from
\bigskip

\noindent
{\bf Lemma 4.9}\qua {\sl Let $C$ be an irreducible, HWZ subvariety
that is invariant under some $1$--parameter subgroup of the group $T$.
Then, the corresponding operator $D$ has trivial cokernel.}
\bigskip

\noindent
{\bf Proof of Lemma 4.9}\qua If $C$ is fixed by a circle subgroup in
$T$, then the corresponding operator $D$ is equivariant and its
analysis can be simplified with the help of a separation of variables
strategy.  This is to say that $D$ preserves the character eigenspaces
of the circle's action on the domain and range, and the restriction of
$D$ to such an eigenspace reduces the partial differential equation
$D^*\eta = 0$ to a first-order ODE.  This last reduction, plus some
timely applications of the maximum principle prove the asserted
triviality of $D$'s cokernel in each of the cases in Section 4a.  The
details here are straightforward and left to the reader.
\bigskip

With this last lemma in hand, the only remaining cases for Proposition
4.8 are the thrice-punctured spheres from Part d of the third point in
Proposition 4.2.  In this regard, remember that these spheres are
immersed; and remember that an immersed HWZ subvariety has a
well-defined normal bundle and that the operator $D$ is a differential
operator, as in (3.5), on the space of sections of said normal bundle.

The preceding understood, the argument for a thrice-punctured sphere
case given below is a generalization of the argument introduced by
Gromov \cite{Gr} when considering pseudoholomorphic planes and spheres, and
it also has antecedents in some of HWZ's work.  The argument begins
with the observation that for such $C$, there exist vectors $\eta\in
\ker (D)$ whose restriction to each end of $C_0$ does not have limit
zero as $|s| \rightarrow\infty$.  Indeed, to find such an element,
note that $C$ is not preserved by any 1--parameter subgroup of $T$, and
thus the infinitesimal version of the $T$ action through its
generators produces a two-dimensional subspace, $V$, in the kernel of
$D$.  Moreover, as only the closed Reeb orbits with $\theta_0 \in \{0,
\pi\}$ are completely $T$ invariant, and as no such orbit appear in
$C$'s limit set, so the generic element in $V$ has the desired
property.
 
In fact, $V$ has the following `universal limit' property: Given
$\eta\in \text{kernel}(D)$ and an end $\cE\subset C_0$, there exists a
unique element $v \in V$ such that $\eta - v$ has limit zero as $|s|
\rightarrow\infty$ on $\cE$.  Indeed, this follows because each closed
Reeb orbit with $\theta_0\notin \{0, \pi\}$ has precisely a
one-dimensional family of deformations and the latter is the orbit of
a 1--parameter subgroup of $T$.
      
The next observation is that a vector $\eta\in \text{kernel}(D)$ with
non-zero limit on all three ends of $C_0$ as $|s| \rightarrow\infty$
must have exactly one zero on $C_0$, and a non-degenerate one at that.
Indeed, such a vector $\eta$ is, \emph{a priori}, a section of the
normal bundle, $N$, of $C$.  Moreover, as $\eta$ can be approximated
at large $|s|$ on each end of $C$ by a vector from $V$, so at large
$|s|$, $\eta$ is both non-zero and homotopic through non-zero vectors
to the section from Part 2 of Section 3a which is used to define the
expression $\la e, [C]\ra$ in Proposition 3.1.  Thus, a count of the
zero's of $\eta$ with the appropriate multiplicities computes the
expression $\la e, [C]\ra - 2 m_C$ which appears in Proposition 3.1.
Moreover, each such zero counts with positive weight, its order of
vanishing, as a consequence of its annihilation by $D$ in (3.5).
Given all of the above and the fact that $\chi(C_0) = -1$, the formula
in Proposition 3.1 can hold if and only if $\eta$ has precisely one
zero, and this zero has multiplicity one.
  
With the preceding understood, suppose now that the dimension of the
cokernel of $D$ is positive.  This implies that the dimension of the
kernel of $D$ is at least $\aleph + 2$, and, as is shown next, such a
condition leads to the absurd conclusion that this kernel has a vector
$\eta$ with non-zero $|s| \rightarrow\infty$ limit on each end of
$C_0$ and with at least two zeros.  To view this vector in the $\aleph
= 3$ case, start with the observation that $D$'s kernel has at least
five linearly independent vectors when its cokernel is nontrivial.
Choose any $\eta$ as in the preceding paragraph and let $z \in C_0$
denote its one zero.  As dim(kernel$(D)) \ge 5$, there is a
two-dimensional subspace in kernel$(D)/(\mr\eta)$ of vectors which
vanish at $z$.  Let $W \subset \text{kernel}(D)$ project
isomorphically onto such a subspace.  As $W$ is two-dimensional, there
is, given any $z' \not= z$, a vector $\eta_{z'}\in W$ such that $\eta
+ \eta_{z'}$ vanishes at both $z$ and $z'$.  So, to avoid the desired
contradiction, each such $\eta + \eta _{z'}$ must have zero $|s|
\rightarrow \infty$ limit on some end of $C_0$.
  
To see that such behavior is absurd, note that a non-vanishing limit
is an open condition on the kernel of $D$ and so there must be an end
of $C_0$ with the property that each such $\eta+ \eta_{z'}$ has zero
$|s| \rightarrow\infty$ limit on this end.  This last conclusion
cannot occur unless $W$ has a basis, $\{\eta_1, \eta_2\}$, such that
$\eta$ can be written as $\eta_1 + \beta \eta_2$ where $\eta_1 - \eta$
and $\eta_2$ have zero $|s| \rightarrow\infty$ limit on this end.
Here, $\beta$ is a smooth, real-valued function on $C_0$.  Moreover,
as the operator $D$ annihilates $\eta$, $\eta_1$ and $\eta_2$, the
function $\beta$ must obey the equation $\beta = 0$ on $C_0$.
However, the only real valued functions with this property are the
constants, and these are ruled out since $\eta\notin W$.

Consider next the case where $\aleph = 2$.  In this case, the
non-triviality of the cokernel of $D$ implies that the kernel of $D$
has dimension 4 and so now the analogous vector space $W$ may only be
one-dimensional.  In any event, choose $\eta$ as before, to have
non-vanishing $|s| \rightarrow\infty$ limit on each end of $C_0$, and
let $\eta'$ be a non-trivial section of $W$.  Let $\cE\subset C_0$
denote the concave side end and the following is true:

\narrower\sl
There exists a set where all points are accumulation points,
$s$ is unbounded, and $\eta$ is proportional to $\eta'$.\autonumm
\endnarrower

Accept (4.45) and the $\aleph = 2$ case follows unless $\eta'$ is a
constant multiple of $\eta$ along this set.  Of course the latter
would imply that $\eta\in W$ since a non-trivial element in the kernel
of $D$ has isolated zeros.  (In fact, the set in question can be shown
to be a piecewise smooth curve.)

The proof of (4.45) begins with the observation that the operator $D$
at large $s$ on $C_0$ can be viewed as an operator on $\mc$--valued
functions of coordinates $\tau$ and $s$ with $\tau$ periodic which has
the form in (3.9) with $\rho = s$ and with $A_0$ the identity matrix.
By viewing $D$ in this way, both $\eta$ and $\eta'$ become
$\mc$--valued functions.  In this regard, $\eta$ at large $s$ has the
form
\begin{equation}
\eta = r_0 e_0 + \kappa ,
\label{(4.46)}
\end{equation} 
where $r_0$ is a constant, non-zero real number and $|\kappa| \le
e^{-\delta s}$ at large $s$ with $\delta$ a positive constant.
  
Meanwhile, $\eta'$ can be written as
\begin{equation}
\eta' = r'_0 \eta + \lambda ,
\label{(4.47)}
\end{equation}
where $r'_0$ is a constant real number and where $|\lambda| \le
e^{-\delta's}$ at large values of $s$ with $d' > 0$ being constant.
By assumption, $\lambda$ is not identically zero.  Now, $\lambda$ can
vanish only on a finite or countable set where $s \ge 0$ since it can
have only a finite number of zeros on any compact set.  Thus, there is
some $s_0 > 0$ and a countable or finite set $U\subset [s_0, \infty)$
such that $\lambda$ maps the constant $s \in [s_0, \infty)-U$ circles
in the $(\tau, s)$ cylinder to $\mc-\{0\}$.  As such, $\lambda$ has a
winding number, the latter defines a locally constant function on
$[s_0, \infty)-U$ and the only way (4.45) can fail is if this winding
number is zero on all components of $[s_0, \infty)-U$ where $s$ is
sufficiently large.  Indeed, if $\lambda$ has non-zero winding number
on a sufficiently large, but constant $s$ circle, then it follows from
(4.46) that $\lambda$ and $\eta$ must be colinear at no fewer than two
points on such a circle.  This last fact implies (4.45).

Thus, to prove (4.45), it is enough to prove that $\lambda$ has
non-zero winding number on all sufficiently large and constant $s \in
[s_0,\infty)-U$ circles.  In this regard, it can be proved that $U$ is
infact finite since $\lambda$ cannot have an infinite number of zeros
without vanishing all together.  However, as the proof of this last
assertion is much longer than the proof that an infinite number of
positive $s$ zeros of $\lambda$ implies (4.45), the proof of the
latter claim follows.  For this purpose, assume for the moment that
$\lambda$ actually has an infinite number of $s \ge 0$ vanishing
points.
  
To see how this last assumption leads to (4.45), note first that there
is in this case at most one connected component of $[s_0, \infty)-U$
where the winding number is zero.  This is because the zeros of
$\lambda$ all occur with positive multiplicity and thus the winding
number changes in a monotonic fashion between consecutive components
of $[s_0, \infty)-U$.
 
In particular, there exists some $s_1 \ge s_0$ such that the winding
number on every constant $s\in [s_1, \infty)-U$ circle is non-zero.

Now assume that $\lambda$ has only a finite set of positive s zeros.
In this case, the winding number of $\lambda$ is defined on all
sufficiently large and constant $s$ circles.  Then, for the sake of
argument, suppose that this winding number is zero.  Under this
assumption, view $\lambda$ as a $\mc-\{0\}$--valued function and
introduce real-valued functions $v$ and $w$, defined at large $s$ on
the $(\tau, s)$ cylinder by writing $\lam = e^{-v-iw}$.  These
functions then obey (4.41) where $\rho\equiv s$ and where $c$ and $c'$
obey (4.42).  In addition, the condition that $|\lam| \ge e^{-\delta
s}$ at large s forces the condition
\begin{equation}
\lim_{\rho\rightarrow\infty} \text{inf}_\tau v(\rho, \tau) = \infty.
\label{(4.48)}
\end{equation}
However, as demonstrated using (4.43), this last condition is
incompatible with (4.41) and (4.42).  Thus, $\lambda$ cannot have zero
winding number on all sufficiently large and constant $s$ circles.

\section{The structure of the $\aleph = 2$, thrice-punctured\nl
sphere moduli space}

The purpose of this section is to give a complete description of the
moduli space of thrice-punctured spheres with two convex side ends
that arise in Part d of the third point in Proposition 4.2.  In
particular, the arguments given here establish Theorem A.2.  The
following proposition provides a restatement of Theorem A.2:
\bigskip

\noindent 
{\bf Proposition 5.1}\sl\qua {The components of the moduli space of
thrice-punctured, $\aleph = 2$ spheres which arise in Part d of the
third point in Proposition 4.2 can be put in $1--1$ correspondence with
the sets of two ordered pairs of integers, $\{(p, p'), (q, q')\}$,
which obey}
\begin{itemize}
\item $\Delta\equiv p q'- q p' > 0$.
\item $q'- p' > 0$ {unless both are non-zero and have the same
sign}.
\item {If $(m, m')\in\{(p, p'), (q, q')\}$ and if $|m'/m| <
\sqrt{3}/\sqrt{2}$, then $m > 0$.  On the other hand, if $m < 0$, then
$|m'/m| > \sqrt{3}/\sqrt{2}$.}
\end{itemize}

\noindent
{Moreover, the component that corresponds to a given set $I
\equiv\{(p, p'), (q, q')\}$ is a smooth manifold which is $\mr\times
T$ equivariantly diffeomorphic to $\mr\times T$.}\rm
\bigskip 

Subsection 5a, below, explains how such sets $I$ of integers are
associated to the moduli space components and derives the constraints
on those sets which arise.  The second subsection proves that each set
$I$ of four integers can be associated to at most one moduli space
component.  A proof is also given for the assertion that the
associated moduli space, for a given $I$ is either empty or
diffeomorphic to $\mr\times T$.  Subsections 5c-g are occupied with
the proof of the assertion that every set of $I = \{(p, p'), (q,
q')\}$ that satisfies the constraints has an associated moduli space
component.
  
With the proof of Proposition 5.1 complete, the final subsection
provides a formula for the number of double points of Proposition
5.1's subvarieties in terms of the corresponding set $I$.  Proposition
5.9 summarizes the latter and its assertions directly imply the part
Theorem A.4 that concerns Theorem A.2's subvarieties.

\sh{(a)\qua Constraints on thrice-punctured spheres with one concave side
end}

A pseudoholomorphic, thrice-punctured sphere $C \subset\mr\times(S^1
\times S^2)$ with one concave side end and no intersections with the
$\theta\in\{0, \pi\}$ locus determines a set of three pair of
integers, $\{(p, p'), (q, q'), (k, k')\}$, in a manner that will now
be described.
  
To start, remember that the large and constant $|s|$ slices of $C$
consist of a disjoint union of three embedded circles, with one on the
concave side of $\mr\times(S^1 \times S^2)$ (where $s$ is positive)
and two on the convex side.  And, as $|s|$ tends to infinity, these
constant $|s|$ circles converge pointwise as multiple covers of the
closed Reeb orbits which comprise $C$'s limit set.  In particular,
\begin{equation}
m \varphi - m' t = \text{constant} + \cO(e^{-|s|/\zeta}) \mod(2\pi)
\label{(5.1)}
\end{equation}
on each such circle; here $\zeta > 0$ is constant while $m$ and $m'$
are integers that are associated to the given end of $C$.  In fact,
the limiting closed Reeb orbit in question is determined in part using
the pair of integers in (1.8) provided by the quotient of $(m, m')$ by
their greatest common divisor.  The multiplicity of covering over the
closed Reeb orbit is then equal to this greatest common divisor of $m$
and $m'$.  Finally, the signs of $m$ and $m'$ are fixed by the
following convention: Take the signs of $m$ and $m'$ (when non-zero)
to equal the signs of the respective restrictions of $f$ and $h$ to
the closed Reeb orbit in question.  In this regard, note that $m = 0$
if and only if $f$ restricts as zero to the closed Reeb orbit, and
likewise $m' = 0$ if and only if $h$ restricts as zero to the closed
Reeb orbit.  Also, note that (2.7) guarantees the compatibility of
this sign determination with (5.1).  The constraint for the third
point in Proposition 5.1 arises from this use of integer pairs $(m,
m')$ to parameterize the closed Reeb orbits in $S^1 \times S^2$.
  
In this way, $C$ determines the set $\{(p, p'), (q, q'), (k, k')\}$ of
three pair of integers.  The convention here is that the third pair
listed, $(k, k')$, comes from the concave side end of $C$.  Meanwhile,
the order of appearance of the first two pair has, as yet, no
intrinsic significance since this order corresponds to an arbitrary
labeling of the convex side ends of $C$.  However, the first two pair
will be ordered shortly.

With the preceding understood, this subsection finds necessary
conditions for a set $\{(p, p'), (q, q'), (k, k')\}$ to arise from a
pseudoholomorphic, thrice-punctured sphere with one concave side end
and no intersections with the $\theta\in \{0, \pi\}$ locus.
   
To begin, note that there is an evident first constraint:
\bigskip

{\bf Constraint 1}\qua {\sl No pair in $\{(p, p'), (q, q'), (k,
k')\}$ can vanish identically.}
\bigskip

\noindent  
There is also a second constraint which comes from (5.1) and the fact
that $d\varphi$ and $dt$ both restrict to $C$ as smooth, closed forms:
\bigskip

{\bf Constraint 2}\sl\qua $k = p + q$ {and} $k' = p' + q'$.\rm
\bigskip

(The reader is left with the task of verifying that this constraint is
consistent with the given sign conventions.)  Note that the
constraints listed in Proposition 5.1 for $\{(p, p'), (q, q')\}$ imply
the third constraint for $(m, m') = (k, k')$.  Thus, no new
constraints appear with the association of $(k, k')$ to an end of an
$\aleph = 2$, thrice-punctured sphere.

There are additional constraints.  The next one involves the integer
\begin{equation}
\Delta\equiv pq' - qp' .
\label{(5.2)}
\end{equation}
and asserts:
\bigskip

{\bf Constraint 3}\qua $\Delta\not= 0$.
\bigskip

Indeed, suppose, to the contrary, that $\Delta= 0$.  Now, both $p$ and
$q$ can't vanish as then all ends of $C$ would have $\cos^2 \theta_0 =
1/3$ closed Reeb orbit limits and an argument from Section 4 proved
this impossible.  In addition, the vanishing of $\Delta$ implies that
$p'/p = q'/q = (p' + q')/(p + q)$ so none of $p$, $q$ or $p + q$ can
vanish.  Moreover, the equality of these ratios implies that there are
at most two values for the $|s|\rightarrow \infty$ limits of $\theta$
on $C$.  However, as seen in Section 4, this cannot happen unless $C$
is a cylinder.
 
By the way, given that $\Delta\not= 0$ and that $\Delta$ changes sign
upon interchanging $(p, p')$ with $(q, q')$, an ordering of these
pairs is unambiguously defined by requiring $\Delta$ to be positive.
This last convention is implicit in all that follows.

Of course, with $\Delta > 0$ there are obvious constraints on the
relative signs between the four integers $p$, $p'$, $q$ and $q'$ that
arise just from the definition of $\Delta$.  However, a less obvious
constraint is:
\bigskip

{\bf Constraint 4}\qua $q' - p' > 0$
{\sl unless both are non-zero and have the same sign.}
\bigskip

\noindent 
The proof of this constraint is quite lengthy, so is broken into nine
steps.
\bigskip

{\bf Step 1}\qua This first step proves the following assertions:
\begin{itemize}\sl
\item {If $p' + q' = 0$, then $p' < 0$.}
\item {If $p' = 0$, then $q' > 0$.}
\item {If $q' = 0$, then $p' < 0$.}
\end{itemize}
\bigskip

\noindent
To see why the first point holds, suppose, to the contrary that $p' >
0$.  Then, as $\Delta$ in this case is $-p' (p + q)$ and is positive,
so $k = p + q < 0$.  This implies that $f < 0$ on the concave end of
$C_0$.  However, if this is the case, then $h$ cannot vanish on the
concave end of $C_0$ because $h$ is zero and $f$ is negative if and
only if $\theta\in\{0, \pi\}$.  Thus, because $h$ tends to zero as $s
\rightarrow\infty$ but is nowhere zero at large $s$, there are, given
$R > 0$, compact components of constant $h$ level sets that lie where
$s > R$.  However, this is absurd since (4.17) guarantees that
$d\varphi$ has non-zero integral over such a level set while the $p' +
q' = 0$ condition guarantees that $d\varphi$ is exact where $s$ is
large on $C$.

To see why the second point above holds, note that $\Delta = p q'$ in
this case.  Thus, positivity of $\Delta$ requires positivity of $q'$
or else $p$ would be negative.  If $p$ were negative, then $f$ would
be negative on the $(p, p')$ end of $C$ and as $p' = 0$, this end of
$C$ would be asymptotic as $|s| \rightarrow\infty$ to a component of
the $\theta = \{0, \pi\}$ locus.
  
A similar argument establishes that $p'$ is negative when $q' \in 0$.
\bigskip

{\bf Step 2}\qua This step proves Constraint 4 with the extra
assumption that $p$ and $q$ cannot both be either strictly positive or
strictly negative.  Indeed, if both are strictly positive, $p' > 0$
and $q' < 0$, then $\Delta < 0$.  On the other hand, if both are
strictly negative, then $f < 0$ on $C_0$ while $h$ changes sign.
Thus, the $h = 0$ locus is non-empty and occurs where $f$ is negative,
which is precluded since $C$ does not intersect the $\theta\in\{0,
\pi\}$ locus.
\bigskip

{\bf Step 3}\qua This step proves Constraint 4 with the added
assumption that one of $p$, $q$ or $p + q$ is zero.  To start, assume
that $p = 0$, $p' > 0$ and $q' < 0$.  Positivity of $\Delta$ then
requires that $q < 0$ and so $f$ is negative at large $s$ on the $(q,
q')$ end of $C$ and also at large $|s|$ on the concave (that is, $(k,
k')$) end of $C$.  In fact, as is argued momentarily, $f$ is strictly
negative on $C$.  Granted this claim, the assumed violation of
Constraint 4 results in the following absurdity: As $h$ is positive on
the $(p, p')$ end of $C_0$ and negative on the $(q, q')$ end, so the
$h = 0$ locus is non-empty.  But, as $f < 0$, this means that $C$
intersects the $\theta\in \{0, \pi\}$ locus.
 
To see that $f < 0$ on $C$, suppose not.  Then $f$ has non-trivial,
positive regular values (by the maximum principle).  Moreover, as $f$
is negative on the $(q, q')$ and concave ends of $C$, there exists,
given $R \ge 1$, such a regular value whose level set sits entirely
where $|s| \ge R$ in the $(p, p')$ end of $C$.  Because $dt$ restricts
as a non-zero form on such a level set, by virtue of (4.16), and
because $dt$ is exact at large $s$ on the $(p, p')$ end of $C$ when $p
= 0$, so this level set cannot be compact.  In fact, it must be a
properly embedded copy of $\mr$ with each end leaving via the $(p,
p')$ end of $C$.  But this last conclusion is absurd as $dt$ is
nowhere zero on such a level set yet (5.1) asserts that $t$ tends to a
constant as $|s|$ tends to infinity on the $(p, p')$ end of $C$.

Now assume that Constraint 4 is violated when $q = 0$.  As $q' < 0$
and $\Delta > 0$, this requires $p < 0$ and the preceding argument
applies with insignificant changes to find $f < 0$ on $C$.  Again, as
$h$ must change sign on $C$, such a conclusion is absurd.

Finally, assume that Constraint 4 is violated when $p + q = 0$.  Now,
the $\Delta > 0$ condition requires that $p$ and $p' + q'$ have the
same sign.  The argument given below takes this sign to be positive.
But for some straightforward sign changes, the same argument also
handles the case where this sign is negative.

To start, note that when $p > 0$, then $f$ is positive on the $(p,
p')$ end of $C$ and negative on the $(q, q')$ end.  Moreover, $f$ has
no finite limits on either end, as it tends uniformly to $\infty$
because $|s| \rightarrow\infty$ on the $(p, p')$ end and to $-\infty$
on the $(q, q')$ end.  Meanwhile, even as $f$ limits to 0 as $s
\rightarrow\infty$ on the concave side end of $C$, this function must
take negative and positive values at arbitrarily large values of $s$.
Indeed, were it strictly negative at very large $s$, then a component
of the level set of some very small, negative regular value of $f$
would be compact and lie entirely where $s$ is very large.  And, as
the integral of $dt$ over such a level set could not be zero (due to
(4.16)), the existence of such a level set would run afoul of (5.1).
An analogous argument explains why $f$ cannot be strictly positive at
large $s$.
  
With the preceding understood, suppose that $-\epsilon$ is a regular
value of $f$ with $\epsilon$ positive and very small.  For tiny
$\epsilon$, a component of this level set will extend far down the
concave side end of $C$.  In particular, as $p' + q' > 0$, when
$\epsilon$ is small, the function $h$ will be positive on some of this
component.  Then, this component sits entirely where $h > 0$ in $C$
because $\theta\notin \{0, \pi\}$ on $C$.  Now, consider increasing
$\epsilon$ and viewing the behavior of the $f = -\epsilon$ level set.
In particular, when $\epsilon\gg 1$, then this level set necessarily
sits entirely in the $(q, q')$ end of $C$ where $h < 0$.  Thus, as no
component of this level set is null-homologous, there is some
intermediate values of $\epsilon$ where a component of the $f =
-\epsilon$ level set intersects both the $h > 0$ region and the $h <
0$ region.  Of course, such an event is absurd for $\theta$ would take
value 0 or $\pi$ on $C$.
\bigskip

{\bf Step 4}\qua The subsequent steps prove the remaining cases of
Constraint 4 by establishing the following:
\narrower\sl
If none of $p, p', q, q', p + q,$ and $p' + q'$
are zero and if the signs of $p'$ and $q'$ are opposite,
then $p' < 0.$
\autonumm\endnarrower

To begin the justification for (5.3), use Step 2 to conclude that if
both $p$ and $q$ are non-zero, then one is positive and the other
negative.  Thus, $f$ changes sign on $C$.  Moreover, where $s$ is very
negative (thus, on the convex side of $\mr \times (S^1 \times S^2)$),
the function $f$ tends to infinity on one end of $C$ and minus
infinity on the other.  In particular, when $R \gg 1$, then the $f =
R$ locus is an embedded circle in the convex side end of $C$ that
corresponds to the positive member of the pair $(p, q)$.  Meanwhile,
the $f = -R$ locus is likewise a circle on the convex side end of $C$
that corresponds to the negative member of the pair $(p, q)$.  Now,
(4.16) implies that $f$ and $t$ have the same critical points, and so
Section 4f proves that there is only one critical point of $f$ on $C$,
and thus only one critical value.  The sign of this critical value
depends on the sign of $p + q$.  In particular, if $p + q > 0$, then
this critical value is positive, and if $p + q < 0$, then this
critical value must be negative.  The maximum principle is involved
here, since each non-critical level set of $f$ is either a single
embedded circle or a pair of embedded circles.  Indeed, when $p + q >
0$, then the $f = 0$ and the very negative and constant $f$ loci must
be isotopic in $C$ to accommodate the maximum principle.  In this
case, the $f \le 0$ portion of $C$ is an infinite half cylinder.  On
the other hand, if $p + q < 0$, then the $f = 0$ locus and the very
positive and constant $f$ loci in $C$ must be isotopic and so now the
$f \ge 0$ portion of $C$ is a half infinite cylinder.  The analogous
conclusions hold for the $h = 0$ locus with $p$ and $q$ replace by
$p'$ and $q'$.
\medskip

{\bf Step 5}\qua With the preceding understood, consider the possible
location of the $h = 0$ locus in a hypothetical case where (5.3) is
violated.  When $p + q < 0$ this locus lies in the half cylinder where
$f \ge 0$, and this implies that either the $h \ge 0$ locus or the $h
\le 0$ locus is a subcylinder of the $f \ge 0$ locus.  The former can
happen only if $p' + q' < 0$ and $p > 0$ and the latter only if $p' +
q' > 0$ and $p < 0$.  

\narrower\sl
If $p' > 0, q' < 0$ and $p'+q < 0,$ while $p' + q', p$ and $q$
are each non-zero then $p' + q'$
and $p$ have opposite signs.\autonumm%(5.4)
\endnarrower

Now suppose that $p + q > 0$ with $p'$ still positive and $q'$
negative.  Assume that neither $p$ nor $q$ is zero so that one is
positive and the other negative.  If $p > 0$, then positivity of
$\Delta$ requires that $p' + q' > 0$.  On the other hand, if $p < 0$,
then the positivity of $\Delta$ requires that $p' + q' < 0$.  To
summarize:

\narrower\sl
If $p' > 0, q' < 0$ and $p'+q > 0,$ while $p' + q', p$ and $q$
are each non-zero then $p' + q'$
and $p$ have the same signs.\autonumm%(5.5)
\endnarrower

\noindent
{\bf Step 6}\qua Steps 6 and 7 argue that there are no cases of (5.4)
with $p« + q« < 0$.  To see why, remark that in this case, both $f$
and $h$ are negative on the concave side end of $C_0$.  With this
understood, let $C_+ \subset C_0$ denote the connected component of
the $h/f > 0$ locus which contains this concave side end.  Moreover,
since $h/f$ has neither local maxima nor minima, $C_+$ must also
contain the $(q, q')$ end of $C_0$.  On the other hand, as $h^{-1}(0)$
must lie where $f > 0$, $C_+$ does not intersect the sufficiently
large $|s|$ portion of the $(p, p')$ side end of $C_0$.  (The boundary
of the closure of $C_+$ is the $f = 0$ locus.)  They are also both
negative on the convex side end corresponding to $(q, q')$ and both
are positive on the $(p, p')$ convex side end.  However, as
$h^{-1}(0)$ occurs where $f > 0$, it follows that $C_+$ contains both
the concave side $(q, q')$ end and the concave side end of $C_0$.

Given the preceding, remark that $h/f$ tends to $(k'/k) \sin^2
\theta_{0K}$ as $s \rightarrow\infty$ on the concave side end of
$C_0$, while it approaches $(q'/q) \sin^2 \theta_{0Q}$ on the $(q,
q')$ side end.  Here, $\theta_{0K}$ is the value of $\theta_0$ on the
closed Reeb orbit that is determined by the concave side end and
$\theta_{0Q}$ is the value of $\theta_0$ that is determined by the
corresponding closed Reeb orbit for the $(q, q')$ end.  With this last
point understood, the next claim is that
\begin{equation}
(k'/k) \sin^2 \theta_{0K} < (q'/q) \sin^2 \theta_{0Q} .
\label{(5.6)}
\end{equation}
By way of justification, note first that $k'/k < q'/q$ because $\Delta
= k q' - k' q$ and $\Delta > 0$.  Meanwhile, the assignment of
$\theta_0\rightarrow\sqrt{6} \cos(\theta_0)/ (1 - 3 \cos^2 \theta_0)$
defines a smooth function on the (connected) subset of $\theta_0\in
(0, \pi)$ where $\cos \theta_0 < -1/\sqrt{3}$.  As the derivative of
this function is negative on this interval, the fact that $k'/k <
q'/q$ implies that $\theta_{0K} > \theta_{0Q}$.  Therefore, $\sin^2
\theta_{0K} < \sin^2 \theta_{0Q}$ and (5.6) follows.
\medskip

{\bf Step 7}\qua As remarked in the previous section, the restriction
to $C_+$ of $h/f$ has neither local maxima nor minima.  It thus
follows from (5.6) that the infimum of $h/f$ is the limiting value,
$(k'/k) \sin^2 \theta_{0K}$, on the concave side end of $C_0$.  To see
that such an event is absurd, introduce coordinates $(\tau, u)$ on
this end as in (2.19) where $\tau$ is periodic and $u$ is identified
with the pullback of $f$ and so ranges through $(-\delta, 0)$ for some
$\delta> 0$.  Then, parameterize this end of $C_0$ as in (2.19) in
terms of functions $(x, y)$ of the variables $(\tau, u)$.  In
particular, the function $y$ must be non-positive if $h/f$ has infimum
$(k'/k) \sin^2 \theta_{0K}$.  The latter constraint is inconsistent
with (2.20) for the following reasons: First, it follows from
Proposition 2.3 that $|x|$, $|u|^{-1}|y|$ and $|y_u|$ all tend to zero
as $|u|$ tends to zero.  This last fact implies that $x$ and $w \equiv
\sin^{-2} \theta y$ obey an equation which has the schematic form
\begin{equation}
x_\tau = w_u + (\sigma_0 + r) u^{-1} w\;,
\label{(5.7)}
\end{equation}
where $\sigma_0 = (2/3)^{1/2} (k'/k) \cos \theta_0 (1 + (k'/k)^2
\sin^2 \theta_0)^{-1}$ and $r = r(u, w)$ is a smooth function with
$\lim_{u\rightarrow 0} r = 0$.  (In fact, $|r| = \cO(|u|^{-1} |w|)$.)
With regard to $\sigma_0$, note that $\sigma_0 < 0$ since its sign is
that of $k$.  However, it is crucial to note that $\sigma_0 > - 1$.
  
With (5.7) in hand, suppose now that $w \le 0$ for all sufficiently
small values of $|u|$.  It then follows that from (5.7) that there
exists a positive constant $\sigma_1 < 1$ such that for all
sufficiently small $|u|$,
\begin{equation}
x_\tau \ge w_u - \sigma_1 u^{-1} w \;.
\label{(5.8)}
\end{equation}
And, with this last point understood, let $\overline{\omega}(u)$
denote the average of $w$ over the constant $u$ circles.  The latter
function is negative on circles where $w$ is not identically zero.
Moreover, by virtue of (5.8), this function obeys $0 \ge
\overline{\omega}_u - \sigma_1 u^{-1}\overline{\omega}$ from which it
follows that $\overline{\omega}\le - c |u|^{\sigma_1}$ where $c > 0$.
(Since $C_0$ is not a cylinder, $w$ cannot vanish identically on any
open set; thus under the assumption that $w \le 0$, the function
$\overline{\omega}$ cannot vanish identically.)  Thus, as $\sigma_1 <
1$, so $|\overline{\omega}| |u|^{-1}\ge c |u|^{\sigma_1-1}$; and so
$|\overline{\omega}| |u|^{-1}$ diverges as $|u| \rightarrow 0$.  This
last conclusion is absurd because the divergence of this ratio is
precluded by Proposition 2.3.
\medskip

{\bf Step 8}\qua This step eliminates the case of (5.4) where $p' + q'
> 0$.  To start the argument, note that in this case, $f$ is negative
on the concave side end of $C_0$ and also negative on the end that
corresponds to $(p, p')$.  On the other hand, $f$ is positive on the
end that corresponds to $(q, q')$.  This is the end where $h$ is
negative, but $h$ is positive on the concave side end and the end that
corresponds to $(p, p')$.  As $h^{-1}(0)$ lies where $f > 0$, there is
a component, $C_-\subset C_0$, of the locus where $h/f < 0$ which
contains both the concave side end and the $(p, p')$ end of $C_0$.
Furthermore, the closure of $C_-$ has the $f = 0$ locus as its
boundary, and $h/f \rightarrow -\infty$ as this locus is approached
from $C_-$.  Meanwhile, $h/f$ converges to $p'/p \sin^2 \theta_{0P}$
as $|s| \rightarrow\infty$ on the $(p, p')$ end of $C_0$, and it
converges to $k'/k \sin^2 \theta_{0K}$ as $s \rightarrow\infty$ on the
concave side end of $C_0$.
  
With the preceding understood, then the argument just completed in
Step 6 adapts to this case with essentially no modifications given
that
\begin{equation}
k'/k \sin^2 \theta_{0K} > p'/p \sin^2 \theta_{0P}\; .
\label{(5.9)}
\end{equation}
To justify this last claim, note first that $k'/k$ is less negative
than $p'/p$ since $\Delta > 0$.  Thus, (5.9) follows directly if
$\sin^2 \theta_{0K} < \sin^2 \theta_{0P}$.  To see the latter
inequality, note first that both $\theta_{0K}$ and $\theta_{0P}$ lie
in the subinterval of $(0, \pi)$ where $1 > \cos \theta > 1/\sqrt{3}$.
In this interval, the assignment to $\theta_0$ of $\sqrt{6}
\cos(\theta_0)/(1 - 3 \cos^2 \theta_0)$ defines a monotonically
decreasing function.  As the assignment to such $\theta_0$ of $\sin^2
\theta_0$ defines an increasing function of $\theta_0$ on this same
interval, the desired conclusion follows.
\medskip

{\bf Step 9}\qua This step rules out any examples of (5.5).  The first
case to consider here has $p' + q' < 0$.  In this case, $p < 0$ but
both $p + q$ and $q$ are positive.  Thus, $f > 0$ on both the concave
side end and the $(q, q')$ side end of $C_0$, but $f < 0$ on the $(p,
p')$ end.  Meanwhile, $h > 0$ on the $(p, p')$ end of $C_0$ and $h <
0$ on the other two ends.  In this case, there is a connected
component, $C_- \subset C_0$ of the $h/f < 0$ locus which contains
both the concave side end and the $(q, q')$ end of $C_0$.  The
boundary of the closure of $C_-$ is the $h = 0$ locus again, and $h/f
\rightarrow 0$ as this locus is approached from the $C_-$ side.  Then,
with little change, the previous arguments apply to rule this case out
given that (5.6) holds.

To see (5.6) in this case, note first that $k'/k < q'/q$ since $\Delta
> 0$.  In addition, note that both $\theta_{0K}$ and $\theta_{0Q}$ lie
in the subinterval of $(0, \pi)$ where $-1/\sqrt{3} < \cos \theta <
0$.  On this interval, the expression $\sqrt{6} \cos(\theta_0)/(1 - 3
\cos^2 \theta_0)$ defines a decreasing function of $\theta_0$ and so
$\theta_{0K} > \theta_{0Q}$.  Now, $k'/k$ and $q'/q$ are both
negative, so the inequality $\sin^2 \theta_{0K} < \sin^2 \theta_{0Q}$
does not imply (5.5).  However, as $(k'/k) \sin^2 \theta_{0K}$ and
$(q'/q) \sin^2 \theta_{0Q}$ are the values of
\begin{equation} 
\sqrt{6} \cos \theta \sin^2 \theta (1 - 3 \cos^2 \theta)^{-1}
\label{(5.10)}
\end{equation}
at $\theta = \theta_{0K}$ and $\theta_{0Q}$, the inequality in (5.6)
does follow from the fact that (5.10) is a decreasing function of
$\theta$ on the interval in question.  Indeed, the $\theta$--derivative
of (5.10) is
\begin{equation}
-\sqrt{6} \sin \theta (1 + 3 \cos^4 \theta) (1 - 3 \cos^2 \theta)^{-2
 }\;.
\label{(5.11)}
\end{equation}

Finally, consider the possibility that (5.5) holds with $p' + q' > 0$.
Now $p > 0$ so if $q < 0$, then $f > 0$ on the concave side end of
$C_0$ and also on the $(p, p')$ side end.  However, as $q < 0$, so $f
< 0$ on the $(q, q')$ side end.  Thus, there is a component, $C_+
\subset C_0$ of the locus where $h/f > 0$ which contains both the
concave side end and the $(p, p')$ end of $C_0$.  The closure of $C_+$
has the $h = 0$ locus for its boundary and $h/f$ tends to zero as this
boundary is approached in $C_+$.  With this point understood, then the
previously used argument applies given that (5.9) holds.
  
To see (5.9) in this case, first note that $k'/k > p'/p$ since $\Delta
> 0$.  Now, as both $\theta_{0K}$ and $\theta_{0P}$ lie where $0 <
\cos \theta < 1/\sqrt{3}$, and as the function $\theta_0 \rightarrow
\sqrt{6} \cos(\theta_0)/(1 - 3 \cos^2 \theta_0)$ is decreasing on this
subinterval, it follows that $\theta_{0K} < \theta_{0P}$.  With this
understood, (5.9) follows from the fact that (5.10) is also
decreasing, by virtue of (5.11), on this same interval.
 
\sh{(b)\qua Moduli space components}

Fix an ordered set $I = \{(p, p'), (q, q')\}$ of integers subject to
the constraints listed in the statement of Proposition 5.1.  As noted
at the outset of the preceding subsection, this set labels those
components of the moduli space of pseudoholomorphic, $\aleph = 2$,
thrice-punctured spheres in $\mr \times (S^1 \times S^2)$ with ends
that are characterized by the set of three integer pairs $\{(p, p')$,
$(q, q')$, $(k = p + q, k' = p' + q')\}$.  Use $\cH_I$ to denote this
subspace of $\fM$.  The question arises as to the number of components
$\cH_I$.  Here is the answer:
\bigskip

\noindent
{\bf Proposition 5.2}\qua {\sl Let $I = \{(p, p'), (q, q')\}$ denote
a set of pairs of integers that obeys the constraints listed in the
statement of Proposition 5.1.  Then, the space $\fM_I$ of
pseudoholomorphic, $\aleph = 2$, thrice-punctured spheres from
Proposition 4.2 with ends characterized by the set $\{(p, p')$, $(q,
q')$, \\ $(p + q, p' + q')\}$ has at most one connected component.
Moreover, if non-empty, the latter is a smooth manifold that is
$\mr\times T$ equivariantly diffeomorphic to $\mr\times T$.}
\bigskip

The remainder of this subsection is occupied with the following proof.
\bigskip

\noindent
{\bf Proof of Proposition 5.2}\qua The subsequent discussion for the
proof of Proposition 5.2 treats the case where the integer $k = p + q$
is non-zero and positive.  The proof when $k < 0$ is identical to that
for $k > 0$ except for some judicious sign changes.  Meanwhile, if $k
= 0$, then $k' \not= 0$ and the discussion below applies after the
roles of the pair $(t, f)$ are interchanged with those of $(\varphi,
h)$.  Thus, assume throughout that $k > 0$.  Also, assume until
further notice that neither $p$ nor $q$ is zero.

Each component of $\fM_I$ is a smooth manifold by virtue of
Propositions 3.2 and 4.8.  Moreover, as $\dim(\fM_1) = 3$, the
subgroup $\mr\times T$ of $\text{Isom}(\mr \times (S^1 \times S^2)$
acts transitively on each component of $\fM_1$, and so each is
$\mr\times T$ equivariantly diffeomorphic to $\mr\times T$.  In this
regard, remember that the $T$ action on $S^1 \times S^2$ is generated
by the vector fields $\partial_t$ and $\partial_\varphi$, while the
action of $\mr$ on $\mr \times (S^1 \times S^2)$ is generated by
$\partial_s$.  By the way, the $T$ action on $\fM_I$ must be a free
action since the Riemann sphere has no complex automorphisms that fix
three given points.

In any event, if $C \in \fM_I$ and a component $\cH'\subset \fM_I$
have been specified, there exists $C' \in \cH'$ with two special
properties:
\begin{itemize}\sl
\item {The $(t, f)$ coordinates of the critical point of $f$'s
restriction to $C$ are the same as those of the analogous critical
point on $C'$.}
\item {The constant term on the right-hand side of the concave
side end version of (5.1) for $C$ is the same as that for $C'$.}
\end{itemize}
\bigskip

With regard to the first point here, remember that the restriction of
$f$ to any thrice-punctured sphere from Proposition 4.2 with each of
$p$, $q$ and $k$ non-zero has precisely one critical point so
precisely one critical value.  Moreover, this critical value is
non-zero with sign that of $k$.  Indeed, the latter conclusions follow
from the maximum principle since $f$'s pullback to $C_0$ obeys (4.16)
and (4.16) implies a second-order equation for $f$ with the schematic
form $d^*df + v\cdot df = 0$.  Finally, given that the critical values
of $f$ on $C$ and $C'$ have the same sign, then a suitable translation
of $C'$ along the $\mr$ factor of $\mr \times (S^1 \times S^2)$ makes
them equal.  Meanwhile, a suitable rotation of the $S^1$ factor moves
$C'$ so that $f$'s critical point on the resulting subvariety has the
same $t$ coordinate as that of $f$'s critical point on $C$.  Such a
rotation does not change the value of $f$ at its critical point.

With regard to the second point above, by virtue of the fact that $k
\not= 0$, there is an equatorial rotation of the $S^2$ factor moves
any $C'$ so that the resulting subvariety obeys the desired condition.
Note that such a rotation will not change the $(t, f)$ coordinates of
$f$'s critical point.

Now let $C_0$ denote the model thrice-punctured sphere.  As noted in
the preceding sections, $C_0$ comes with a pseudoholomorphic
immersion, $\phi$, into $\mr \times (S^1 \times S^2)$ whose image is
$C$.  There is a similar immersion with image $C'$.  The latter is
denoted by $\phi'_0$ since a subsequent modification, $\phi'$, is
needed for later arguments.  The following lemma describes the salient
features of $\phi'$:
\bigskip

\noindent
{\bf Lemma 5.3}\qua {\sl There exists an immersion $\phi'$: $C_0
\rightarrow \mr \times (S^1\times S^2)$ with image $C'$ such that
$\phi^{'*}(t, f) = \phi^*(t, f)$.}
\bigskip

What follows is a digression for the proof of this lemma.
\bigskip

\noindent
{\bf Proof of Lemma 5.3}\qua Let $f_0$ denote the critical value of
$f$'s restrictions to $C$ and to $C'$.  Now, let $C_f \subset C'_0$
denote the portion where $\phi^*f \not= f_0$.  Likewise, define $C'_f
\subset C'_0$ as the portion where $\phi_0^{'*}f \not= f_0$. There are
three components of $C_f$, each is a cylinder and each corresponds to
an end of $C_0$ and hence a pair from $\{(p, p'), (q, q'), (p + q, p'
+ q')\}$.  As the same assertions hold for $C'_f$, there is a
canonical 1--1 correspondence between the components of $C_f$ and those
of $C'_f$: Components correspond if they correspond to the same pair
from $\{(p, p'), (q, q'), (p + q, p' + q')\}$.
  
Meanwhile, the assignment of $\phi^*(t, f)$ to each component of $C_f$
defines a proper covering map to some subcylinder of the $(t, f)$
coordinate cylinder.  The same is true for the components of $C'_f$,
and corresponding components have the same image.  With a component of
$C_*\subset C_f$ fixed, the covering map to the appropriate $(t, f)$
cylinder is a cyclic covering which is determined by the pair from
$\{(p, p'), (q, q'), (p + q, p' + q')\}$.  Thus, the analogous
covering map from its partner component, $C'_*\subset C'_f$, is
isomorphic.  As a consequence, there exists a diffeomorphism
$\psi'_*$: $C_* \subset C'_*$ which intertwines the projection maps to
the relevant $(t, f)$ sub-cylinder.  Note that this diffeomorphism is
determined up to composition with the group of deck transformations.
  
This freedom with the group of deck transformations can be used to
insure that the three versions of $\psi_*$ from the components of
$C_f$ patch together along $(\phi^*f)^{-1}(f_0)$ to define a
diffeormorphism, $\psi$, from $C_0$ to $C_0$.  To explain, let $t_0$
now denote the value of the $t$--coordinate of $f$'s critical points on
$C$ and $C'$.  Consider some $(t_1, f_0)$ in the $(t, f)$ cylinder
that with $t_1 \not= t_0$.  Take a small disk about this point that is
disjoint from the images of the critical point of $\phi^*f$ and of
$\phi_0^{'*}f$.  This done, then each component of $C_*$ can be
extended by adding the $\phi$--inverse image disks.  Meanwhile, each
component of $C'_*$ can similarly be extended with the addition of the
$\phi'_0$--inverse image disks.  As these extensions remain proper
covering maps over their images in the $(t, f)$ cylinder, so the
corresponding maps $\psi$ can be extended as well.  This understood,
consider a point $z \in C_0$ that lies in a component, say $C_1$, of
$C_*$ and also in the extension of another component, $C_2$.  Then
$\psi_1^*$ is defined near $z$ and so is the extended $\psi_2$.  As
both compose with $\phi_0$ to give the same $(t, f)$ values, so
$\psi_2$ differs from $\psi_1$ on a neighborhood of $z$ by at most a
deck transformation of $C_1$.
  
Now, local agreement of one $\psi_*$ with the extended version of
another implies global agreement for the three.  Indeed, this all
follows from the geometry of the critical locus.  In particular,
because the critical point of $\phi^*f$ is non-degenerate, the
closures of two components of $C_*$ have piece-wise smooth circle
boundary, and that of one component has a figure eight boundary.
Here, the bad point in the figure eight is the critical point of
$\phi^*f$.  Moreover, the critical point is the only point where the
two circle boundary components intersect.  These last points
understood, the deck transformations of the two components of $C_f$
with circle boundary can be used independently to create a smooth map
$\psi$ that is defined on the whole complement in $C_0$ of the
$\psi^*f$ critical point.

Now, some further checking should be done to insure that all is well
with this map $\psi$ near the critical point of $f$.  The latter task
is left to the reader save for the following remark: Let $g_0$ denote
the value of the function $g$ at the critical point of $f$.  Then, the
behavior of the complex function $\lambda\equiv\phi^*(g_0^{-1} f - i
t)$ near the critical point of $f$ can be analyzed using (4.16).  In
particular, $\lambda = \lambda_0 + a z^2 + \cO(|z|^3)$ with respect to
a holomorphic coordinate $z$ centered at the critical point of $f$.
Here, $\lambda_0$ is a constant and $a$ is a non-zero constant.  This
local form for $\lambda$ follows from (4.16).  See, Appendix A in \cite{T4}
where a completely analogous assertion is proved.  Also, note that the
first-order vanishing of $d\lambda$ follows from the observations in
Section 4f that the standard algebraic count of $dt$'s zeros is $-1$,
and that all zero's of $dt$ count with negative weight.

In any event, with $\psi$ in hand, the lemma's map $\psi'$ is the
composition of $\phi_0$ with $\psi$.

With Lemma 5.3 proved, the digression is over.  To continue the proof
of Proposition 5.2, introduce the 1--form $d\uvar\equiv\phi
^*d\varphi-\phi^{'*}d\varphi$.  This is a smooth, closed 1--form on
$C_0$.  It is also exact since $C'$ and $C$ determine the same set
$I$.  Thus, the difference $\uvar=\phi^*\varphi-\phi^{'*}\varphi$ can
be viewed as a bonafide function on $C_0$.  Note that in principle,
there is a choice involved in so viewing $\uvar$, but any two choices
differ by an integer multiple of $2\pi$.  In any event, by virtue of
the fact that $C$ and $C'$ have the same constant term on the concave
end version of (5.1), there is a unique choice for, $\uvar$ that
limits to zero on the concave end of $C_0$.  This said, then $\uvar$
has finite limits on the two convex ends of $C_0$ by virtue of (5.1).
Thus, $\uvar$ is a bounded function on $C_0$.

Now let $\uh\equiv\phi^*h -\phi^{'*}h$, which is automatically a
smooth function on $C_0$.  It then follows from (4.16) and (4.17) that
the pair $(\uvar, \uh)$ obeys an elliptic, first-order differential
equation which has the schematic form:

\begin{itemize}
\item $g \sin^2 \theta\underline{\varphi}_1 = \uh_2 + \sigma \uh$,
\item $g \sin^2 \theta\underline{\varphi}_2 = -\uh_1 + \sigma' \uh$.
\autonum%(5.12)
\end{itemize}

Here, $g$ and $\theta$ are identified with their $\phi$--pullbacks,
while $\sigma$ and $\sigma'$ are smooth functions.
 
Equation (5.12) is employed to justify certain 
remarks that follow about the locus, 
$\underline{G}\subset C_0$, where $\uh = 0$.  
In particular, either $\uh$ is 
identically zero, in which case so is $\varphi$
and 
$C' = C$, or else $\underline{G}$ has the structure of an 'embedded graph' 
as defined in Step 7 of the 
proof in Section 2 of Proposition 2.2.  In this regard, 
the vertices of $\underline{G}$ are the $\uh = 0$ 
critical points of $\uh$.  As with the directed graph 
which appears in Proposition 2.2's proof, 
this graph is naturally oriented.  Its orientation is 
defined by the pullback of $d\underline{\varphi}$ to each 
edge; the latter is non-zero by virtue of (5.12).  
It is also a fact that $\GF$,  as with 
its Section 2 counterpart, has a non-zero, even number 
of incident edges impinging on each 
vertex; and of these, half are point towards the 
vertex and half point away.  (All of these last 
remarks are proved by copying the arguments 
in Steps 3--5 of Part b of the Appendix to 
\cite{T4}.)

Now, an argument along the lines of that used in 
Step 8 of Proposition 2.2's proof 
establishes the following:  If $C' \not= C$ and $\GF\not=\emptyset$, and 
if the extreme values of $\uvar$'s restriction 
to $\GF$ are not  its limiting values on 
the ends of $\GF$ 
then the aforementioned properties 
of $\GF$ cannot hold.  Thus, 
Proposition 5.2 follows via \emph{reductio ad 
absurdem} with an argument that proves
when $C' \not=C$, then $\GF$ is non-empty and neither of 
$\uvar$'s extreme values on $\GF$ are its limiting values 
on the ends of $\GF$.   

For the purposes of establishing that $\GF \not=\emptyset$, consider 
(5.12) where $s \gg 1$, thus, far down the concave 
side end of $C_0$.  Now parameterize this 
portion of $C_0$ by coordinates $(\tau, u)$ where $\tau\in \mr/(2\pi
|k| \mz)$ 
and $u \in (0,\delta)$ or $u \in(-\delta, 0)$ 
depending on whether $k > 0$ or $k < 0$.  Here, $\delta > 0$ is 
very small.  The argument that follows 
considers first the case where $k > 0$. Thus, with $u > 0$ 
understood, the relevant portion of $C$ 
is parameterized as in (2.19) in terms of functions $x$ and $y$.  
In particular, the pair $(x, y)$
obey  (2.20) and of particular interest here is 
the second equation in (2.20).  In this regard,
view $h$  as a function of $u$ and $y$ and thus view $\theta$
as function of $u$ and $y$.  Now, introduce the 
function $w \equiv y \sin^{-2}\theta$ and then the second equation 
in (2.20) implies an equation for $x$ and 
$w$ that has the schematic form as in (5.7).   
Note that in this version of (5.7), the constant $\sigma_0 
> 0$ since $k > 0$ and $\sigma_0$ has the same sign as $k$.
  
Meanwhile, the analogous part of $C'$ also has a 
parameterization as in (2.19) in terms 
of functions $(x', y')$.  Then, $x'$ and the primed analog, $w'$, 
of $w$ obey the analog of (5.7).  
With this understood, subtract the primed version of (5.7) 
from the original to obtain the 
following equation for $\uvar= x - x'$ and $\underline{w} = w - w'$:
\begin{equation}
\uvar_\tau = \underline{w}_u + (\sigma_0 + \underline{r}) u^{-1}\underline{w}\; ,
\label{(5.13)}
\end{equation}
where $r$ is a smooth function with $\lim_{ u\rightarrow 0} 
|\underline{r}| = 0$

An equation such as (5.13) for $\uw$ is useful for two reasons.  
First, when $u > 0$ but 
very small, then, as is demonstrated below, $\uw = 0$ if and only if 
$\uh = 0$.  Given that such is the 
case, it is sufficient for the purposes of proving Proposition 5.1 
to establish the existence of 
a zero of $\uw$ along each sufficiently constant but small 
$u$ circle.  In this regard, remember that 
$\lim_{u\rightarrow 0}\uvar = 0$ on the concave side end of $C_0$.  
Second, (5.13) does indeed imply
the existence  of zeros of $\uw$.  To see (5.13) lead to this last 
conclusion, suppose to the
contrary that $\uw > 0$  for all sufficiently small $u$.  
It then follows from (5.13) that the
average, $\oww\equiv\oww(u)$ of $\uw$  around all constant but 
small and positive $u$ circles obeys the
differential inequality
\begin{equation}
\oww_u + 2^{-1}\sigma_0 u^{-1}\oww < 0\; .
\label{(5.14)}
\end{equation}
The latter implies that $\oww > c u^{-\sigma_0/2}$ as 
$u \rightarrow 0$ with $c$ a positive constant.  This
conclusion is  ludicrous as both $y$ and $y'$ tend to zero as $u$ 
tends to zero.  Likewise, if $\uw < 0$
for all  sufficiently small $u$, then (5.14) implies that $\oww < -c u^{-\sigma_0/2}$ 
as $u \rightarrow 0$ with $c$
a positive constant,  which is an equally ludicrous 
conclusion.  Thus, (5.14) is consistent only
with the  conclusion that $\uw = 0$ at some point 
on each constant, but small $u$ circle.

With the preceding understood, return to the claim 
that $\uw = 0$ if and only if $\uh = 0$ 
when $u$ is very small.  In this regard, note first that $\uh = 0$ 
requires $y = y'$ and $\theta = \theta'$, and
so  $\uw = 0$.  To prove the converse, note that the dependence 
on the coordinates $f$ and $h$ of $\theta$ is 
such that \\
$\theta(f, h = k'/k\sin^2 \theta_0 f + \eta) = 
\theta_0 + \gamma(\eta/f)$ when $|\eta|/f$ is small.  Moreover,
$\gamma$ is $\cO(\eta/f)$.   Thus, as both $y/u$ and $y'/u$ 
tend to zero as $u \rightarrow 0$ (as attested by
Proposition 2.3), so both $w = y \sin^{-2}\theta_0 (1 + \cO(y/u))$ 
and $w' = y' \sin^{-2} \theta_0 (1 +
\cO(y'/u))$ when $u$ is small.   Therefore, when $u$ is small, 
$w = w'$ forces $y = y'$ and so $h = h'$.

With the $k > 0$ argument understood, it can be said 
that the argument for the $k < 0$ 
case is similar although not identical.  In particular, 
the only substantive difference arises in 
the argument for the vanishing of $\oww$ because the 
constant $\sigma_0$ which appears in (5.7) is 
negative when $k < 0$.  To argue that $\oww= 0$ when $k < 0$, 
remark first that though negative,
$\sigma_0>  -1$.  Thus, if $\uw > 0$ where $|u| > 0$ 
is sufficiently small, then $\oww$ obeys the
following analog of  (5.14):
\begin{equation}
\oww_v - \sigma_1 v^{-1}\oww  < 0 \;,
\label{(5.15)}
\end{equation}
where $v \equiv -u > 0$ is small.  Here, $\sigma_1$ is positive, 
but $\sigma_1 < 1$.  This last equation
implies that $\oww\ge c |u|^{\sigma_1}$ where $c$ is a positive 
constant.  In particular, $|u|^{-1}\oww$ is
unbounded as $|u|\rightarrow 0$ which  is impossible since, 
as previouly noted, both $y/|u|$ and
$y'/|u|$ tend to zero as $|u|$ tends to zero.   With some 
judicial sign changes, the preceding
argument also rules out the possibility that  $\uw$ is strictly 
negative where $|u|$ is sufficiently
small.
 
With it now established that $\GF \not=\emptyset$, turn now to the 
question of whether $\uvar$'s extreme 
values on $\GF$ are its limiting values on $\GF$'s ends.  For this 
purpose, suppose that $\GF$ has non-compact
intersection with the closure, $C_1 \subset C_0$, 
of one of the three components of the 
complement of the $f = f_0$ locus.  Now, fix some 
$s_0 \gg 1$, a regular value of $|s|$ on $C_1$ and such 
that the $|s| = s_1$ locus in $C_1$ is a circle having 
transversal intersection with $\GF$.  As $\GF$
divides  the $C_1$ into the portion where $\uh > 0$ and where 
$\uh < 0$, so $\GF$ must have an even number of 
intersections with the $|s| = s_0$ circle in $C_1$ 
and these points alternate upon a
circumnavigation  of this circle between points 
where $\varphi$ is increasing and decreasing in the
direction of  increasing $|s|$. 
 
Remark now that a compact, oriented path in the $|s| \ge s_0$ part of $\GF$ with both 
boundary points on the $|s| = s_0$ circle pairs up two 
$|s| = s_0$ boundary points.  This 
understood, remove from $\GF$ a maximal set of such paths, 
no two sharing edges, to obtain a 
new embedded and directed graph, $\GF_1$, in the $|s| \ge s_0$ 
portion of $C_1$.  Note that $\GF_1$ cannot
be  empty (by assumption).  In addition, each $|s| > s_0$ 
vertex has an even number of impinging 
edges, half oriented by $d\uvar$ to point outward and 
half to point inward.  Moreover, $\GF_1$ must 
intersect the $|s| = s_0$ locus.  Indeed, otherwise there 
would be a properly embedded path in 
$\GF_1$ with $|s|$ unbounded at both ends and this is ruled out 
by the fact that $\uvar$ has a unique, $|s| 
\rightarrow\infty$ limit on $C_1$.  Now, given that 
$\GF_1$ interstects the $s = s_0$ locus, it does so in
an even  number of points, and again, these alternate 
between those where $\uvar$ is increasing with $|s|$
and  those where $\uvar$ is decreasing.  Now, by assumption, 
no pair of these points comprise the 
boundary of a compact, oriented path in $\GF_1$, and so each is the 
sole boundary point of a path 
in $\GF_1$ on which $|s|$ is unbounded.  This the case, 
then $\uvar$ is increasing on half of these paths 
in the unbounded direction and decreasing on the other 
half.  In particular, $\uvar$'s limit as $|s| 
\rightarrow\infty$ on $\GF \cap C_1$ is not an extreme value of $\uvar$ on $\GF$.
  
Now consider the case where one of $p$, $q$ vanishes. In this regard, the argument 
below considers the case where $q = 0$, $p > 0$.  
The argument for the $p = 0$, $q < 0$ is virtually 
identical and is left to the reader. 
The following lemma is needed:
\bigskip

{\bf Lemma 5.4}\qua  {\sl Let $r$ be a smooth function of the 
angle $\theta$, defined on some open interval 
$\uZ \subset [0, \pi]$ where it is everywhere distinct from 
$-f/h$. Then there exists a function $\zeta$ on
$\uZ$ and  a function $u$ of the variables $(f, h)$, 
defined where $\theta\in\uZ$, with $du = e^\zeta(df + r
dh)$.}
\bigskip

{\bf Proof of Lemma 5.4}\qua  First, the $f$ and $h$ derivatives of $
\theta$ are related via $\theta_f = -
h/f \theta_h$.  With  this understood, then $\zeta$ is 
determined up to an additive constant by the
requirement that
\begin{equation}
\zeta_\theta(1 + r h/f) + r_\theta h/f = 0 \;.
\label{(5.16)}
\end{equation}
As $h/f = \sqrt{6} \cos \theta \sin^2 \theta (1 - 3 \cos^2 
\theta)^{-1}$, this last equation reads
\begin{equation}
\zeta_\theta = - r_\theta \sqrt{6} 
\cos \theta \sin^2 \theta [(1 - 3 \cos^2 \theta)+ r 
\sqrt{6} \cos \theta \sin^2 \theta]^{-1}\;.
\label{(5.17)}
\end{equation}
As $r \not= -f/h$ on its domain of definition, the 
latter equation can be integrated to obtain $\zeta$.

To employ this lemma, remark that if $C \in \fM_I$, then 
$q' > 0$ as both $p$ and $\Delta$ are positive 
and $q = 0$.  Thus, one of the convex side ends of $C$ corresponds 
to a closed Reeb 
orbit with $\cos \theta_0 = 1/\sqrt{3}$.  Meanwhile, $f > 0$ 
on the other two ends of $C$.  As argued
in the  previous subsection, this implies 
that $f > 0$ on the whole of $C$.
 
With $f$ everywhere positive on $C$, it follows 
that the range in $\mr$ of the restriction to $C$ 
of the function $-\sin^2 \theta f/h$ is disjoint from an 
interval of the form $[0, \delta)$.  Given such $\delta$
that  is less than the minimum of $p/|p'|$ and 
$p/|p'+ q'|$, let $\nu\in(0, \delta/2)$ be a rational number.  
Now employ Lemma 5.4 using $r = \nu \sin^{-2}\theta$.  
It then follows that any choice for the 
resulting function $u$ is defined on a neighborhood of 
every $C \in \fM_I$.  By the way, as any two 
choices for $u$ differ by an additive constant, so it 
follows from (1.3) that there is a unique 
choice which limits to zero as $s \rightarrow\infty$ along every 
$C \in \fM_I$.  This particular choice for
$u$  should be taken in what follows.

Concerning this function $u$, note first that $|u|$ tends 
uniformly to infinity as $s \rightarrow-\infty$ 
along each of the convex side ends of any $C \subset \fM_I$.  
Moreover, with $C \in \fM_I$ specified, the 
function $u$ pulls back via the defining pseudoholomorphic 
immersion from $C_0$ as a function 
with only one critical point, the latter being 
non-degenerate and hyperbolic.  Indeed, this 
follows from the following two observations:  
First, $J$ maps du to a nowhere zero multiple of 
$dt + \nu d\varphi$.  Meanwhile, as $\nu$ is neither 
0, $p/p'$ nor $p/(p' + q')$, the arguments from the
proof  of Proposition 4.7 apply to prove that 
$dt + \nu d\varphi$ pulls back to $C_0$ with but one zero, which
is  hyperbolic.
 
With the preceding as background, suppose that $C \in \fM_I$ 
has been specified as well 
as a component $\cH'\subset \fM_I$.  Once again, $\cH'$ is a 
smooth manifold and  
equivariantly diffeomorphic to $\mr\times T \subset \text{Isom}
(\mr \times (S^1 \times S^2)$).  In particular, 
this implies that there exists $C'\in \cH'$ with the following 
properties:  First, the pair $(t + \nu
\varphi,  u)$ at the  critical point of $u$'s 
restriction to $C$ is identical to that
at the  critical point of $u$'s restriction to $C'$.      Second, the constant term on the
right-hand side of the concave version of (5.1)  for $C$ 
is the same as that for $C'$.
  
To proceed, let $\phi$: $C_0 \rightarrow \mr \times 
(S^1 \times S^2)$ denote the defining
pseudoholomorphic  immersion with image $C$, and let $\phi'_0$ 
denote the corresponding immersion
with image $C'$.   The arguments given above for the 
proof of Lemma 5.2 can be modified in a minor
way to  find a diffeomorphism $\psi$: $C_0 \rightarrow C_0$ 
for which the pullbacks of $(t + \nu\varphi, u)$
by $\phi'\equiv\phi'_0\circ\psi$ and $\phi$ are identical.  
This understood, introduce $\uvar\equiv\phi^*\varphi-\phi^{'*}\varphi$ 
and $\uh =
\phi^*h-\phi^{'*}h$ as  functions on $C_0$.  In this 
regard, note that $\uvar$ can be viewed as a
bonafide $\mr$ valued function  on $C_0$ that is bounded, 
limits to zero as $s \rightarrow\infty$ 
on $C_0$ and has finite  limits on the other ends 
of $C_0$.  Moreover, as demonstrated momentarily,
the functions $\uvar$  and $\uh$ enjoy the following properties:
\bigskip

\narrower\sl
{If $ C \not= C'$, then the $\uh = 0$
level set is an embedded graph,
$\underline{G}\subset C_0,$
whose vertices are the
$\uh = 0$ critical points of  $\uh$. Moreover,
$d\uvar$ pulls back
without zeros to the
edges of
$\GF$ and,  provided $\nu>0$ is sufficiently small,
$\GF $ has a non-empty set
of edges.}
\autonumm%(5.18)
\endnarrower

As before, these properties are inconsistent 
and so $\uh = 0$, $\uvar = 0$ and $C = C'$.
  
The validity of (5.18) is easiest seen (perhaps) by 
using the $\phi$ pullbacks of the pair 
$(\tau\equiv t+\nu\varphi,u)$ as local coordinates away 
from the critical point of $u$.  In terms of
these  coordinates, the $\phi$ pull back of $(\varphi, h)$ 
obeys, by virtue of $C$ being
pseudoholomorphic, an  equation with the schematic form
\begin{equation}
\varphi_\tau = \sin^{-2} \theta\; e^\zeta h_u\quad
 and\quad  \varphi_u = - \sin^{-2 }\theta \;g^{-2} 
e^{-\zeta} h_\tau\;.
\label{(5.19)}
\end{equation}
Here, $\theta$ and $\zeta$ are the pullbacks via $\phi$ 
of their namesakes on $\mr \times (S^1 \times S^2)$. 
Meanwhile, $\varphi'\equiv\phi^{'*}\varphi$ and $h'\equiv\phi^{'*}h$ 
obey (5.19) but with $\theta$ and $\zeta$
replaced by $\phi^{'*}\theta$ and $\phi^{'*}\zeta$, respectively.   
As $\varphi$ and $\zeta$ are implicit
functions of $f$ and $h$, it follows from (5.19) and its 
primed analog  that the pair $(\uvar, \uh)$ obey an
equation with the schematic form of (5.12).  
And, with this last  point understood, the arguments
given about $\GF$ in the case where none of $p$, $q$, and $p + q$  
vanish can be used with only
minor modifications to prove (5.18).  As these modifications  
are slight, their details are left to
the reader.

\sh{(c)\qua  Thrice-punctured spheres in $\mc^*\times \mc^*$}

The argument for the existence of Proposition 5.1's 
thrice-punctured spheres starts 
here and runs through Subsection 5g.  By way of a 
beginning, this subsection first 
considers holomorphic, triply-punctured spheres in the 
complex manifold $\mc^*\times \mc^*$ where 
$\mc^* = \mc-\{0\}$.  The punctured spheres are then used 
to construct symplectic, thrice-punctured 
spheres in $\mr \times (S^1 \times S^2)$ with prescribed 
asymptotics.  The subsequent subsections
explain  how to deform the latter to obtain 
those predicted by Proposition 5.1.
  
To start the $\mc^* \times \mc^*$ discussion, identify 
the thrice-punctured sphere, $C_0$, as a 
complex manifold with the complement of the points 0, 1 and 
$\infty$ in the Riemann sphere, $\mP^1 
= \mc \cup\infty$.  Alternately, $C_0 = \mc-\{0, 1\}$.  
Note that $C_0$ comes with the order six subgroup
$G$  of the complex automorphism group 
$PSL(2, \mc)$ which permutes the points $\{0, 1,\infty\}$.  This 
group is generated by the automorphisms that send $z$ to $1-z$ and $z$ to $1/z$.

Now, consider a holomorphic map $\phi$: $C_0 \rightarrow 
\mc^* \times \mc^*$ which has the form
\begin{equation}
\phi = (a r^{p+q} z^{-p} (1 - z)^{-q}, a'
r^{p'+q'} z^{-p'} (1 - z)^{-q'})\;,
\label{(5.20)}
\end{equation}
where  $p'$, $q$,  and $q'$ are integers, $a$ and $a'$ 
are unit length, complex numbers, and $r
\gg 1$  is a constant whose lower bound will be specified 
shortly.  Here,  $I\equiv 
\{(p,  p'), (q, q')\}$ is constrained so that $\Delta\equiv 
p q' - q p'$ is non-zero.  In this
regard, note that $\phi$  factors through a map to 
$\mc^*$ if and only if $\Delta = 0$, so this is a reasonable
constraint on $I$.

It is important to point out that with $r$ fixed, each 
$\phi$ in (5.20) has five companions 
with the same image in $\mc^* \times \mc^*$.  These companions are 
obtained by composing $\phi$ with the 
non-trivial elements of $G$.  However, in the subsequent 
constructions, the puncture at $\infty$ 
plays a special role and only the involution $z \rightarrow (1 - z)$ 
in $G$ preserves this puncture.  
As  this involution preserves the set of maps 
having the form in (5.20), it sends the four integers 
$\{(p, p'), (q, q')\}$ to another such set, namely 
$\{(q, q'), (p, p')\}$
and  thus changes the sign of $\Delta$.  In particular, 
no generality is lost in the subsequent
discussions  by restricting to the $\Delta > 0$ case.
	
For future reference, make a note that the map $\phi$ 
is an immersion.  Indeed, a singular 
point of $\phi$ can occur only where both $(p/z - q/(1 - z))$ 
and \\
$(p'/z - q'/(1 - z))$ vanish.  Since
$\Delta\not=0$, this happens nowhere on the thrice-punctured 
disk.  It is also important to note
that $\phi$  has only a finite number of double points 
and thus embeds neighborhoods in $C_0$ of the 
punctures.

\sh{(d)\qua  The map $\Phi$}
  
Fix a set $I = \{(p, p'), (q, q')\}$ of integers which 
obey the constraints in the three 
points of Proposition 5.1.   An embedding in $\mr \times (S^1 \times S^2)$ 
of a neighborhood of the 
image in $\mc^* \times \mc^*$(5.20) that is proper on this 
image defines a properly immersed,
triply-punctured sphere in $\mr \times (S^1 \times S^2)$.   
This subsection describes a map that
provides an  embedding of this sort and that sends $C_0$ 
to a subvariety, $C'$, with certain
favorable  properties.  In particular, $C'$ is 
pseudoholomorphic for a $T$--invariant almost complex 
structure, $J'$, that tames $\omega$ and is asymptotic to 
(1.5)'s fiducial almost complex structure
$J$ as  $|s| \rightarrow\infty$ on $\mr \times (S^1 \times S^2)$.  
Also, $C'$ is constructed with
prescribed asymptotics.  This is to  say that the 
components of the large and constant $|s|$ slices
of $C'$ converge to multiple  covers of closed Reeb orbits 
where the multiplicities and
the closed  Reeb orbits are suitably determined by the set $I$.  
Finally, $C'$ avoids the $\theta\in
\{0, \pi\}$ locus  precisely when the set $I$ obeys 
Proposition 5.1's constraints.  The preceding
properties play  key roles in  subsequent arguments.
  
The desired embedding, $\Phi$, is described in this 
subsection simply as a map from a 
neighborhood of $\phi(C_0)$ into $(S^1 \times \mr) \times (S^1 \times \mr)$.  
Subsequently $\Phi$ is
interpreted as a map  into $\mr \times (S^1 \times S^2)$ 
by the identification the former space with
the complement of the $\theta\in
\{0,  \pi\}$ locus in the latter via the coordinates 
$(t, f, \varphi, h)$.  Of course, such an
interpretation is  possible only when the initial map into 
$(S^1 \times \mr) \times (S^1 \times \mr)$ avoids
$(S^1 \times 0)\times (S^1 \times 0)$.  The  next subsection 
provides the verification that such is
the case.   

To begin the discussion, take $(t, f, \varphi, h)$ as 
coordinates on $(S^1 \times \mr) \times (S^1 \times \mr)$. 
Here,  $t$ and $\varphi$ are $\mr/(2\pi \mz)$ valued while $f$ 
and $h$ are $\mr$ valued.  Meanwhile, introduce
$\mr/(2\pi \mz)$  valued functions $(t, \varphi)$ and $\mr$--valued 
functions $(u, v)$ for $\mc^* \times \mc^*$ by writing
the complex  coordinates $\lambda$ and $\lambda'$ as 
$\lambda = e^{u-it}$ and $\lambda' = e^{v-i\varphi}$.    
The map $\Phi$, 
sends a neighborhood of  $\phi$'s image into 
$(S^1\times \mr) \times (S^1 \times \mr)$, so that both the
coordinates $t$ and $\varphi$ on $(S^1 \times \mr) \times  (S^1 \times \mr)$ 
are respectfully identified
via $\Phi$'s pullback  with their namesake  functions on $\mc^*\times \mc^*$.  
Meanwhile, $\Phi$ respectively
identifies $f$ and $h$ with $u$ and $v$ except  near the 
$\phi$ image of neighborhoods of 0, 1 and $\infty$. 
Near the $\phi$ image of neighborhoods of  the punctures, 
the definition of the map $\phi$ is completed below by writing 
$f$ and $h$ as  suitable functions of $u$
and $v$.
     
To see what is involved here, observe that $\phi$ near $z = 0$ has the form
\begin{equation}
\phi(z) = (a r^{p+q} z^{-p} (1 + \cO(|z|)), 
a' r^{p'+q'} z^{-p'} (1 + \cO(|z|)))\; .
\label{5.21}
\end{equation}
Thus, the $\phi$ image of a neighborhood of $z = 0$ has 
\begin{itemize}
\item $p \varphi - p' t  =  p'\;\text{ arg}(a) - p\; 
\text{arg}(a') +  \cO(|z|)    \mod(2\pi)$,
\item $p v - p' u  = \Delta \ln r +  \cO(|z|)$.\autonum%(5.22)
\end{itemize}
Even with the $\cO(z)$ and $\ln r$ terms absent, the identification 
$(f = u, h = v)$ does not 
make the locus in (5.22) close to a $J$--pseudoholomorphic 
subvariety when both $p$ and $p'$ are 
non-zero.  In any event, the plan is to define $\Phi$ near the 
$\phi$--image of a neighborhood of 0 so 
that the composition $\Phi\circ\phi$ sends a neighborhood of 
$z = 0$ close to the $J$--pseudoholomorphic 
locus that is defined by
\begin{itemize}
\item $p \varphi - p' t  =  p' \arg(a) - p \arg(a')\;    \mod(2\pi)$,
\item $p h - p' \sin^2 \theta_0 f  =  0$.
\autonum%(5.23)
\end{itemize}

\noindent
Here, $\theta_0$ is determined as follows:  First, $\theta_0$ 
obeys $p'(1 - 3 \cos^2 \theta_0) = p \sqrt{6} \cos
\theta_0$.  In this  regard, remember that there is one 
solution when $|p'/p| < (3/2)^{1/2}$ and two
otherwise.  In the  former case, the condition 
$p'(1 - 3 \cos^2 \theta_0) = p \sqrt{6} \cos \theta_0$ completely
determines $\theta_0$.   Meanwhile, when $|p'/p| > (3/2)^{1/2}$, 
the sign of $\cos \theta_0$ differs for
the two solutions; and with  this understood, the angle $\theta_0$ 
is chosen so that the signs of $p'$ and
$\cos\theta_0$ agree.  The third  constraint in Proposition~5.1 
guarantees the existence of such a $\theta_0$.
 
The subsequent definition of $u$ and $v$ in terms of $f$ 
and $h$ depends on whether or not 
one of $p$ and $p'$ is zero, and when both are 
non-zero, it depends on their signs.  The case 
where both are non-zero and positive is presented 
immediately below; the other cases where 
both are non-zero can be obtained from this one 
by suitable notational changes.  The 
case where one of $p$ or $p'$ is zero is presented afterwards.
    
In the case where both $p$ and $p'$ are positive, $u$ 
is chosen to be an increasing function 
of $f$ and $v$ an increasing function of $h$.  
In this regard, keep in mind that with $p$ positive, 
small $|z|$ means large $u$ and thus large $f$; and with 
$p'$ positive, the same conclusion holds for 
the pair $v$ and $h$.  In particular, the small $|z|$ 
region of $D$ is to be mapped by $\Phi\circ\phi$ 
to the convex end of $\mr \times (S^1 \times S^2)$, that where $s \rightarrow-\infty$.
    
To define $\Phi$ near the $\phi$ image of a neighborhood of 
$z = 0$, first introduce
\begin{equation}
\sigma_0 \equiv\sqrt{6} (p^2 + p^{'2} \sin^2 \theta_0)^{1/2} .
%(5.24)
\end{equation}
Next, select a function $\chi$: $\mr\rightarrow[0, 1]$ 
that equals 1 on $(-\infty, 1]$ and 0 on $[2, \infty)$. 
With $\chi$ chosen and for each $R > 1$, introduce $\alpha_R
(\cdot) \equiv\chi((\cdot)/(4R))$ and $\beta_R(\cdot) \equiv 1 -
\chi((\cdot)/2R)$.

Now, take $R_0 \gg 1 + r$, set $R \equiv R_0$ to 
define $\alpha_R$ and $\beta_R$, and declare the $\Phi$
pullback  of the function $f$ to be
\begin{equation}
f(u) = \alpha_R(p^{-1} u) u + \beta_R(p^{-1}u) 
(p' \sin^2 \theta_0)^{-1}
e^{\sigma_0(u-q\ln(r))/p}   .
%(5.25)
\end{equation}
Meanwhile, declare the $\Phi$ pullback of the $h$ to be
\begin{equation}
h(v) =\alpha_R(p^{'-1} v) v + \beta_R(p^{'-1} v) p^{-1}
e^{\sigma_0(v-q'\ln(r))/p'}  .
%(5.26)
\end{equation}

Concerning this definition of these $\Phi$ pullbacks, 
note in particular that both $f_u > 0$ 
and $h_v > 0$ when $R_0$ is sufficiently large as 
$\alpha_R + \beta_R \ge 1$ everywhere, $|d\beta_R$ and $|d\alpha_R|$ are 
$\cO(1/R)$, and $\beta_R = 1$ where $d\alpha_R \not= 0$.
   
Also, as $f(u) = u$ when $u < 2 pR_0$ and $f(v) = v$ 
when $v < 2 p' R_0$, the definition of $\Phi$ 
via (5.25)--(5.26) is consistent (when $R_0$ is 
large) with the previous definition for $\Phi$'s 
restriction to a neighborhood of a particular compact 
subset of $C_0$.  Moreover, by virtue of 
(5.22), the fact that $f(u) = (p' \sin^2 \theta_0)^{-1}
e^{\sigma_0(u-q\ln(r))/p}$
when $u > 8 pR_0$ and the fact that $h(v) = 
p^{-1}
e^{\sigma_0(v-q'\ln(r))/p'}$  when $v > 8 p'R_0$, the composition
$\Phi\circ\phi$  embeds a neighborhood in $C_0$
of  the puncture 0 as a locus in 
$\mr\times(S^1\times S^2)$ which is defined by equations of the following 
sort:
\begin{itemize}
\item $p \varphi - p't = p' \arg(a) - 
p \arg(a') + r_1   \mod(2\pi)$ .
\item $p h - p' \sin^2 \theta_0 f = r_2 $.
\autonum %(5.27)
\end{itemize}

\noindent
Here, the terms $r_1$ and $r_2$ come from the 
$\cO(z)$ terms in (5.22).  In particular, they can be 
viewed as functions which are defined on an open 
set which contains the $\Phi\circ\phi$ image of a 
neighborhood of $z = 0$ in $D$.  In this regard, such an 
open set can be chosen so that each of 
$r_1$, $f^{-1} r_2$ and their derivatives are $\cO(e^{-\sqrt{6}|s|/\sigma_0})$.

By the way, as the image of $\phi$ in $\mc^*\times \mc^*$ 
is holomorphic, it is pseudoholomorphic 
with respect to the integrable almost complex structure 
which sends $\partial_t$ to $\partial_u$ and $\partial_\varphi$
to $\partial_h$.  As  a consequence, the image of $C_0$ in 
$\mr\times(S^1\times S^2)$ by $\varphi$'s composition
with $\Phi$, as defined  above by $t$, $\varphi$ and (5.25--26), 
is pseudoholomorphic for the almost
complex structure $J'$ that  is defined near the $\Phi\circ\phi$ 
image of a neighborhood of $z = 0$ by
\begin{equation}
J'\partial_t = f_u \partial_f\quad and  \quad
J'\partial\varphi = h_v \partial_h .
\end{equation}
%(5.28)
Here, $f_u$ is evaluated at $u(f)$ and likewise, $h_v$ at $v(h)$.  
Note specifically that (5.25),
(5.26) and  (5.27) imply that
\begin{itemize}
\item  $f_u = p^{-1} \sigma_0 f = g (1 + \cO(e^{-\sqrt{6}|s|/\sigma_0} ))$ ,
\item $h_v = p^{'-1} \sigma_0 h = g \sin^2  (1 + \cO(e^{-\sqrt{6}|s|/\sigma_0}))$ 
\autonum%(5.29)
\end{itemize}

\noindent
near the $\Phi\circ\phi$ image of a neighborhood of $z = 0$.  
Here, $g = \sqrt{6} (f^2 + h^2 \sin^{-2} \theta)^{1/2}$. 
Thus, the  almost complex structure $J'$ and the 
fiducial almost complex structure $J$ in (1.5) are 
exponentially close,
\begin{equation}
|J' - J| = \cO(e^{-\sqrt{6}|s|/\sigma_0} ) , 
%(5.30)
\end{equation} 
on an open set which contains the $\Phi\circ\phi$ 
image of a neighborhood of $z = 0$ in $D$.

Now consider $\Phi$ when $p = 0$ and $p' < 0$.  
The $p'> 0$ case is obtained from this one 
by switching various signs and is left to the reader.  
In the $p = 0$ and $p' < 0$ case, $v$ tends to 
$-\infty$ as $|z|$ tends to zero.  With this 
understood, again take $R \equiv R_0$ very large and choose $h
=  h(v)$ to be a favorite function of $v$ with the following properties:
\begin{itemize}
\item $h_v > 0$ .
\item  $h = v$    {	where}  $v \le -R$ .
\item $h = -e^{-\sqrt{6} v}$  { where}  $v < -R^2$.  
\autonum%(5.31)
\end{itemize}

\noindent
With $h$ now chosen, take $f$ to be a favorite function of 
both $u$ and $v$ that satisfies
\begin{itemize}
\item $f_u > 0$ .
\item $f = u$  \quad 		{	where}  $v \le -R$ .
\item $f = - 6^{1/2} h(v) (u - q \ln r)$ \quad 	{where}   $v \le -R^2$ .
\autonum% (5.32)
\end{itemize}

Given the preceding, it follows that a 
neighborhood in $C_0$ of $z = 0$ is mapped by 
$\Phi\circ\phi$ as the large $-s$ portion of a 
locus that is defined by equations of the form
\begin{itemize}
\item $\varphi = -\arg(a') +  r_1 \;   \mod(2\pi)$ ,
\item $f = r_2$ ,
\autonum%(5.33)
\end{itemize}

\noindent
where $r_1$ and $r_2$ are functions on a 
neighborhood of this portion of the image of
$\Phi\circ\phi$.  Here,  as with their namesakes in 
(5.27), each of $r_1$, $|h|^{-1} r_2$, and their
derivatives are $\cO(e^{-|s|/\zeta})$ where  $\zeta\ge 1$ is a constant.

In this case, the locus image of $\Phi\circ\phi$ 
is pseudoholomorphic for the almost complex 
structure $J'$ which is defined so that
\begin{equation}
J'\partial_t = f_u \partial_f\quad \hbox{\rm and}\quad        
J' \partial_\varphi=f_v \partial_f
+h_v\partial_h . 
%(5.34)
\end{equation}
In particular, this last equation implies that $J'$ and 
$J$ again obey (5.12) on an open set which 
contains the $\Phi\circ\phi$  image of a 
neighborhood of $z = 0$ in $C_0$.

As remarked, there is an analogous description of $\Phi$ 
in the $p = 0$, $p' > 0$ case, and 
also in the cases where $p' = 0$.  The details 
here differ from those just given only in the 
notation so they won't be given.
          
Of course, a corresponding description of $\Phi$ 
can be simultaneously made near the $\phi$ 
image of a neighborhood of $z = 1$ as well.  The discussion 
is identical to that just ended save 
for some signs and the interchange of $p$ with $q$ and 
$p'$ with $q'$.  By the way, as the definition 
of $\Phi$ near the $\phi$ image of a neighborhood in 
$C_0$ of 0 required the choice of the constant
$R_0$,  so the definition near the $\phi$ image of a 
neighborhood in $C_0$ of 1 requires the choice of
a  constant, $R_1$.  The freedom to make these choices 
is exploited in a subsequent argument.

The final order of business in this subsection 
completes $\Phi$'s definition by 
describing this map near the $\phi$ image of a neighborhood in 
$C_0$ of $z = \infty$.  Here, $\Phi$ is defined 
so that the composition $\Phi\circ\phi$  properly 
embeds the complement of a large radius disk
in $C$ in  the concave end of $\mr\times(S^1\times S^2)$.  
In particular, with $\Phi$ as specified below,
the constant $s$  slices of the ($\Phi\circ\phi$)-image 
converge as $s \rightarrow-\infty$ to the closed Reeb
orbit  whose $\theta_0$ parameter obeys $(p' + q') (1 - 3
\cos^2 \theta_0) = (p + q) \sqrt{6} \cos \theta_0$
and is such that
$\cos  \theta_0$ and $p' + q'$ have the same 
sign when the latter are non-zero.
 
The details of the definition of $\Phi$ near the 
$\phi$ image of a neighborhood of $\infty$ depends 
here on whether $k \equiv p + q$ or $k'\equiv p' + q'$ 
are zero and on their signs when they are
non-zero.  The discussion that immediately follows 
assumes that both are positive.  This 
discussion also covers the other possibilities after 
some straightforward sign changes.  The 
cases where one of $k$ or $k'$ is zero is discussed 
subsequently.  In what follows, $\sigma_0$ denotes 
$\sqrt{6} (k^2 + k'^{\,2} \sin^{2} \theta_0)^{1/2}$.
 
To see what is involved in the $k$ and $k'$ both 
positive case, note that when $k > 0$, then 
the function $u = \ln|\lambda|$ on $\mc^*\times \mc^*$ 
tends to $-\infty$ as $|z| \rightarrow 0$.  
Likewise, so does $v =
\ln|\lambda'|$  when $k' > 0$.  Moreover, the $\phi$ 
image of a neighborhood of $\infty$ has the form
\begin{itemize}
\item $k \varphi - k' t = k' \arg(a) - k \arg(a') + 
\Delta\pi + r_1 \;   \mod(2\pi)$,
\item $k v - k' u = r_2 $.
\autonum%(5.35)
\end{itemize}

\noindent
Here, $r_1$, $r_2$ and their derivatives are 
$\cO(|z|^{-1})$.  Of course, this is by virtue of the fact
that  where $|z| > 2$,
\begin{equation} 
\phi = (a r^k z^{-k} (1 + \cO(z^{-1})), 
a' r^{k'} z^{-k'} (1 + \cO(z^{-1})) \;.
%(5.36)
\end{equation}

To proceed, assume that $r \gg  e^{40} \sigma_0^4$, 
take $R \equiv R_\infty \gg 1$, and then take $f$ to be a 
favorite function of the variable $u$ which has the 
following properties:  First, $f(u)$ is as 
previously defined near the $\phi$ image of the $|z| < 10$ 
portion of $C_0$.  On the remaining
portion  of $C_0$, require that
\begin{itemize}
\item $f_u > 0$.
\item $f = u$ \qquad 	\qquad \qquad	{where}   $u > 2 \sin^{-2} \theta_0 R^{-1}$.
\item $f = (k' \sin^2 \theta_0)^{-1} R^{-1}
e^{\sigma_0u/k}$   	\qquad {where}   $u \le 0$.
\autonum%(5.37)
\end{itemize}

\noindent
Meanwhile, take $h$ to be a favorite function of 
the variable $v$ which is as previously defined 
near the $\phi$ image of the $|z| < 10$ portion of $C_0$ 
and is contrained near the remainder of
$\phi(C_0)$  to obey
\begin{itemize}
\item $h_v > 0$.
\item $h = v$   	\qquad	\qquad		{where}   $v > 2 R^{-1}$.
\item $h = k^{-1} R^{-1}$  \qquad 	{where}   $v < 0$.
\autonum%(5.38)
\end{itemize}

Note that the constraints in (5.37--38) are consistent 
with the previously specified constraints 
on $u$ and $v$ because $|u|$ and $|v|$ are 
$o(\ln r)$ where $10 < |z| < \sqrt{r}$.
 
By virtue of (5.35), these definitions of $f$ and $h$ 
imply that the $\Phi\circ\phi$--image 
of the $|z| > 10$ portion of $C_0$ has the form
\begin{itemize}
\item $k\varphi - k' t = k'\arg(a) - k \;\arg(a') + 
\Delta\pi + r_1 \;   \mod(2\pi)$,
\item $k \sin^{-2} \theta_0 h - k' \sin^2 \theta_0 f = r_2$ .  
\autonum%(5.39)
\end{itemize}

Here, $r_1$ and $r_2$ can be viewed as functions where 
are defined in a neighborhood of the $\Phi\circ\phi$--image
of the $|z| > 10$ portion of $C_0$.  In this regard, 
such a neighborhood can be defined so 
that $r_1$, $f^{-1} r_2$ and their derivatives are 
$o(e^{-|s|/\zeta})$ for some positive constant $\zeta$.  
(Note
that $r_2$ is  not the same function as that which appears in (5.35).)

As before, since the $\phi$ image of $C_0$ in $\mc^*\times \mc^*$ 
is holomorphic for the almost 
complex structure which maps $\partial_t$ to $\partial_u$ 
and $\partial_\varphi$ to $\partial_v$, so the $\Phi\circ\phi$  image
of the $|z| > 10$  portion of $C_0$ is pseudoholomorphic 
for an almost complex structure which has
the same  form as that depicted in (5.28).  
Moreover, as (5.29) holds for this almost complex 
structure, so does (5.30).
 
Now consider the definition of $\Phi$ near the 
$\phi$ image of a neighborhood of $z = \infty$ in 
the case where $k > 0$ and $k' = 0$.  The definition 
of $\Phi$ when $k < 0$ and $k' = 0$, or when $k = 0$ 
and $k'\not= 0$ is omitted since the definition in 
these cases is obtained from the $k > 0$, $k' = 0$ 
definition by changing notation.
  
To begin, remember that $\sin^2 \theta_0 = 1$ and 
$\sigma_0 = \sqrt{6} k$ because $k' = 0$.  Also, note that 
the function $u$ tends to $-\infty$ as $|z|$ 
tends to $\infty$ on $C_0$ because $k > 0$.  However,
when $r$ is very  large, then $u$ is on the order of  
$2^{-1} k \ln r$ on the circle where $|z| =
\sqrt{r}$.  With the preceding  understood, define the 
embedding $\Phi$ near the $\phi$ image of the $|z| > 10$
portion of $C_0$ by taking  $f$ to be any function of 
the coordinate $u$ with the following
properties:
\begin{itemize}
\item $f_u > 0$.
\item $f = u$   \quad 	where  $u > 2 R^{-1}$.
\item $f = R^{-1} e^{\sqrt{6} u}$  \quad	where $u < 0
$.\autonum%(5.40)
\end{itemize}

Here, $R \equiv R_\infty\gg 1$.  Then, with $f(u)$ chosen, 
take $h$ to be any function of both $u$ and $v$
which  obeys
\begin{itemize}
\item $h_v > 0$.
\item $h = v$  \quad  		where  $u > 2 R^{-1}$.
\item $h = 6^{1/2} f(u)v$  \quad	where  $u < 0 $.
\autonum%(5.41)
\end{itemize}

These definitions of $f$ and $h$ complete the 
definition of $\Phi$ near the $\phi$--image of a 
neighborhood of $z = \infty$ in $C_0$.  
Note that the $\Phi\circ\phi$ image of a neighborhood in $C_0$ of $\infty$ 
obeys the $k' = 0$ version of (5.39).  Moreover, 
the image of such a neighborhood is 
pseudoholomorphic as defined by an almost complex 
structure $J'$ that sends $\partial_t$ to the vector 
field $f_u\partial_f + h_u\partial_h$ and 
$\partial_\varphi$ to $h_v\partial_h$.  In this regard, 
note that (5.35), (5.40) and
(5.41) insure that  $J$ and $J'$ are nearly identical 
where $s$ is large  near $C_0$. 
Indeed, the just  referenced equations imply that 
$|J' - J| \le e^{-ss/\zeta}$ at large $s$ where $\zeta\ge 1$
is a constant.

\sh{(e)\qua  $\Phi$ as a map into $\mr\times(S^1\times S^2)$}
 
The first task for this subsection is to verify 
that $\Phi$ can be constructed so that its 
image avoids the locus where both $f$ and $h$ vanish.  
Having verified that such is the case, $\Phi$ is 
then interpreted as a map into $\mr\times(S^1\times S^2)$.  
The subsection next summarizes the
salient  features of $\Phi$ as a map into the latter 
space. In particular, as $\Delta > 0$ and $q' - p' > 0$
unless both  $p'$ and $q'$ are non-zero and have 
the same sign, the map $\Phi$ can be defined as above so
that it  and the resulting $C'\equiv (\Phi\circ\phi)(C_0)\subset 
\mr\times(S^1\times S^2)$ have the
properties listed below:

\begin{itemize}
\item $\Phi$ is an embedding on a neighborhood of $\phi(C_0)$.
\item There exists $\delta > 0$ such that $\sin \theta > \delta$ on $C'$.
\item $C'$ has two ends on the convex side of 
$\mr\times(S^1\times S^2)$ and one on the concave
side.   Moreover, one of the convex side ends is 
described by (5.27) where $|s|$ is large, and the 
other by (5.27) with $(q, q')$ replacing $(p, p')$.  
Meanwhile, the concave side end of $C'$ is 
described by (5.39) where $s$ is large.
\item $C'$ is symplectic and it is pseudoholomorphic 
for an almost complex structure,  $J'$, on
$\mr\times(S^1\times S^2)$ with the following list of features:
 
\begin{itemize}
\item[\rm(a)]  	$J'$ tames $\omega$.
\item[\rm(b)]  	$J' = J$ where  $\sin \theta < \delta/10$.
\item[\rm(c)]   $J'\partial_t = a_t \partial_f + b_t \partial_h$  
and   $J'\partial_\varphi = a_\varphi \partial_f + b_\varphi
\partial_h$; moreover, $J'$ is $T$--invariant so the  
coefficients $(a_r, b_t, a_\varphi, b_\varphi)$ depend only on
the coordinates $f$ and $h$. 
\item[\rm(d)]  There exists $\zeta \ge 1$ such that  $|J' - J| \le \zeta 
e^{-|s|/\zeta}$  on 
$\mr\times(S^1\times S^2)$.\break
\hbox{}\autonum %(5.42)
\end{itemize}
\end{itemize}

By the way, to say that $J'$ tames $\omega$ is to 
assert that $\omega(\cdot, J'(\cdot))$ is a positive definite,
bilinear  form on $T(\mr\times(S^1\times S^2))$.
  
The proof of the assertions of (5.42) rounds out 
the discussion in this subsection.  
In this regard, the arguments for (5.42) and 
also the argument that justifies $\Phi$'s 
interpretation as a map into $\mr\times(S^1\times S^2)$ 
make certain requirements on the parameters $r$,
$R_0$,  $R_1$ and $R_\infty$, the first being that they should all be large.

The first task is to prove that the image in $(S^1 \times \mr)
\times (S^1 \times \mr)$ via $\Phi$ of a 
neighborhood of $\phi(C_0)$ in  avoids the $f = h = 0$ locus.  
For this purpose, and for use in the 
discussion on (5.42), it proves useful to describe $C_0$ 
as the union of four open sets, $U$, $U_0$, 
$U_1$ and $U_\infty $.  Here, $U = \{z \in C_0: |z| > r^{-1/2}, 
|z - 1| > r^{-1/2}\; \text{and}\;|z| < r^{1/2}\}$. 
Meanwhile, $U_0$ is the  subset of $C_0$ on which $|z| 
< 2 r^{-1/2}$, $U_1$ is where $|z - 1| < 2
r^{-1/2}$ and $U_\infty$  is where $|z| > 2^{-1} r^{1/2}$.
  
To continue, note that one or both of
the $\Phi$--pullbacks $h$ and $f$ have large absolute  
value in $\phi(U_0)$ and in $\phi(U_1)$ so the $\Phi\circ\phi$
image of both of these sets avoids the $h = f = 0$  locus.  
Meanwhile, on $\phi(U)$, the $\Phi$--pullbacks of 
$f$ and $h$ are identified  with $u$  and $v$.  
Since $\Delta\not=  0$, these are both zero only
when $|z| = |1 - z| = r$.  In particular, when $r >  100$, 
the pullbacks of $f$ and $h$ by $\Phi$ are not
simultaneously zero on $\phi(U)$.  Finally, the  
definition of $\Phi$ near $\phi(U_\infty)$ 
insures that the $\Phi$
pullbacks of $f$ and $h$ are non-zero.  For  example, 
with the case where both $k$ and $k'$ are
positive, $f$ is defined by (5.37) to be positive  
as long as $u \le 2 \sin^{-2}\theta_0 R_\infty ^{-1}$ and $f$ is
defined to equal $u$ when $u > 2
\sin^{-2}\theta_0 R^{-1}$.  In this  regard, note 
that the $u = 2 \sin^{-2} \theta_0 R_\infty^{-1}$ locus is an
embedded circle where $|z| = r \exp(-2 (k 
\sin^{-2} \theta_0 R_\infty )^{-1})(1 + \cO(1/r))$.  
Moreover, as the derivative of $u$ with respect to
$|z|$ on this circle  is equal to $-k/|z| + \cO(1/|z|^2)$, 
the function $u$ and thus $f$ only increase on
the inside of this  circle.

With the preceding understood, $\Phi$ should henceforth be viewed as a map from a 
neighborhood of $\phi(C_0)$ into $\mr\times(S^1\times S^2)$ and 
then the next order of business is to
verify  that the assertions in (5.42) are correct.  
For this purpose, note that the third point of
(5.42)  follows directly from the definition of $\Phi$.  
Likewise, Parts a) and d) of the fourth point 
follow directly from the definition as does 
Part c) if the asserted $T$--invariance of $J'$ is 
disregarded.  Meanwhile, Part b) of the fourth point 
follows by standard arguments from 
the second point of (5.42).  Thus, the subsequent 
discussion concerns only the first and 
second points of (5.42) and the assertion of $T$--invariance 
in part c) of the fourth point.  In 
this regard, note that the first point in (5.42) and 
the $T$--invariance assertion are both 
consequences of the following lemma:
\medskip

\noindent
{\bf Lemma 5.5}\qua{\sl When $r$ is large, the parameters 
$R_0$, $R_1$ and $R_\infty$  can be chosen for
the  definition of $\Phi$ with arbitrarily large minimum 
and so that the set $U$ contains all distinct 
pairs $z$, $z' \in C_0$ that are mapped to the 
same point in $\mr^2$ by the $\Phi\circ\phi$ pullback of $(f, h)$.}
\medskip

Note that the ability to make the minimum of $R_0$, 
$R_1$ and $R_\infty$  large ensures that $\Phi$
defines a  map into $\mr\times(S^1\times S^2)$.
  
The proof of Lemma 5.5 is supplied momentarily.  
Granted this lemma, the first 
point in (5.42) can be argued as follows:  
First, as the map $\Phi$ is both proper on $C_0$ and a 
local diffeomorphism, it is enough to demonstrate 
that all distinct pairs $\{z, z'\} \in C_0$ which 
are sent to the same point by $\Phi\circ\phi$ 
lie in $U$ and are thus already identified by $\phi$.  This
last  property is, of course, guaranteed by Lemma 5.5.
  
To argue for the $T$--invariance assertion in part c) 
of the fourth point of (5.42), note 
first that one strategy for the construction of 
the required $J'$ would take $J'$ as given near $C'$ 
by (5.28) or (5.34) and rotate the latter via the $T$--action on 
$\mr \times (S^1 \times S^2)$.  In particular,
this  strategy would exploit the explicit lack of 
$(t, \varphi)$ dependence of the functions which multiply 
the vector fields on the right-hand side of (5.28) 
or (5.34).  However, as the functions in 
(5.28) and (5.34) are defined \emph{a priori} only near $C'$, 
such a strategy has the following 
prerequisite for success:  When a $T$ orbit intersects $C'$ 
more than once, then the versions of 
$J'$ given by (5.28) or (5.34) at the various 
intersection points must all agree.  Lemma 5.5 
provides just this prerequisite.
\medskip

\noindent
{\bf Proof of Lemma 5.5}\qua  Let $\cP\equiv\Phi^*(f, h)\co \phi(C_0) 
\rightarrow \mr \times \mr$.  
By virtue of its very
definition,  this map $\cP$ embeds the subset of of $\phi(C_0)$ 
where $\Phi^*(f, h) = (u, v)$.  This means, in 
particular, that $\cP$ embeds $\phi(U)$.  Also, $\cP$ 
separately embeds the $\phi$ images of $U_0$, $U_1$ and
$U_\infty$   since the pullbacks via $\Phi$ of the 
functions $f$ and $h$ have been constructed to change 
monotonically with $u$ and $v$, near these images.  
Thus, after the introduction of $V_{0,1,\infty} \subset
U_{0,1,\infty}$  to denote the subset where $(\Phi\circ\phi)^*(f,
h)$ differs from $(u, v)$, the question at issue here  
is whether any of $\cP\circ\phi(V_{0,1,\infty })$ intersect
$\cP\circ\phi(U)$, whether $\cP\circ\phi(V_\infty )$ intersects any 
of  $\cP\circ\phi(U_{0,1})$ and whether $\cP\circ\phi(V_0)$ intersects $\cP\circ\phi
(U_1)$ or $\cP\circ\phi(V_1)$ intersects $\cP\circ\phi(U_0)$.          

To address this question, note first that with
$R_\infty$  very large, both $f$ and $h$ will be  
nearly zero on $\phi(V_\infty )$.  Meanwhile, as $\Phi^*(f^2
+ h^2)$ is bounded away from zero on $\phi(U)$,  $\phi(U_0)$ and $\phi(U_1)$ 
by an $r$--dependent, but $R_\infty$--independent
constant.  Thus, when $r$ is large and  then $R_\infty$  
is very large, the $\cP$ images of
$\phi(V_\infty )$ and each of $\phi(U)$, $\phi(U_0)$ or $\phi(U_1)$ 
are  necessarily disjoint.  Likewise, when
$R_0$ and $R_1$ are very large, then at least 
one of the $\Phi$  pullbacks of $f$ or $h$ has large absolute
value on $\phi(V_0 )$ and $\phi(V_1)$ where these pullbacks  
differ from $(u, v)$.  Meanwhile, both $f$ and $h$
are uniformly bounded on $\phi(U)$ once $r$ has been  
specified.  Thus, when $R_0$ and $R_1$ are
sufficiently large, there are no intersections 
between  the $\cP$ images of $\phi(U)$ and those of $\phi(V_0 )$ and
$\phi(V_1)$.  

The final possibility for trouble lies with 
the $\cP$ images of $\phi(V_0 )$ and $\phi(V_1)$.  
However, this issue is nontrivial only when $p$ and $q$ 
are both non-zero and have the same 
sign while $p'$ and $q'$ are also both non-zero 
and have the same sign.  With these sign 
equalities now understood, the $\Delta > 0$ condition 
implies that $q'/q > p'/p$.  The argument when 
$q'/q < 0$ is the same but for cosmetic changes as that 
for the case when $q'/q > 0$, so only the 
$q'/q > 0$ case is discussed below.  In any event, 
given that $q'/q>0$, (5.22), (5.25) and (5.26) 
imply that the $\Phi$ pullback of the ratio $h/f$ obeys
\begin{itemize}
\item $\Phi^*(h/f) = p'/p + \cO(r^{1/2})$   on  $\phi(U_0 - V_0 )$.
\item $p'/p + \cO(R_0^{-\delta}) \ge  \Phi^*(h/f)  \ge  p'/p \sin^2 \theta_
{0P} - \cO(R_0^{-\delta})$ 
on $\phi(V_0 )$.
\autonum%(5.43)
\end{itemize}
Here, $\theta_{0P}$ is the value of $\theta$ on the closed Reeb 
orbit that is associated to the end 
$\Phi\circ\phi(V_0 )$ of $C'$.  Also, $\delta > 0$ is determined by the 
integers $\{(p, p'), (q,
q')\}$. Meanwhile, change $p$ to $q$, $p'$ to 
$q'$ and $P$ to $Q$ in (5.43) and the resulting
inequalities then  describes the $\Phi$ 
pullback of $h/f$ to $\phi(U_1)$.
 
Now, as argued in Step~9 of the proof of Constraint~4 in 
Subsection~5b, when $q'/q$ 
and $p'/p$ are both positive and $q'/q$ is greater 
than $p'/p$, then $q'/q \sin^2 \theta_{ 0Q} > p'/p \sin^2
\theta_{0P}$.   With this point understood, and given both 
(5.43) and its $(q, q')$ analog, the $\cP$ images of 
both the pairs $\{\phi(V_0 ), \phi(U_1)\}$ and $\{\phi(V_1), \phi(U_0)\}$ 
are necessarily disjoint when $r$ is first 
chosen to be very large, then $R_0$ is chosen to be much 
larger than $r$ and finally, $R_1$ is chosen 
to obey $R_1 \gg  256 (p/q) R_0$.  Indeed, such a choice 
insures that the ratios of $\Phi^*(h/f)$ differ on 
$\phi(V_0 )$ and $\phi(U_1)$, as do the analogous ratios 
on $\phi(V_1)$ and $\phi(U_0)$.  
\medskip

With the proof of Lemma 5.5 complete, only the second 
point of (5.42) is unspoken 
for.  Here is the strategy for proving that such $\delta$ 
exists:  The first step proves that $\delta > 0$ 
exists so that $\sin \theta > \delta$ on the $\Phi\circ\phi$ 
images of $U_0$, $U_1$ and $U_\infty$.  Having
done so, the  remaining step considers whether 
$\sin \theta$ is ever zero on the $\Phi\circ\phi$ 
image of $U$.  As $U$ has 
compact closure, these two steps suffice to 
establish the second point of (5.42).  By the way, 
in order to prove that $\sin \theta \not=  0$ on a given set, 
it is enough to prove that $f$ is positive
where $h$  vanishes on the set in question.

To begin, consider the case with $U_\infty$.  
Here, it follows from (5.35) and (5.39) that 
when both $k$ and $k'$ are non-zero and both $r$ 
and $R_\infty$  are large, then the required $\delta$
exists.   Thus, the only troublesome case has $p' + q' = 0$ 
and $p + q < 0$.  However, as $\Delta > 0$, this
case  requires $q' < 0$ and $p' > 0$ 
and so is excluded from consideration.
  
Now consider whether the required $\delta$ exists 
for $U_0$ when $R_0$ is large.  For this 
purpose, note that the case where $p' = 0$ and $p < 0$ 
can be ignored as the conditions $\Delta > 0$ 
and $q' - p' > 0$ would otherwise be violated.  
Thus, one can assume that $p' \not=  0$ or else $p' =
0$  and $p > 0$.  As explained momentarily, 
the existence of the required $\delta$ then follows from 
(5.22), (5.25) and (5.26).  Indeed, to obtain $\delta$, 
note first that when $R_0$ is large, (5.26) 
provides a non-zero lower bound for $\sin \theta$ on the 
$\Phi\circ\phi$ image of $V_0  \subset U_0$ when $p' >
0$.   Meanwhile, if $p' = 0$, then $p > 0$ and so (5.25) 
insures that $\Phi^*f > 0$ on $\phi(V_0 )$.  Thus, it 
remains only to verify that $\sin \theta > 0$ on the closure 
in $C_0$ of the subset of $U_0$ where $\Phi^*f = u$ 
and $\Phi^*h = v$.  For this purpose, note that 
$v = 0$ on $\phi(C_0)$ only on the image of points $z$ 
where
\begin{equation}
r^{p'+q'} = |z|^{p'} |1 - z|^{q}
%(5.44)
\end{equation}
and so such a zero occurs with $|z| \ll 1$ if and 
only if $p'$ and $p' + q'$ have opposite signs.  
This requires $q'$ and $p'$ to have opposite 
signs and thus $p'$ is negative since $q' - p'$ must be 
positive.  But, if $r$ is large and $p' < 0$, 
then (5.22) requires that $u > 0$ where $v = 0$.

An analogous argument  finds a non-zero 
lower bound for $\sin \theta$ on the $\Phi\circ\phi$ 
image of $U_1$.

With the preceding understood, it remains 
only to verify that the zeros of the $\Phi$ pullback
of $h$ on $\phi(U)$ occur where the $\Phi$ pullback of $f$ 
is positive.  In this regard, the 
straightforward case to verify occurs when either 
one of $p'$ or $q'$ is zero, or when both $p'$ 
and $q'$ are non-zero and they have the same sign; 
for in this case, $\Phi^*h$ has no zeros at all on 
$\phi(U)$.  Indeed, as $\Phi^*h = v$ on $\phi(U)$ so at 
issue are the zeros of $v$.  As $p'$ and $q'$ have the 
same sign if one is not zero, the latter locus is 
placed by (5.44) where $|z| = r + O(1)$ which is 
not in $U$.

The next case to verify has $p' + q' = 0$.  In this case, 
$\Phi^*h = v$ is zero on $\phi(U)$ only 
on the $\phi$ image of the line where $|z| = |z - 1|$.  On this 
line, $\Phi^*f = u = (p + q) \ln (r/|z|)$.  
Now, the positivity of $\Delta$ and the vanishing of 
$p' + q' = 0$ and the positivity of $q'$ requires the 
positivity of $p + q$.  Thus, as $|z| \le r^{1/2}$ on 
$U$, so $\Phi^*f > 0$ on the zero locus in $\phi(U$) of
$\Phi^*h$. 

Next, assume that neither $p'$, $q'$ or 
$p' + q'$ are zero but $p'$ and $q'$ have oppositive 
signs.  Thus, $p' < 0$ and $q' > 0$.  In this case, 
the locus in (5.44) has two components.  One 
component has $|z| = r + \cO(1)$ so doesn't sit in $U$.  
The other component has  $|z|$ or $|1 
- z|$ small and can lie in $U$.   In particular, if 
$p' + q' > 0$, this second component occurs 
where $|z| = r^{(p'+q')/p'}$; and if $p' + q' < 0$, 
this second component occurs where $|1 - z| = 
r^{(p'+q')/q'}$.  In any event, the sign of $u = \Phi^*f$ 
on the $\phi$--image of this other component can be 
determined with the help of one or the other of 
two identities that follow from (5.20):
\begin{equation}
q v - q' u = - \Delta \ln(r/|z|)\quad \hbox{\rm and}\quad     
p v - p' u = \Delta\ln(r/|1 - z|) .
%(5.45)
\end{equation}
In particular, as $\Delta > 0$, it follows from (5.45) that 
$u > 0$ on these small $|z|$ or $|1 - z|$ 
components of the $v = 0$ locus.

\sh{(f)\qua  Deformations}

Suppose that $I \equiv \{(p, p), 
(q, q')\}$ are integers such that $\Delta > 0$ and $q' - 
p' > 0$ unless both $p'$ and $q'$ are non-zero and 
have the same sign.  Take $r \gg  1$ to be large, 
and choose large $R_0$, $R_1$ and $R_\infty$  so that 
the map $\Phi$ is defined as in the preceding two 
subsections near $\phi(C_0)$ as an embedding into 
$\mr\times(S^1\times S^2)$ as described in (5.42).  As
in  (5.42), use $C'$ to denote $\phi(C_0)$.  This $C'$ 
is the starting member of a set of symplectic 
subvarieties which is parameterized by a non-trivial 
subinterval $[0, T_1] \subset [0, 1]$ whose end 
member, $C$ is in the space $\fM_I$ of Proposition 5.2.  
Thus, $C$ is a $J$--pseudoholomorphic, thrice-punctured 
sphere that is parameterized by the same set $I$ as is 
$C'$, and whose existence is 
asserted in Proposition 5.1.  This subsection describes 
the relevant parameterized set of 
symplectic subvarieties and it provides a proof 
that the end member of the set is a thrice-punctured
sphere whose existence is predicted by Proposition 5.1.
 
The promised parameterized set of symplectic 
subvarieties is constructed with the 
help of a set of almost complex structures on $\mr\times(S^1\times S^2)$.  
If $\hat{J}$  is an almost complex 
structure in this set, then $\hat{J}$  is constrained to 
obey the following four conditions:
\begin{itemize}
\item[\rm(a)]  $\hat{J}$  tames  $\omega$.
\item[\rm(b)]    $\hat{J}= J$   where  $\sin \theta < \epsilon$.
\item[\rm(c)]  $|\hat{J}  - J| \le \epsilon^{-1} 
e^{-\epsilon|s|}$  on $X$.
\item[\rm(d)]  $\hat{J} \partial_t = \alpha_t \partial_f + \beta_t 
\partial_h$   and     
$\hat{J}\partial_\varphi = \alpha_\varphi \partial_f + \beta_\varphi
\partial_h$  where the coefficients $(\alpha_r, \beta_t, \alpha_\varphi, 
\beta_\varphi)$ 
depend only on the coordinates $f$ and $h$; thus $\hat{J}$  is $T$--invariant.
\autonum%(5.46)
\end{itemize}
Here, $\epsilon > 0$ is determined by $I$. The particular set of 
almost complex structures under 
consideration is a certain, continuous, 1--parameter 
family $\{J_r\}_{r\in[0,1]}$ with $J_1 = J$.  Note that 
each   $\hat{J}\in \{J_ r\}_{r\in[0,1]}$ obeys
\begin{equation}
\hat{J}  = J\quad  \hbox{\rm where} \quad   |s| > 1/\epsilon ,
%(5.47)
\end{equation}
which is a stronger condition than (5.46c).  (The 
condition in (5.47) is not imposed by 
necessity but to simplify subsequent arguments.)

The complex structure $J_r$ is defined by the requirement that
\begin{equation}
J_r v  = (1 - r) (\chi_\epsilon J' v + (1 - \chi_\epsilon) J v) + r J v
%(5.48)
\end{equation}
when $v = \partial_t$ or $\partial_\varphi$.  Here, $\chi_\epsilon$ is a smooth 
function of $|s|$ that equals 1 where $|s| <
1/(2\epsilon)$,  vanishes where $|s| > 1/\epsilon$ and has 
derivative bounded by $4 \epsilon$.  The equality in (5.47)
follows  from the presence of $\chi_\epsilon$ in (5.48).  In any 
event, $J_r$ satisfies (5.46) because both $J'$ and
$J$ do.
 
As stated above, the parameter, $r$, for the parameterized 
set, $\{\Sigma_r\}$, of subvarieties takes 
values in $[0, 1]$ where $T \in (0, 1]$.  For each 
such $r$, the corresponding $\Sigma_r$ is a $J_r$--pseudoholomorphic, 
immersed, thrice-punctured sphere
that avoids the $\theta\in\{0, \pi \}$ locus and  whose 
ends are constrained by the set $I$ exactly as
those of $C'$ are via the third point in  (5.42).  
The construction of $\Sigma_r$ is facilitated by the
following:
\bigskip

{\bf Lemma 5.6}\qua  {\sl Let $\hat{J}$  denote an almost complex 
structure on $\mr\times(S^1\times S^2)$ which
obeys  (5.46).  Suppose that $I = \{((p, p')(q, q')\}$ 
is a set of pairs of integers  
with $pq'  - qp' > 0$ and such that $q' - p' > 0$ 
when both $p'$ and $q'$ are non-zero and do not
have the  same sign.  Given $I$, let $\Sigma\subset \mr\times(S^1\times S^2)$ 
denote a subvariety that is the
image via a  $\hat{J}$--pseudoholomorphic map of a thrice-punctured 
sphere whose ends are constrained by $I$
as  those of $C'$ are by the third point of (5.42).  
In addition, require that $\Sigma$ is disjoint from
the  $\{0, \pi \}$ locus.  
Then $\Sigma$ has the following properties:}
\begin{itemize}
\item $\Sigma$ {\sl is immersed and so the deformation 
operator $D$ as described by (2.6) is well defined.}
\item {\sl This operator $D$ has  three-dimensional 
kernel and trivial cokernel.}
\end{itemize}
\bigskip

This lemma is proved below, so accept its validity for now.  
In particular, given the 
first two points of the lemma, a straightforward 
extension of the arguments in Section~3c for 
Proposition 3.2 establish the following two key facts:  
\begin{itemize}\sl
\item There exists $\epsilon_0 > 0$ such that $J'$ and $J$ 
obey (5.46) with $\epsilon\le \epsilon_0$; and for each such  $\epsilon$,
the  almost complex structure $J_0$ from (5.48) 
admits a $J_0$--pseudoholomorphic map from $C_0$ 
into $\mr\times(S^1\times S^2)$ whose image, $\Sigma_0$ is 
disjoint from the $\theta\in \{0, \pi \}$ locus and has
its  ends constrained by $I$ as those of $C'$ are by 
the third point of (5.42).  
\item Fix  $\epsilon < \epsilon_0$ so that each $J_r$ in (5.48) obeys 
the constraints in (5.46) and (5.47).  Let $\Lambda\subset
[0, 1]$ denote the set of points $r$ such that 
there exists a $J_r$--pseudoholomorphic map $C_0$ 
into $\mr\times(S^1\times S^2)$ whose image, $\Sigma_r$ is 
disjoint from the $\theta\in \{0, \pi \}$ locus and has
its  ends constrained by $I$ as those of $C'$ are by 
the third point of (5.42).  This set is non-empty
and open. 
\autonum%(5.49)
\end{itemize} 
Indeed, the assertions in (5.49) are simply perturbation 
theoretic consequences of the 
vanishing of the cokernel of $D$ for the subvariety $C'$.  
The arguments are essentially 
identical to those for Proposition 3.2 that 
appear in Section 3c.

By the way,  note that $J'$ in (5.42) obeys 
(5.46) but not (5.47).  However, a small 
perturbation at very large values of $|s|$ on 
$\mr\times(S^1\times S^2)$ gives an almost complex
structure $J_0$  that obeys both (5.46) and (5.47) 
with $\epsilon$ as small as desired.  This understood,
perturbation  theory using the invertibility of the 
$C'$ version of $D$ produces the required $\Sigma_0$ as a
small  deformation of $C'$.

Now, let $T \in [0, 1]$ denote the least 
upper bound of those $r \in\Lambda$. Of prime 
importance is the limiting behavior as $r \rightarrow T$ 
of the set $\{\Sigma_r: r \in\Lambda\}$. The nature of this limit, 
as given in the next lemma, explains how this 
parameterized set determines a subvariety $C$ in 
Proposition 5.2's moduli space $\fM_I$. 
\bigskip

\noindent
{\bf Lemma 5.7}\qua  {\sl Let  $\Lambda$ be as just described 
and let $T$ denote the least upper bound of $\Lambda$.
Then  there exists a  $J$--pseudoholomorphic, triply punctured 
sphere $C$ in the moduli space $\fM_I$ of 
Proposition 5.2, a countable, increasing set 
$\{r(i)\}_{i=1,2, \ldots}\subset\Lambda $ and a corresponding set 
$\{w(i)\}_{i=1,2,\ldots} \subset \mr$ with the following properties:} 
\begin{itemize}
\item $\lim_{i\rightarrow\infty}  r(i) = T$,
\item {\sl Let $\Sigma_i$ denote the push-forward of 
$\Sigma_{r(i)}$ via the diffeomorphism of $\mr\times(S^1\times
S^2)$ defined  as translation by $w(i)$ on the $\mr$ factor.  Then }
\[
\lim_{i\rightarrow\infty}  (\text{sup}_{x\in(C\cap K)}
{\rm dist}(x, \Sigma_i) + 
\text{sup}_{x\in(\Sigma_i\cap k)}
{\rm dist}(C, x)) = 0
\]
{\sl for all compact sets  $K \subset \mr\times(S^1\times S^2)$.}
\end{itemize}
\bigskip

By the way, the arguments below imply that the 
sequence $\{w(i)\}$ has no convergent 
subsequences when $T < 1$.

Note that the assertion in Proposition 5.1 that $\fM_I$ 
is non-trivial follows directly from 
this last lemma.

The remainder of this subsection is occupied with the 
proof of Lemma 5.6, while the 
next subsection gives that of Lemma 5.7.
\bigskip

\noindent
{\bf Proof of Lemma 5.6}\qua   Modulo two observations, 
the argument for the assertion that $\Sigma$ is 
immersed is identical to the argument in Section 4f 
for the proof of Proposition 4.7.  The 
first observation asserts that the form $\kappa_t dt - \kappa_\varphi d\varphi$ 
pulls back to $\Sigma$ with at most a single 
non-degenerate, hyperbolic zero when $(\kappa_t, \kappa_\varphi)$, are 
constants that define a vector which is not a 
multiple of either $(p', p)$, $(q', q)$ or $(k', k)$.  
Indeed, to see why this must be the case,
note  first that the form of $\hat{J}$  in (5.46) implies that 
the  $\hat{J}$--pseudoholomorphic map sending $C_0$ to 
$\Sigma$ pulls back $\kappa_t t - \kappa_\varphi\varphi$ 
as a multivalued function, $\nu $, 
that obeys an equation with the 
schematic form $\Delta\nu + r\cdot d\nu = 0$.  Here, $\Delta$ is a 
Laplacian on $C_0$ and $r$ is an appropriate
vector  field.  The maximum principle applies to this 
last equation and rules out local maxima and 
minima for $\nu$.  Thus, non-degenerate zeros of $d\nu$ 
are hyperbolic.  Moreover, this same 
equation for $\nu$ implies that degenerate zeros of $d\nu$ 
provide a count of $-2$ or more to the Euler 
characteristic of $C_0$.  Thus, $d\nu$ has exactly one 
zero, which is non-degenerate and hyperbolic.

The second observation needed concerns the 
part of the proof of Proposition 4.7 
that follows (4.44).  In particular, it follows from (5.47) 
that this portion of the argument 
applies here with no changes provided that it is 
applied to an end of $\Sigma$ where the 
corresponding closed Reeb orbit has $\cos^2 \theta_0$ different from both 1/3 and 
1.  Since $\Delta > 0$, such an end is always present.

With the first point of Lemma 5.6 understood, the 
argument for the second point is 
identical in all essential aspects to that 
used in Section 4g to prove the $\aleph= 2$ case of 
Proposition 4.8.  

\sh{(g)\qua  The limits of the deformations}

The purpose of this subsection is to provide the 
proof of Lemma 5.7.  Being 
lengthy, the proof is divided into eleven steps.  
The first nine steps comprise what might be 
termed Part 1 of the proof.  These steps give the 
proof when none of $p$, $q$ and $p + q$ is zero.  
The last three steps provide the proof in the 
remaining cases.  In any event, keep in mind 
that Steps 1--9 implicitly assume that each 
of $p$, $q$ and $p + q$ is non-zero.
\bigskip

{\bf Step 1}\qua  Choose an increasing sequence 
$\{r(i)\}_{i=1,2, \ldots} \subset\Lambda$ with limit $T$.  This step 
defines a corresponding sequence $\{w(i)\}_{i=1,2,\ldots }\in \mr$.  
To start, remember that the proof
of  Lemma 5.6 established that $dt$'s pull back to $C_0$ 
via the $\Sigma_r$--defining, $J_r$ pseudoholomorphic 
immersion has a single zero, one that is non-degenerate 
and hyperbolic.  As neither $p$, $q$ nor 
$p + q$ is zero, the almost complex structure $J'$ 
from Section 5f that is described in (5.42) has 
$J'$ $dt$ proportional to $df$.  Therefore, this is 
also the case for $J_r$ as defined in (5.48) and so
$df$  also pulls back to $C_0$ by the $\Sigma_i$--defining 
map as a 1--form with but one zero, also non-degenerate 
and hyperbolic.  Let $x_r$ denote the image of
this zero on $\Sigma_r$.  With $x_r$ understood,  take $w(i)$ 
to equal minus the value of the coordinate $s$
at $x_{ri}$.  Thus, the action of $w(i)$ on $\mr\times(S^1\times S^2)$ 
by translation of the $\mr$ factor
sends $x_{r(i)}$ to a point, $\ux(i)$, which lies on the zero 
locus  of the function $s$.  (Although
certain arguments in this Part 1 of Lemma 5.7's proof rely on  
this critical point structure of $df$,
there are alternative arguments that do not.  In fact, some of  
the latter are used for Part 2 of
Lemma 5.7's proof.)

If the sequence $\{r(i)\}$ has a subsequence 
for which the corresponding sequence 
$\{w(i)\}$ is unbounded, then replace $\{r(i)\}$ by a 
subsequence for which the corresponding 
$\{w(i)\}$ sequence is unbounded and either strictly 
increasing or strictly decreasing.  Agree to 
relable this new $\{r(i)\}$ sequence and its 
corresponding $\{w(i)\}$ sequence by consecutive 
integers starting at 1.
 
Let $\underline{\Sigma}_i$ denote the image of $\Sigma_{r(i)}$ under the 
translation isometry defined by this same 
number $w(i)$.  Likewise, let $\underline{J}_i$ denote the 
push-forward by this isometry of the almost 
complex structure $J_{r(i)}$.  Thus, $\underline{\Sigma}_i$ is $\underline{J}_i$--pseudoholomorphic.  
Also, $\underline{J}_i$ obeys the conditions in 
(5.46), albeit with the functions $(\alpha_t, \beta_r,\alpha_\varphi, 
\beta_\varphi)$ 
replaced by their translated versions.  
\bigskip

{\bf Step 2}\qua This step and Step 3 argue that 
the integral of $\omega$ over the intersection of $\underline{\Sigma}_i$ 
with a given compact set enjoys an upper bound which is 
independent of the index $i$.  In 
particular, this step considers the case when neither 
$p'$ nor $q'$ are zero.  In this regard, 
remember that $p$ and $q$ are already assumed to 
be non-zero.  The subsequent argument 
makes no use of the assumption that $p + q \not=  0$.

In any event, when $p'$ and $q'$ are non-zero, 
then a compact set $K$ intersects $\underline{\Sigma}_i$ where 
$-f' \le f \le f'$ for suitable $f' > 0$.  Likewise, 
a compact set $K$ intersects $\underline{\Sigma}_i$ where $-h' \le h
\le h'$  for suitable $h'$.  Having said this, consider 
first the subcase where the sequence
$\{w(i)\}$ is  bounded.  Under this extra assumption, 
the constants $f'$ and $h'$ can be large to
guarantee that  $\underline{J}_i = J$ where $|f| \ge f'$ and where 
$|h| \ge h'$.  On this same portion of $C$, $dt \wedge
df$ and $d\varphi\wedge dh$ are  non-negative and so
\begin{equation}
\int_{K\cap\underline{\Sigma}_i}
\omega\le\int_{\{-f'\le f\le f'\}
\cap\underline{\Sigma}_i}
dt\wedge df+\int_{\{-h'\le h\le h'\}
\cap\underline{\Sigma}_i}
d\varphi\wedge dh\;.
%(5.50)
\end{equation}
With this last inequality in hand, apply Stokes' theorem to find
\begin{align}
\int_{K\cap\underline{\Sigma}_i}
\omega\le f'&
\left(\int_{\{f=-f'\}\cap\underline{\Sigma}_i}
dt-\int_{\{f=f'\}\cap\underline{\Sigma}_i}dt
\right)\notag\\
& +h'
\left(\int_{\{h=-h'\}\cap\underline{\Sigma}_i}
d\varphi-\int_{\{h=-h'\}\cap\underline{\Sigma}_i}
d\varphi\right).
%(5.51)
\end{align}
\noindent
With regard to the derivation of (5.51), note that 
the concave side end of $\underline{\Sigma}_i$ makes no 
contribution to the boundary terms in Stokes' theorem 
since both $f$ and $h$ limit to zero on 
this end as $s$ limits to infinity.  In any event, the 
inequality in (5.51) provides the promised 
upper bound for the integral of $\omega$ since the four path 
integrals which appear above are 
determined {\sl a priori} by the set I.

In the case where $\{w(i)\}$ has no upper bound, 
the preceding argument still applies 
with no extra comments required.  With one small 
revision, this same argument also applies 
to in the case where $\{w(i)\}$ has no lower bound.  
As before, there exists such $f'$ and $h'$ so 
that the set where $|f| < f'$ as well as that where 
$|h| < h'$ contain $K$.  Moreover, given that 
$\{w(i)\}$ has no lower bound, the part of $\underline{\Sigma}_i$ where 
either $|f| < f'$ or where $|h| < h'$ has $\underline{J}_i = J$ 
when $i$ is sufficiently large.  Thus, the inequalities 
in (5.50) and (5.51) still hold provided 
that $i$ is sufficiently large.  In particular, (5.51) 
again provides an $i$--independent bound for 
$\omega$'s integral over $K \cap\underline{\Sigma}_i$.
\bigskip

{\bf Step 3}\qua  Now consider the finiteness of the integral 
of $\omega$ in the case where either $p'$ 
or $q'$ is zero.  Here, the argument is very similar 
to the one just given in Step 2.  Indeed, the 
argument in this case requires but one substantial 
novelty.  The novel part of the argument 
exploits the fact that the function $h$ has precisely 
one non-compact level set when either $p'$ 
or $q'$ is zero.  Only the $p' = 0$ case is considered 
below as the considerations for the $q' = 0$ 
case are essentially the same.

As the existence of a single non-compact $h$--level 
set is proved momentarily, accept it 
for now and let $h_0$ denote the value of $h$ on this set.  
Then, there exists $h' > 0$ which depends 
only on $K$, and given $i$, there exists $\epsilon > 0$ 
such that the integral of $\omega$ over $K \cap\underline{\Sigma}_i$ is no more 
that twice that of $\omega$ over the subset of 
points in $\underline{\Sigma}_i$ where either $-h' < h \le h_0 - 
\epsilon$ or $h_0 + \epsilon 
\le h  < h'$.  This last conclusion constitutes the novel 
input.  (Even with $p'$ or $q'$ zero, there
exists  $f' > 0$ such that $K$ is contained in 
the subset where $|f| \le f'$.)
 
With the preceding understood, turn to the 
promised justification of the assertion 
that $h$ has but one non-compact level set.  For 
this purpose, it is necessary to return to $\Sigma_{r(i)}$.  
Having done so, the first point is that $h$ is bounded 
on the $(p, p' = 0)$ end of $\Sigma_{r(i)}$ because $\varphi$ 
could not have an $|s| \rightarrow \infty$  limit otherwise.  
Indeed, as $h$ has no limits as $|s| \rightarrow \infty$  on
$(q, q')$  end of $\Sigma_{r(i)}$ and limits to zero on 
the concave side end, were $h$ unbounded on the end
where  $(p, p' = 0)$, this end would contain whole 
components of regular value level sets of $h$
where  $J_r = J$ and where $|s|$ is everywhere large.  
But, $d\varphi$ pulls back without zeros to such a
level set,  so these properties preclude the 
existence of $\varphi$'s $|s| \rightarrow \infty$  limit.

Note that the preceding argument also shows 
that $h$ must have at least one 
non-compact level set, for otherwise there would be 
a component of a compact level set which sat 
entirely at large $|s|$ in the $(p, p' = 0)$ end of $\Sigma_{r(i)}$.
  
Now, as $h$ tends to zero on the concave side of $\Sigma_{r(i)}$, 
tends uniformly to infinity on the 
$(q, q')$ end of $\Sigma_{r(i)}$ and is bounded on the $(p, p' = 0)$ 
end of $\Sigma_{r(i)}$, it follows that the 
non-compact level set of $h$ must have their ends in 
the $(p, p' = 0)$ end of $\Sigma_{r(i)}$.  Thus, to prove
that  there is at most one non-compact level set of 
$h$, it is sufficient to establish the following 
lemma:
\bigskip

\noindent
{\bf Lemma 5.8}\qua {\sl The function $h$ has a unique $|s| \rightarrow \infty$  
limit on the $(p, p' = 0)$ end of
$\Sigma_{r(i)}$.  In  fact the restrictions of $h$ to the large 
but constant $|s|$ circles in this end
converge as $|s| \rightarrow
\infty$   in the $C^\infty$  topology to the constant function. 
}
\bigskip

\noindent
{\bf Proof of Lemma 5.8}\qua  Where $|s|$ is large on the 
$(p, p' = 0)$ end of $\Sigma_{r(i)}$, $J_r = J$.  Thus, by 
virtue of the asymptotics asserted in Proposition 2.3, 
this part of $\Sigma_{r(i)}$ can be parameterized as 
in (2.13) so that the column vector $\lambda$ with top 
component $x$ and bottom component $w$ obeys 
the analog,
\begin{equation}
\partial_\rho\lambda + L_0 \lambda = w \sigma ,
%(5.52)
\end{equation}
of (2.15).  Here, $L_0 \equiv\left(
\begin{smallmatrix}
0&\partial_\tau\\
-\partial_\tau&\sqrt{6}
\end{smallmatrix}\right)$   and $\sigma$ is a vector with 
$|\sigma| e^{\rho/\delta}$ bounded as 
$\rho\rightarrow \infty$  for 
some positive constant $\delta$.  (The entries of $\sigma$ are 
functions of $w$, $w_\rho$ and $w_\tau$ with affine 
dependence on the latter two.)  Note that it is 
permissable to assume that $x$ limits to zero 
here as $\rho\rightarrow \infty$.
 
Now, the point of writing (5.52) is to employ a 
simplified version of the analysis 
that was used in Steps 2--5 of the proof of Proposition 2.3 
as given in Section 2.  For this 
purpose, note that the operator $L_0$ has a zero 
eigenvalue with multiplicity 1, the eigenvector 
has constant top component $x$ and bottom component $w = 0$.  
Meanwhile, $L_0$'s smallest 
positive eigenvalue is $\sqrt{6}$ and the eigenvector 
has $x = 0$ and constant $w$.  Let $E_+ > \sqrt{6}$
denote  the next smallest positive eigenvalue of $L_0$ and 
let $-E_-$ denote the largest negative
eigenvalue  of $L_0$.  Next, let $f^+(r)$ denote the $L^2$ 
norm, defined by integration over the circle
parameterized  by $\tau$, of the $L^2$ orthogonal projection of 
$\lambda( \cdot, r)$ onto the span of the
eigenvectors of $L_0$ with  eigenvalue greater 
than $\sqrt{6}$.  Define $f^-(r$) in an analogous fashion
using the span of the  eigenvectors of $L_0$ with 
negative eigenvalue.  Likewise, define $f^{\sqrt{6}}(r)$
to denote the $L^2$ norm of  the $L^2$ orthogonal 
projection of $\lambda$ onto the span of the eigenvector of
$L_0$ with eigenvalue
$\sqrt{6}$.   According to Lemma 2.5, $\lambda$ has zero 
component along the zero eigenvalue eigenvector
of  $L_0$.
  
Here is a key point with regard to the 
functions $f^{\pm,\sqrt{6}}$:  As $y$ is bounded, so each of 
$f^{\pm,\sqrt{6}}$ is bounded by $\zeta e^{-\sqrt{6}\rho}$.
 
With the definitions of 
$f^{\pm,\sqrt{6}}$ in hand, then (5.52) 
implies that for all sufficiently large 
$\rho$,  
\begin{itemize}
\item $\partial_\rho f^+ +E_+f^+ \le \epsilon (f^- + f^{\sqrt{6}})$ ,
\item $\partial_\rho f^- -E_-f^- \ge- \epsilon (f^+ + f^{\sqrt{6}})$ ,
\item  $|\partial_\rho f^{\sqrt{6}}+ \sqrt{6}f^{\sqrt{6}}|
\le \epsilon (f^+ + f^-+f^{\sqrt{6}})$ ,
\autonum%(5.53)
\end{itemize}

\noindent
where $\epsilon\equiv\epsilon(\rho)$ is a positive, integrable function 
of $\rho$ on a domain of the form $[\rho_0, \infty$ ). 
Now,  given the bounds $f^{\pm,\sqrt{6}} < \zeta e^{-\sqrt{6}\rho}$ 
and the integrability of $\epsilon(\cdot)$, the
inequalities in (5.53) can be  integrated to yield 
\begin{equation}
f_+ + f_- \le \epsilon_1(\rho) e^{-\sqrt{6}\rho} \quad   \hbox{\rm and }\quad     
|e^{\sqrt{6}\rho} f^{\sqrt{6}} - c| \le \epsilon_1(\rho)
%(5.54)
\end{equation}
where $c$ is a constant and $\epsilon_1(\cdot)$ is a 
positive function of $\rho$ on a domain of the form $[\rho_1,
\infty $)  which limits to zero as $\rho\rightarrow \infty$ .
  
These last bounds imply that the restrictions, 
$\{h(\cdot,\rho)\}_{\rho\gg 1}$, of $h$ to the large but 
constant $\rho$ circles in the $(p, p' = 0)$ end of $\Sigma_r$ 
limit in the $L^2$ sense to a constant as $\rho\rightarrow 
\infty $.   This $L^2$ statement can be readily 
bootstrapped as in Steps 4 and 5 of Section 2's proof
of  Proposition 2.3 to establish that the restrictions, 
$h(\cdot,\rho)$ of $h$ to the large and constant $\rho$ 
circles limit pointwise and uniformly to a constant 
when $ \rho\rightarrow\infty$.  Moreover, these same 
arguments prove that the derivatives of $h(\cdot,\rho)$ 
to all orders converge to zero as $\rho\rightarrow \infty$.
\bigskip  

{\bf Step 4}\qua  With the results of the previous 
steps in hand, Proposition 3.3 in \cite{T4} can 
be invoked to conclude the following:   
First, there is a finite, non-empty set $\{(S_k, m_k)\}$ 
where $\{S_k\}$ is a set of distinct, irreducible, 
$J_T$--pseudoholomorphic subvarieties and each $m_k$ 
is a positive integer.  Second, after passing 
to a subsequence of $\{\uSigma_i\}$ (and renumbering 
consecutively from 1), this set converges pointwise to $C_T 
\equiv\cup S_k$ in that 
\begin{equation}
\lim_{i\rightarrow\infty}  ({\rm sup}_{x\in(C_T\cap K)}
 \dist(x, \Sigma_i) 
+ {\rm sup}_{ x\in(\underline{\Sigma}_i\cap K)} 
\dist(C_T, x)) = 0 .
%(5.55)
\end{equation}
for all compact sets $K \subset \mr\times(S^1\times S^2)$.  Third, and 
here is where the integers $\{m_k\}$
enter,  the sequence of currents defined by $\{\uSigma_i\}$ converges 
to the current $\sum_k m_k S_k$ in the sense
that  if $\kappa$ is any smooth, compactly supported 
2--form on $\mr\times(S^1\times S^2)$, then
\begin{equation}
\lim_{i\rightarrow\infty}\int_{\uSigma_i}
\kappa=\sum_k m_k\int_{S_k}\kappa.
%(5.56)
\end{equation}

By the way, with regard to the conclusion that $C_T \not= \emptyset$, 
this follows from (5.55) since 
$\uSigma_i$ contains a point on the fixed, compact submanifold 
in $\mr\times(S^1\times S^2)$ where $s = 0$.   
\bigskip

{\bf Step 5}\qua  Steps 5--9 establish that $C_T$ is the image 
of the thrice-punctured sphere $C_0$ by 
a $J_T$--pseudoholomorphic map, is disjoint from the 
$\theta\in \{0, \pi \}$ locus, and has its ends 
constrained by $I$ as those of $C'$ are by the third 
point of (5.42).  Given the above, then 
Lemma 5.7 follows immediately when the sequence $\{w(i)\}$ 
is unbounded as then $J_T = J$.  In 
the case where $\{w(i)\}$ is bounded, then $T = 1$ since 
$\Lambda$ is open; thus $J_T = J$ and Lemma 5.7 
again follows.
 
This step starts the task by verifying that $C_T$ is an 
HWZ subvariety.  In this regard, 
remember that an HWZ subvariety is characterized by the 
requirement that the exterior 
derivative of the contact form for $S^1 \times S^2$ has 
finite integral over both the concave side
and  the convex side ends.  Here, the contact form is 
$-(1 - 3 \cos^2 \theta) dt - \sqrt{6} \cos \sin^2 \theta
d\varphi$, and its  exterior derivative is integrable 
over each end of $\uSigma_i$.  Indeed, the value of
this integral is  bounded independent of the index $i$, for an 
application of Stokes' theorem identifies it
only as a  function of the integers in the set $I$.  
Now, as the exterior derivative of the contact
form  restricts to any $J$--pseudoholomorphic subvariety 
as a non-negative form, it follows that
this  exterior derivative is a non-negative form on the ends of $C_T$.  
Thus, (5.56) and the
dominated  convergence theorem guarantee the 
finiteness of the integral of this 2--form over each
end of  $C_T$.
\bigskip
	
{\bf Step 6}\qua  This step verifies that $C_T$ avoids the 
$\theta\in\{0, \pi \}$ locus.  In particular, it 
follows from (5.55) that $\sin \theta > \epsilon$ on 
$C_T$ granted that such is the case on each $\uSigma_i$. 
To find  such an $\epsilon$, remember that each $\underline{J}_i$ 
agrees with $J$ where $\sin \theta < \delta$. 
With this understood, take $\epsilon$  to be less 
than the minimum of $\delta, (2/3)^{1/2}$ and the values of
$\sin \theta_0$ on the closed  Reeb orbits that 
correspond to the integers in the set $I$. 
The set of points in $\uSigma_i$ where $\sin \theta < \epsilon$ 
is then compact and
$J$--pseudoholomorphic.  Were it non-empty, then the function $h/f$ 
on  each component would have a
finite minimum or maximum depending on the sign of
$\cos\theta$.   However, as previously established, 
this function doesn't have finite local maxima
or  minima on a region in a $J$--pseudoholomorphic 
subvariety which avoids the $\theta\in \{0, \pi \}$ 
locus.
\bigskip

{\bf Step 7}\qua  The topology of $C_T$ is the subject of 
Steps 7 and 8.  The story starts with a 
brief digression to recall that the critical 
point of $df$ on each $\uSigma_i$ occurs on 
the $s = 0$ locus.  As 
this locus is compact, so the critical value, 
$f_i$, of $f$ on $\uSigma_i$ takes values in some compact and 
index $i$ independent interval of $\mr$.  
Thus, there is an infinite subset of indices $i$ for which
the  corresponding set $\{f_i\}$ converges.  
Pass to this subsequence when considering $\{\uSigma_i\}$ and 
then renumber consecutively from 1.  With this now 
understood, let $f_0$ denote the limit of 
the set $\{f_i\}$.
  
To continue, remark that each end of $C_T$ has an 
associated pair of integers $(m, m')$ 
which is characterized in part by the condition that 
$m d\varphi - m' dt$ restricts to the end in 
question as an exact form.  In this regard, the 
large and constant $|s|$ slices on the end are 
circles which converge as $|s| \rightarrow \infty$  
to a multiple wrapping of a closed, 
Reeb  orbit,  whose $\theta_0$ values is 
determined by the quotients of $m$ and $m'$
by their  greatest common divisor.  Up to the action 
of the group $T$, the sign of $m$ or $m'$ finishes
the  specification of this closed Reeb orbit.  
Meanwhile, the greatest common divisor 
of $m$ and $m'$ is the multiplicity by which the 
constant $|s|$ circles wrap the limit  
closed Reeb orbit.
  
With the preceding understood, suppose that a 
given convex side end of $C_T$ is 
characterized by $(m, m')$ with $m \not=  0$.  
Then, for all $c \gg  1$, the $|f| = c$ level set in $C_T$
has a  component in this end.  It then follows from 
(5.55) that the $|f| = c$ level set of $\uSigma_i$ lies in
a  tubular neighborhood of the $|f| = c$ level set 
of $C_T$ when $i$ is large, and so the integral of
the  1--form $m d\varphi - m' dt$ about the $|f| = c$ 
level set in $\uSigma_i$ is zero.  However, the integral
of the  form $m d\varphi - m' dt$ about any component 
of a large and constant level set of $|f|$ in $\uSigma_i$
must be  a non-zero multiple of either $m p' - m' p$ 
or else $m q' - m q$.  Indeed, this can be seen
by  using Stokes theorem to compute these integrals 
for $\uSigma_i$ in the limit that $c \rightarrow \infty$.  This
implies  that $m'/m$ is either $p'/p$ or $q'/q$.  
Moreover, it then follows from (5.56) that the
signs of $m$  and $p$ agree in the former case, while 
those of $m$ and $q$ agree in the latter. 
Finally, as there  are precisely two components 
to the $|f| = c$ locus in $\uSigma_i$ when $c > |f_0|$, so (5.55)
and (5.56)  have the following additional consequence:  
There are two $m \not= 0$ convex side ends to
$C_T$,  one with $(m, m') = a^{-1} (p, p')$, 
the other with $(m, m') = b^{-1} (q, q')$.  Here, $a$ and $b$
are positive  integers.
  
Slight modifications of the preceding arguments 
also prove that $C_T$ has a single $m \not= 
0$, convex side end; and for this end, 
$(m, m') = c^{-1} (p + q, p' + q')$ where $c$ is a positive 
integer.
\bigskip

{\bf Step 8}\qua  Now consider the possibility that $C_T$ 
has an end with $m = 0$.  If such an end 
were on the convex side of $\mr\times(S^1\times S^2)$ 
and if $h > 0$ on this end, then such an end
would  contain a circle component of the $|h| = \kappa$ 
level set for all $\kappa \gg  1$.  Moreover, $dt$
would restrict  to this circle and thus to a neighborhood 
as an exact 1--form.  The absurdity of
this  conclusion in the case where neither $p'$, 
$q'$ nor $p' + q'$ is zero can be seen as follows:  
Under this last assumption, the 1--form dh also has 
exactly one critical point on each $\uSigma_i$.  
This follows from from the fact that $d\varphi$ 
has but one zero, and that $J'$ maps $d\varphi$ into a
multiple  of $dh$.  (The latter point follows from 
(5.28), (5.37) and (5.38).)   Thus, as long as
$\kappa$ is not  the absolute value of the critical 
value of $h$ on a given $\uSigma_i$, there are two
components to the $|h|  = \kappa$ level set, both 
circles and $dt$ is not exact on either one since
neither $p$, $q$, or $(p + q)$ is  zero.  Thus, 
(5.55) forbids a convex side $m = 0$ end in this
case.   A completely analogous  argument forbids 
a concave side $m = 0$ end when neither $p'$, $q'$
nor $p' + q'$ is zero.
  
On the other hand, if $p'$ were zero, then $h$ has  
precisely one non-compact level set 
on $\uSigma_i$ because, as argued in Step 3, it has but 
one such on $\Sigma_{r(i)}$.  Now, either $\kappa$ is less
than or  greater than the value of $h$ on its non-compact 
level set.  If either, then the $h = \kappa$
level set has  a single component and $q d\varphi - q' dt$ 
is zero on this component, not $dt$.

Modulo some straightforward notational changes, 
this last argument rules out the 
case $q' = 0$ too.  Thus, the only case left to consider 
has $p' + q' = 0$.  In this case, the
locus  on $\uSigma_i$ where $|h| = \kappa > 0$ consists 
of a pair of circles, and $p d\varphi - p' dt$ is
not exact on both, so  $dt$ can't be exact on either.  
\bigskip
	
{\bf Step 9}\qua  With the preceding understood, it follows 
that $C_T$ has two convex side ends 
and one concave side end, with integers $a^{-1} (p, p')$ 
and $b^{-1}(q, q')$ associated to the convex 
side ends, while $c^{-1} (p + q, p' + q')$ is 
associated to the concave side end.  Here, $a$, $b$ and
$c$  are positive integers.  Moreover, the 
integers $a$, $b$ and $c$ are constrained by the
requirement  that $p/a + q/b = (p + q)/c$ and its 
analog using $p'$ and $q'$.  In particular, because
$\Delta\not=  0$, these  constraints require $a = b = c$.  
The next paragraph establishes that
$C_T$ is the image of $C_0$ via a  $J_T$--pseudoholomorphic 
immersion, and the subsequent
paragraph then argues that $a = b = c  = 1$ and thus 
ends the proof of Lemma 5.7 where neither $p$,
$q$ nor $ p + q$ is zero.
 
To argue that $C_T$ is the image of $C_0$ via a 
$J_T$--pseudoholomorphic immersion, first 
observe that the function $f$ has precisely one critical 
value on $C_T$, this being $f_0$.  Indeed, a 
critical value at some $f = f_1$ must change the 
connectivity of the constant $f$ level sets as
this  critical value is crossed.  When the index $i$ is large, 
such a change must also occur near
the $f  = f_1$ level set on $\uSigma_i$.  In particular, such 
a change happens only at the
critical value of $f$ on $\uSigma_i$,  and as these critical 
values tend to $f_0$ as $i$ tends to
infinity, so it follows that $f_1 = f_0$.   Moreover, 
the fact that there are at most two points
in $\uSigma_i$ which share any given $(t, f)$ value  implies 
that there is precisely one critical
point of $f$ on $C_T$, and that the latter is non-degenerate.  
Then, Euler characteristic
considerations imply that $C_T$ is the image of $C_0$ 
via a  $J_T$--pseudoholomorphic map and Lemma
5.6 implies that this map is an immersion.
   
Now consider the assertion that $a = b = c = 1$.  For this 
purpose, let $R$ be very large 
and let $C_T^R \subset C_T$ denote the subset 
where two constraints are obeyed:  First, $|f| \le
R$; and  second, $|f| \ge 1/R$ at the points where 
$s \ge (2 \sqrt{6})^{-1} \ln R$.   For large $R$, this
is a smooth  manifold with boundary, where the boundary 
consists of three circles, one on each
end of  $C_T$.  Let $\uSigma_i^R \subset \uSigma_i$ 
denote the analogous set.  Being immersed, $C_T$
has a well defined normal  bundle, and then when $i$ 
is large with $R$ fixed, $\uSigma_i^R$ lies in the image
of a disk subbundle via  the exponential map.  
Moreover, the normal bundle lifts to $C_0$, and the
restriction of its  projection to the inverse image 
of $\uSigma_i^R$ can be used to define a proper
covering map from $C_0$  to itself whose degree 
is the integer $a$.  However, this means that $a = 1$
as Euler  characteristics behave multiplicatively under coverings.	    
\bigskip
	  
{\bf Step 10}\qua  This step and Step 11 discuss the 
proof of Lemma 5.7 in the case where 
one of $p$, $q$ or $p + q$ is zero.  In all of these 
three cases, the arguments are much like
those  used in the preceding steps to prove 
the lemma when none of $p$, $q$ and $p + q$ is zero.  In 
particular, only the significant differences 
are highlighted and the details are left for the most 
part with the reader.  Only the $p = 0$ case is given 
below as the $q = 0$ and $p + q = 0$ 
discussions are completely analogous to the $p = 0$ one.

As indicated above, this step and the next assume 
that $p = 0$.  For this case, as in 
Step 1, the argument begins with a definition 
of the sequence $\{w(i)\}_{i=1,2,\ldots}$.  However, the 
task here requires a preliminary digression.  
To start the digression, choose some small, 
constant $\kappa$ which is neither 0 nor $q/q'$.  
For such a choice, the form $dt - \kappa d\varphi$
restricts  without zeros on the closed Reeb orbits 
that are defined by the ends of
$\uSigma_i$.  As  argued in the proof of Lemma 5.6, 
this form, $d\nu$, has but one zero on $\uSigma_i$, with
the latter non-degenerate and of hyperbolic type.  
Meanwhile, it follows from (5.28) and
(5.31--32) that the  form $J'\cdot(dt - \kappa d\varphi)$ 
has can be written as $\alpha_0 df + \beta_0 dh$,
where $\alpha_0$ and $\beta_0$ are functions of $f$  
and $h$, and where $\alpha_0$ is positive.  Thus, each
$J_r$ in (5.48) also sends $dt - \kappa d\varphi$ 
to a 1--form $\alpha_r  df + \beta_r dh$ with both $\alpha_r$ and $\beta_r$
functions of $(f, h)$ and with $\alpha_r > 0$.  
In particular, as $\alpha_r > 0$,  there exist smooth
functions $u_r$ and $\zeta_r$ of $f$ and $h$, defined where
$\sin \theta > 0$ and such that $(\alpha_r df  + \beta dh) =  
e^{\zeta_r}du_r$.  Indeed, with $\alpha_r > 0$, the
vector field $\alpha_r
\partial_h - \beta_r
\partial_f$ is nowhere zero and  nowhere tangent to 
the constant $h$ level sets.   As a
consequence, $u_r$ can be viewed as a  measure of 
distance along the integral curves of this
vector field starting from the $f > 0$  portion 
of the $h = 0$ level set.
  
Note, by the way, that $u_r$ varies smoothly with 
$r \in [0, 1]$.  Moreover, up to additive 
constants, $u_r$ is $r$--independent on the two components 
of $\mr\times(S^1\times S^2)$ where $|s|$ is
so large  that $J_r = J$.
 
With the digression now over, choose a sequence 
$\{r(i)\}_{i=1,2,\ldots} \subset\Lambda$ with limit $T$.  As 
$du_{r(i)}$ is proportional to $J_{r(i)} (dt - \kappa d\varphi)$, 
so the function $u_{r(i)}$ has a
single critical point on $\Sigma_{r(i)}$,  and the latter 
is a non-degenerate saddle point.  In
any event, take $w(i)$ to equal minus the  value of 
the coordinate $s$ at this critical point.  As
in Part 1, if the resulting sequence $\{w(i)\}$  is not 
bounded, then pass to subsequence of $\{1,
2, \ldots \}$ for which the corresponding  sequence $\{w(i)\}$ 
is non-convergent and either
strictly positive and increasing, or strictly  
negative decreasing.  Then, re-index the
subsequence so the labels $i$ are consecutive and start  at 1.
  
Define, as before, the corresponding translated 
sequence $\{\uSigma_i\}$.  With but minor 
modifications (one such switches the roles of $f$ and $h$), 
the argument in Step 3 establish an 
index $i$ independent bound for the integral of $\omega$ 
over the intersection of $\uSigma_i$ with any
given  compact subset $K \subset \mr\times(S^1\times S^2)$.  
A repeat of Step 4 now finds a
subsequence of the  indexing set (hence renumbered 
consecutively from 1) and data $\{(S_k, m_k)\}$
such that (5.55)  and (5.56) hold.  Here again, 
$\{S_k\}$ is a finite set of distinct, irreducible,
$J_T$--pseudoholomorphic subvarieties in $\mr\times(S^1\times S^2)$ 
and each $m_k$ is a positive
integer.  Once  again, set $C_T \equiv \cup_k S_k$.
\bigskip
	
{\bf Step 11}\qua  This step discusses the verification 
that $C_T$ is the image of the thrice 
punctured sphere $C_0$ by a $J_T$--pseudoholomorphic 
map, that $C_T$ avoids the $\theta\in \{0, \pi \}$
locus,  and that its ends are constrained by $I$ as 
those of $C'$ are by the third point of
(5.42).  Given  the above, then Lemma 5.7 follows 
as argued at the beginning of Step 5.

To start, note that the arguments in Steps 5 and 6 
work here to prove that $C_T$ is an 
HWZ subvariety on which $\sin \theta$ enjoys a positive lower bound.
   
To continue the verification, pass to a subsequence 
of the index set (hence 
renumbered consecutively from 1) so that the 
resulting sequence of the critical values of the 
functions $u_{r(i)}$ on $\uSigma_i$ converge, and let $u_0$ 
denote the limiting value.  In this
regard, remember  that the sole critical point of 
each $u_{r(i)}$ on $\uSigma_i$ sits on the compact, $s = 0$
locus while $\{u_{r(i)}\}$  converges uniformly to $u_T$.  
Thus, the associated sequence of
critical values is bounded.
  
Given the preceding, then the arguments in Steps 7 and 9 
apply in this case after 
some minor modifications and establish that $C_T$ 
has precisely three ends, two on the convex 
side of $\mr\times(S^1\times S^2)$ and one on the 
concave side.  Moreover, the convex side ends
are  characterized by integers $(0, p')$ and $(q, q')$, 
while the concave side end is
characterized by  $(q, p' + q')$.  Furthermore, 
the final arguments in Step 9 also apply here
after a minor  change to prove that $C_T$ is the 
image of $C_0$ via a $J_T$--pseudoholomorphic
immersion.  In this  regard, the only substantive 
modification to the arguments in Steps 7 and 9
involves the  replacement of the function $f$ by $u_T$ 
when level sets are considered on $C_T$, and
the  replacement of $f$ by $u_{r(i)}$ when considering level sets on $\uSigma_i$. 

\sh{(h)\qua  Theorem A.4 and the number of double points}

This last section explains how Proposition 5.1's 
classifying set $I$ determines 
Proposition 3.1's double point number, $m_C$, for any 
subvariety in $I$'s component of the 
moduli space of $\aleph = 2$, thrice-punctured spheres.  
The following proposition summarizes 
the story and the remainder of this subsection 
is then occupied with its proof.  Note that 
all assertions in Theorem A.4 that concern 
the subvarieties from Theorem A.2 follow 
directly from the next proposition.
\medskip

\noindent
{\bf Proposition 5.9}\qua  {\sl Suppose that $I = \{(p, p'), (q, q')\}$ 
obeys the constraints
listed in  Proposition 5.1.  If  $C\in \fM$ is parameterized 
by $I$ as in Proposition 5.1  then $C$
has only  transversal double point singularities.  
Moreover, the double point number $m_C$ from 
Proposition 3.1 is  one half of the 
number of ordered pairs $(\eta, \eta) \in S^1 \times S^1$ with 
$\eta\not=\eta'$, neither equal to 1 and such that \
$\eta^p\eta'^{\,q} = \eta^{p'}\eta'^{\,q'} = 1$.  For
example, $m_C  = 0$ if and only if  one of the following conditions hold }
  
\begin{itemize} 
\item $\Delta  = 1$,
\item $\Delta = 2$,
\item {\sl There exists a pair  $(m, m') \in \{(p, p'), 
(q, q'), (p + q, p' + q')\}$  with both
$m$ and  $m'$ divisible by $\Delta$.}
\end{itemize}  
{\sl For a second example, if $\Delta\ge 3$ is prime and 
the third condition above does not
hold, then  $m_C = 2^{-1} (\Delta - 1)$.}
\bigskip

By the way, 
the condition in the third 
point of this proposition is equivalent to the following:  
\narrower\sl
There exist distinct pairs  $(a, a'), (b, b')$ 
 in the indicated set 
 and integers, $c$ 
 and $c'$, 
 such that
$ac' - a'c = 1$   
 and     
$bc' - b'c = -1$. 
\autonumm%(5.57)
\endnarrower

Thus, if $m_C = 0$ by virtue of the third 
condition in Proposition 5.9, then both pairs in
the  indicated set that are not mentioned consist of 
relatively prime integers.  In any event,
an  example where the third point holds has 
$p = 2$, $p' = 1$, $q = 1$, $q' = 2$ and so $\Delta = 3$.  By 
comparison, the condition in the third point of the 
proposition does not hold in the $\Delta = 3$ 
case where $p = 4$, $p' = 1$, $q = 1$ and $q' = 1$. 
\bigskip

\noindent
{\bf Proof of Proposition 5.9}\qua  The proof starts with a 
computation of the double point number 
for the spheres in the image of the composition 
$\Phi\circ\phi$ as described in Subsections 5c-e.  The 
proof then explains why this double point number 
is preserved through the deformations of 
Subsections 5f and 5g.  The details are presented in six steps.
\bigskip

{\bf Step 1}\qua  This step computes the corresponding 
double point number, $m_\phi$, for the 
image in $\mc^*\times \mc^*$ of the map $\phi$ that is 
defined in Section 5c via (5.20) from
the set  $I$.  In particular, the discussion that follows 
proves the singularities of the $\phi$ image
of $C_0$ are  purely transversal double points and 
that each assertion about $m_C$ in Proposition
5.9 also  hold for $m_\phi$.  In particular, $m_\phi = 0$ 
if and only if one of the three points
in Proposition 5.9  hold, and $m_\phi  = 2^{-1} (\Delta - 1)$ 
when $\Delta \ge 3$ is prime and the third
point in Proposition 5.9 does  not hold.
 
With the preceding   goal in mind, 
note that $z \not= w \in C_0$ are 
sent to the same point by $\phi$ if and only if
\begin{equation}
z^{p} (1 - z)^{q} = w^{p} (1 - w)^{q}  
\quad  \hbox{\rm and } \quad     z^{p'} (1 - z)^{q'} = 
w^{p'} (1 - w)^{q'} .
%(5.58)
\end{equation}
The first observation here is that these 
equalities admit no non-trivial first-order 
deformations by virtue of the fact that $\Delta \not= 0$.  
Thus, the intersections of the $\phi$ image
of small  open disks in $C_0$ are either empty or transversal.
 
Here is the second observation:  The two equations in (5.58) imply 
that $z^\Delta = w^\Delta$ and also $(1 -
z)^\Delta = (1 - w)^\Delta$ .  Thus, 
\begin{equation}
w = \eta z  \quad   \hbox{\rm and }\quad      
(1 - w) = \eta' (1 - z) 
%(5.59)
\end{equation}
where $\eta$ and $\eta'$ are distinct complex 
numbers, neither are equal to 1,  and $\eta^\Delta = 
\eta^{-\Delta'} = 1$.  The fact that $m_\phi = 0$ when $\Delta$ 
is either 1 or 2 follows
immediately from this last  point.    
This noted, assume that $\Delta >  2$ in the remainder of this step.
  
With the pair $ (\eta, \eta')$ given, subject to the 
constraints just stated, then the 
solutions $z$ and $w$ to (5.59) are given by
\begin{equation}
z = (1 - \eta^{'-1})/(1 - \eta\eta^{'-1})\quad    
\hbox{\rm and }\quad  w = (1 - \eta')/(1 - \eta^{-1}\eta') .
%(5.60)
\end{equation}
In this regard, note that there are $(\Delta - 1)(\Delta - 2)$ 
choices for $\eta$ and $\eta'$ that
satisfy the  constraints and no two choices for such a 
pair determine the same pair $(z, w)$. 
Moreover, as  $z$ is the complex conjugate of $ w$, 
no two pair of $\eta$, $\eta'$ determine either the
same $z$ or the same  $w$ in (5.60).  This last point 
implies the claim that the singularities of
$\phi$ are purely transversal  double points.  
Moreover, it implies that the number $m_\phi$ is half
the number of pairs $(\eta, \eta')$,  neither 1, both 
$\Delta$--roots of unity, distinct, and such
that $(z, w)$ in (5.60) satisfies (5.58).     
Thus, $m_\phi$ it  is at most $2^{-1}
(\Delta - 1)(\Delta - 2)$.
  
Now, not all pairs $(\eta, \eta')$ with neither 1, 
distinct and both $\Delta$--roots 
of unity produce, via (5.60), a solution to (5.58).  
In particular,  (5.58).  
places  the following additional constraints on $ (\eta, \eta')$.   
\begin{equation}
\eta^p\eta^{'q} = 1 \quad    
\hbox{\rm and }\quad  \eta^{p'}\eta^{'q'} = 1 .
%(5.61)
\end{equation}
As the solutions to (5.61) are pairs of $\Delta$--roots 
of unity, the number $m_\phi$ 
is thus seen equal to one half of the number  solutions $(\eta,
\eta')$ to (5.61) which are distinct and such that neither is equal to
1.   Thus, the general characterization of $m_\phi$ is 
identical to Proposition 5.9's
characterization of  $m_C$.

Now turn to the characterization of the $m_\phi = 0$ cases.  
The first part of this task 
verifies that $m_\phi \not= 0$ unless one of the three 
points in the proposition hold.  For this
purpose,  suppose $\Delta > 2$.  
Now, fix a primitive $\Delta$--root of unity, $\beta$.  
This done, then the
pair $(\eta= \beta^{q'}, \eta' = \beta^{-p'})$ 
solves (5.61) as does $(\eta = \beta^q, \eta' = \beta^{-p})$.  
Note that
the former obeys the condition $\eta\not=\eta'$ 
if $p' + q' \not= 0 \mod(\Delta)$ and the latter if $p +
q \not=  0 \mod(\Delta)$.  Thus, as the third point of  
the proposition is not operative, one of these
pair, say $(\beta^q, \beta^{-p})$, is still a viable 
candidate for  producing a double point.  In
particular, the latter fails to produce a double point 
only if one  of $q$ or $p$ is divisible by
$\Delta$.  For the sake of the argument, suppose 
that $p$ is divisible by $\Delta$.   This understood,
then $q$ is not divisible by $\Delta$ since $p + q$ is
not.  However, as $\Delta = p q' - p' q$  and $p$ is
divisible by $\Delta$ and $q$ is not, so $p'$ is 
divisible by $\Delta$.  But this conclusion is
nonsense  as the divisibility of both $p$ and $p'$ by $\Delta$ 
invokes the third point of the
proposition.  Thus, $p$  is not divisible by $\Delta$, 
and by symmetry, neither is $q$.  Thus, $m_\phi$
is non-zero when the third  point of the proposition is not operative.

As $m_\phi = 0$ if either of the first two points in 
the proposition hold, consider whether 
$m_\phi$ is necessarily zero when the third point 
holds.  For this purpose, suppose, for the sake
of  argument, that both $q$ and $q'$ are divisible by $\Delta$.  
This the case, then a solution
$(\eta, \eta')$ to  (5.61) with $\eta\not= 1$ 
requires $\eta^{p} = 1$ and $\eta^{p'} = 1$.  However, there is no
$\eta\not= 1$ solution to these  last two constraints 
because, as noted in (5.57), $p$ and $p'$ are
relatively prime.  (To see why  such is the case, 
fix integers $c$ and $c'$ such that $pc' - p'c =
1$.  Then, given that both $\eta^p = 1$  and $\eta^{p'} = 1$, 
so $\eta^{pc'} = \eta^{p'c}$ and thus 
$\eta^{pc'-p'c} = 1$.)   

Analogous arguments deal with the 
cases where either $p$ and $p'$ or else $p +
q$ and $p'  + q'$ are divisible by $\Delta$.

To finish the story for $m_\phi$, consider now the value of 
$m_\phi$ when $\Delta\ge 3$ is prime and the 
third point in the proposition is not operative.  
For this purpose, note first that $q$ and $q'$ are 
relatively prime.  Indeed, as $\Delta = p q' - p' q$, 
any common divisor of $q$ and $q'$ must divide
$\Delta$.   Thus, as $\Delta$ is prime and not a common 
divisior of $q$ and $q'$, so only 1 divides
both  simultaneously.  Likewise, $p$ and $p'$ are relatively 
prime.  In any event, as $q$ and $q'$
are  relatively prime, there exist a pair of 
integers, $c$ and $c'$, such that 
\begin{equation}
q c' - q' c = 1.  
%(5.62)
\end{equation}
Note that the pair $c$ and $c'$ are unique up to 
adding to $c$ an integer multiple of $q$ while 
simultaneously adding the same multiple of $q'$ to $c'$.

Now, given such a pair, define a homomorphism from the 
group of $\Delta$--roots of unity 
to itself via the formula
\begin{equation}
\eta\rightarrow\eta'[\eta] \equiv\eta^{-(pc'-p'c)} .
%(5.63)
\end{equation}
Note that this homomorphism is insensitive to the 
choice of the pair $(c, c')$ in (5.62).
  
By virtue of (5.62), the homomorphism in (5.63) 
has the property that both $(\eta'[\eta])^q 
= \eta^{-p}$ and $(\eta'[\eta])^{q'} = \eta^{-p'}$.  
Thus, the pair $(\eta, \eta'[\eta])$ solves (5.61). 
Moreover, (5.61)  together with the failure of the 
third point in the proposition imply that this
homomorphism  has trivial kernel and no fixed elements.  
As a direct consequence, the set of
solutions $(\eta,  \eta')$ to (5.61) with $\eta\not= 1$, 
$\eta'\not=1$ and $\eta \not=\eta'$ are parameterized by
the non-trival $\Delta$--roots of  unity via $\eta
\rightarrow(\eta, \eta'[\eta])$.  In this
regard, note that one element of a solution pair to (5.61)  
determines the other member in the
case where $\Delta$ is prime.  In any event, as there 
are $\Delta - 1$  non-trivial $\Delta$--roots of unity,
so $m_\phi = 2^{-1} (\Delta - 1)$ as claimed.
\bigskip

{\bf Step 2}\qua  This step counts the number of intersections 
between the $\phi$ image of $C_0$ and 
the image of the map that has the form of (5.20) but
with $a'$ different from the choice for $\phi$.  
For this purpose, fix $(a, a')$ for $\phi$, some non-negative 
real number $\delta \ll 1$, and another
small  number $\epsilon > 0$.  Note that when 
neither $p$, $q$ nor $p + q$ is zero, then $\delta = 0$ is a
permissable  choice, but not so otherwise.  In any 
event, with these parameters chosen, let $\phi_\epsilon$
denote the  version of (5.20) where $(a, a')$ are replaced 
by $(e^{i\delta\epsilon} a, e^{i\epsilon} a')$.  A
straightforward argument  along the lines of the derivation 
of (5.20) finds that the images of $\phi$
and $\phi_\epsilon $ intersect in
\begin{equation}
1 + 2 m_\phi + (gcd(p, p') - 1) + (gcd(q, q') - 1) + (gcd(p + q, p' + q') - 1)
%(5.64)
\end{equation}
points.  Moreover, all of these intersection points count with weight +1 to the intersection 
number between $\phi(C_0)$ and $\phi_\epsilon(C_0)$.
  
By way of an explanation of the various terms in (5.64), 
note that for small $\epsilon$, the left 
most term, 1, is contributed by an intersection very near 
the single critical point of $\phi$.  Of 
course, the $2m_\phi$ term counts the intersections 
near the double points of $\phi$.  The other
terms  count intersections that are, for small $\epsilon$, in the 
$\phi$ images of points that have
distance $d$ on the  Riemann sphere that is $o(\epsilon)$ 
from one or another of the punctures. 
In particular, for each  such point,
\begin{equation}
\zeta^{-1} |\epsilon| \le d \le \zeta |\epsilon|
%(5.65)
\end{equation}
where $\zeta \ge 1$ is independent of the point, $\epsilon$ 
(if less than $10^{-3}$) and the
parameters $r$, $a$ and $a'$  that appear in (5.20).  
However, $\zeta$ depends on $\delta$ if one of $p$, $q$ or $p
+ q$ vanishes.

For example, the contribution of the term that is 
third to the left, 1 less than the 
greatest common divisor of $p$ and $p'$, counts the 
number of intersection points whose $\phi$--inverse
image have distance $d$ as in (5.65) from the puncture 
at the origin in $\mc$.  The next 
term counts points whose $\phi$--inverse image has 
distance $d$ from $1 \in \mc$; and the right most 
term counts those whose $\phi$--inverse image has 
distance $1/d$ from the origin in $\mc$.
\bigskip

{\bf Step 3}\qua   With $\delta$ and $\epsilon$ 
chosen, make the parameter $r$ in (5.20) very large
and make  the parameters $R_0$, $R_1$ and $R_\infty$  
in the definition of Section 5d's map $\Phi$
extremely large as  well.  With such choices, 
the composition $\Phi\circ\phi$ immerses $C_0$ in
$\mr\times(S^1\times S^2)$ as described  in Section 5e.  
Note that $\Phi\circ\phi$ also immerses $C_0$ in
$\mr\times(S^1\times S^2)$.  Moreover, if the  parameters 
$r$, $R_0$, $R_1$ and $R_\infty$  are sufficiently
large, then the intersection points between the  $\Phi\circ\phi$ 
and $\Phi\circ\phi$ images of $C_0$ are the
images under $\Phi$ of the corresponding $\phi$ and 
$\phi_\epsilon$ versions  in $\mc^*\times \mc^*$.  Indeed,
this conclusion follows from (5.65) and Lemma~5.5.
  
With the preceding understood, choose one of the $\mr$'s worth of 
subvarieties $C \in \fM_I$ 
that approaches $\Phi\circ\phi(C_0)$ as $|s| \rightarrow\infty$.  
This done, let $C_\epsilon \in \fM_I$ denote
the image of $C$ under  the action of the subgroup of $T$ 
that sends the coordinate $\varphi$ to $\varphi + \epsilon$
and the coordinate $t$ to $t  + \delta \epsilon$.  
Thus, the large $|s|$ asymptotics of $C_\epsilon$ are
identical to those of $\Phi\circ\phi_\epsilon$.

With $C$ and 
$C_\epsilon$ defined, here is the key point to take from this step:
\begin{equation}  
\hbox{\sl The intersection number of}\; C \;\hbox{\rm with }\; 
C_\epsilon\; \hbox{\sl is given by }\; (5.64).
%(5.66)
\end{equation}
As explained next, (5.66) follows from homological 
considerations.  Indeed, fix some small 
neighborhood, $U \subset S^1 \times S^2$ of the closed 
Reeb orbits in $C$'s limit set
that is  disjoint from its rotated version, $U_\epsilon$.  
This done, then there exists $s_0 \gg  1$
so that the constant  $|s| > s_0$ slices of $C$ and 
also $\Phi\circ\phi(C_0)$ lie in $U$ while those of $C_\epsilon$ and
$\Phi\circ\phi_\epsilon(C_0)$ lie in $U_\epsilon$.   
Now, remark that for sufficiently large $s_0$,
the homology classes defined by $C$ and $\Phi\circ\phi(C_0)$  
are homologous rel the
set of points $(s, w) \in \mr\times(S^1\times S^2)$ with $|s|
\ge s_0$ and $w \in U$.   Likewise, the classes defined 
by $C_\epsilon$ and $\Phi\circ\phi_\epsilon(C_0)$ are homologous rel
those points $(s, w)$  where $|s| \ge s_0$ and $w \in U_\epsilon$.  
These last two points imply
(5.66).
\bigskip

{\bf Step 4}\qua   For each large $s_0$, view the number of $|s| < s_0$ 
intersections between $C$ and 
$C_\epsilon$ as a piece-wise constant function, $f(\epsilon; s_0)$, 
of the parameter $\epsilon$.  As
explained below, the  $\epsilon\rightarrow 0$ limit of 
$f(\epsilon, s_0)$ exists for all $s_0$ sufficiently
large and this limit is equal to $2m_C + 1$.   Granted this last 
conclusion, then the value of $m_C$
can be found by counting the the large $|s|$  intersections of 
the small $\epsilon$ versions of $C$ and $C_\epsilon$. 
The latter task occupies Steps 5 and 6  below.
  
To verify that asserted $\epsilon\rightarrow 0$ behavior of 
$f(\epsilon, s_0)$, reintroduce the $C$--version
of the  normal disk bundle $N_0 \rightarrow C_0$ and its 
exponential map $q\co N_0 \rightarrow
\mr\times(S^1\times S^2)$ as described  in (3.12).  Thus, 
$C$ is the image via $q$ of the zero section
of $N_0$.  Moreover, given $s_0$, then the  $|s| \le s_0$ 
portion of $C_\epsilon$ is obtained as the
composition of $q$ with a section, $\lambda_\epsilon$, 
of the normal  bundle $N \rightarrow C_0$ that lies in
$N_0$ over the $|s|
\le s_0$ portion of $C$ and whose norm is everywhere  
bounded by $|\epsilon|$.  This then implies that the
number of $|s| \le s_0$ intersections of $C$ and $C_\epsilon$ 
is  equal to the sum of twice $m_C$
with the algebraic count of the zeros of $\lambda_\epsilon$ 
that lie where $|s| \le  s_0$.

To prove that the latter count equals 1, remember that 
$\{C_\epsilon\}$ is constructed by the 
action of the  1--parameter subgroup of $T$ generated by 
$\delta \partial_t + \partial_\varphi$.  This implies
that $|\epsilon|^{-1} \lambda_\epsilon$  converges pointwise 
on $C$ as $\epsilon\rightarrow 0$ to the section of $N$
that is obtained by projecting the  vector field 
$\delta \partial_t + \partial_\varphi$ onto $C_0$'s normal
bundle.  In this regard, the latter vanishes only at the  
critical point of the 1--form $dt - \delta
d\varphi$, for then $\delta\partial_t + \partial_\varphi$ 
is tangent to $C$.  Now, as argued in 
Section 4f, the generic $\delta$ version of the 1--form 
$dt - \delta d\varphi$ has a single, non-degenerate
zero  on $C$.  In fact, unless one of $p$, $q$ or $p + q$ is 
zero, such is the case for the $\delta = 0$
version.  In  any event, the preceding implies that the 
generic $\delta$ and sufficiently small $\epsilon$
versions of $\lambda_\epsilon$  have one zero, this with 
multiplicity one sitting very close to
the zero of $dt$.  
\bigskip
 
{\bf Step 5}\qua  The proof of Proposition 5.9 is completed here and in 
Step 6 with the 
verification of the following claim:
\narrower\sl
{Let $(m, m')$ denote either $(p, p')$, $(q, q')$ or $(p + q, p' + q')$.  
There are}
$gcd(m, m') -
1$ {intersections between $C$ and the small $\epsilon$ versions of 
$C_\epsilon$}
{ at very large values of $|s|$
on the  corresponding $(m, m')$ end of $C$.}
\autonumm %(5.67)
\endnarrower

With regards to this claim, note that it follows 
immediately when $gcd(m, m') = 1$ from the 
identification in the previous step of the small $\epsilon$ limit 
of $\epsilon^{-1}\lambda_\epsilon$.  This understood,
suppose  that $gcd(m, m') = n > 1$.  Also, suppose 
that $m \not= 0$.  The $m = 0$ cases are treated by 
analogous arguments and so left to the reader.
  
To begin, use (2.13) to parameterize the large $|s|$ 
portion of the end in question 
using coordinates $(\tau, \rho) \in \mz/(2\pi |m|\mz) \times [\rho_0, \infty)$ 
and functions $x(\tau, \rho)$
and $w(\tau, \rho)$.  In this  regard, no generality is lost 
by assuming that both $|x|$ and $|w|$
tend to zero as $\rho\rightarrow
\infty$, and as  indicated in Proposition 2.3, their sizes 
are $\cO(e^{-\sigma \rho})$ at large $\rho$ where $\sigma > 0$ is
constant.  This  said, then the analogous end of the 
small $C$ version of $C_\epsilon$ is parameterized
using coordinates  $(\tau, \rho)$ with the pair $(x, w)$ 
replaced by $(x + \nu \epsilon, w)$ where $\nu\equiv 1 + \delta
m'/m$ is positive when $\delta$ is very small.
  
Given the preceding, it then follows that the 
intersections at large $\rho$ on the $(m, m')$ 
end of $C$ occur at and only at the values of $(\tau, \rho)$ where
\begin{equation}
x(\tau, \rho) = x(\tau + 2\pi  r m/n, \rho) + \nu\epsilon\;\;    
and\;\;    w(\tau, \rho) = w(\tau + 2\pi  r m/n, \rho)
%(5.68)
\end{equation}
for some $r \in \{1, \ldots , n - 1\}$.  Thus, (5.67) follows 
by demonstrating that there exists,
for  each such $r$ and for small $\epsilon$, precisely one pair 
$(\tau, \rho)$ where (5.68) holds.  The
remainder of  this step argues that there is at 
least one pair $(\tau, \rho)$ where (5.68) holds for any
given $r$.  The  subsequent step proves that 
there is at most one such pair for each $r$.

To establish that (5.68) can be solved, note first that 
the zero locus of the difference, 
$\uw \equiv w(\tau + 2\pi  r m/n, \rho) - w(\tau, \rho)$ 
intersects each large and constant $\rho$
circle.  This follows  from the fact that the average of $\uw$ 
over such a circle is zero.  Let $\ux$
denote the corresponding  difference of the values of $x$.  
Then, as a direct consequence of the
discussion in Step 7 of  Section 2a, the $\uw = 0$ locus 
has the structure of an oriented,
embedded graph.  Here,  the orientation is given by 
the restriction of the 1--form $d\ux$.  To
elaborate, this graph has the  following properties:
\begin{itemize}
\item  {Each vertex  is a $\underline{w} = 0$ 
critical point of $\underline{w}$.}
\item  {Each edge is an embedded, open arc, disjoint 
from the vertices and from all other 
edges.}
\item  {The 1--form $d\ux$ has nowhere zero pull-back to the edges.} 
\item  {The intersection of the graph with some 
open neighborhood of each vertex is a finite 
union of embedded, half open arcs with 
endpoints lying on the vertex, but disjoint 
otherwise.  Moreover, the tangent lines to 
the arcs at the vertex are well defined and 
disjoint. The interior of each arc is 
part of an edge of the graph.  The number of such 
arcs is non-zero and even. Exactly half of the 
arcs are oriented by $d\ux$ so $\ux$ increases 
towards the vertex while half are oriented by $d\ux$ so that $\ux$ 
decreases towards the vertex.}
\item  {If $\rho_1 \gg \rho_0$ is chosen to be sufficiently 
generic, then the $\rho=\rho_1$ locus
intersects the graph  only in its edges and this 
intersection is transversal.}\break
\hbox{}\autonum%(5.69)
\end{itemize}

It is important to note that the  conclusions in the preceding 
paragraph and the assertions in (5.69) hold for the $\underline{w}  = 0$ 
locus and for the same reasons, 
when $\underline{w}$ and $\ux \equiv (x(\tau + 2\pi  r m/n,\rho) - x(\tau, \rho)$ 
are defined using any
value of $r \in (0,  n)$.
   
The fact that  any $r \in (0, n)$ version of the $\underline{w} = 0$ 
locus has the structure 
just described and the fact that $\ux$ is a bonafide function implies a 
great deal about the large $\rho$ portion of these 
loci.  The picture of the $\underline{w} = 0$
locus  drawn below is based on two consequences of 
this structure:  First,   there
are no  non-trivial, closed, oriented paths in the 
$\uw= 0$ locus.  Second, there are no
non-compact,  oriented paths in the  $\uw = 0$ locus with 
$\rho$ unbounded at both ends.
  
These last two points imply that the $\uw = 0$ locus can be 
described as follows:    Fix  sufficiently large $\rho_1 > \rho_0$ 
making sure that 
the $\rho = \rho_1$ circle misses all vertices of the $\uw = 0$ 
locus and such that $\rho_1$ is a
regular value of  $\rho$ on each edge.  This done, 
there exist three finite sets, $\vartheta_0, \vartheta_+$ and $\vartheta_-$,
of oriented paths in  the $\rho \ge \rho_1$ portion of 
the $\uw = 0$ locus.  The union of all of
these paths is the whole $\rho \ge \rho_1$  portion 
of this locus.  Meanwhile, no two paths in
this union share any edge, although two  paths can 
intersect at a vertex.  The set $\vartheta_0$ is
distinguished by the fact that each of its paths  have both 
ends on the $\rho = \rho_1$ locus.  Meanwhile
$\rho$ is unbounded on paths in either $\vartheta_+$  or 
$\vartheta_-$; each such path has but one end on
the $\rho = \rho_1$ circle.  The paths in $\vartheta_+$ 
are distinguished
from those in $\vartheta_-$ by the fact that $\rho$ is  
unbounded in the oriented direction on a
path in $\vartheta_+$, while $\rho\rightarrow\rho_1$ in the 
oriented direction  on a path in $\vartheta_-$. 
By the way, as $\ux \rightarrow 0$ as $\rho\rightarrow\infty$,
these last points imply that $\ux$ is negative  and 
increasing in the unbounded direction
on the paths in $\vartheta_+$, but positive and decreasing 
in  the unbounded direction on the
paths in $\vartheta_-$. Finally, note that as $\rho$ 
is unbounded on the $\uw =  0$ locus, neither $\varphi_+$
nor $\varphi_-$ is empty and these sets must have the 
same number of elements.  

The preceding
picture of the $\uw = 0$ locus directly implies the following: 

\narrower\sl
With $\rho_2 > \rho_1$ and $r \in (0, n)$ fixed, there exists 
$\epsilon_1$ such that   $0 < \epsilon < \epsilon 1$
is  sufficiently generic, then each $\rho > \rho_2$ 
point where the $\epsilon$--version of (5.68) holds
lies on  an edge of a path in $\varphi_-$.  Conversely, each path 
in $\vartheta_-$ contains precisely one $\rho
> \rho_2$ point  where this same version of 
(5.68) holds. 
\autonumm%(5.70)
\endnarrower

{\bf Step 6}\qua  This step completes the proof of Proposition 5.9 
by demonstrating that the 
set $\vartheta_-$ contains at most one element.  
To start the demonstration, remark that the
number of  elements of $\vartheta_-$ is constant as $r$ 
varies in $(0, n)$.  Indeed, this number is
locally constant  because $\uw$ vanishes transversely 
on the edges of the graph and therefore
constant on $(0, n)$  as the latter is connected.
 
As an aside, note that the number of paths in $\vartheta_0$ 
can change as $r$ varies.  For 
example, a path in $\vartheta_0$ can concatonate at 
some value of $r$ with one in either $\vartheta_+$ or
$\vartheta_-$.   Conversely, a piece of a path in either 
$\vartheta_+$ or $\vartheta_-$ can become a path in 
$\vartheta_0$ at some
value of $r$.   The fact that the number of elements in 
$\vartheta_0$ can change is related to the fact that
the number  of elements in $\vartheta_0$ depends on the 
choice of the parameter $\rho_1$.
 
In any event, consider the $\uw = 0$ locus when $r$ 
is positive but very small.  In this case, 
$\uw(\tau,\rho) = (2\pi  r n/m) w_\tau(\tau,\rho) + \cO(r^2)$ 
and so the $\uw = 0$ locus converges as $r \rightarrow 0$
to the locus  where $w_\tau = 0$.  In fact, fix a 
generic $\rho_3 \gg  \rho_1$ and then take $r$ very small
but positive.  Thus, if  $\uw(\cdot, \rho_3) = 0$ at some $\tau = \tau_3$, 
then there is precisely
one $\tau\in (\tau_3, \tau_3 + 2\pi  r n/m)$ where $w_\tau, (\tau,\rho_3)  = 0$.  
Thus, the $\rho = \rho_3$
circle has at least as many zeros of $w_\tau$ as elements in 
$\vartheta_- \cup \vartheta_+$.   

Meanwhile, as $\ux = (2\pi 
r n/m) x_\tau(\tau,\rho) + \cO(r^2)$, it follows that 
$x_\tau < 0$ on at least as  many of the $w_\tau = 0$
points on the $\rho = \rho_3$ circle as there are elements in 
$\vartheta$ and $x_\tau > 0$ on at  least as many
as there are elements in $\vartheta_-$.  This said, note 
that by virtue of (2.15) or, 
equivalently, the first point in (2.19), the 
$w_\tau = 0$ locus is the $x_u = 0$ locus, where $u =
\text{sign}(p)  e^{\kappa\rho}$ with $\kappa$ defined as 
in (2.15) from the pair $ m$ and $m'$ via the
corresponding angle $\theta_0$ for  the associated closed 
Reeb orbit.  All of this understood,
then the argument subsequent to  (4.44) in Section 4g finds 
the absolute value of the Euler
class of $C_0$ no less than the  number of elements in $\vartheta_-$.  
Thus, $\vartheta_-$ has at most one
element as claimed.

\section{The structure of the $\aleph = 3$, thrice-punctured\nl
sphere moduli space}

This section describes the moduli spaces of the triply punctured 
spheres with three convex side ends that arise in Part d of the 
third point of Proposition 4.2 and in so doing provides a proof 
of Theorem A.3.  The following proposition restates the latter for use here:
\bigskip

\noindent
{\bf Proposition 6.1}\qua  \sl{The components of the moduli space of thrice-punctured, 
$\aleph = 3$ spheres which arise in Part (d) of the third point 
in Proposition 4.2 can be described as follows: }
\begin{itemize}
\item[\rm(A)]	 {The components are in 1--1 correspondence with the unordered 
sets of three pair of integers that are constrained in the 
following way:  The set in question, I, can be ordered as 
$\{(p, p')$,$(q, q')$, $(k, k')\}$ with}
\begin{itemize}
\item[$\bullet$]  $p + q + k = 0$    { and}   $ p' + q' + k' = 0$.
\item[$\bullet$]  $|k/k'| > \sqrt{3}/\sqrt{2}$.
\item[$\bullet$]  $\{(p, p'), (q, q')\}$  {obey the three constraints in Proposition 5.1.}
\end{itemize}
\item[]  {In this regard, a set $ I$ with an ordering that satisfies 
these three conditions has precisely two distinct orderings 
that satisfy the conditions.} 
\item[\rm(B)]  { This 1--1 correspondence has the following properties }
\begin{itemize}
\item[$\bullet$]  {The component of the moduli space that is labeled by  
$I \equiv  \{(p, p'),$ $(q, q'), (k, k')\}$ is a smooth manifold that 
is $\mr \times   T$ equivariantly diffeomorphic to $(0, 1) \times   \mr \times   T$.  
Here, $\mr \times   T$ acts on the moduli space via  its isometric 
action on $\mr \times   (S^1  \times  S^2 )$ and it acts on $(0, 1) \times   \mr \times   T$ by 
ignoring the $(0, 1)$ factor and acting as itself on the $\mr \times   T$ factor.}  
\item[$\bullet$]  { Thus, the quotient of the moduli space by the $\mr \times   T$ 
action is $(0, 1)$, and an identification is provided by 
composing the preceding diffeomorphism with the projection map from  
$(0, 1) \times   \mr \times   T$.}  
\item[$\bullet$]   {This quotient has a natural compactification as $[0, 1]$ 
where the two added points label the $\mr \times   T$ quotient of 
two components of the moduli space of thrice-punctured, 
$\aleph = 2$, spheres that arise in Part d of the third point 
in Proposition 4.2.  Moreover, the relevant components of 
the $\aleph = 2$ moduli space are labeled as in Proposition 5.1 
by the first two pairs from the two orderings of $I$ that 
obey the three constraints listed above in Part A.} 
\end{itemize}
\end{itemize}
\bigskip
\rm
The proof of this proposition is summarized in 
Subsection 6d, below.  The results in Subsections 
6a-d provide various pieces of the proof.  The final 
subsection explains via Proposition 6.9 how the set $I$
determines the number of double points of any immersed 
sphere in its moduli space component.  All assertions of 
Theorem A.4 about the subvarieties that appear in Theorem 
A.3 follows directly from Proposition 6.9.

\sh{(a)\qua  Initial constraints on the $\aleph = 3$, thrice-punctured 
sphere moduli spaces}

	  Let $C \subset  \mr \times   (S^1  \times S^2 )$ denote one of the thrice-punctured 
spheres under consideration.  As described 
in Section 4a, each end of $C$ determines a unique, 
ordered pair $(m, m')$ of integers.  However, the 
three pair, $\{(p, p'), (q, q'), (k, k')\}$ determined by 
$C$ are not completely independent.  Some initial constraints on this set are:
\bigskip

{\bf Constraint 1}\qua No pair in $\{(p, p'), (q, q'), (k, k')$ 
can vanish identically.
\bigskip

{\bf Constraint 2}\qua $p + q + k = 0$  {and}   $p' + q' + k' = 0$.
\bigskip

{\bf Constraint 3}\qua  $\Delta \equiv  p q' - q p' \not=   0$.
\bigskip

{\bf Constraint 4}\qua If $(m, m') \in  \{(p, p'), (q, q'), (k, k')\}$ 
and $m < 0$, then $m' \not=   0$.
\bigskip

With regard to the third constraint, note that $\Delta$ is also 
equal to $p k' - k p'$ and $k q' - q k'$ and thus the $\Delta \not=   0$ 
assertion does not require the ordering of the set $\{(p, p'), (q, q'), (k, k')\}$.  
The arguments for the first three constraints are analogous 
to those given in Section 5a for the triply punctured spheres 
with $\aleph = 2$, and so omitted except for the remark that the 
second constraint here differs by a sign from the second 
constraint in Section 5a.  The fourth constraint follows 
from the observation that an $m < 0$ and $m' = 0$ end would be 
asymptotic to the $\theta \in  \{0, \pi  \}$ locus.  

\sh{(b)\qua  The structure of each component}

Let $I \equiv    \{(p, p'), (q, q'), (k, k')\}$ denote an unordered 
set of integer pairs that obeys the three constraints listed.  
Given such $ I$, let $\fM_I$ denote the space of pseudoholomorphic, 
$\aleph = 3$, thrice-punctured spheres from Part d of the third 
point of Proposition 4.2 whose ends are described as in (5.1) by $I$.
 
An argument that is essentially identical to the one used 
in Section 5b's proof of Proposition 5.2 establishes that 
each component of $\fM_I$ is a smooth manifold, $ \mr \times   T$ equivariantly 
diffeomorphic to $\cN \times   \mr \times   T$; here $\cN$ is either a circle or an 
interval.  As seen below in this subsection $\cN$ is, in all cases, 
an interval.  The next proposition makes the official assertion:
\bigskip

\noindent
{\bf Proposition 6.2}\qua {\sl Let $I = \{(p, p'), (q, q'), (k, k')\}$ denote an 
unordered set of pairs of integers that obeys the constraints 
listed in the statement of Proposition 6.1.  Then, every 
component of  $\fM_I$ is diffeomorphic to  $(0, 1) \times   \mr \times   T$.  
Moreover, there exists such a diffeomorphism that is $\mr \times   T$
equivariant when $ \mr \times   T$ acts on $\fM_I$ via its isometric action 
on $\mr \times   (S^1  \times S^2 )$, and on $(0, 1) \times  \mr \times  T $ via 
its natural action on the factor $\mr \times   T$.}
\bigskip  

The remainder of this subsection contains the following proof.
\bigskip

\noindent
{\bf Proof of Proposition 6.2}\qua  The proof is lengthy and so divided into seven steps.
\bigskip
 
{\bf Step 1}\qua  Let $\cH \subset  \fM_I$ denote a component.  Suppose $C \subset  \cH$ 
and focus attention on an end of $C$ characterized by integers 
$(m, m')$ as in (5.1).  In the case where $m \not=   0$, this end is 
parameterized at large values of $|s|$ by coordinates $(\tau, u)$ 
where $\tau \in  \mr/(2\pi   |m| \mz)$ and sign$(m)u \in  [R, \infty)$ for some suitably
large constant $R$.  This parameterization is given by (2.19) with $p = m$ and $p' = m'$.  
Thus, $x$ and $y$ are functions of $\tau$ and $u$, which obey (2.20).  In the $m = 0$ case, the $(0, m')$ 
end of $C$ is parameterized by $(\tau, u)$ where $\tau \in  \mr/(2\pi  |m'| \mz)$ and 
sign$(m')u \in  [R,\infty)$ 
via $(t = t_0 + x, f = y, \varphi = \tau, h = u)$.  Here, $x$ and $y$ are functions of 
$\tau$ and $u$ which obey 
the version of (2.20) where the role of $(t, f)$ are switched with that of $(\varphi, h)$.
 
The following lemma describes the $|u| \rightarrow\infty$ behavior of $y$ 
for any pair $(m, m')$.  For this purpose, introduce 
\begin{equation}
\sigma_0  \equiv  4 |\cos \theta_0| |m'/m| (1 + (m'/m)^2 \sin^2 \theta_0)^{-1}.  
%(6.1)
\end{equation}
Note that $\sigma_0$  is determined by the pair $(m, m')$.
\bigskip  

\noindent
{\bf Lemma 6.3}\qua {\sl The restrictions  $|u|^{\sigma_0} 
y(\cdot, u)$ to the constant $u$ 
circles converge in the $C^\infty$ topology as $|u| \rightarrow  \infty$ to the 
constant function.  Moreover, the assignment to each 
$C \in \cH$ of the resulting constant defines a smooth function on $\cH$.} 
\bigskip

\noindent 
{\bf Proof of Lemma 6.3}\qua  The convergence of $y(\cdot, u)$ in the 
case where $m' = 0$ is stated and proved as Lemma 5.8 in 
the previous section.  Save for notational changes, the
$ m ' \not=   0$ cases have the identical argument.  For example, 
in the cases where $m \not=   0$, one should set $w \equiv  6^{-1/2} 
(1 + (m'/m) \sin^2 \theta_0)^{-1/2}$ \\ 
$y/ (u \sin^2 \theta)$ and $\rho\equiv  (1 + 
(m'/m)^2 \sin^2 \theta_0)^{-1/2} 6^{-1/2} \ln(|u|)$.  This done, then 
(5.52) holds for $\lambda$ provided that the $\sqrt{6}$ in the lower 
right entry in $L_0$ is replaced $\sigma_1  = (1 + (m'/m) \sin^2 
\theta_0)^{1/2} (\sqrt{6} + \sigma_0 )$.  The modifications to the argument 
given subsequent to (5.52) are of a similar nature.
  
With the preceding understood, let $Y\co \cH \rightarrow \mr$ denote 
the function which assigns to $C$ the $|u| \rightarrow  \infty$ limit of 
the $C$ version of  $|u|^{\sigma_0}y(\cdot, u)$.  To prove that $Y$ has bounded 
first derivative on $\cH$, consider that a tangent vector to 
$C$ in $\cH$  can be represented where $|s|$ is large on the 
$(m, m')$ end of $C$ by functions $(x', y')$ of $\tau$ and $u$ which 
obey the linearized version of (2.20).  Moreover, $|x'|$ 
and $|y'/u|$ are bounded as $|u| \rightarrow  \infty$.  With this understood,
essentially the same argument that proved Lemma 5.8 
proves that  $|u|^{\sigma_0}y'$ also has a unique limit as $|u| \rightarrow  \infty$.  
The latter conclusion implies that the function $Y$ is 
at least Lipschitz on $\cH$.  The behavior of the higher 
order derivatives of $Y$ can be analyzed in a similar 
manner using the parameterization of a neighborhood of 
$C$ in $\cH$ that is provided by Proposition 3.2 and Steps 5 
and 6 of its proof.  The details of this analysis is 
straightforward and so omitted.
\bigskip

{\bf Step 2}\qua  The purpose of this step is to highlight 
and then prove the following:
\bigskip

\noindent
{\bf Lemma 6.4}\qua  {\sl The set $I = \{(p, p'), (q, q'), (k, k')\}$ 
has a canonical element for which the corresponding 
function $Y\co \cH \rightarrow  \mr$ described by Lemma 6.3 is never zero.  
In this regard, if $(k, k')$ denotes the canonical element, 
then $k < 0$ and $k' \not=   0$.}
\bigskip

\noindent
{\bf Proof of Lemma 6.4}\qua  As is demonstrated below, this 
constraint is a consequence of $C$'s avoidance of the 
$\theta \in  \{0, \pi  \}$ locus.  To start the argument, consider 
an end of $C$ where the integer $m$ is not zero.  Then 
the ratio $h/f$ is defined at large $|s|$ on this end 
and for large $|s|$, behaves as
\[
h/f = (m'/m) \sin^2 \theta_0 + \text{sign}(m)|u|^
{-\sigma_0-1}  (Y + o(1))\;.
%(6.2)
\]
In particular, when $m'$ and $Y$ have the same sign, then 
$|h/f| > |m'/m| \sin^2 \theta_0$ at large $|s|$ on this end; and
when their signs differ, then $|h/f| < |m'/m| \sin^2 \theta_0$ 
at large $|s|$ on the end in question.
  
To see where this leads, note that at least one of $p$, 
$q$, and $k$ is negative.  If one, denote it by $k$.  If two, 
then, as argued in the eighth step of Section 5's proof 
of Constraint 4, the number $|m'/m| \sin^2 \theta_0$ is smallest 
for precisely one, so denote the latter by $k$.  In either 
case, let $C_- \subset  C$ denote the component of the subset of 
points where $f < 0$ which contains the $(k, k')$ end.  
Note that $k' \not=   0$ because of the fourth constraint in 
Subsection 6a, above.

With the preceding understood, suppose that the sign 
of the $(k, k')$ end's version of $Y$ is opposite that of 
$k'$.  Then $|h/f|$ at large $|s|$ on this end is less than 
its $|s| \rightarrow  \infty$ limit.  Meanwhile, $h/f$ cannot change sign 
in the closure of $C_-$ since a zero of $ h$ on $C$ where $f \le 0$ 
is precluded when $C$ is disjoint from the $\theta \in  \{0, \pi  \}$ 
locus.  Thus, $|h/f|$ is nowhere zero in $C_-$, diverges 
as the boundary of $C_-$ is approached, and its infimum is 
not its $|s| \rightarrow  \infty$ limit on $C_-$.  This can happen only if $h/f$ 
has a local extreme point in $C_-$ which is a forbidden event
as (4.21) makes $h/f$ subject to the maximum principle.  In 
this regard, a point must be added in the case where the 
$f = 0$ locus is non-compact.  For such to happen, then one 
of $p$ and $q$ must be zero, say $p$.  This implies that the 
large $|s|$ portion of the $f = 0$ locus sits on the $(p = 0, p')$ 
end of $C$.  However, $|h|$ diverges at large $|s|$ on this end, 
while Lemma 6.3 guarantees that $f$ limits to zero as $|s| \rightarrow  \infty$
on this end.  Therefore, $|h/f|$ diverges on the $(p = 0, p')$
end of $C$ as well.

Now consider the possibility that $Y = 0$ for this same $(k, k')$ 
end of $C$.  Were such the case, the argument given below establishes 
that there are still points at large $|s|$ on this end where 
$|h/f| < |k/k'| \sin^2 \theta_0$.  Thus, the argument from the previous 
paragraph applies and yields the same contradictory conclusions.  
With the preceding understood, suppose that $Y = 0$ on an end of $C$
characterized by integers $(m, m')$ with $m \not=   0$.
    
For this purpose, return to the $(m, m')$ version of (2.20).  
Set $v = \text{sign}(m) u$ and view $w 
\equiv  \sin^{-2}  \theta y$ as a function of $\tau$ 
and $v$.  With these changes, the second of the equations in 
(2.20) has the schematic form
\begin{equation}
x_\tau = w_v + \sigma_0  v^{-1} w + r (w),
%(6.2)
\end{equation}
where $|r(w)| \le \zeta v^{-1} w^2$ with $\zeta$ some constant.  This last 
equation yields an impossible conclusion under 
the assumption that $y$ and thus $w$ have a definite 
sign at large values of $v$.  To see how, introduce $\uw(v) \equiv  \int w(\tau, v) d\tau$ 
to.  Then (6.2) implies that
\begin{equation}
|(v^{\sigma_0} \uw)_v| < \zeta  v^{\sigma_0}  v^{-1} \int w^2(\tau, v)\; d\tau\; .
\end{equation}
%(6.3)

Now, as $Y = 0$, so $\lim_{v\rightarrow  \infty}  v^{\sigma_0} w = 0$ 
and so integration of (6.3) 
yields the inequality 
\begin{equation}
v^{\sigma_0} |\uw(v)| < \zeta \sup_{v'\ge v}(v^{'\sigma_0}
|w(v')| ) \; | \int_{ v'\ge v} \uw(v') v^{'-1}\; dv'|\; .
%(6.4)
\end{equation}
Note that the derivation of (6.4) exploits the lack of 
a sign change in $\uw$.  In any event, recall that $\sigma_0  > 0$, 
and thus (6.4) implies that
\begin{equation}
v^{\sigma_0}  |\uw(v)| < \zeta' v^{-\sigma_0}  \sup_{v'\ge v}
(v^{'\sigma_0}  |\uw(v')|) 
%(6.5)
\end{equation}
where $\zeta' = 2 \zeta \sigma_0^{-1} \sup (v^{'\sigma_0} |w(v')| )$.  Of course 
(6.5) asserts the ridiculous when $v$ is large.

An argument that is almost identical to that just 
given also establishes the following:  
\bigskip

\noindent
{\bf Lemma 6.5}\qua  {\sl If there is a single $(m, m') \in  I$ with $m > 0$, 
denote the latter by $(p, p')$.  If there are two $m > 0$ 
pairs in $I$ and if $m'$ has the same sign for both, denote 
by $(p, p')$ that for which $ |m/m'| \sin^2 \theta_0$ is greatest.  
Finally, if there are two $m > 0$ pairs and the corresponding 
$m'$ are opposite, denote by $(p, p')$ that for which  $|m'|$ is 
greatest.  Then, the $(p, p')$ version of $Y$ is never zero and 
its sign is opposite that of $p'$.}
\bigskip

{\bf Step 3}\qua  An embedded copy of $\cN$ in $\cH$ is obtained by the 
choice of a slice across $\cH$ for the action $\mr \times   T$.  For the 
purpose of defining such a slice, remember that $\mr \times   T$ is 
the subgroup in Isom$(\mr \times   (S^1  \times S^2 ))$ where the action of 
the $\mr$ factor is generated by $\partial_s$ while $\partial_t$ and $\partial_\varphi$ generate 
the actions of the respective $S^1$  factors of $T = S^1  \times   S^1$.   
This said, the definition of the slice involves three 
conditions.  Before these conditions are stated, remark 
that at least one of the integers $p$, $q$ and $k$ must be 
positive, and with this understood, take $(p, p')$ as in 
Lemma 6.5.  Then, the first condition assert the vanishing 
of the constant term on the right hand side of (5.1) for 
$(p, p')$ end of $C$, and the second asserts the vanishing of 
the $(k, k')$ version of this same constant.  As was the case 
with the $\aleph = 2$ triply punctured spheres, these first two 
conditions force the vanishing of the $(q, q')$ version of 
(5.1)'s right hand side constant term.  The final condition 
on the slice involves the $(k, k')$ version of the function $Y$ 
from Lemma 6.3.  According to Lemma 6.4, this function is 
nowhere zero, and with this understood, a slice of the $\mr$ 
factor in $\mr \times   T$ is obtained by requiring that $|Y| = 1$.  
Here, it is important to note that $a \in  \mr \subset  \mr \times   T$ 
sends $Y$ to $e^{-\sqrt{6}(1+\sigma_0)a}  Y$.
\bigskip

{\bf Step 4}\qua  With $\cN \subset  \cH$ so identified, define a function, 
$F\co\cN \rightarrow  \mr$ 
by assigning to each $C \in  \cN$ the value of the $(p, p')$ end version 
of the function $Y$.  With $F$ understood, then Proposition 6.2 is 
an immediate consequence of  
\bigskip

\noindent
{\bf Lemma 6.6}\qua  {\sl  So defined, the function $F$ has no critical points on $\cN$.}  
\bigskip

The proof of Lemma 6.6 occupies this step and the next two 
steps in Proposition 6.2's proof.
\bigskip

\noindent
{\bf Proof of Lemma 6.6}\qua   To begin the proof, suppose, for the 
sake of argument, that $C \in  \cN$ were a critical point of this 
function $F$.  To study the implications for $C$, it is 
important to remember that the tangent space to $C$ in $\cH$ 
consists of the 4--dimensional space of bounded sections 
of $C$'s normal bundle that are annihilated by the $C$ 
version of the operator in (3.5).  In this regard, 
remember that $C$ is the image of $C_0$ via a pseudoholomorphic 
immersion and so the pullback of $T(\mr \times   (S^1  \times S^2 ))$ via this 
immersion splits orthogonally as a direct sum of oriented, 
$J$--invariant 2--plane bundles.  The differential of the 
immersion canonically identifies one of these subbundles 
with $TC_0$; the orthogonal complement of the latter is the 
`normal bundle to $C$' in $\mr \times   (S^1  \times S^2 )$.    Thus, with a
pseudoholomorphic immersion of $C_0$ into $\mr \times   (S^1  \times S^2 )$ 
understood, vector fields pullback to $C_0$ as sections 
of the pullback of $T(\mr \times   (S^1  \times S^2 ))$.  For example, the 
vectors $\partial_s$, $\partial_t$ and $\partial_\varphi$ 
pullback to sections over $C_0$ of 
the pullback of $T(\mr \times   (S^1  \times S^2 ))$.

With these reminiscences now ended, note that $T\cH|_C$ has a 
three-dimensional subspace which is spanned by the 
projections onto $C$'s normal bundle of the vector fields 
$\partial_s$, $\partial_t$ and $\partial_\varphi$.  Denote 
the latter by $\eta^s$, $\eta^t$ and $\eta^\varphi$,
respectively.  The tangent vector, $\eta$, to $\cN$ at $C$ provides 
a fourth basis element for $T\cH|_C$.  In particular, to say 
that $\eta$ is tangent to $\cN$ at $C$ places three constraints on 
the behavior of $\eta$, whether or not $C$ is a critical point 
of $F$.  The assumption that $C$ is a critical point places 
an additional constraint on $\eta$.  As argued 
subsequently, these constraints are not mutually self-consistent.  
\bigskip

{\bf Step 5}\qua   This step describes the aforementioned four 
constraints on $\eta$.  In particular, these constraints 
involve the asymptotics of $\eta$ on both the $(p, p')$ and 
$(k, k')$ ends of $C$. In particular, fix attention first 
on the $(k, k')$ end of $C$.  Reintroduce the coordinates 
$(\tau, u)$ on this end and then any element of $T\cH|_C$ at 
large $|u|$ can be viewed as a pair $(x', y')$ of functions 
of $\tau$ and $u$ which obey the linearized version of (2.20) 
subject to the constraint that $|x'|$ and $|y'/u|$ are 
bounded as $|u| \rightarrow  \infty$.  As indicated in the proof of
Lemma 6.3, both $x'$ and $|u|^{\sigma_0} y'$ have limits as $|u| \rightarrow  \infty$.  
With this point understood, here are the first and 
second requirements for $\eta$:   The $\eta$ version of the 
functions $x'$ and $y'$ obey
\begin{itemize}
\item  $\lim_{|u|\rightarrow \infty} x' = 0$
\item  $\lim_{|u|\rightarrow \infty} |u|^{\sigma_0} y' = 0$
\autonum%(6.6)
\end{itemize}

\noindent
Indeed, these constraints are required for movement 
along $\cN$ in the direction of $\eta$ to preserve the defining 
slice conditions that come from the $(k, k')$ end.
  
By the same reasoning, the $(p, p')$ version of the 
function $x'$ must also satisfy the first point (6.6) if $\eta$ 
is to be tangent to $\cN$.  This is the third requirement for $\eta$.  
Finally, the $(p, p')$ version of the function $y'$ must 
satisfy the second point in (6.6) if $C$ is to be a critical
point of the function $F$.  This is the fourth constraint on $\eta$.  
\bigskip

{\bf Step 6}\qua    As is demonstrated momentarily, it is 
impossible for (6.6) to hold at both the $(k, k')$ 
and $(p, p')$ ends of $C$.  The proximate cause of 
this incompatibility stems from the second point in 
(6.6) which requires both the $(k, k')$ and $(p, p')$ 
versions of the function $y'$ to have zeros at arbitrarily 
large values of the coordinate $|u|$.  Indeed, the proof 
that such zeros exist amounts to a linearized version of 
the argument in the proof of Lemma 6.4 that surrounds 
(6.2)--(6.5).  However, more must be said about these 
zeros of $y'$.  In particular, both the $(k, k')$ and $(p, p')$ 
versions of $x'$ and $y'$ have the following property: 

\narrower\sl
Given $R > 0$, there exists $\epsilon > 0$ such that the function  $x'$ 
take every value in $(-\epsilon, \epsilon)$ at some $|u| > R$ zero of $y'$.
\autonumm%(6.7)
\endnarrower

The validity of (6.7) is justified below, so accept it for now.

This large $|u|$ behavior of the $\eta$ version of $y'$ on the 
two ends of $C$ is used below to find a constant $r$ and 
two distinct points in $C_0$ where $\eta^\varphi + r \eta$ is zero.  But 
$\eta^\varphi + r \eta$ can have at most one zero for the following 
reason:  All of its zeros count with the positive weight 
to any Euler class computation.  Thus, if $\eta^\varphi + r \eta$ has 
two or more zeros, then a standard perturbation argument 
provides at least two zeros for $\eta^\varphi + r \eta + \epsilon \eta^t$ for any 
small in absolute value constant $\epsilon$.  Of course, the latter 
is also annihilated by the operator $D$, so all of its zeros 
also count with positive weight.  Thus any Euler number 
calculation that uses $\eta^\varphi + r \eta + \epsilon \eta^t$ gives 2 or more for 
an answer.  On the other hand, for a suitably generic 
choice of $\epsilon$, the Euler count of the zeros of 
$\eta^\varphi + r \eta + \epsilon \eta^t$ 
can be used to compute the expression $\la e, [C]\ra - 2 m_C$ that 
appears in Proposition 3.2.  And, as this expression equals 
1, so $\eta^\varphi + r \eta + \epsilon \eta^t$ has but one zero on $C_0$.
  
To find $r$ which makes $\eta^\varphi + r \eta$ vanish at two points, first
note that on both $(k, k')$ and $(p, p')$ ends of $C$, the $\eta^\varphi$ 
version of the functions $x'$ and $y'$ has $x' \equiv  1$ and $y' \equiv  0$.  
Thus, $\eta^\varphi$ and $\eta$ are colinear at all $(\tau, u)$ where the $\eta$ 
version of $y'$ is zero.  With the preceding understood, 
choose $R \gg 1$ and then take $r$ so that $1/r$ is less than 
the smallest of the two $\epsilon$ values provided by the $(k, k')$ 
and $(p, p')$ versions of (6.7).  Then (6.7) guarantees a 
zero for $\eta^\varphi + r \eta = 0$ on the $(k, k')$ end of $C$ 
and another on the $(p, p')$ end.
\bigskip
	
{\bf Step 7}\qua   It remains now only to justify (6.7).  For this 
purpose, note that the pair $(x', y')$ obey the linearized 
version of (2.20) and so 
\begin{equation}
x'_u = - g^{-2} \sin^{-2}    \theta \; y'_\tau + \sigma_1  y' \quad     
\hbox{\rm and }\quad      x'_\tau = \sin^{-2}  \theta \; y'_u
+\sigma_2  y' ,
%(6.8)
\end{equation}
where $\sigma_{1,2}$ are bounded functions of $\tau$ and $u$.  Note that 
(6.8) asserts that $(x', y')$ obeys a Cauchy Riemann equation 
to order $|y'|$.  In particular, arguments akin to those 
used in Appendix A of \cite{T4} can be brought to bear and 
establish the following:  First, the simultaneous zeros of 
$y'$ and $dy'$ are isolated points.  Second, the zero set of 
$y'$ constitutes an oriented graph, $G$, whose vertices are 
these simultaneous zeros.  Here, the orientation on the 
edges is defined by the pullback of $dx'$.  Finally, each 
vertex has a non-zero and even number of incident edges, 
with precisely half oriented by $dx'$ to point into the vertex.
  
The preceding properties of $G$ have (6.7) as a consequence.  
To see that such is the case, note first that $G$ has no closed, 
oriented loops since its orientation form, $dx'$, is the 
differential of a bonafide function.  Thus, as there are
points in $G$ where $|s|$ is arbitrarily large, so there is
an oriented path in $G$, starting on some large, but constant
$|s|$ slice and on which $|s|$ is unbounded.  Indeed, such a
path is constructed by following an oriented edge until it
hits a vertex, and then continuing out from the vertex along
another edge with the outward pointing orientation.  Only 
finitely many paths of this sort which start at some finite 
$|s|$ slice will have $|s|$ bounded since none are closed 
loops and only finitely many can intersect any given 
constant $|s|$ locus.  Use up the finitely many with both
ends hitting a chosen constant $|s|$ locus and then start 
another at a very large $|s|$ zero of $y'$.  The latter must
have one end where $|s|$ is unbounded.

Next, remark that $y'$ is negative on one side of any 
edge, and positive on the other; and which side has 
which sign is determined by the orientation because 
(6.8) asserts that $dx' \wedge dy' > 0$ where $y' = 0$ and $dy' \not=   0$.  
In particular, as $y'$ is also a bonafide function, there are 
two distinct paths in $G$ which start on some large, but 
constant $|s|$ slice on which $|s|$ is unbounded.  More to 
the point, $|s| \rightarrow  \infty$ in the direction oriented by $dx'$ on 
one of them, while $|s| \rightarrow  \infty$ in the direction oriented 
by $Ðdx'$ on the other.   However, since $|x'|$ limits to 
zero as $|s| \rightarrow  \infty$, so $x'$ must be negative on the first of 
these paths, and $x'$ must be positive on the second.   
Given that $|x'| \rightarrow  0$ as $|s| \rightarrow  \infty$, these last conclusions 
directly imply (6.7).

\sh{(c)\qua  The limits on a component of $\fM_I$}

Supposing that $\fM_I$ is non-empty, focus attention on a 
component $\cH \subset  \fM_I$.  As was just proved, $\cH/(\mr \times   T)$ is an open 
interval, and so any slice in $\cH$ of the $\mr \subset  \mr \times   T$ action 
is non-compact.  This said, the focus here is on the 
limiting behavior of non-convergent sequences in such a 
slice.  In this regard, the discussion is simplest when the 
integer $q$ from $I$ is non-zero, for in this case, none of $p$, $q$ 
or $k$ is zero and $dt$ has exactly one zero upon pullback to 
any $C \in  \cH$.  In particular, with $q \not=   0$, a slice, $\cH^s$, of the $\mr$ 
action on $\cH$ is obtained by requiring that $s$ equal 0 at the 
zero of $dt$.  In the case where one of $p$, $q$ or $k$ is equal to 
zero, choose some $\kappa \not=   0$ and not equal to one of the finite 
ratios $p'/p$, $q'/q$ or $k'/k$.  For such $\epsilon$, the form $dt - \kappa d\varphi$
pulls back to any $C \in  \cH$ with but one zero.  In this case, a 
slice, $\cH^s$, of the $\mr$ action on $\cH$ is obtained by the condition 
that $s = 0$ at the zero of $dt -\kappa d\varphi$.

The assignment to each $C \in  \cH^s$ of the maximum of $s$ on $C$ defines 
a continuous function $S\co \cH^s \rightarrow  [0,\infty)$.  
This function $S$ is proper; 
as can be shown with an essentially verbatim repetition of 
arguments from Section 5g.  Thus, a sequence $\{C_i\} \in  \cH^s$ has 
no convergent subsequences provided that $\{S(C_i)\}$ is unbounded.  
In particular, fix attention on such a sequence where $S(C_i) > 
S(C_{i+1})$.  For this sequence, minimal modifications of arguments 
from Section 5g find a subsequence (hence renumbered consecutively
from 1) of $\{C_i\}$ which converges as described 
in Step 4 of Section 5g to one of the $\aleph = 2$, triply punctured 
spheres from Part d of the third point of Proposition 4.2.

Let $C$ denote this limit triply punctured sphere.  Up to the 
action of $T$, this sphere is described by a set $I_C = \{(a, a'), 
(b, b')\}$ of pairs of integers with $\Delta_C = a b' - b a'$ positive.  
In this regard, arguments from Section 5g can be employed to 
deduce that $I_C$ is a subset of $I$.  As is argued below, the 
pair $(m, m')$ from $I$ that is missing from $I_C$ obeys
\begin{equation}
\sqrt{3}/\sqrt{2} < |m'/m| <\infty .
%(6.9)
\end{equation}
These last conclusions have the following consequence:  
\narrower\sl
{Let $I$ denote an unordered set of pairs of integers that 
corresponds to a}
{component of the moduli space of $\aleph = 3$, 
thrice-punctured spheres from}
{Part~d of the third point of 
Proposition 4.2.  Then $I$ satisfies}
{ the three constraints in 
Part A of Proposition 6.1.}
\endnarrower 
Before discussing (6.9), agree on the notational 
convention that orders the integers in $I$ so that 
$I_C = \{(p, p'), (q, q')\}$; thus $(k, k')$ is the pair from $I$ 
that is missing from $I_C$.  Now, to prove (6.9), translate 
the elements in the convergent sequence $\{C_i\}$ by the $\mr$--factor 
in $\mr \times   T$ so that $S(C_i) = 0$.  Use 
$\{\underline{C}_i\} \subset  \cH$ to denote the 
resulting sequence.  The zero of $dt$ on $\underline{C}_i$ now occurs 
where $s = -S(C_i)$ and so the corresponding sequence of 
points has no convergent subsequence in any part of $\mr \times   
(S^1  \times S^2 )$ where $s$ is bounded from below.   Arguments 
such as found in Section 5g apply now to the sequence 
$\{\underline{C}i\}$ and find a subsequence which converges as 
described in Step 4 of Section 5g to a cylinder from
Example 6 in Section 4a, parameterized by the pair 
whose respective components are the quotients of $k$ 
and $k'$ by their greatest common divisor.  This implies 
that $|k'/k| > \sqrt{3}/\sqrt{2}$ because a cylinder with $|k'/k| < \sqrt{3}/\sqrt{2}$ 
is ruled out by the arguments from Step 6 of Section 5g 
which find $\delta_I > 0$ such that $\sin \theta > \delta_I$ on every 
$C \in  \fM_I$.

Being an open interval, $\cH/(\mr \times   T)$ has two ends and so two 
distinct, triply punctured, $\aleph = 2$ spheres can be expected 
as limits of non-convergent sequences in $\cH^s$.  As argued in 
the next subsection, this is indeed the case; and the next 
lemma describes the relation between resulting two versions 
of Proposition 5.1's data set $I_C$:
\bigskip

\noindent
{\bf Lemma 6.7}\qua   {\sl Suppose that $I_C = \{(p, p'), (q, q')\}$ 
satisfies 
Proposition 5.1's constraints with the additional requirement 
that $k \equiv  p + q$ and $k' \equiv  p' + q'$ obey 
$|k'/k| > \sqrt{3}/\sqrt{2}$.  Then $I_C$ 
contains a unique pair $(m, m')$ such that the following is true:}
\begin{itemize}
\item  $|m'/m| > \sqrt{3}/2$.
\item  {\sl When $(m, m') = (p, p')$, define $I_{C'} \equiv  
\{(q, q'), (-k, -k')\}$; 
and when $(m, m') = (q, q')$, define $I_{C'} \equiv  \{(-k, -k'), (p, p')\}$.  
Then $I_{C'}$ also satisfies Proposition 5.1's constraints.}
\end{itemize}
Note here that $I_C$ is never the same as $I_{C'}$.
\bigskip

\noindent
{\bf Proof of Lemma 6.7}\qua   The discussion here considers various 
cases in turn, with the sign$(p) = \text{sign}(q)$ case treated first.  
Consider first the existence of $(m, m')$ in this case.  To 
start the story, note that $\Delta\equiv  p q' - q p' > 0$ is both $k' p - k p'$ 
and $k q' - q k'$, and as the signs of $p$, $q$ and $k$ agree, so
\begin{equation}
q'/q > k'/k > p'/p \;.
%(6.10)
\end{equation}
Thus, at least one of $q'/q$ or $p'/p$ has absolute value greater 
than $\sqrt{3}/\sqrt{2}$.   Note that the first and third constraints in 
Proposition 5.1 hold automatically for $I_{C'}$ no matter which of 
the pair $(p, p')$ or $(q, q')$ ends up as $(m, m')$.  Thus, only 
the second constraint is open to debate.
  
For the sake of argument, suppose that $k'/k > 0$ so that $q'/q$ 
is guaranteed larger than $\sqrt{3}/\sqrt{2}$.   Either the obvious choice, 
$\{(-k, -k'), (p, p')\}$, for $I_{C'}$ obeys the constraint of the 
second point of Proposition 5.1 or not.   
If not, then $|p'/p|$ 
had better be greater than $\{\sqrt{3}/\sqrt{2}\}$, and it is.  Indeed, if the 
obvious $\{(-k, -k'), (p, p')\}$ does not obey the constraint in 
Proposition 5.1's second point, then $p' + k'$ must be negative 
and both $p'$ and $k'$ must also have the same sign unless $p' = 0$. 
Now, $p' \not=   0$ is precluded by the second point of Proposition 
5.1; the latter requires $q' > 0$ when $p' = 0$, so $k' = q' + p'$ 
makes $k' > 0$ when $p' = 0$.  Thus, both $k'$ and $p'$ are negative.  
Since $k' < 0$, so $k < 0$ and then the positivity of $\Delta = p k' - p' k$ 
requires $p < 0$.  But then $|p'/p| > \sqrt{3}/\sqrt{2}$ by virtue of the 
third constraint in Proposition 5.1.

Thus, $(m, m') = (p, p')$ satisfies the first constraint of 
Lemma 6.7 and so the set  $\{(q, q'), (-k, -k')\}$ is also a 
candidate for $I_{C'}$.  To verify the second constraint of 
Proposition 5.1, consider the consequences of its violation.  
This occurs when $q' > 0$ and $k' > 0$.  But, as previously argued, $k' < 0$ 
or else the set $\{(-k, -k'), (p, p')\}$ would serve for $I_{C'}$.
 
Now consider the uniqueness question for this first case 
when $k'/k > 0$.  In particular, in the light of the preceding 
discussion, it is enough to verify that  $\{(q, q'), (-k, -k')\}$ 
cannot satisfy the requirements of Lemma 6.7 when 
$\{(-k,$ $-k'), 
(p, p')\}$ does.  For this purpose, note that when $-k'$ and $p'$ 
have opposite signs then the second 
point of Proposition 5.1 applies to $\{(-k, -k'), (p, p')\}$ 
to force $k' > 0$ and $p' > 0$.  If the second point is to
apply to $\{(q, q'), (-k, -k')\}$ as well, then $q'$ must be 
negative.  However, if $k'$ is positive, then so is $k$ and 
thus so is $q$, and then (6.10) is violated.  On the other
hand, if $Ðk'$ and $p'$ have the same sign, then $k'$ and $p'$ 
have opposite signs.  This requires $q'$ and $p'$ to have 
opposite signs so $q' > 0$ and $p' < 0$ by an application of 
the second point in Proposition 5.1 to $\{(p, p'), (q, q')\}$.  
But now the fact that $q' > 0$ and $Ðk' < 0$ violates this 
same point when applied to $\{(q, q'), (-k, -k')\}$.

The argument for the case where the signs of $p$ and $q$ agree 
and $k'/k < 0$ is obtained from the preceding argument by 
changing various signs.
 
Consider next the case for Lemma 6.7 when $p$ and $q$ have 
different signs.  Again,  the existence question is treated 
first.  To start, suppose that $p$ is negative.  Thus $(p, p')$ 
obeys the first constraint of Lemma 6.7 because of the third
constraint in Proposition 5.1.  In this case, $I_{C'} = \{(q, q'), 
(-k, -k')\}$ does not also satisfy the constraints from 
Proposition 5.1 only if both $k' > 0$ and $q' \ge  0$.  In this 
regard, consider first the case that $k > 0$ too.  As $k$ and 
$q$ have the same sign, the positivity of $\Delta$ demands $q'/q > 
k'/k > \sqrt{3}/\sqrt{2}$, so the $(q, q')$ pair all obeys the first 
constraint in Lemma 6.7.  Moreover, now $I_{C'} = \{(-k, -k'),
(p, p')\}$ has the same sign for both primed pair and so 
obeys the second constraint in Proposition 5.1.
  
Suppose next that $p < 0$ and $k < 0$.  As noted, the second 
constraint of Proposition 5.1 fails for the corresponding 
$I_{C'}$ only when both $k' > 0$ and $q' \ge  0$.  However, this last
possibility is precluded by virtue of the fact that 
$\Delta  = k q' - q k'$ is positive.

Now suppose that $p > 0$ and $q < 0$.  Here, the $(q, q')$ pair 
obeys the first constraint in Lemma 6.7.  The corresponding 
$I_{C'}$ does not obey the second constraint of Proposition 5.1 
only when $k' < 0$ and $p' > 0$.  However, this requires $q' < 0$ 
and so is precluded by the $\{(p, p'), (q, q')\}$ version of the 
second constraint of Proposition 5.1.
  
The argument for the uniqueness assertion in Lemma 6.7 for 
the case where $p$ and $q$ have different signs is straightforward 
and left to the reader.

Finally, consider the case where one of $p$ or $q$ vanishes.  
The existence argument where $p = 0$ is given below; the $q = 0$
existence argument and the uniqueness arguments are left to 
the reader as both are reasonably straightforward.

In the case where $p = 0$, then $p' \not=   0$ and so the $(p, p')$ pair
obeys the first condition in Lemma 6.7.  The second condition 
of Proposition 5.1 is violated for the corresponding $I_{C'}$ only 
if $q' > 0$ and $k' > 0$.  Since $p = 0$, the positivity of $\Delta$ requires
that $k$ and $q$ have opposite sign to $p'$.  If $p' < 0$, then both
have positive sign and $q'/q > k'/k > \sqrt{3}/\sqrt{2}$; thus $(q, q')$ obeys
the first constraint in Lemma 6.7.  Moreover, the corresponding 
$I_{C'}$ obeys the second constraint as its primed entries have the 
same sign.  Meanwhile, if $p' > 0$, then $k$ and $q$ are negative and
now $(q, q')$ obeys the first constraint in Lemma 6.7 by virtue of 
the fact that $\{(p, p'), (q, q')\}$ obeys the third constraint in 
Proposition 5.1.  In this case, the corresponding $I_{C'} = \{(-k, -k'),
(p, p')\}$ obeys the second constraint in Proposition 5.1 
because its first primed entry is negative and the second is positive.

\sh{(d)\qua  The existence of $\aleph = 3$, 
thrice-punctured sphere moduli spaces}

This subsection gives a construction for points in the moduli 
spaces that are described by Proposition 6.1.  For this purpose, 
fix a set $I$  as described in Part A of the proposition, and fix 
one of the two ways to order $I$ as $\{(p, p'), (q, q'), (k, k')\}$ 
so that the three points of this same Part A hold.  The subset 
$\{(p, p'), (q, q')\} \subset  I$ then labels a component of the moduli 
space of $\aleph = 2$, thrice-punctured spheres from Proposition 5.1. 
There is a unique point, $C'$, in this same moduli space where
the following two constraints hold:  The first constraint 
requires the vanishing of the constant term on the right 
side of (5.1) for each end of $C'$.  The statement of the 
second constraint depends on whether or not one of $p$, $q$ or 
$p + q$ is zero.  If none vanish, then $C'$ is constrained so 
that the sole point where $dt$ is zero on $TC'$ occurs where 
$s = 0$.  In the case where one of these integers equals zero, 
fix $\epsilon \not=   0$ nor equal to any $m'/m$ for $(m, m') \in  I$.  
With $\epsilon$ 
fixed, the second constraint requires that $s = 0$ at the 
sole point in $C$ where $dt + e d\varphi$ is zero on $TC'$.

Now consider that the pair $(k', k)$ labels a component of 
the moduli space of pseudoholomorphic cylinders from Example 6 
of Section 4a.  To be precise here, introduce $m$ to denote the 
greatest, common (positive) divisor of $k$ and $k'$, and then 
introduce $\underline{k} \equiv  k/m$ and $\underline{k'} \equiv  k'/m$.  
It is the pair $(\uk, \uk')$ 
that determines the moduli space.  In any event, with $S \in  \mr$ 
fixed, the latter moduli space contains a unique point, $C_S$, 
where the following two conditions hold:  First, the constant
term on the right-hand side of (5.1) is zero for each end of 
$C_S$.  Second, the maximum value of $s$ is $C_S$ is equal to $S$.

When $S \gg 1$, these subvarieties $C'$ and $C_S$ are used below to 
construct an element, $C$, in the $I$--labeled component, $\fM_I$, of 
Proposition 6.1's moduli space of thrice-punctured.  Of course, 
the construction of $C$ verifies the existence assertion of
Proposition 6.1.  Moreover, certain properties of the 
construction are used to prove some of the other assertions (4.15)
of Proposition 6.1.  In particular, the latter are described in the
following result.
\bigskip

\noindent
{\bf Proposition 6.8}\qua   {\sl Fix a set, $I$, of three pair of integers 
that satisfy the constraints in Part A of Proposition 6.1, 
and fix one of the two orderings $I$ as $\{(p, p'), (q, q'), (k, k')\}$ 
so as to satisfy these constraints.  Given $\epsilon > 0$, there exists 
$S_0 > 1$ and an embedding, $\Phi$, from $[S_0,\infty)$ into $\cH_I$ with the 
following properties:}  
\begin{itemize}
\item  {\sl If $C = \Phi(S)$, then $C = C_- \cup   C_+$ 
where $C_\pm$ are open sets, 
the $(k, k')$ end of $C$ is in $C_-$, the $(p, p')$ and $(q, q')$ 
ends are in $C_+$, and}
\begin{equation}
\text{sup}_{z\in C_-}    \dist(z, C_S) + \text{sup}_{z\in C_-}    
\dist(z, C') < \epsilon
%(6.11)
\end{equation}
\item {\sl  The projection of the image of $\Phi$ to $\fM_I/(\mr \times   T)$ defines 
a proper embedding of the closed half line $[S_0,\infty)$ into $\fM_I/(\mr \times   T)$.}
\end{itemize}

{\sl Moreover, if $C \in \fM_I$ and $S(C) > S_0 + 1$, then there exists 
$\gamma \in  \mr \times   T$ such that $\gamma(C)$ is in the image of the version 
of the map $\Phi$ as defined by one or the other of the two 
orderings of $I$ that satisfy the constraints in Part A of Proposition 6.1.} 
\bigskip
  
This proof of this proposition is given momentarily.  
Accept it for now, and here is the completion of the following proof.
\bigskip

\noindent
{\bf Proof of Proposition 6.1}\qua    Consider first the assertions 
of Part A:  The fact that $\fM_I = \emptyset$ unless $I$ obeys the 
constraints in Part A was established in Subsection 6c, 
above, while the sufficiency of these constraints follows 
directly from the asserted existence of Proposition 6.8's 
map $\Phi$.  (Remember that Lemma 6.9 justifies the claimed 
existence of two orderings for $I$ that obey Part A's conditions.)

As for the assertions in Part B, the claimed structure of 
a component of $\fM_I$ as the product $(0, 1) \times   (\mr \times   T/\Gamma_I)$ is 
proved as Proposition 6.2.  The fact that $\fM_I$ has but a 
single component follows from the third point of 
Proposition 6.8.  Finally, the asserted properties of the 
two point compactification of $\fM_I/(\mr \times   T)$ are direct 
consequences of the first and second points of Proposition 6.8.
\bigskip

The remainder of this subsection contains the following proof.
\bigskip

\noindent
{\bf Proof of Proposition 6.8}\qua   The construction of $C$ from $C'$ 
and $C_S$ constitutes a by now standard `gluing construction' 
for pseudoholomorphic subvarieties.  Such gluing 
constructions were first introduced by Floer \cite{F} as analogies
to similar constructions that are used to construct self-dual 
connections on 4--manifolds.  In any event, the details of the
construction are left to the reader save for the outline 
that follows.

To begin the outline, remark that all of the gluing 
construction can be construed in the following way:  
The subvariety $C$ is obtained from a small normed section, 
$\eta$, of certain complex line bundle over $C_0$.  Here, this line 
bundle is the normal bundle to a symplectic immersion of $C_0$
whose image is very close to $C$ on the $s < 3S/4$ portion of 
the latter and whose image is very close to $C_S$ on the 
complement of the $s < -S/4$ portion of this cylinder's $(-\uk, -\uk')$ 
end.  To be more explicit, the image of this immersion agrees 
with $C'$ on the $s < -S/4$ portion of the latter and it sits very 
close to $C_S$ on the complement of the $s < 3S/4$ portion of this 
cylinder's $(-\uk, -\uk')$ end.  In particular, near this part of 
$C_S$, the immersion is an embedding obtained by composing the
exponential map from (3.12) with a suitably chosen, m-multivalued 
section of the cylinder's normal bundle.   Then, where $s \in  [-S/4, 3S/4]$, 
these two portions of the immersion are extended and sutured 
together using cut-off functions.

In any event, the resulting immersion of $C_0$ is everywhere
symplectic and pseudoholomorphic for an almost complex 
structure that is pointwise close to $J$, even where $s \in  [S/4, 3S/4]$.  
(Proposition 2.3 is
needed for this last conclusion.)  With the  immersion of 
$C_0$ constructed, the subvariety $C$ is
the image of the  composition of a section, $\zeta$, of $C_0$'s normal bundle with an 
exponential map such as that used in (3.12).

Defined as it is by $\eta$, the subvariety $C$ is pseudoholomorphic 
provided $\eta$ solves an inhomogeneous version of the $C_0$ version
of the equation in the third point of (3.12).  In this regard, 
fix $n \ge  0$ and then the relevant inhomogeneous term is $\cO(e^{S/\zeta})$ 
as measured with the $L_{n,\delta}^2$ norm in (3.14).  Here, $\zeta \ge  1$ is 
independent of $S$ when the latter is large.  (Proposition 2.3 is 
required to obtain these bounds.)  The implicit function theorem
can be employed as in Step 1 of the proof of Proposition 3.2 to 
find a small pointwise norm solution to this inhomogeneous 
equation which is $L^2$ orthogonal to the kernel of $D$.  The existence of 
such a solution verifies the first point of Proposition 6.8.

Of course, the application here of the implicit function theorem 
to the inhomogeneous version of (3.12) requires some sort of $S$ 
independent bound on the norm of the inverse of the operator $D$ 
that appears in this version of (3.12).  In the case at hand, 
such a bound exists when $S$ is large; its existence follows from 
the fact that the analogous $D$ has trivial cokernel on both $C'$ 
and on the $m$--fold covering space of $C_S$.  The fact that the $C'$ 
version of $D$ has trivial cokernel is proved as Proposition 4.6.  
Analogous arguments, or a separation of variables analysis can
be used to establish this same conclusion for the operator on 
the $m$--fold cover of $C_S$.  (Remember that $C_S$ is invariant under 
an $S^1$  subgroup of $T$.)
  
The use made here of the implicit function theorem has the 
second point of Proposition 6.8 as a straightforward consequence.

Implicit function theorem conclusions typically have 
uniqueness as well as existence ramifications, and the 
former play a key role in the proof of the final assertion 
of Proposition 6.8.  In particular, as applied here, the 
implicit function theorem can be used to prove that $\eta$ is 
the only small pointwise norm solution to this inhomogeneous 
equation which is $L^2$ orthogonal to the kernel of $D$.  This 
last conclusion has the following consequence:  

\narrower\sl
{There exists $\epsilon> 0$ such that if $S$ is large and if a 
given  $C \in  \fM_I$ has all points distance $\epsilon$ or less from 
$C_S \cup   C'$, then there exists $\gamma \in  \mr \times   T$ such that 
$\gamma(C)$ is in the image of $\Phi$.}
\autonumm%(6.12)
\endnarrower\rm

With (6.12) understood, the proof of the final assertion 
of Proposition 6.8 begins by noting that the two orderings 
of $I$ that satisfy the constraints of Part A of Proposition 6.1
each define, for $S \ge  1$, a pair $(C_S, C')$.  Then, the arguments 
from Section 6c lead to the following conclusion:  Given $\epsilon > 0$, 
a good portion of the points in $C$ have distance $\epsilon$ or less
from one of these versions of $C_S \cup  C'$ when $S(C)$ is large.  
To be more precise, fix $R > 1$ and there exists $S_0$ such that 
when $S(C) > S_0$, then the portion of $C$ with distance less than 
$\epsilon$ from $C_S \cup   C'$ contains all points except possibly those in 
the set, $U$, of points that both lie in the bounded component 
of the locus in $C$'s $(k, k')$ end where $s \le S(C) - R$, and lie
where $s \ge  R$. 

This said, turn attention to the subset $U$.  For this purpose,
reintroduce the cylinder $C$ where both $\uk \varphi - \uk' t = 0$ and $h/f = 
(\uk'/\uk) \sin^2 \theta_0$.  According to Proposition 1.2, when $R$ is large, 
then all of the $s \le S(C) - R$ half cylinder in $C_S$ whose boundary 
is nearest to $U$ has distance less than $e^{-R/\zeta}$ from $\uC$.  The same 
is true for the $s \ge  R$ half cylinder in $C'$.  Here, $\zeta \ge  1$ can be 
taken to be independent of $R$.  (Remember that $|\uk'/\uk| > \sqrt{3}/\sqrt{2}$ and
that $\sin^2 \theta_0$ is bounded away from zero.)  With these points taken, 
the manner of convergence described in Section 6c implies the 
following:  When $\epsilon$ is fixed, $R$ is large, and then $S(C)$ is very 
large, the points in $U$ where $s = R$ or $s = S(C) Ð $R obey
\begin{equation}
|\uk \varphi - \uk' t| < \epsilon^2 + 2 e^{-R/\zeta}\quad    
\hbox{\rm and }\quad    |(h/f) - (\uk'/\uk) \sin^2 \theta_0| 
\ll \epsilon^2 + 2 e^{-R/\zeta}. 
%(6.13)
\end{equation}

Now, as remarked previously on many occasions, both $\uk \varphi - \uk' t$ 
and $h/f$ are subject to the mercy of the maximum principle.  
This implies that $|\uk \varphi -\uk' t|$ is bounded by 
$\epsilon^2 + 2 e^{-R/\zeta}$ on 
the whole of $U$, as is $|(h/f) - (\uk'/\uk) \sin^2 \theta_0|$ provided 
that $f \not=   0$ on $U$.  In this regard, the convergence described 
in Section 6c implies that the sign of $f$ is that of $\uk$ on 
both boundaries of $U$, and as $f$ is also subject to the 
maximum priniciple, so $f$ is, indeed, nowhere zero on $U$.
  
Next, remark that as $U$ lies in the $(k, k')$ end of $C$ and 
is far from a critical point of $f$ (where $s = 0$), so it 
can be parameterized as in (2.19) by variables $(\tau, u)$ 
using the functions $(x, y)$ of $(\tau, u)$.  The conclusions of 
the previous paragraph now translate as the assertion that 
both $|x|$ and $|u|^{-1} |y|$ are bounded by 
$\epsilon^2 + 2 e^{-R/\zeta}$ on $U$.
  
The preceding point is relevant because the metric distance 
of a point in $U$ from the cylinder $\uC$ is bounded by $\zeta (|x| + |u|^{-1}
|y|)$ where $\zeta$ is independent of $R$ and $\epsilon$ when $S(C)$ is large.  Thus,
small $\epsilon$, large $R$ and very large $S(C)$ implies that $C$ obeys the 
assumptions in (6.12) for one or the other of the two possible 
versions of $C_S \cup   C'$ that are provided by $I$.  

\sh{(e)\qua  The number of double points}

The next proposition describes the double point number 
$m_C$ of a subvariety $C$ that comes from any of the moduli 
space components that appear in Proposition 6.1.  In 
particular, this proposition directly implies all assertions of 
Theorem A.4 that concern subvarieties from Theorem A.3.
\bigskip

\noindent
{\bf Proposition 6.9}\qua   {\sl Let $I = \{(p, p'), (q, q'), (k, k')\}$ 
satisfy the requirements in Proposition 6.1 so as to 
describe a component of the moduli space of $\aleph = 3$, 
thrice-punctured spheres that appear in Part d of the 
third point of Proposition 4.2.  Let $C$ denote some 
subvariety in this component of the moduli space.  
Then, $m_C$ is equal to one half of the number of ordered
pairs  $(\eta, \eta') \in  S^1  \times   S^1$  with $\eta \not=   \eta'$, neither 
equal to 1 and such that $\eta^p\eta^{'q} = \eta^{p'}\eta^{'q'} = 1$.}
\bigskip  

Note that the count for $m_C$ is the same as that given in
Proposition 5.9 for the case when $\{(p, p'), (q, q')\}$ 
determine the latterÕs moduli space component.  This 
said, it follows that $m_C = 0$ here if and only if either 
$|p q' - p' q|$ is one or two, or if both integers in at 
least one of the pair $(p, p')$, $(q, q')$ and $(k, k')$ is 
evenly divisible by $p q' - q p'$.

By the way, in spite of the appearance to the contrary, 
the preceding and Proposition 6.9's count of $m_C$ is 
insensitive to the ordering of the set $I$.
  
The remainder of this subsection contains the following proof.
\bigskip

\noindent
{\bf Proof of Proposition 6.9}\qua   Let $C$ be as described in the
proposition, and, with $\delta$ small and sufficiently generic
and $\epsilon > 0$ but small, let $C_\epsilon$ denote 
the translate of $C$ 
via the element in $T$ that moves $t$ to $t + \delta \epsilon$ and $\varphi$ 
to $\varphi + \epsilon$.  
Here, $\delta$ is chosen so that the Reeb orbits that are 
determined by $C$ are pairwise distinct from those deterimined 
by $C_\epsilon$.  This understood, then it follows from the gluing 
construction of $C$'s moduli space in Section 6d that the 
intersection number between $C$ and $C_\epsilon$ is given by the 
formula in (5.64) where $m_\phi$ is the number that is described 
in Step 1 of the proof of Proposition 5.9.  In this regard, 
note that $m_\phi$ is equal to Proposition 6.9's asserted value for $m_C$.

The preceding understood, then the relationship between 
the afore-mentioned $C-C_\epsilon$ intersection number and the number 
$m_C$ is as described in Step 4 of the proof of Proposition 5.9;
the same argument gives the justification.  Moreover, 
the claim in (5.67) is valid here as well, and with an 
identical proof.  Of course, this implies that $m_C$ has the value 
asserted by Proposition~6.9.

\rk{Acknowledgements}

The author thanks Jerrel Mast for pointing out some errors in a
preprint version of this article.
 
The author is supported in part by the National Science Foundation.

\end{document}

%% file: gtoutput.tex
\def\ifplaintex{\expandafter\ifx\csname documentclass\endcsname\relax}

\ifplaintex 
\else
\fi

\expandafter\ifx\csname beginpicture\endcsname\relax
\expandafter\ifx\csname documentclass\endcsname\relax
\input pictex \else
\input prepictex \input pictex \input postpictex \fi\fi

\def\gt{{\mathsurround=0pt\it $\cal G\mskip-2mu$eometry \&\ 
$\cal T\!\!$opology}}        %  journal title in recommended style

\def\gtp{{\mathsurround=0pt\it $\cal G\mskip-2mu$eometry \&\ 
$\cal T\!\!$opology $\cal P\!$ublications}}  % GT publications

\def\lognumber#1{\def\thelognumber{#1}}
\def\volumenumber#1{\def\thevolumenumber{#1}}
\def\papernumber#1{\def\thepapernumber{#1}}
\def\volumeyear#1{\def\thevolumeyear{#1}}

\def\pagenumbers#1#2{\def\startpage{#1}\def\finishpage{#2}}
\def\published#1{\def\publishdate{#1}}
\def\proposed#1{\def\theproposer{#1}}
\def\seconded#1{\def\theseconders{#1}}
\def\received#1{\def\receiveddate{#1}}

\def\accepted#1{\def\accepteddate{#1}}
\def\asciititle#1{\def\theasciititle{#1}}
\def\covertitle#1{\def\thecovertitle{#1}}

\long\def\asciiabstract#1{\long\def\theasciiabstract{#1}}

\let\\\par\let\thelognumber\relax
\let\thevolumenumber\relax\let\thepapernumber\relax
\let\thevolumeyear\relax\let\thesamplenumber\relax\let\startpage\relax
\let\finishpage\relax\let\publishdate\relax\let\receiveddate\relax
\let\reviseddate\relax\let\accepteddate\relax\let\theasciititle\relax
\let\thecovertitle\relax\let\theasciiauthors\relax
\let\theasciiabstract\relax
\let\theasciiemail\relax\let\theshortauthors\relax\let\theshorttitle\relax

\long\def\maketitlep{   % start of definition of \maketitlep

\count0=\startpage

\gt\hfill      %   Journal title (top left) 
\beginpicture
\setcoordinatesystem units <0.33truein, 0.33truein> point at 2.2 0.9
\setplotsymbol ({$\cal G$})
\plotsymbolspacing=9truept
\circulararc 315 degrees from 0 1 center at 0 0
\setplotsymbol ({$\cal T$})
\circulararc 315 degrees from 1 -1 center at 1 0
\endpicture
\break
{\small\ifx\thesamplenumber\relax % sample?  
Volume \else Sample
\fi\thevolumenumber\ (\thevolumeyear)
\startpage--\finishpage\nl
Published: \publishdate}
\vglue 0.5truein plus 0.4fil minus 0.1truein

{\parskip=0pt\leftskip 0pt plus 1fil\def\\{\par\smallskip}{\ifplaintex\large
\else\Large\fi\bf\thetitle}\par\medskip}   

\vglue 0pt plus 0.1fil 

{\parskip=0pt\leftskip 0pt plus 1fil\def\\{\par}{\sc\theauthors}
\par\medskip}

\vglue 0pt plus 0.1fil 

{\small\parskip=0pt\let\newline\\
{\leftskip 0pt plus 1fil\def\\{\par}{\sl\theaddress}\par}
\expandafter\ifx\theemail\relax    % email address?
\relax\else\vglue 5pt plus 0.02fil minus 2pt\def\\{\stdspace{\rm 
and}\stdspace} 
\cl{Email:\stdspace\tt\theemail}\fi
\ifx\theurl\relax                  % URL given?
\relax\else\vglue 5pt plus 0.02fil minus 2pt\def\\{\stdspace{\rm 
and}\stdspace}
\cl{URL:\stdspace\tt\theurl}\fi\par}

\vglue 7pt plus 0.3fil minus 3pt

{\bf Abstract}
\vglue 5pt plus 0.1fil minus 2pt

\theabstract

\vglue 7pt plus 0.3fil minus 3pt

{\bf AMS Classification numbers}\quad Primary:\quad \theprimaryclass

Secondary:\quad \thesecondaryclass

\vglue 5pt plus 0.3fil minus 2pt

{\bf Keywords}\quad \thekeywords

\vglue 10pt plus 0.5fil minus 5pt

{\small  Proposed: \theproposer\hfill Received: \receiveddate\nl
Seconded: \theseconders\hfill 
\ifx\reviseddate\relax                         % paper revised?
Accepted: \accepteddate                        % no
\else
Revised: \reviseddate                          % yes
\fi}
\eject
}       %  end of definition of \maketitlep

\let\maketitlepage\maketitlep
\let\maketitle\maketitlepage

\font\phead=cmsl9 scaled 950
\font=cmsl9 scaled 1050
\font\pnum=cmbx10 scaled 913
\font\lnum=cmbx10 
\font\pfoot=cmsl9 scaled 950
\font=cmsl9 scaled 1050
\ifplaintex
\headline{\vbox to 0pt{\vskip -4.5mm\line{\small\phead\ifnum
\count0=\startpage ISSN 1364-0380 (on line)
1465-3060 (printed) \hfill {\pnum\folio}\else\ifodd\count0\def\\{ }% 
\ifx\theshorttitle\relax\thetitle\else\theshorttitle\fi\hfill{\pnum\folio}
\else\def\\{ and }{\pnum\folio}\hfill\ifx\theshortauthors\relax\theauthors
\else\theshortauthors\fi\fi\fi}\vss}}
\footline{\vbox to 0pt{\vglue 0mm\line{\small\pfoot\ifnum\count0=\startpage
\copyright\ \gtp\hfill\else
\gt, Volume \thevolumenumber\ (\thevolumeyear)\hfill\fi}\vss
}}
\else
\makeatletter
\def\@oddhead{{\small\lhead\ifnum\count0=\startpage ISSN 1364-0380 (on line)
1465-3060 (printed) \hfill {\lnum\number\count0}\else\ifodd\count0
\def\\{ }\ifx\theshorttitle\relax \thetitle \else\theshorttitle\fi\hfill
{\lnum\number\count0}\else\def\\{ and }{\lnum\number\count0}
\hfill\ifx\theshortauthors\relax 
\theauthors\else\theshortauthors\fi\fi\fi}}\def\@evenhead{\@oddhead}
\def\@oddfoot{\small\lfoot\ifnum\count0=\startpage\copyright\ \gtp\hfill\else
\gt, Volume \thevolumenumber\ (\thevolumeyear)\hfill\fi}
\def\@evenfoot{\@oddfoot}
\makeatother
\fi

\newwrite\gtoutfile
\long\gdef\makeheadfile{  %%% start of definition of \makeheadfile
{\def\\{, }\def\s{ }
\immediate\openout\gtoutfile head.xxx
\immediate\write\gtoutfile{To: math@arxiv.org}
\immediate\write\gtoutfile{Subject: put or rep NNNNN:pppp}
\immediate\write\gtoutfile{--text follows this line--}
\immediate\write\gtoutfile{Proxy-for: \ifx\theasciiauthors\relax
\theauthors\else\theasciiauthors\fi\s<\ifx\theasciiemail\relax\theemail\else\theasciiemail\fi>}
\immediate\write\gtoutfile{\noexpand\\}
\immediate\write\gtoutfile{Authors: \ifx\theasciiauthors\relax
\theauthors\else\theasciiauthors\fi}
{\def\\{ }\immediate\write\gtoutfile{Title: \ifx\theasciititle\relax
\thetitle\else\theasciititle\fi}}
\immediate\write\gtoutfile{Subj-class: GT or SG or MG etc}
\immediate\write\gtoutfile{MSC-class: \theprimaryclass\ifx\thesecondaryclass\relax\else, \thesecondaryclass\fi}
\immediate\write\gtoutfile{Journal-ref: Geom. Topol. \thevolumenumber
(\thevolumeyear) \startpage-\finishpage}
\immediate\write\gtoutfile{Comments: Published by Geometry and Topology at}
\immediate\write\gtoutfile{\s\s http://www.maths.warwick.ac.uk/gt/GTVol\thevolumenumber/paper\thepapernumber.abs.html}
\immediate\write\gtoutfile{\noexpand\\}
\immediate\write\gtoutfile{}
\ifx\theasciiabstract\relax
\immediate\write\gtoutfile{\theabstract}\else
\immediate\write\gtoutfile{\theasciiabstract}\fi
\immediate\write\gtoutfile{}
\immediate\write\gtoutfile{\noexpand\\}
\immediate\write\gtoutfile{}
\immediate\closeout\gtoutfile}}  %%% end of definition of \makeheadfile

\def\maketitlepage{\maketitlep\makeheadfile}
\let\maketitle\maketitlepage